\documentclass[10pt]{amsart}
\usepackage{amsmath,amssymb}

\usepackage[dvipsnames,svgnames]{xcolor}
\definecolor{myblue}{RGB}{47,79,127}

\usepackage{tikz}
\usetikzlibrary{calc,arrows.meta}
\usepackage{subcaption}

\usepackage{mathrsfs}

\usepackage[active]{srcltx}
\usepackage{graphicx,color}

\usepackage{mdframed}

\usepackage{epic}
\usepackage{pstricks}
\usepackage{comment}
\usepackage{amsthm}

\newtheorem{theorem}{Theorem}[section]
\newtheorem{proposition}[theorem]{Proposition}

\newtheorem{lemma}[theorem]{Lemma}

\newtheorem{corollary}[theorem]{Corollary}
\theoremstyle{definition}

\newtheorem{remark}[theorem]{Remark}

\newtheorem{hypothesis}{Hypothesis}[section]

\numberwithin{equation}{section}

\newcommand{\ii}{\mathrm i}
\def\abs{\mathrm{abs}}

\def\R{\mathbb{R}}

\def\N{\mathbb{N}}
\def\C{\mathbb{C}}
\def\ZZ{\mathbb{Z}}

\def\la{\lambda}

\def \sgn{\operatorname{sgn}}

\def \supp{\operatorname{supp}}

\def\dist{\textrm{dist}\,}
\newcommand{\wt}{\widetilde}
\newcommand{\lan}{\langle}
\newcommand{\ran}{\rangle}

\def\Re{{\mathrm{Re\,}}}
\def\Im{{\mathrm{Im\,}}}
\def\cA{{\mathcal A}}

\def\cD{{\mathcal D}}
\def\cX{{\mathcal X}}
\def\cF{{\mathcal F}}

\def\cK{{\mathcal K}}
\def\cL{{\mathcal L}}
\def\cM{{\mathcal M}}

\def\cR{{\mathcal R}}
\def\cS{{\mathcal S}}
\def\cU{{\mathcal U}}

\def\cW{{\mathcal W}}

\def\cO{{\mathcal O}}

\def\cX{{\mathcal X}}
\def\cG{{\mathcal G}}
\def\cT{{\mathcal T}}
\def\cM{{\mathcal M}}
\def\cB{{\mathcal B}}
\def\cP{{\mathcal P}}

\def\cQ{\mathcal{Q}}

\def\m{\mathcal{M}}

\newcommand{\ud}{\mathrm{d}}

\renewcommand{\t}[1]{\widetilde{#1}}

\DeclareMathOperator{\Dom}{Dom}

\DeclareMathOperator{\Ran}{Ran}
\def\Ker{{\mathrm{Ker}\,}}

\definecolor{darkred}{rgb}{0.7,0.1,0.1}

\newcommand{\ip}[2]{\langle #1,#2\rangle}
\newcommand{\norm}[1]{\left\|#1\right\|}

\title[Diffusion phenomena for evolution equations]{Diffusion phenomena for evolution equations}

\author{Sebastian Kr\'{o}l}

\address{Sebastian Kr{\'o}l, Faculty of Mathematics and Computer Science, Adam Mickiewicz University in Pozna{\'n}, ul. Uniwersytetu Pozna{\'n}skiego 4, 61-614 Pozna{\'n}, Poland}
\email{sebastian.krol@amu.edu.pl}

\author{Yuri Latushkin}

\address{Yuri Latushkin, University of Missouri Department of Mathematics 810 East Rollins Street Columbia, MO 65211, USA}
\email{latushkiny@missouri.edu}

\author{Yuri Tomilov}

\address{Yuri Tomilov, Institute of Mathematics, Polish Academy of Sciences, Śniadeckich 8, 00-
656 Warsaw, Poland, and 
Faculty of Mathematics and Computer Science, Nicolaus Copernicus University, Chopina 12/18, 87-100 Toru{\'n}, Poland}
\email{ytomilov@impan.pl}

\thanks{
The first and the third authors were supported by the NAWA/NSF grant BPN/NSF/2023/1/00001. They were also partially supported by the NCN grant Weave Unisono, 2024/06/Y/ST1/00044.  
The second author was supported by the NSF grant DMS-2106157, and would like to thank the Courant Institute of Mathematical Sciences at NYU and especially Prof. Lai-Sang Young for their hospitality.
}

\keywords{hyperbolic equation, parabolic equation, damped wave equation, heat equation, duffusion, resolvent bounds, energy decay, holomorphic semigroups}

\subjclass{35B40, 47B44, 35B27, 47A10}

\begin{document}

\maketitle

\begin{abstract}
We offer a novel perspective on the classical problem of shadowing a solution $u$ of the abstract hyperbolic equation $u''(t)+Bu'(t)+Au(t)=0$, $t\ge 0$, on a Banach space, by a solution $v$ of its possibly degenerate parabolic counterpart $Bv'(t)+Av(t)=0$. The operator $B$ is assumed to be bounded and the damped wave equation is assumed to be well-posed; no further assumptions are imposed on $A$ or $B$ being 
symmetric, injective, or commuting. The shadowing we study is understood in a very robust sense and exhibits the \emph{diffusion phenomenon}: the norms of the solutions $u$ and $v$ decay at infinity more slowly than the norm of their difference, so that $\|u(t)-v(t)\|=o(\|v(t)\|)$ and hence $\|u(t)\|/\|v(t)\|\to1$ as $t\to\infty$. 

Our first core result concerns the Fourier-analytic properties of two resolvent expressions associated with the quadratic and linear operator pencils of (HE) and (PE), for suitably matched initial data. It establishes decay estimates for the inverse Fourier transform of their difference and, under appropriate well-posedness assumptions, translates these estimates into decay bounds for the difference of the corresponding hyperbolic and parabolic mild solutions.

This leads to our next major result, a fairly complete description of the initial data for which the traditional upper estimates of $\|u(t)-v(t)\|$ of orders
$t^{-1}$ and $t^{-2}$ hold, for the solutions and their time derivatives separately, in the quite general setting of the paper. 
Next, given an arbitrary slowly decaying non-negative function $r(t)\to0$, we prove the much less expected \emph{lower} estimates
$\|u(t)\|\gtrsim r(t)$ and $\|v(t)\|\gtrsim r(t)$ for a large set of initial data, which is described explicitly in the paper. In turn, this yields the existence of genuine diffusive profiles for data satisfying certain high- and low-frequency resolvent estimates, together with a new condition on the position of the ranges of $A$ and $B$. This is done through a careful spectral analysis of a solution semigroup associated with the parabolic equation, holomorphic on a sector but not necessarily strongly continuous, together with a thorough discussion of the well-posedness of the hyperbolic equation, which may be of independent interest.
 
For the classical wave equation on $L^2(\Omega)$, $\Omega\subseteq\mathbb{R}^n$, with \(A\) the Dirichlet Laplacian and \(B\) the
multiplication operator induced by a bounded non-negative function that may degenerate on large sets, we check the direct applicability of the abstract results above whenever possible, and outline a clearly defined circle of open problems in the area. We establish effective criteria
for the required resolvent and range conditions, verify the abstract
hypotheses on the whole space and on exterior domains, and construct explicit
diffusive data by means of Weyl sequences for the Laplacian.
\end{abstract}

{\small
\tableofcontents
}

\renewcommand{\subjclassname}{\textup{2020} Mathematics Subject Classification}

\section{Introduction}

A recurrent phenomenon in dissipative wave propagation is that damping suppresses the genuinely hyperbolic component of the motion before it suppresses the low-frequency, mass-carrying part.  In the simplest constant-damping model on $L^2(\mathbb R^n$), 
\begin{equation}\tag{DW$_0$}\label{eq:DW0}
\begin{aligned}
 &u_{tt}+u_t-\Delta u=0, \qquad \quad t\ge0, \\
   & u(0) = x_0 \in W^{1,2}(\mathbb R^n),\quad u'(0) = x_1 \in L^2(\mathbb R^n),
 \end{aligned}
\end{equation}
solutions behave for large times like solutions of the heat equation
\begin{equation}\tag{P$_0$}\label{eq:P0}
 v_t-\Delta v=0, \qquad t\ge0,  \quad \, v(0) = x_0+x_1.
\end{equation} 
In this form the phenomenon is classical. What becomes much less classical is the following. Suppose that the damping is spatially distributed and may vanish on large sets, that the damping operator does not commute with the conservative part, that explicit formulas are unavailable, and that the first-order comparison dynamics does not naturally evolve on the whole ambient space. How should one formulate the correct parabolic profile, and how should one prove that it really gives the leading term of the hyperbolic asymptotics rather than merely their convenient comparison?

The present paper addresses these questions. We study the abstract damped wave equation
\[
u''(t)+ Bu'(t) + Au(t) =0,\,  t\ge0, \quad \, u(0) = x_0,\, u'(0) = x_1, \tag{HE}
\]
and the associated first-order problem
\[
B v'(t) + Av(t) =0,\, t\ge0,  \quad \, v(0) = y_0, \tag{PE}
\]
on a Banach space $X$, where $B\in \mathcal L(X)$ is bounded and $A$ is a linear operator. The guiding model is the damped wave equation 
\begin{equation*}\label{eq:DWEintro}
 u_{tt}(t,\xi)+b(\xi)u_t(t,\xi)-\Delta u(t,\xi)=0,
\end{equation*}
where $A=-\Delta$ and $B=M_b$ is multiplication by a non-negative function $b$. The  main interest is precisely when $B$ is not injective and does not commute with $A$.
In this setting, methods based on a common spectral representation
are no longer available, and the associated parabolic
problem need not generate a strongly continuous semigroup
on $X$. The present paper develops an abstract theory accommodating both difficulties without assuming selfadjointness, commutativity, or injectivity of the damping operator.

The object of the paper is the following diffusion phenomenon in the strong asymptotic sense: Given a solution $u$ of (HE), one seeks a solution $v$ of (PE) such that
\begin{equation}\label{eq:diffintro}
 \|u(t)-v(t)\|=o(\|v(t)\|), \qquad t\to\infty.
\end{equation}
When \eqref{eq:diffintro} holds, the parabolic orbit $v$ is not merely close to $u$, it is truly the first non-vanishing term of its large-time asymptotics. In particular,
\begin{equation}\label{eq:diffintro_equi}
 \Bigl|\frac{\|u(t)\|}{\|v(t)\|}-1\Bigr|\le \frac{\|u(t)-v(t)\|}{\|v(t)\|}\to 0, \qquad t \to \infty,
\end{equation}
so the decay rate and asymptotic size of $u$ are governed by $v$.

This distinction is crucial. A great many estimates in the literature control the norm of $u-v$ without proving that $v$ is asymptotically visible. Such estimates are often strong and deep, but by themselves they do not identify the leading term of the motion. In the present paper this issue is not treated as a secondary refinement. It becomes a structural part of the analysis. We first derive decay estimates for the inverse Fourier transform of a multiplier and, under suitable assumptions, identify this transform with the difference of actual solutions. We next construct the reduced semigroup carrying the parabolic profile and investigate whether its orbits decay slowly enough to dominate that difference.

The separation of tasks reflects a genuine gap between the classical reference points and the setting considered here. In the case of commuting selfadjoint $A$ and $B$ on a Hilbert space one can often write down the associated parabolic evolution by functional calculus, compare the two dynamics spectrally
in a direct manner, and read off the asymptotics. Once commutativity is lost and the damping becomes degenerate, none of these simplifications is available. Instead, one must first determine the appropriate space for the parabolic dynamics, then prove its well-posedness,
then estimate the difference between the hyperbolic and parabolic solutions, and finally show that the resulting parabolic orbit is large enough to dominate this difference. The present work proposes a unified framework for these steps.

A further distinctive feature is that we obtain separate estimates for solutions and their time derivatives, rather than controlling only their combined energy. This matters because, when \(A\) is non-negative selfadjoint, an energy bound involving \(\|A^{1/2}u\|\) and \(\|u'\|\) need not control \(\|u(t)\|\). Our theory combines separate hyperbolic–parabolic comparison estimates with lower bounds for individual solutions, thus establishing genuinely dominant diffusive profiles.

\subsection{Literature and related works}

The literature on diffusion phenomena for damped waves is broad, but only a few lines are directly relevant to the present paper. For our purposes, the main distinction is not between linear and nonlinear equations, or between bounded and unbounded domains, but between results that identify an actual parabolic asymptotic profile and results that only estimate the difference between hyperbolic and parabolic solutions. 
The latter are often strong and important, but they do not by themselves imply the strong asymptotic relation \eqref{eq:diffintro}, which is the notion of diffusion phenomenon adopted in the present paper.

Hsiao and Liu \cite{HsLiu92,HsLiu93} showed that, after suitable scaling, solutions of certain damped hyperbolic systems converged to nonlinear diffusion waves. Yang and Milani \cite{YangMilani} subsequently compared solutions $u$ of the quasilinear equation
$u_{tt}+u_t-\operatorname{div}(a(\nabla u)\nabla u)=0$ with solutions $v$ of the corresponding nonlinear parabolic equation
$v_t-\operatorname{div}(a(\nabla v)\nabla v)=0$, with $v(0)=u(0)+u_t(0)$. Assuming $a(x)=1+O(\|x\|^\alpha)$ near zero and suitably small, smooth initial data, they proved
\[
\|u(\cdot,t)-v(\cdot,t)\|_{L^\infty(\mathbb R^n)} = O\bigl(t^{-(n+1)/2}\bigr).
\]
This rate is faster than the $O(t^{-n/2})$ upper bounds obtained for the individual solutions. Although Yang and Milani described this as the difference decaying faster than either solution, those upper bounds do not establish \eqref{eq:diffintro}; without a lower bound for $v$, the latter can fail even for $a\equiv1$.

For the constant-damping model \eqref{eq:DW0} on $\mathbb R^n$, Karch \cite{Karch00} established an explicit Gaussian asymptotic profile. For initial data $u(0)=u_0$ and $u_t(0)=u_1$ satisfying the assumptions of \cite{Karch00}, let
$G(t,x)=(4\pi t)^{-n/2}e^{-\|x\|^2/(4t)}$ be the Gauss--Weierstrass kernel and put
$M_0=\int_{\mathbb R^n}(u_0+u_1)\,dx$. The result gives
\begin{equation}\label{eq:karch-intro}
 \|u(t)-M_0G(t)\|_{L^2}=o\bigl(t^{-n/4}\bigr) \quad \textrm{as } t\to \infty.
\end{equation}
If $M_0\ne0$, then $\|M_0G(t)\|_{L^2}\asymp t^{-n/4}$. Thus $M_0G$, which solves \eqref{eq:P0} for $t>0$, is the leading asymptotic term. For the semilinear source $-|u|^{\rho-1}u$ also treated in \cite{Karch00}, the Gaussian coefficient includes the space-time integral of the source.

Hyperbolic--parabolic comparison also arises in singular-perturbation problems. Gallay and Raugel \cite{GallayRaugel} 
obtained, for sufficiently small initial data, a second-order large-time
expansion for a nonlinear variable-coefficient equation of the form
\[
  \varepsilon u_{tt}+u_t=(a(x)u_x)_x+N(u,u_x,u_t),
\]
where $\varepsilon>0$ need not be small, and showed that the expansion is
governed by a linear parabolic equation. In an abstract setting, Chill and
Haraux \cite{ChillHarauxSingular} proved sharp estimates for the
singular perturbation limit as $\varepsilon\to0$.

Further work on \eqref{eq:DW0} concerned exterior domains, decompositions of solutions, and estimates for the difference from the heat evolution. Building on the work of Ikehata and Matsuyama \cite{IkMa2002} on $L^2$-bounds and decay for heat and wave equations in exterior domains, Ikehata \cite{Ikehata2002} obtained comparison estimates for the Dirichlet problem in this setting. Nishihara \cite{Nish03} subsequently decomposed 
solutions of \eqref{eq:DW0} on $\mathbb R^3$ into an exponentially decaying wave part and a diffusive part approximated by the heat semigroup. This decomposition did not imply that the heat solution was the leading term for every choice of data. Indeed, \eqref{eq:diffintro} fails for suitable data whose Fourier transforms are supported away from the origin. In this setting, the matched solution of \eqref{eq:P0} is $v(t)=e^{t\Delta}(u_0+u_1)$. For instance, for $u_0\in H^1(\mathbb R^n)$ and $u_1\in L^2(\mathbb R^n)$ one has
\begin{equation}\label{eq:const-damp-intro}
 \|u(t)-e^{t\Delta}(u_0+u_1)\|_{L^2(\mathbb R^n)}
 \lesssim t^{-1}\bigl(\|u_0\|_{H^1}+\|u_1\|_{L^2(\mathbb R^n)}\bigr),\qquad t\ge1.
\end{equation}
These arguments used explicit kernels, Fourier analysis, or other properties of constant damping. They do not directly extend to degenerate spatial damping, for which the parabolic evolution may have to be defined on a smaller space.

Chill and Haraux \cite{Chill03} extended the constant-damping comparison to (HE) and (PE) with $B=I$ and $A$ an arbitrary non-negative selfadjoint operator on a Hilbert space.  Their uniform $O(t^{-1})$ estimate for the difference of solutions, in the norm of $D(A^{1/2})$, sharpened earlier results of Ikehata and Nishihara \cite{IkNi2003}. When $0$ is a nonisolated point of the spectrum of $A$, the rate cannot be improved to $o(t^{-1})$ uniformly over normalized initial data. Nevertheless, for each fixed pair of initial data $(x_0,x_1)$, the difference is $o(t^{-1})$. Thus the sharpness of the remainder estimate does not provide a lower bound for any fixed parabolic orbit, as required in the analysis of \eqref{eq:diffintro}.

A different historical perspective on the comparison between
hyperbolic and parabolic dynamics originates in the study of the telegraph equation and its probabilistic interpretation pioneered by S.~Goldstein and
M.~Kac \cite{SGoldstein,Kac}. This viewpoint was developed in an
abstract selfadjoint setting in
\cite{EcksteinGoldsteinLeggas,ClarkeEcksteinGoldstein,
FragnelliGoldsteinGoldsteinRomanelli}, where the solutions of
(HE) and (PE) are considered with
$A=S^2$ and $B=2F(S)$, for a non-negative selfadjoint operator $S$ on a Hilbert space. In the earlier works $F$ was the positive constant function $F\equiv a$, while \cite{FragnelliGoldsteinGoldsteinRomanelli} allowed a suitable positive function $F$, possibly unbounded. The initial value of $v$ was prescribed by an explicit spectral formula in terms of the initial data for $u$ such that \eqref{eq:diffintro} holds.
This orbitwise formulation is closely related to the present paper, although the common
spectral representation of $S$ and $F(S)$ used in those works excludes
noncommuting damping. On the probabilistic side, Mascia showed in
\cite{Mascia} that a multidimensional generalization of the
Goldstein--Kac model has parabolic large-time behaviour.

A different PDE line of research, important for variable damping, is represented by the work of Yamazaki and later by Sobajima and Wakasugi. Yamazaki \cite{Ya07} extended the abstract hyperbolic--parabolic comparison to time-dependent scalar dissipation, replacing the constant damping coefficient by a positive function $b(t)$ which decays in time. Complementary results on effective time-dependent dissipation were obtained by Wirth \cite{Wirth07}, who also established quantitative parabolic approximation results for partially dissipative hyperbolic systems satisfying a uniform Kalman rank condition \cite{Wirth14}. The spatially dependent case calls for different techniques. Sobajima and Wakasugi \cite{SoWa16} studied effective damping in exterior domains, comparing the damped wave equation with the weighted parabolic problem $a(x)v_t-\Delta v=0$ by means of heat-semigroup and weighted energy estimates. For $a(x)=\|x\|^{-\alpha}$, $0\le\alpha<1$, they also established a lower bound for the parabolic solution corresponding to suitable nonzero compactly supported data. Together with their remainder estimate, this gave a leading parabolic profile in a weighted $L^2$ norm. Their later work \cite{SoWa18} treated, in particular, damping growing at infinity, while Sobajima \cite{So21} obtained higher-order asymptotic expansions. A related spectral treatment for non-negative selfadjoint operators on general measure spaces, based on decay properties of the corresponding heat semigroup, can be found in \cite{IkedaTaniguchiWakasugi}. Diffusion phenomena for equations with space- and time-dependent coefficients were studied by Taylor \cite{Taylor2018}, and subsequently by Taylor and Todorova \cite{TaylorTodorova} for time-dependent selfadjoint operators associated with Dirichlet forms on metric measure spaces. These works provided comparison estimates and, in certain concrete models, leading parabolic profiles. They did not address general degenerate, noncommuting operator pairs $(A,B)$.

Among the abstract predecessors were the Hilbert-space works of Radu, Todorova, and Yordanov \cite{RaToYo11,RaToYo16},
Ikehata, Todorova, and Yordanov \cite{IkToYo2013} and of Nishiyama \cite{Nishiyama2016}.
Radu, Todorova, and Yordanov established diffusion-type estimates and higher-order approximations for abstract damped wave equations, including equations with noncommuting selfadjoint coefficients. In a different setting, Ikehata, Todorova, and Yordanov studied strongly damped wave equations on Hilbert spaces and described the large-time behavior in terms of parabolic and wave evolutions. Nishiyama considered (HE) with noncommuting selfadjoint operators $A$ and $B$, but assumed $B\ge cI>0$. His weighted comparison estimates therefore did not cover degenerate damping. For strongly damped wave equations on Hilbert spaces, Ikehata, Todorova, and Yordanov [new reference] obtained a generalized diffusion asymptotic involving both parabolic and wave components.

Also, there is a substantial literature on damping, stabilization, and low-frequency asymptotics in unbounded geometries which is directly relevant to the application side of the present paper. This includes the work of Burq and Joly \cite{BuJo16} on decay and stabilization in wave guides and unbounded domains, the paper \cite{AlIbKh15} by Aloui, Ibrahim, and Khenissi  on exterior domain energy decay, the research of Joly and Royer in \cite{JoRo18} on low-frequency asymptotics in asymptotically periodic media, 
and the recent preprints \cite{Aranal,Fahs} on unbounded damping.
 These papers are not profile theorems in the sense of \eqref{eq:diffintro}, but they provide a highly relevant background for the PDE part of the present work: they clarify the role of geometric control, the distinction between low and high frequencies, and the extent to which one can expect effective parabolic behavior in unbounded media. For additional closely related in spirit material, we refer to \cite{FrSiTr18} and \cite{MaRo18}.

Finally, it is worth noting that, independently of diffusion asymptotics, abstract damped wave equations
have long been studied through their energy spaces, phase space generators,
operator matrices, quadratic pencils, analyticity, and energy decay; see, for
example,
\cite{AnLa14,BanksItoWang, ChillPa, EngelDissipative,HrynivShkalikov,GesztesyGoldsteinHoldenTeschl, KlWa}. 
In particular, \cite{AnLa14} develops an abstract
resolvent framework for polynomial stability of damped
wave equations, relates energy decay to observability
of the associated Schrödinger group, and applies
these results to the damped wave equation on the torus.
More recent contributions to asymptotics for abstract second-order
equations includes \cite{ChillPa},
where observability conditions yield nonuniform energy
decay for wave, fractional Klein--Gordon, and beam
equations, and \cite{KlWa}, which
uses control-theoretic estimates to
derive decay rates depending explicitly on the degree
of unboundedness of the damping.
The comprehensive study in \cite{GesztesyGoldsteinHoldenTeschl} relates the phase space generator to a Dirac-type operator and to the associated quadratic operator pencil, and applies this, in particular, to the study of equipartition of energy.
We also note that form-based well-posedness results
for damped second-order equations on Hilbert triples
were obtained by Banks, Ito, and Wang \cite{BanksItoWang}; see also \cite{BanksRemarks}.  Their framework allows unbounded damping and input operators and relates variational solutions to semigroup solutions on the associated energy space.

Taken together, these works leave a rather clear gap. Genuine profile results are known in some important concrete models, but much of the literature still uses the phrase ``diffusion phenomenon'' for estimates of the  remainder alone. Abstract operator-theoretic treatments exist, but they are usually tied to selfadjointness, commutation, invertibility of the damping, or energy-space formulations in which the norms of $u-v$ and $(u-v)'$ are not separated. On the other hand, PDE works that treat more realistic damping geometries rely on methods too specific to produce a general abstract approach. The aim of the present paper is to bridge the gap: to isolate an operator-theoretic framework that can handle degenerate noncommuting damping, still distinguish genuine profiles from mere remainders, and nevertheless lead to a substantial concrete theory for the damped wave equation with multiplier damping.

\subsection{Overview of the main abstract results}
Our primary goal is to control the size of the difference between solutions of (HE) and (PE) in terms of the operators \(A\), \(B\), and the initial data. Formally taking Laplace transforms $\mathcal L$ in (HE) and (PE), subtracting the resulting expressions, and using the resolvent identity, we obtain
\[
\begin{aligned}
(\mathcal Lu)(\lambda)-(\mathcal Lv)(\lambda)
={}&(\lambda^2+\lambda B+A)^{-1}(Bx_0+x_1-By_0)\\
&+\lambda(\lambda^2+\lambda B+A)^{-1}
\bigl(x_0-\lambda(\lambda B+A)^{-1}By_0\bigr).
\end{aligned}
\]
The resolvent assumptions used below control
\(
is(-s^2+isB+A)^{-1},
\)
$ s\in \mathbb R\setminus\{0\},$
as \(s\to0\), but need not control \((-s^2+isB+A)^{-1}\) itself. We therefore match the initial data by setting \(Bx_0+x_1=By_0\), which eliminates the first term. For \(B=I\), this reduces to the familiar choice \(y_0=x_0+x_1\) in (P0). Evaluating the remaining resolvent expression at \(\lambda=is\), we arrive at
\begin{equation}\label{eq:mintro}
 m(s)= (-s^2 +is B +A)^{-1} (is x_0+By_0)-(is B + A)^{-1} By_0,
 \qquad s\in\mathbb R\setminus\{0\}.
\end{equation}
The function \(m\) compares the second-order and first-order dynamics directly on the frequency side, before the existence of either propagator has been assumed.

Here and in the sequel, for $\lambda\in\mathbb C$, whenever the notation
$(\lambda^2+\lambda B+A)^{-1}$ or $(\lambda B+A)^{-1}$ is used,  we
implicitly assume that the corresponding operator is boundedly
invertible. The basic resolvent hypotheses take the
form
\begin{equation}\label{eq:theta-resolvent-intro}
\sup_{s\ne0}|s|\bigl\|(-s^2+is B +A)^{-1}\bigr\|<\infty,
 \quad
\sup_{s\ne0} \min (|s|,|s|^\theta)\bigl\|(isB+A)^{-1}B\bigr\|<\infty,
\end{equation}
with a parameter $\theta\in[0,1]$.  The parameter $\theta$ measures the high-frequency loss allowed in the parabolic resolvent estimate. 
At high frequencies it permits the parabolic resolvent to decay only as $|s|^{-\theta}$. When $\theta<1$, this loss allows the associated parabolic solution family to be singular at the initial time.  Thus the case $\theta=1$ is the regular endpoint of a genuinely larger scale of degenerate problems.

The high-frequency regularity imposed on the data is encoded by the spaces
\begin{align*}
 X_2 & :=\Bigl\{x\in X:\sup_{|s|\ge1} s^2 \bigr\| (-s^2+is B +A)^{-1}x \bigr\|<\infty\Bigr\},\\
 X_D^\eta &:=\Bigl\{x\in X:\sup_{|s|\ge1}
 |s|^\eta\bigl\|(is B +A)^{-1} Bx\bigr\|<\infty\Bigr\},
 \qquad \eta\le1.
\end{align*}
Depending on $\theta$, Theorem~\ref{thm 1} involves three additional paces $Y_D^\theta$, $Z_D^\theta$, and
$\widetilde Z_D^\theta$, defined in the next subsection
in terms of the resolvent decay spaces $X_D^\eta$ and
the operators $A$ and $B$.

Theorem \ref{thm 1} is the central analytic result. Its proof establishes Mihlin-type estimates for $m$  from the resolvent hypotheses and uses a vector-valued estimate from Fourier multiplier theory to show that $\cF^{-1}m$, initially a distribution, is represented away from zero by a function with the decay stated below. For $x_0\in\Dom(A)\cap B^{-1}(X_2)$ and $y_0\in Y_D^\theta\cap B^{-1}(X_2)$, this function satisfies
\begin{equation}\label{eq:multiplier_intro}
\|\mathcal F^{-1}m(t)\|\lesssim |t|^{-1},
\qquad t\ne0.
\end{equation}
If, in addition,
\(
x_0\in\Ran(A)\)
and
\(
y_0\in B^{-1}(\Ran(A))\cap Z_D^\theta,
\)
then the estimate improves to 
$ \|\mathcal F^{-1}m(t)\|\lesssim |t|^{-2}$.
Under the additional regularity conditions of Theorem 1.1(iii),
the corresponding estimate
$\|(\mathcal F^{-1}m)'(t)\|\lesssim |t|^{-2}$
holds for the derivative. These conclusions are orbitwise Banach-space statements: they do not require commutativity,
selfadjointness, or the prior existence of global solution operators.  The
orders $t^{-1}$ and $t^{-2}$ are those reached by the multiplier argument
developed here; the obstacles to extending the same method to higher orders
are examined separately in Section~\ref{sect3}.

Theorem~\ref{thm:ind-solu} supplies the bridge from frequency-side estimates to actual
solutions.  If mild solutions of the two problems exist, have at most
polynomial growth, and the initial data are matched by
\begin{equation}\label{matching}
 x_1=By_0-Bx_0,
\end{equation}
then
\[
 \mathcal F^{-1}m(t)=u(t;x_0,x_1)-v(t;y_0),
 \quad t>0.
\]
Thus solvability and decay are deliberately separated: 
Theorem~\ref{thm 1} provides estimates for $\cF^{-1}m,$ while Theorem~\ref{thm:ind-solu} identifies this
function with the difference of the corresponding hyperbolic and parabolic orbits whenever
those orbits are available.
The construction of the parabolic dynamics is the second main component of
the abstract theory.  Under a sectorial version of the second estimate in
\eqref{eq:theta-resolvent-intro}, in which $i\mathbb R\setminus\{0\}$ is
replaced by a sector $\Sigma_{\omega+\pi/2} = \{\la \in \C\setminus \{0\}: |\arg \la|<\omega+\pi/2\}$ ($\omega\in (0,\pi/2)$), for every $z\in \Sigma_\omega$, $\theta \in (0,1]$, and $\beta \in (1-\theta, 1)$, we set
\begin{equation*}\label{eq:S-F-intro}
 \cT(z):=\frac{1}{2\pi i}\int_\Gamma
 e^{\lambda z}(\lambda B+A)^{-1}B\,d\lambda,
 \quad
 J_\beta = \frac{\sin(\pi\beta)}{\pi} \int_0^\infty t^{-\beta} ( (1+ t)B + A)^{-1}B\, \ud t,
\end{equation*} where $\Gamma$ is an appropriate sectorial contour; see Subsection~\ref{subs:main_results} for its precise definition.
Theorem~\ref{thm:mild_DPE} proves that $\cT$ is a holomorphic semigroup on $X$ and that
\begin{equation}\label{eq:S-singular-intro}
 \|\cT(z)\|=O(|z|^{\theta-1})
 \quad\text{as }z\to0
\end{equation}
inside smaller sectors, while $\cT$ is uniformly bounded for large $|z|$ there.
Moreover,
\[
(\lambda B + A)^{-1}B=\int_0^\infty e^{-\lambda t}\cT(t)\,dt,
\qquad \operatorname{Re}\lambda>0.
\]
When $\theta<1$, the bound \eqref{eq:S-singular-intro} allows singular
behavior at the origin, and $\cT$ need not be a $C_0$-semigroup on the whole of
$X$.  The theorem therefore
identifies two natural classes of admissible initial values.  If
\[
 y_0\in A^{-1}(BX),
\]
or if $y_0=J_\beta x$ for some $x\in X$, then $\cT(t)y_0\to y_0$ as
the initial time is approached, and $v(t)=\cT(t)y_0$ is a mild solution of the
degenerate parabolic problem with a bounded holomorphic extension to a sector.
The operator $J_\beta$ maps $X$ into $X_D^1$ and produces initial values at
which $\cT(t)$ is strongly continuous at $t=0$.
It therefore supplies, alongside $A^{-1}(BX)$, a second natural class of
admissible initial values in the possibly singular regime $\theta<1$.  At the
endpoint $\theta=1$, the reduced $C_0$-semigroup description below becomes
available.

The same theorem also provides the slowly decaying parabolic orbits needed
for a genuine diffusion phenomenon. We define recursively $D_1:=A^{-1}(BX)$, $D_\ell=A^{-1}(BD_{\ell-1})$, $\ell\ge2$. If
\begin{equation}\label{eq:range-condition-intro}
 \operatorname{Ran}(B)\setminus\operatorname{Ran}(A)\neq\emptyset, \quad (\textrm{equivalently, \quad $\Ran(B)\not\subset \Ran(A)$),} 
\end{equation}
then the spectral radius of every $\cT(t)$, $t>0$, is at least one.  Consequently,
for every non-increasing function $r$ on $\R_+$ so that $r(t)\downarrow0$ there is a vector $y_0\in\bigcap_{\ell\ge1}D_\ell$ whose orbit satisfies $\|\cT(t)y_0\|\gtrsim r(t)$ for all $t>0$.  This argument
does not require an infinitesimal generator for $\cT$ and is therefore available
also when $\theta<1$.

At $\theta=1$ the structure becomes more classical and more precise.  On a
reflexive space one obtains the decomposition
\[
 X=Y\oplus\operatorname{Ker}(B) \quad \textrm{and} \quad  \cT(t)=T(t)\oplus0,\quad t>0,
 \]
where $Y:=\overline{A^{-1}(BX)}^{\,X}$ and $T$ is a sectorially bounded holomorphic $C_0$-semigroup on $Y$.  Its
generator is the operator
\[
 L=-(B|_Y)^{-1}\circ A
\]
on the maximal domain
\[
 \operatorname{Dom}(L)=\{x\in\operatorname{Dom}(A)\cap Y:Ax\in BY\}.
\]
This endpoint theory gives a sharper spectral description of slowly
decaying orbits and is used in Theorem~\ref{thm:diffusion_theta=1} for the refined $t^{-2}$ and
derivative statements.  It should, however, be viewed as the regular endpoint
of the general $\theta$-theory rather than as the framework for the whole
paper.

For the second-order equation we impose well-posedness only at the stage where
actual hyperbolic solutions are needed.  If $V\times X$ is its phase space, the
auxiliary results of Section~4 show that $V\subset X_2$.  When $X$ is a UMD
space, one has $V=\operatorname{Dom}((A+\omega)^{1/2})$ for every $\omega\ge0$ for which $A+\omega$ is sectorial.  Assuming in addition
that the phase space semigroup is polynomially bounded, Theorem~\ref{thm: diffusion intro} combines
this fact with Theorems~\ref{thm 1}--\ref{thm:mild_DPE}.  For matched data \eqref{matching} it gives
\begin{equation}\label{eq:solution-rates-intro}
 \|u(t)-v(t)\|\lesssim t^{-1}.
\end{equation}
Under the additional range assumptions $x_0,By_0\in\operatorname{Ran}(A)$ the
rate improves to
\[
 \|u(t)-v(t)\|\lesssim t^{-2},
\]
and, under the kinetic regularity assumptions, we obtain
\[
 \|u'(t)-v'(t)\|\lesssim t^{-2}.
\]
For the first estimate, the parabolic data $y_0\in Y_D^\theta$ may be
chosen either from $A^{-1}(BX)$ or from $\Ran(J_\beta)$, where
$\beta\in(1-\theta,1)$, with the further conditions imposed for the
$t^{-2}$ and derivative estimates. At $\theta=1$, all data in
$\overline{A^{-1}(BX)}$ are available.

The estimates in \eqref{eq:solution-rates-intro} still do not by themselves
produce the diffusion phenomenon in the sense of \eqref{eq:diffintro}.
The parabolic orbit $v$ must decay more slowly than the difference $u-v$.
Theorem~\ref{thm: diffusion} combines Theorem~\ref{thm: diffusion intro} with the slow orbit part of Theorem~\ref{thm:mild_DPE}.  Under
\eqref{eq:range-condition-intro} and the compatibility condition
\[
 B(\operatorname{Dom}(A))\subset
 \operatorname{Dom}((A+\omega)^{1/2}),
\]
one can choose a single datum
\[
 y_0\in D_2=A^{-1}\bigl(BA^{-1}(BX)\bigr)
\]
which works for every $x_0\in\operatorname{Dom}(A)$ after setting $x_1=By_0-Bx_0$. Given any bounded non-increasing function $r$ with
\[
 r(t)\to0,
 \qquad
 t^{-1}=o(r(t)),
\]
the data can be chosen so that, as $t\to\infty$,
\[
 \|u(t)-v(t)\|=O(t^{-1})=o(r(t)),
 \qquad
 \|v(t)\|\gtrsim r(t),
 \qquad
 \|u(t)\|\gtrsim r(t),
\]
and therefore \eqref{eq:diffintro_equi} holds. Theorem~\ref{thm:diffusion_theta=1} exploits the additional endpoint structure at $\theta=1$.  In a
Hilbert space it gives a weak scalar lower bound of prescribed slow order; under
the finer spectral condition $0\in\sigma(L)\setminus\sigma_p(L)$ it combines a
$t^{-2}$ remainder with a limsup lower estimate for $u$, and it gives an
analogous statement for the derivatives $u'$ and $v'$.
We also obtain several results which may be of independent interest. In
particular, Theorem~\ref{thm: DWE well-posed} and
Corollary~\ref{cor:DWE_well_posed} show that bounded damping does not affect
the well-posedness criterion or the phase space of the abstract second-order
equation. Namely, (HE) is well posed for some, and hence for every,
$B\in\mathcal L(X)$ if and only if $-A$ generates a cosine function on $X$,
and in this case its phase space coincides, with equivalent norms, with that
of the undamped equation.

The same phase space analysis also leads to an energy decay statement for
abstract damped waves. Proposition~\ref{prop:total_energy} shows that, if
$A$ is an injective non-negative selfadjoint operator, $B$ is bounded and
non-negative, and
\[
 i\mathbb R\setminus\{0\}\subset\rho(P),\qquad
 \sup_{s\neq0}\min\bigl(|s|^{1+\beta},|s|\bigr)\|R(is)\|<\infty
\]
for some $\beta\geq0$, then, for
\(
 y\in\Dom(A^{1/2})\) and  \(x\in\Ran(A)-By,\)
the corresponding solution of (HE) satisfies
\[
 \|A^{1/2}u(t)\|+\|u'(t)\|
 =O_{y,x}\bigl(t^{-1/(2+\beta)}\bigr),
 \qquad t\to\infty.
\]
Finally, on the semigroup side, Proposition~\ref{cor:fixed-scale}
characterizes polynomial decay of an individual orbit of the holomorphic
solution family $\cT$ associated with (PE):
\[
 \|\cT(t)y\|=O(t^{-\eta}),\qquad 0<\eta<1,
\]
as $t\to\infty$, if and only if the corresponding local resolvent
\(
 \lambda\longmapsto(\lambda B+A)^{-1}By
\)
has the prescribed behaviour near zero. At the endpoint $\theta=1$, where
$\cT(t)=T(t)\oplus0$ and $T$ is a $C_0$-semigroup on the reduced space,
this gives in particular an orbitwise resolvent criterion for the
$C_0$-semigroup $T$.

The results reviewed above separate several issues that are usually
intertwined. Theorem \ref{thm 1} gives estimates for
$\cF^{-1}m$ without a generation hypothesis.
 Under sectorial parabolic resolvent estimates, Theorem \ref{thm:mild_DPE} constructs the holomorphic solution family \(\cT\), which may be singular at \(t=0\) when \(0<\theta<1\). After imposing
well-posedness and choosing admissible parabolic initial values, these
operator estimates become estimates for the difference $u-v$ of actual
solutions. Finally, the existence of sufficiently slowly decaying parabolic
orbits makes this difference negligible relative to $v$ and hence yields a
genuine asymptotic profile. The classical commuting selfadjoint case is
contained in this scheme, but neither commutativity nor selfadjointness is
built into the general multiplier framework.

\subsection{Applications to the damped wave equation with Dirichlet Laplacian}

The abstract theory is complemented by a substantial analysis of the model
problem with the Dirichlet Laplacian,
\begin{equation*}\label{model}
 u_{tt}(t)+b u_t(t)-\Delta u(t)=0,
 \qquad
 b v_t(t)-\Delta v(t)=0,
\end{equation*}
on $L^2(\Omega)$, where $b\in L^\infty_+(\Omega)$ and $\Omega\subset\mathbb R^n$
is an arbitrary open set unless additional geometry is explicitly imposed.
This part is not merely an illustration of the abstract results.  It gives
concrete criteria for the abstract hypotheses, separates the different
hyperbolic, parabolic, and slow orbit requirements, and supplies both broad
classes of coefficients and explicit examples.

For $A=-\Delta_D^\Omega$ and $B=M_b$, the model analysis uses three classes
of damping coefficients,
\[
 \cB(\Omega),\qquad
 \cB_D^\theta(\Omega),\qquad
 \cR(\Omega).
\]
The class $\cB(\Omega)$ encodes the resolvent assumptions for the
hyperbolic pencil $-s^2+isM_b-\Delta_D^\Omega$,
$\cB_D^\theta(\Omega)$ encodes the sectorial resolvent assumptions for the
degenerate parabolic pencil $\lambda M_b-\Delta_D^\Omega$, and
$\cR(\Omega)$ encodes the range condition
\[
 \Ran(M_b)\setminus\Ran(-\Delta_D^\Omega)\neq\emptyset
\]
used, together with the parabolic resolvent assumptions, to obtain slowly
decaying parabolic orbits. Their precise definitions are given below in
Subsection~\ref{subs:applications}; Section~\ref{sec:examples} develops
criteria for membership in these classes.

The hyperbolic problem is well posed for every bounded non-negative $b$.  Its phase space is $W^{1,2}_0(\Omega)\times L^2(\Omega)$, 
and the associated first-order semigroup has at most linear growth.  The truly delicate hyperbolic input is the high-frequency resolvent estimate.  Rather than reproducing the microlocal theory behind it, the paper isolates the class $\cB_\infty(\Omega)$ of non-zero functions $b\in L_+^\infty(\Omega)$ such that $\sup_{|s|\ge1} |s| \bigl\|(-s^2 + is M_b - \Delta)^{-1}\bigr\|<\infty$,
and imports available geometric control results in the whole space and
exterior domain settings.  The low-frequency relation between hyperbolic pencil and parabolic pencil is treated abstractly, so that this external PDE input can be joined to the intrinsic parabolic analysis.

The parabolic class $\cB_D^\theta(\Omega)$ is developed directly.  Theorem~2.2 shows, in particular, that $b\in \cB_{D}^{1/2}(\Omega)$ when there are sets $U'\subset U\subset\Omega$ such that $U$ has the Poincar\'e property, $U'$ admits a cut-off supported in $U$, and $b$ is uniformly positive on $\Omega\setminus U'$.
Thus the damping may vanish on a set, provided this possible zero region is
localized inside a Poincar\'e set.  Besides the sectorial resolvent inclusion, one obtains the characteristic high-frequency bounds
\[
 \sup_{|\lambda|\ge1}\|(\lambda M_b - \Delta)^{-1}\|<\infty,
 \qquad
 \sup_{|\lambda|\ge1}
 |\lambda|^{1/2}\|(\lambda M_b - \Delta)^{-1}M_b\|<\infty
\]
in every fixed proper sector.  Proposition~2.4 then covers the full scale
$\theta\in[1/2,1]$.  Under the assumptions
\begin{equation}\label{eq:gradient-theta-intro}
 b\in W^{1,\infty}(\Omega),
 \qquad
 \|\nabla b(\xi)\|\le c\,b(\xi)^{\,2-1/\theta},
\end{equation}
the $|\lambda|^{-1/2}$ decay improves to $|\lambda|^{-\theta}$.  Thus the
parameter $\theta$ has a concrete PDE meaning: it quantifies how the
differential control of $b$ near its zero set translates into improved
high-frequency decay of the parabolic resolvent.  At the endpoint $\theta=1$,
condition \eqref{eq:gradient-theta-intro} becomes
$\|\nabla b\|\lesssim b$.

The positive results are complemented by the analysis of converse phenomena.  If the damping becomes too small along the supports of suitable Weyl sequences for $0\in \sigma_{\mathrm{ap}}((-\Delta)^{1/2})$, the resolvent set of the parabolic pencil may even be empty. Simple abstract examples also show that the high-frequency behavior of the hyperbolic and parabolic pencils can differ markedly. These facts explain why $\cB(\Omega)$ and $\cB_{D}^\theta(\Omega)$ must be kept separated and why the parabolic estimates cannot simply be inferred from the hyperbolic ones.

The class $\cR(\Omega)$ addresses a different question: whether the
parabolic family has orbits decaying sufficiently slowly to dominate a remainder of order $t^{-1}$. A useful sufficient condition is provided by the Weyl sequence class $\cB_{\cW}(\Omega)$ defined below in Subsection~\ref{subs:applications}. Roughly speaking, it consists of damping coefficients which stay bounded away from zero on the supports of a suitable sequence of pairwise disjoint Weyl functions for $0\in\sigma_{\mathrm{ap}}((-\Delta)^{1/2})$. We prove that
$\cB_{\cW}(\Omega)\subset\cR(\Omega)$. Geometrically, for $b\in \cB_\cW(\Omega)$ the functions $\varphi_j$ in the Weyl sequence are concentrated on disjoint regions of increasingly low Dirichlet energy on which the damping remains effective. 

The analysis of the Dirichlet Laplacian also yields geometric criteria
which may be of separate interest.  
In particular,
Lemma~\ref{lem:steiner_formula}, in conjunction with Propositions \ref{cor:diff_copy} and \ref{prop:correction_b},
shows that diffeomorphic copies in $\Omega$ of domains $G_j$ with uniform Lipschitz character and
$|\partial G_j|/|G_j|\to0$ as $j\to\infty$,
obtained through diffeomorphisms with uniformly bounded
derivatives and Jacobian determinants bounded away from zero,
give rise to Weyl sequences with normalized $L^2$-mass and
vanishing gradient, and hence to explicit damping coefficients
in $\mathcal R(\Omega)$. This links the geometry of large
subsets of $\Omega$ to the range condition $\mathcal R(\Omega)$,
and in turn gives explicit Weyl sequence constructions of
slowly decaying parabolic initial data.

Section~2 develops several ways of constructing such sequences.  In the whole space, the range condition, $b\in \cR(\R^n)$, automatically holds for all nonzero $b\in L^\infty_+(\R^n)$ provided $n\le4$; in more general situations the class $\cB_{\cW}(\Omega)$ gives a flexible sufficient condition for $b\in \cR(\Omega)$.

For uniformly continuous damping on $\mathbb R^n$, (AGCC) supplies the
high-frequency estimate for the damped waves. Combined with the
low-frequency and parabolic criteria mentioned above, this gives the
resolvent assumptions needed for the wave--parabolic comparison and for
the construction of the parabolic family. Under the same assumption,
(AGCC) also implies the range condition needed for slowly decaying
parabolic orbits. In dimensions \(n\leq4\), the latter condition 
automatically holds for every non-zero damping coefficient.

Our approach is sufficiently concrete to produce nontrivial examples,
rather than only abstract membership criteria. In dimension two,
Proposition~\ref{prop:example_b} gives a smooth AGCC coefficient $b$ for which the damped-wave resolvent estimate, the full scale of parabolic resolvent
estimates, and the range condition hold simultaneously, even though its
damping geometry is highly sparse. More precisely, for every $r>1$ and $\varepsilon\in(0,1)$, one may choose a coefficient $b_1$ for which $|B(x,r)\cap\operatorname{supp}b_1|<\varepsilon$ for every $x$, and another coefficient $b_2$ for which $|B(x,r)\setminus\operatorname{supp}b_2|<\varepsilon$ for every $x$.

The culmination of the model analysis is Theorem~\ref{thm:diffusion_for_Delta}. On $\mathbb R^n$, let $b\in C^1(\mathbb R^n)\cap L^\infty_+(\mathbb R^n)$ satisfy the averaged geometric control condition and have bounded gradient.  Then the conclusions
of Theorems~\ref{thm:mild_DPE} and~\ref{thm: diffusion intro} hold with $\theta=1/2$, $X=L^2(\mathbb R^n)$, $V=W^{1,2}(\mathbb R^n)$.
At this value of the parameter the abstract data spaces simplify to
$Y_D^{1/2}=X_D^1$, $Z_D^{1/2}=X$, $\widetilde Z_D^{1/2}=A^{-1}\big(BX\big)$.  
For example, if $u_0\in W^{2,2}(\mathbb R^n)$ and
$v_0\in W^{2,2}(\mathbb R^n)$ satisfies $\Delta v_0=bw_0$ for some
$w_0\in L^2(\mathbb R^n)$, then, with $u_1=bv_0-bu_0$, the corresponding solutions satisfy \eqref{eq:solution-rates-intro}.
 
Moreover, since $b$ is uniformly continuous, (AGCC) implies
$b\in\cR(\mathbb R^n)$; see Corollary~\ref{cor:agcc_implies_r}. Hence
Theorem~\ref{thm: diffusion} also applies and yields genuine diffusive profiles whose norms are bounded below by any prescribed rate $r(t)$ satisfying $t^{-1}=o(r(t))$ as $t\to\infty$.

For a smooth exterior domain, the corresponding statement holds under the
exterior geometric control condition, boundedness of $\nabla b$, and uniform
positivity of $b$ near infinity.  In this case Theorems~\ref{thm:mild_DPE}--\ref{thm: diffusion} apply with  $\theta=1/2$, $X=L^2(\Omega)$, $V=W^{1,2}_0(\Omega)$. 
In both the whole space and exterior domain settings, the geometric control
results provide the high-frequency hyperbolic estimate, while the
low-frequency comparison and the parabolic resolvent estimates are supplied
by the abstract and elliptic arguments developed in this paper.

Finally, Section~\ref{subs:explict} gives a constructive counterpart to the slow orbit
existence theorem.  Starting from a Weyl sequence $(\varphi_j)$ for $0\in \sigma_{\mathrm{ap}}(-\Delta)$ with pairwise disjoint supports on which $b$ stays separated from zero, one forms data of the form
\[
 y_0=\sum_{j\ge1}a_j\varphi_j.
\]
An individual resolvent criterion then gives explicit choices of $y_0$ for which
\[
 \limsup_{t\to\infty}t^\eta\|\cT(t)y_0\|=\infty
\]
for prescribed, or even all, $\eta\in(0,1)$. Combined with the \(O(t^{-1})\) estimate for the difference of solutions, this shows that the parabolic solution \(v\) dominates the difference \(u-v\) along suitable sequences of times. The conclusion is weaker than the full asymptotic relation
\eqref{eq:diffintro}, but the initial data are explicit. The model part of the paper therefore does more than verify the abstract hypotheses: it exhibits the geometric mechanisms, the coefficient classes, and concrete initial data that make the diffusive component visible.
\smallskip

\noindent{\bf Convention.}
 For simplicity of presentation, throughout the paper we assume \(n\ge2\) and \(b\in L^\infty_+(\Omega)\setminus\{0\}\) whenever we consider applications to the Dirichlet Laplacian on \(\Omega\subset\mathbb R^n\).
\smallskip

\noindent{\bf Notation.} Throughout, $X$ denotes a complex Banach space with norm $\|\cdot\|$.
If $X$ is a Hilbert space, $\langle\cdot,\cdot\rangle$ denotes its inner
product. We write $\mathcal L(X,Y)$ for the space of bounded linear
operators from a Banach space $X$ to a Banach space $Y$, and
$\mathcal L(X):=\mathcal L(X,X)$. The continuous dual of $X$ is denoted
by $X^*$, and, for a densely defined operator $A$ on $X$, $A^*$ denotes
its adjoint.

We use $\mathbb N=\{1,2,\ldots\}$ and
$\mathbb N_0=\mathbb N\cup\{0\}$. If $n\in\mathbb N$ and
$E\subset\mathbb R^n$ is Lebesgue measurable, then $|E|$ denotes its $n$-dimensional Lebesgue measure. For a Lipschitz hypersurface $S\subset\mathbb R^n$, the same notation $|S|$ denotes its $(n-1)$-dimensional Hausdorff measure. We write $\mathbf 1_E$ for the
characteristic function of a set $E.$ If $f\in L^1_{\mathrm{loc}}(U)$, where $U\subset\mathbb R^n$ is open, then we write $\operatorname{supp}f$ for the support in $U$ of the distribution induced by $f$.  For continuous functions, this clearly  agrees with the usual support. For $\xi \in \mathbb R^n$ and $r>0$, $B(\xi;r)$ denotes the open Euclidean ball of radius $r$ centred at $x$. For non-empty sets $E,F\subset\mathbb R^n$, we write
\[
\operatorname{dist}(\xi,E):=\inf_{\zeta\in E}\|\xi- \zeta\|, \qquad \operatorname{dist}(E,F):=
\inf\{\|\xi - \zeta\|:\xi \in E,\ \zeta\in F\}.
\]

For an open set $\Omega\subset\mathbb R^n$,
$L^\infty_+(\Omega)$ denotes the cone of all non-negative essentially
bounded measurable functions on $\Omega$. Throughout, when we write
$b\in L^\infty_+(\Omega)$, we mean, with a slight abuse of notation,
that $b$ is a fixed non-negative, bounded, Lebesgue measurable
representative of an element of $L^\infty(\Omega)$. 
The the space of smooth functions with compact support in 
\(\Omega\) is denoted by \(C_c^\infty(\Omega)\). 
We write
$W^{k,2}(\Omega)$, $k\in\mathbb N$, for the Sobolev space of functions
with weak derivatives up to order $k$ in $L^2(\Omega)$, and
$W^{1,2}_0(\Omega)$ for the closure of $C_c^\infty(\Omega)$ in
$W^{1,2}(\Omega)$. The notation $\Delta=\Delta_D^\Omega$ stands for the
Dirichlet Laplacian on $L^2(\Omega)$, whereas
$\Delta_{\cD'(\Omega)}$ denotes the distributional Laplacian.

For real-valued functions $f$ and $g$ defined on a set $D$, we write
$f\lesssim g$ if there exists a constant $c>0$ such that $f\leq cg$ on
$D$. For a linear operator $A$ on $X$, its domain, range, kernel,
resolvent set, spectrum, point spectrum, and approximate point spectrum
are denoted by
$\Dom(A)$, $\Ran(A)$, $\Ker(A)$, $\rho(A)$, $\sigma(A)$,
$\sigma_p(A)$, and $\sigma_{\mathrm{ap}}(A)$, respectively. 
If $X$ is a Hilbert space, we call a selfadjoint operator $A$ on $X$
\emph{positive} if
\[
 \langle Ax,x\rangle>0,
 \qquad x\in\Dom(A)\setminus\{0\},
\]
and we call it non-negative if $\langle Ax,x\rangle\geq0$ for every $x\in\Dom(A)$. We denote by $\cF$ and $\cF^{-1}$ the Fourier transform and its inverse on scalar- or Banach-space-valued tempered distributions, defined by duality and normalized so that, for $f\in L^1(\R;X)$,
\[
(\cF f)(s)=\int_{\R}e^{-ist}f(t)\,dt,
\qquad
(\cF^{-1}f)(t)=\frac1{2\pi}\int_{\R}e^{ist}f(s)\,ds.
\]
For $f\in L^1_{\mathrm{loc}}(\R_+;X)$ we write $\cL f$ for its Laplace transform, when it exists:
\[
(\cL f)(\lambda)=\int_0^\infty e^{-\lambda t}f(t)\,dt.
\] 
For $\gamma \in (0,\pi)$ we set $\Sigma_\gamma := \{ \lambda \in \C\setminus \{0\} :  \, |\arg \lambda| < \gamma \}$. Finally,
$\mathbb R_+:=[0,\infty)$ and $\mathbb C_+:=\{\lambda\in\mathbb C:\Re\lambda>0\}$.

\subsection{Formulations of the main abstract results: Theorems \ref{thm 1}--\ref{thm:diffusion_theta=1}}\label{subs:main_results}

For linear operators $A$ and $B$ on a Banach space $X$ we set 
\begin{equation}\label{def:P_PD}
P(\la) := \la^2 + \la B + A \quad \textrm{ and } \quad P_D(\la): = \la B + A, \qquad \la \in \C.
\end{equation}
We write \(\rho(P)\) (resp. \(\rho(P_D)\)) for the set of \(\lambda \in \mathbb{C}\) such that \(P(\lambda)\) (resp. \(P_D(\lambda)\)) is boundedly invertible on $X$
and put
\begin{equation}\label{def:R_RD}
R(\la): = P(\la)^{-1} \qquad (\textrm{resp. } R_D(\la) := P_D(\la)^{-1} ).
\end{equation}
We return to the function $m$ from \eqref{eq:mintro}, now using the notation \eqref{def:P_PD}--\eqref{def:R_RD}. If $i\R\setminus\{0\}\subset\rho(P)\cap\rho(P_D)$, then for $x_0,y_0\in X$ it is given by 
\begin{equation}\label{def of m}
m(s) = m(s; x_0,y_0) : = R(is) (is x_0 +  B  y_0) - R_D (is)  B y_0, \qquad s\neq 0.
\end{equation}
Furthermore, we set 
\begin{equation}\label{set_X2}
\begin{split}
X_2 & := \left\{ x\in X: \sup_{|s|\geq 1} s^2 \| R(is) x\| < \infty\right\},\\
 X_D^{\eta} & :=  \left\{ x\in X: \sup_{|s|\geq 1} |s|^\eta \| R_D(is) B x\| < \infty\right\}, \qquad \eta \leq 1.
\end{split}
\end{equation}

\begin{hypothesis}\label{hyp:main} Let $\theta\in[0,1]$ be fixed. Let $A$ be a linear operator on a Banach space $X$, and let $B$ be a nonzero bounded linear operator on $X$ such that   
\begin{align}
\label{pencils est}
i\R\setminus \{0\} \subset \rho(P), & \qquad [R] := \sup_{s\neq 0} |s| \| R(is) \|<\infty,\\ 
\label{hyp:bound_2}
i\R\setminus \{0\} \subset \rho(P_D), & \qquad [R_D B]_\theta := \sup_{s\neq 0} \min(|s|, |s|^\theta) \| R_D (is)  B  \|<\infty.
\end{align}
\end{hypothesis}
For a subset $D\subset X$, we denote by $A^{-1}(D):=\{x\in\Dom(A): Ax\in D\}$ the preimage of $D$ under $A$. Moreover, we write $BD:=\{Bx:x\in D\}$ for the image of $D$ under $B$. We often write $BX=\Ran(B)$. No injectivity assumption is imposed on either $A$ or $B$.

To formulate assumptions on the initial data at once for all values of $\theta \in [0,1]$, we introduce the notation $Y_D^\theta$, $Z_D^\theta$, and $\widetilde Z_D^\theta$ as follows: 
\begin{equation}\label{def_YD2}
Y_D^{\theta} : = \begin{cases} X_D^{2-2\theta}, & \theta \in [\frac12,1],\\  
A^{-1}\big(BX_D^{1-2\theta}\big),&  \theta \in [0, \frac12),
\end{cases}  
\qquad 
Z_D^\theta :=
\begin{cases}
X_D^{2-3\theta}, & \theta\in[\frac13, 1],\\
A^{-1}\big(BX_D^{1-3\theta}\big), & \theta\in[0,\frac13).
\end{cases}
\end{equation}

\begin{align}\label{def_wtZD}
\wt Z_D^\theta &:=
\begin{cases}
X_D^{3-3\theta}, & \theta \in[\frac23,1],\\
A^{-1}\big(BX_D^{2-3\theta}\big), & \theta\in[\frac13,\frac23),\\
A^{-1}\big(BA^{-1}\big(BX_D^{1-3\theta}\big)\big), & \theta\in[0,\frac13).
\end{cases}
\end{align}

In particular, if $A$ and $B$ satisfy Hypothesis \ref{hyp:main} for $\theta \in [{\frac34},1]$, then  
\begin{equation}\label{Y_D}
X = Y_D^\theta = Z_D^\theta = \wt Z_D^\theta. 
\end{equation}

The spaces introduced above encode the high-frequency
regularity needed for the estimates in Theorem \ref{thm 1} below.
Their definitions reflect the decay of the resolvents
and the identities relating $A$, $B$, and $R_D$,
rather than any particular choice of differential
operators. In particular, the exponents and range
conditions in \eqref{def_YD2}--\eqref{def_wtZD} arise from the estimates
for derivatives of $R_D(is)B$; see e.g. Lemma \ref{lem:m3_in_Mikhlin}.
The piecewise definitions account for the different
amounts of regularity required as $\theta$ varies. 

\begin{theorem}\label{thm 1} Under Hypothesis \ref{hyp:main} the following assertions hold. 
\begin{itemize}
\item [(i)]  For all $x_0, y_0 \in X$ such that 
\begin{equation}\label{assump_thm1(i)}
x_0 \in \Dom(A), \, y_0\in Y_D^{\theta} \quad \textrm{ and } \quad B x_0, \,  B  y_0 \in X_2
\end{equation}
the inverse Fourier transform $\cF^{-1} m$ of $m=m(\cdot; x_0, y_0)$  
is in $C(\R\setminus \{0\}; X)$ and satisfies the estimate
\begin{equation}\label{est_thm1(i)}
\|(\cF^{-1}m)(t)\| \lesssim {1}/{|t|},\qquad t\neq 0.
\end{equation}
In particular, $\Dom(A) \subset X_2$. 
Moreover, if, in addition to Hypothesis \ref{hyp:main}, $\theta \in [1/2, 1]$ and  $\sup_{|s|\geq 1}\|R_D(is)\|<\infty$, then $\Dom(A) \subset X_D^{1}$ and  \eqref{est_thm1(i)} holds whenever $x_0, y_0, Bx_0,  B y_0 \in \Dom(A)$.

\medskip

\item [(ii)] If, in addition to \eqref{assump_thm1(i)},  $x_0, B y_0  \in \Ran (A)$, and $y_0\in Z_D^{\theta}$, then 
\[
\|(\cF^{-1}m)(t)\| \lesssim 1/t^2, \qquad t\neq 0. 
\]

\item [(iii)] If, in addition to \eqref{assump_thm1(i)}, $y_0\in \wt Z_D^{\theta}$ and   
\begin{equation}\label{cond_kinetic}
  Bx_0\in \Dom(A), \, Ax_0, \, B^2x_0 \in X_2, \quad\textrm{and}\quad  By_0 \in \Dom(A), \, B^2 y_0 \in X_2, 
\end{equation}
then $\cF^{-1}m\in C^{1}(\R\setminus \{0\};X)$ and satisfies
\[
\|(\cF^{-1} m)'(t)\|\lesssim 1/t^2, \qquad t\neq 0.
\]
\end{itemize}

\end{theorem}

The proof of Theorem~\ref{thm 1} is sufficiently flexible to accommodate asymptotic behaviours of $R$ and $R_D$ beyond those covered by Hypothesis \ref{hyp:main}. However, it also has certain limitations.   
The reason for restricting our analysis of the decay of $\mathcal{F}^{-1} m$ to the order $O(1/t^2)$ is explained in Subsection~\ref{subs:limit-of-decay}; see Proposition~\ref{prop: sharp}. 
The condition $\sup_{|s|\geq 1}\|R_D(is)\|<\infty$
appearing in the last statement of part~(i) of Theorem~\ref{thm 1} naturally
arises in our model example of (DWE) with the Dirichlet Laplacian; see the first estimate in \eqref{est: 1/sqrt at infty} of Theorem~2.2(ii). 
Analogous statements for the components of the function $m = m_H - m_P$ corresponding to the equations (HE) and (PE), namely $m_H(s): = R(is)(is x_0 + By_0)$ and $m_P(s): = R_D(is)By_0$, $s\neq 0$, can be extracted from the proof of Theorem~\ref{thm 1}. 
The corresponding estimates for $\mathcal{F}^{-1}m_H$ and
$\mathcal{F}^{-1}m_P$ are recorded in Corollary~\ref{cor:decay_of_mH}.

\medskip
 
The following result establishes a connection between
the inverse Fourier transform of $m$ and solutions to {(HE)} and {(PE)}. We refer the reader to Subsection \ref{subsec:well-posed-DE}  for the definition of mild solutions to (PE) and (HE). We provide in Propositions \ref{mild-DE} and \ref{mild_DWE}, respectively, their characterizations in terms of the Laplace transform.

\begin{theorem}\label{thm:ind-solu}  Assume Hypothesis \ref{hyp:main} and, in addition, that $\{\Re\la >0\}\subset \rho(P)\cap \rho(P_D)$. 
 Let $x_0, y_0 \in X$ and denote $x_1 := By_0 - Bx_0$.  Assume that there exists a mild solution $u = u(t)$ to \emph{(HE)} with the initial data $x_0$, $x_1$, and a mild solution $v=v(t)$ to \emph{(PE)} with the initial data $y_0$. Assume that $u$ and $v$ have at most polynomial growth at infinity.
Let $m=m(\cdot;x_0,y_0)$ be the function introduced in
\eqref{def of m}, and let $\wt u$ and $\wt v$ denote the extensions of
$u$ and $v$, respectively, by zero to $(-\infty,0)$. Then
\[
\cF^{-1}m=\wt u-\wt v
\quad\textrm{in }\cS'(\R;X).
\]

In particular, the restriction of $\cF^{-1}m$ to $(0,\infty)$ is represented by the continuous function $u-v$.
\end{theorem}

Our next main result, Theorem~\ref{thm:mild_DPE}, establishes the existence of mild solutions to (PE) that admit bounded holomorphic extensions to a sector in $\C$, under the following parabolic type estimate \eqref{DE-parab_est}, which is stronger than~\eqref{hyp:bound_2}.
Recall the notation $\Sigma_\gamma= \{ \lambda \in \C\setminus \{0\} :  \, |\arg \lambda| < \gamma \},$ $\gamma \in (0,\pi).$
\begin{hypothesis}\label{hyp:main_2} Let $\theta\in (0,1]$ be fixed.
We assume that $A$ is a linear operator and $B$ is a nonzero bounded linear operator on a Banach space $X$ such that for some $\omega \in (0,\pi/2)$
\begin{equation}\label{DE-parab_est}
\Sigma_{\omega+\pi/2}\subset \rho(P_D) \quad \textrm{and} \quad  \sup_{\la \in \Sigma_{\omega+\pi/2}} \min\big(|\la|, |\la|^{\theta}\big) \| R_D(\la)B\| <\infty. 
\end{equation}
\end{hypothesis}

Condition \eqref{DE-parab_est} naturally arises in our model example considered in Section~\ref{sec:examples} (where $A$ is the negative Dirichlet Laplacian and $B$ is a multiplication operator); cf. Theorem~\ref{thm: pencil est} and Proposition \ref{cond: 1/t at infty}.

Assuming \eqref{DE-parab_est}, we extend the approach developed by Favini and Yagi in~\cite[Chapter III]{FaYa99} for studying degenerate parabolic type equations by means of semigroups generated by multivalued operators. For related treatments of degenerate analytic semigroups associated with operator pairs, see also \cite{Sv94}, \cite{BC02}, and \cite{Ba08}.
Rather than considering arbitrary multivalued linear operators, we restrict our attention to operators of the form $B^{-1}A$,
where $A$ is a closed operator and $B$ is a bounded, not necessarily injective operator on the Banach space $X$. The operator $B^{-1} A$ is defined by 
\[
\Dom(B^{-1}A):= A^{-1}\big(BX\big), \qquad  
(B^{-1}A)x := y + \Ker(B),
\]
where $y\in X$ is any vector satisfying $Ax=By$. If $B$ is injective, then $B^{-1}A$ is a single-valued operator on $X$. Unlike Favini and Yagi, we allow $B^{-1}A$ to be non-invertible. This case mainly requires a change of contours in the integral formula defining the semigroup generated by $B^{-1}A$. However, to define the fractional powers of such an operator, we rely on the Balakrishnan formula rather than on the Dunford--Riesz calculus. 
Moreover, we prefer to formulate all assumptions directly in terms of the pencil $P_D(\lambda)=\lambda B+A$ and its inverse $R_D(\lambda)$, rather than using the terminology of multivalued operators. We point out that there is no direct relation between the modified resolvent $R_D$ and the resolvent of $B^{-1}A$; cf. \cite[Chapter~I, Theorems~1.14 and~1.15]{FaYa99}.

We now define the solution family associated with (PE) for a pair of
operators $A$ and $B$ satisfying condition~\eqref{DE-parab_est} for some
$\theta\in(0,1]$ and $\omega\in(0,\pi/2)$. 
For $\rho>0$ and $\omega' \in (0,\omega)$, let $\Gamma = \Gamma(\rho, \omega')$ be the oriented path consisting of
\begin{equation*}
\begin{split}
\Gamma_\pm  & = \Gamma_\pm(\rho, \omega'):=\left\{ r e^{\pm i(\omega' + \pi/2)}: \, r\geq \rho\right\},\\
\Gamma_0 & = \Gamma_0(\rho, \omega') : = \left\{\rho e^{i\gamma}: |\gamma|\leq \omega'+\pi/2\right\}.
\end{split}
\end{equation*}
For $z\in \Sigma_\omega$ and  $\omega' \in (|\arg z |,  \omega)$ we set
\begin{equation}\label{formula:semig}
\cT(z): = \frac{1}{2\pi i} \int_{\Gamma} e^{\lambda z} (\lambda B + A)^{-1}B\,  \ud \la. 
\end{equation}
For any $\beta \in (1 - \theta,1)$ let 
\begin{equation}\label{formula:frac_powers}
J_\beta = \frac{\sin(\pi\beta)}{\pi} \int_0^\infty t^{-\beta} ( (1+ t)B + A)^{-1}B\, \ud t
\end{equation}
(convergence is in the operator norm). We define the family of spaces $D_\ell$, $\ell\in\mathbb N$, inductively by
\begin{equation*}
D_1: = A^{-1}(BX), \qquad D_\ell : = A^{-1}(BD_{\ell -1}) \quad \textrm{for } \ell\geq 2. 
\end{equation*}

\begin{theorem}\label{thm:mild_DPE} Assume Hypothesis \ref{hyp:main_2}. 
Let $\wt \omega \in (0,\omega)$ and $\beta \in (1-\theta,1)$. Then the following assertions hold.
\begin{itemize}
\item [(i)] The family $\cT(z)$, $z\in \Sigma_\omega$, is a holomorphic semigroup on $X$ such that 
\[
\sup_{|\arg z|\leq \wt \omega \atop |z|\ge 1}\|\cT(z)\|<\infty \, \textrm{ and } \, \|\cT(z)\| = O(|z|^{\theta - 1}) \textrm{ as } |z|\to 0, |\arg z|< \wt \omega.
\]
Moreover, for any $\la \in \C$ with $\Re \la >0$ we have
\begin{equation}\label{eq:laplace_int}
R_D(\la) B = \int_0^\infty e^{-\la t} \cT(t)  \ud t 
\end{equation}
(convergence is in the operator norm).  

\smallskip

\item [(ii)] The inclusion $A^{-1}\big(BX\big) \subset X_D^1$ holds. Moreover, for every $y\in A^{-1}\big(BX\big)$,
\[
\sup_{\lambda\in\Sigma_{\omega+\pi/2}}\|\lambda R_D(\lambda)By\|<\infty.
\]
For every $y\in A^{-1}\big(BX\big)$, the function $v:[0,\infty)\to X$ defined by
\[
v(t):=\cT(t)y,\,\, t>0,\quad\text{and}\quad v(0):=y,
\]
is a mild solution to \emph{(PE)} with initial value $y$. Furthermore, $v$ admits a bounded holomorphic extension from $(0,\infty)$ to the sector $\Sigma_{\widetilde{\omega}}$.

\smallskip

\item [(ii')]  $\cT(t)J_\beta\to J_\beta$ in operator norm as $t\to0+$. In particular, the function $v:[0,\infty)\to X$ defined by
\[
v(t):=\cT(t) J_\beta x,\quad t>0,\qquad\text{and}\qquad v(0):= J_\beta x,
\]
is a mild solution to \emph{(PE)} with initial value $J_\beta x$. Moreover, $v$ admits a~bounded holomorphic extension to the sector $\Sigma_{\wt\omega}$. Furthermore,  $\Ran(J_\beta)\subset X_D^1$; more precisely, 
\[
\sup_{\lambda\in\Sigma_{\wt \omega+\pi/2}}\|\lambda R_D(\lambda)B J_\beta\|<\infty.
\]
\smallskip

\item [(iii)] If $\Ran(B)\not\subset\Ran(A)$, then for any $t>0$ the spectral radius of the operator $\cT(t)$ is greater than or equal to $1$, and, consequently, for every non-increasing function $r:\R_+\to \R_+$ with $r(t)\to 0$ as $t\to \infty$ there exists a vector $y\in \bigcap_{\ell\in \N} D_\ell$ such that $\|\cT(t) y\|\gtrsim r(t)$, $t>0$. 
\end{itemize}
\end{theorem}

Parts (ii) and (ii') of Theorem~\ref{thm:mild_DPE} identify two classes of initial values in $X$ for which the corresponding {(PE)} admits mild solutions with bounded holomorphic extensions to some sector in $\mathbb{C}$. For $\theta\in(1/2,1]$ and $\beta\in(1-\theta,\theta)$, we are able to prove that $A^{-1}(\operatorname{Ran}(B))\subset\operatorname{Ran}(J_\beta)$. Whether this inclusion remains valid for other values of $\theta$ and $\beta < 1$ is not known; see Remark~\ref{rem:range_F_kappa}. The proofs of (ii) and (ii') rely on different arguments.

Part (iii) connects the spectral properties of \(\cT\) with the existence of slowly decaying orbits, which are needed for the diffusion phenomenon. For \(C_0\)-semigroups, such properties are often studied through the generator. In our setting, the corresponding operator \(-B^{-1}A\) may be multivalued when \(\theta<1\). Nevertheless, the spectral properties of \(\cT\) can be studied directly, and part (iii) shows that the simple condition \(\operatorname{Ran}(B)\not\subset\operatorname{Ran}(A)\) guarantees the existence of slowly decaying orbits for $\cT$.

For $\theta=1$, a more precise result is available, see Theorem~\ref{thm: mild to DE}. In particular, part~(iv) of Theorem~\ref{thm: mild to DE}, combined with the theory of sectorial operators, provides a more complete characterization of the existence of weakly decaying orbits than the range condition given in part~(iii) of Theorem~\ref{thm:mild_DPE} with $\theta =1$; see the discussion in Section~\ref{subsec:weak-orbits} and Proposition \ref{thm:decying-orbits-of-T} with Corollary~\ref{cor:decay-orbits-L}.

For equation (HE), outside the Hilbert space setting, where the Lumer--Phillips generation theorem is often applicable, establishing the existence of mild or classical solutions for a given class of initial data is a delicate issue. To address this issue, throughout the paper we ensure the existence of \emph{mild solutions} to (HE) for a given class of pairs of initial data (see the definition in Subsection~\ref{subsec:well-posed-DE}) by imposing the assumption of \emph{well-posedness} of (HE), which we describe in Section~\ref{subsec well posed}. We refer the reader to Lemma~\ref{lem: mild s to DWE} for the corresponding result. Moreover, in Theorem~\ref{thm: DWE well-posed} and its corollary, we show that well-posedness of (HE) is equivalent to $-A$ being the generator of an operator cosine function; see, e.g., \cite[Section 3]{ABHN11} for a discussion of cosine function theory.

 Theorem~\ref{thm 1} provides decay estimates for \(\mathcal F^{-1}m\) under assumptions involving the space \(X_2\), while Theorem~\ref{thm:ind-solu} identifies \(\mathcal F^{-1}m\) with the difference \(u-v\) whenever the corresponding solutions exist.

In the UMD setting, assuming that (HE) is well-posed, Theorem~\ref{thm on equiv cond} yields
\[
\Dom((A+\omega)^{1/2})\subset X_2
\]
for every $\omega\geq 0$ such that $A+\omega$ is sectorial. Hence, the space $X_2$ in Theorem~\ref{thm 1} can be replaced by $\Dom((A+\omega)^{1/2})$, which has a more explicit description in terms of the operator $A$. Moreover,  Proposition~\ref{prop:linear_growth} (see also Remark~\ref{Goldstein example}) explains how the polynomial boundedness assumption on the corresponding $C_0$-semigroup $\cU$ on the phase space $V\times X$, imposed in Theorem~\ref{thm: diffusion intro}, arises in concrete applications.

\begin{theorem}\label{thm: diffusion intro} Assume Hypotheses~\ref{hyp:main}~and~\ref{hyp:main_2} and, in addition, that $X$ is a UMD space, that the equation~\emph{(HE)} associated with $A$ and $B$ is well posed, and that the corresponding $C_0$-semigroup $\mathcal{U}$ is polynomially bounded.
Let  $\omega \ge 0$ be such that $A + \omega$ is sectorial on $X$. 

Then, $V = D\bigl((A+\omega)^{1/2}\bigr)$. Moreover, the following assertions hold.

\begin{itemize}
\item [(i)] For every $x_0 \in \Dom(A)$ and every $y_0\in Y_D^\theta$ that belongs to either $A^{-1}(BX)$ or $\Ran(J_\beta)$ for some $\beta \in( 1-\theta, 1)$ (or even for every $y_0\in \overline{A^{-1}\big(BX\big)}$ when $\theta =1$) such that  $B x_0, \, B y_0 \in V$, if  $x_1 : =  B y_0 - Bx_0$, then 
\begin{equation}\label{decay of u-v}
\|u(t; x_0,x_1) - v(t; y_0)\| \lesssim 1/t, \quad t>0,  
\end{equation}
where  $u(\cdot; x_0, x_1)$ and $v(\cdot; y_0)$ denote the corresponding mild solutions to \emph{(HE)} and \emph{(PE)}, respectively.
\smallskip

\item [(ii)] If, in addition to the properties of $x_0$ and $y_0$ listed in \emph{(i)}, $x_0, B y_0 \in \Ran(A)$, and $y_0 \in  Z_D^\theta$, then 
\begin{equation}\label{decay of u-v t2}
\|u(t; x_0,x_1) - v(t; y_0)\| \lesssim 1/t^2,  \quad t>0.  
\end{equation}
Moreover, 
\begin{equation}\label{est:thm4(ii)_u}
\|u(t;x_0,x_1)\| \lesssim 1/t \quad\textrm{ and } \quad \|v(t;y_0)\| \lesssim 1/t, \quad t>0.
\end{equation}

\item [(iii)] If, in addition to the properties of $x_0$ and $y_0$ listed in \emph{(i)}, $y_0\in \wt Z_D^\theta$ and
\begin{equation}\label{cond_kinetic_V}
 Bx_0\in \Dom(A), \, Ax_0, \, B^2 x_0 \in V, \quad\textrm{and}\quad  By_0 \in \Dom(A), \, B^2 y_0 \in V,
\end{equation}
 then $u(\cdot; x_0, x_1)$ and $v(\cdot; y_0)$ are in $C^1((0, \infty);X)$ and 
\[
\|u'(t; x_0, x_1) - v'(t; y_0)\| \lesssim 1/t^2,  \quad t>0.
\]
Moreover,
\begin{equation*}
 \|u'(t;x_0,x_1)\| \lesssim 1/t \quad \textrm{ and } \quad  \|v'(t;y_0)\| \lesssim 1/t,\quad t\geq 1.
\end{equation*}
\end{itemize}
\end{theorem}

The decay of the difference of solutions to (HE) and (PE) (even for a large class of initial data as established in Theorem \ref{thm: diffusion intro}), does not automatically imply the presence of the diffusion effect. To verify that the diffusion effect does occur one also needs to show how to choose initial data for which one of the solutions, $u$ or $v$, decays slower than their difference. The existence of initial data such that the corresponding solutions satisfy both requirements is not obvious.
In particular, additional assumptions (to those  in (iii) of Theorem \ref{thm:mild_DPE}) must be imposed on the operator $B$ for the diffusion effect to occur. 
Regarding estimate \eqref{est:thm4(ii)_u}, clearly \eqref{decay of u-v} yields 
\[
\|u(t; x_0,x_1)\|\lesssim \max\big(1/t, \|v(t,y_0)\|\big).
\]
Therefore, for every $x_0$ and $y_0$ satisfying the assumptions of (i)
and such that $\|v(t;y_0)\|\lesssim t^{-1}$ for $t>0$, we obtain
\eqref{est:thm4(ii)_u}; cf. Corollary \ref{cor:decay_of_mH} and  Remark~\ref{rem:a_statements}.
A~similar observation applies to parts~(ii) and (iii).

As a further development of Theorems~\ref{thm:mild_DPE} and~\ref{thm: diffusion intro}, we present in the next Theorems \ref{thm: diffusion} and \ref{thm:diffusion_theta=1} our main result on the diffusion phenomenon for (HE), which in particular guarantees the existence of the initial data $y_0$ such that the corresponding solution $v(\cdot; y_0)$ to (PE) serves as an asymptotic diffusive profile for solutions $u(\cdot; x_0,x_1)$ to (HE) with arbitrary $x_0\in\Dom(A)$ and a suitably chosen $x_1$, depending on $y_0$ and~$x_0$.

To motivate the formulation of our diffusion results, we note that $D_2\subset  Y_D^\theta$ for every $\theta\in(0,1]$. Moreover, Part~(iii) of Theorem~\ref{thm:mild_DPE} ensures the existence of slowly decaying mild solutions to (PE) with initial data $y_0\in\bigcap_{\ell\in\N}D_\ell$. In particular, such initial data $y_0\in D_2$ satisfy the assumptions of Theorem~\ref{thm: diffusion intro}(i) for every $\theta\in(0,1]$. 
Moreover, since the case $\theta=1$ allows for a significantly more detailed analysis of the spectral properties of the associated semigroup $\cT$ (see Theorem~\ref{thm: mild to DE}), Theorem~\ref{thm:diffusion_theta=1} provides additional diffusion statements with refined decay rates for the solutions and their derivatives, corresponding to parts~(ii) and~(iii) of Theorem~\ref{thm: diffusion intro}.

\begin{theorem}\label{thm: diffusion}
Let $A$, $B$ and $X$ satisfy the assumptions of Theorem \ref{thm: diffusion intro} for some $\theta \in (0,1]$. In~addition, we assume that $\Ran(B)\setminus \Ran(A) \neq \emptyset$ and $B(\Dom(A)) \subset \Dom((A+\omega)^{1/2})$. 
Let $r: \R_+ \rightarrow \R_+$ be an arbitrary non-increasing bounded  function such that $r(t)\rightarrow 0$ and 
$1/t = o(r(t))$ as $t\to \infty$.

 There exists $y_0\in D_2$ such that, for every $x_0\in\Dom(A)$ and $x_1:= B y_0-Bx_0$, the corresponding solutions $u(\cdot; x_0,x_1)$ to \emph{(HE)} and $v(\cdot; y_0)$ to \emph{(PE)} satisfy the estimates
\begin{equation}\label{dp}
\|u(t; x_0,x_1) - v(t; y_0)\| \lesssim 1/t, \,\,\,\, \|v(t;y_0)\|\gtrsim r(t), \, \textrm{ and } \,\,  \|u(t;x_0,x_1)\| \gtrsim r(t)
\end{equation} for all sufficiently large $t>0$,  
and thus the diffusion phenomenon occurs for $u(\cdot; x_0,x_1)$. As a result, the solutions $u(\cdot; x_0,x_1)$ and $v(\cdot; y_0)$ have the same rate of decay, that is, 
\begin{equation}\label{decay u/v}
\|u(t, x_0,x_1)\|/\|v(t;y_0)\| \rightarrow 1 \quad \textrm{as } t\rightarrow \infty.
\end{equation}
\end{theorem}

Theorem~\ref{thm: mild to DE} presented in Section~\ref{sec:proofs} asserts that, if $X$ is reflexive and $A$ and $B$ satisfy Hypothesis~\ref{hyp:main_2} with $\theta=1$, then $X$ and the corresponding semigroup $\cT$ admit the decompositions
\begin{equation}\label{eq:decomp_of_cT}
X=\overline{A^{-1}(BX)}\oplus\Ker(B),
\qquad
\cT(t)=T(t)\oplus {0}, \quad t > 0, 
\end{equation}
where $T$ is a $C_0$-semigroup on $Y:=\overline{A^{-1}(BX)}$, 
and ${0}$ denotes the zero operator on $\Ker(B)$. The infinitesimal generator of $T$ is denoted by $L$; it is computed in \eqref{def:L}.

\begin{theorem}\label{thm:diffusion_theta=1} Under the assumptions of Theorem~\ref{thm: diffusion} with $\theta=1$, the following additional assertions hold.
\begin{itemize} 
\item [(i)] If $X$ is a Hilbert space, then the third estimate in \eqref{dp} can be strengthened as follows: there exists a vector $x\in X$ and $t_0>0$ such that
\begin{equation}\label{T>r}
\left|\langle x, u(t; x_0,x_1)\rangle\right| \gtrsim r(t),
\qquad t\ge t_0.
\end{equation}

\item [(ii)] If $0\in \sigma(L)\setminus \sigma_p(L)$, then there exists $y_0\in Y$ such that for every $x_0\in \Dom(A)\cap \Ran(A)$ and $x_1:=By_0-Bx_0$ the following estimates hold:
\begin{equation}\label{dp t2}
\|u(t; x_0,x_1) - v(t; y_0)\| \lesssim 1/t^2, \qquad  
\|u(t; x_0,x_1)\|\lesssim 1/t, \quad  t>0,
\end{equation} 
\[
\limsup_{t\rightarrow \infty}\, {t}{r(t)}^{-1} \|u(t;x_0,x_1)\| >0.
\]

\item [(iii)] If, in addition, $B(\Dom(A))\subset \Dom(A)$, then there exists $y_0\in \Dom(A)$ such that, for every $x_0\in \Dom(A)$ with $Ax_0\in \Dom((A+\omega)^{1/2})$ and  $x_1: = By_0 - Bx_0$, the corresponding solutions $u(\cdot; x_0, x_1)$ and $v(\cdot; y_0)$ are in $C^1((0, \infty);X)$ and 
\[
\|u'(t; x_0, x_1) - v'(t; y_0)\| \lesssim 1/t^2, \qquad \|u'(t;x_0,x_1)\| \lesssim 1/t,  \quad t>0, 
\]
\[
\limsup_{t\rightarrow \infty}\, t r(t)^{-1} \left\| u'(t; x_0, x_1)\right\|>0.
\]
\end{itemize}
\end{theorem}

\subsection{Overview of  the diffusion phenomenon for (DWE); Theorem~\ref{thm:diffusion_for_Delta}}\label{subs:applications}
We now focus on the model case of (HE), namely the damped wave equation (DWE) on $L^2(\Omega)$ with the Dirichlet boundary conditions for general domains $\Omega \subseteq \mathbb{R}^n$:
\[
u''(t) + bu'(t) - \Delta u(t) = 0, \quad t\geq 0\, \textrm{ with } u(0) = x_0, \, u'(0) =x_1.\tag{DWE}
\]
Here, $b$ denotes a non-negative, bounded, Lebesgue measurable function on $\Omega$, not identically equal to zero (in short, $b\in L^\infty_+(\Omega)\setminus\{0\}$), which we refer to as the \emph{damping coefficient}. We recall our convention: when we write $b\in L^\infty_+(\Omega)$, we mean, with a slight abuse of notation, that $b$ is a fixed non-negative, bounded, Lebesgue measurable representative of an element of $L^\infty(\Omega)$.
The associated problem (PE) is of the form 
\[
 bv'(t) - \Delta v(t) = 0, \quad t\geq 0\, \textrm{ with } v(0) = y_0. \tag{DE}
\]
We denote by $M_b$ the multiplication operator induced by $b \in L^\infty_+(\Omega)$, that is,
\[
M_b f := b f, \quad f \in L^2(\Omega).
\]
We write $\Delta_D^\Omega$, or simply $\Delta$, for the Dirichlet Laplacian on $L^2(\Omega)$, and denote the distributional Laplacian by $\Delta_{\cD'(\Omega)}$. 
Our aim is to formulate conditions on $\Omega$ and on $b \in L^\infty_+(\Omega)\setminus\{0\}$ under which the assumptions of the abstract results, namely Theorems~\ref{thm 1}--\ref{thm:diffusion_theta=1}, are satisfied with $A := -\Delta$, $B := M_b$, and $X = L^2(\Omega)$. In what follows, the pencils $P, P_D$ and the associated resolvent families $R, R_D$ (see \eqref{def:P_PD} and \eqref{def:R_RD}) correspond to $A = -\Delta$ and $B = M_b$.
To make our analysis systematic, we introduce three admissible classes 
of damping coefficients, 
\begin{align}
\cB(\Omega)
& := 
\left\{ b \in L^\infty_+(\Omega)\setminus \{0\}: \textrm{ \eqref{pencils est} in Hypothesis \ref{hyp:main} holds} \right\},
\label{def:B_class}
\\
\cB_D^{\theta} (\Omega)
& :=
\left\{b \in L^\infty_+(\Omega)\setminus \{0\}:  \textrm{ \eqref{DE-parab_est} in Hypothesis \ref{hyp:main_2} holds} 
\right\},
\label{def:B_D_class}
\\
\cR(\Omega)
& := \left\{b \in L^\infty_+(\Omega)\setminus \{0\}: \Ran(B) \setminus \Ran(A) \neq \emptyset \right\}.
\label{def:R_class}
\end{align}

We restrict our attention to the classes $\cB_D^{\theta}(\Omega)$ with $\theta\in\left[\frac12,1\right]$. The classes $\mathcal{B}(\Omega)$, $\cB_D^{\theta}(\Omega)$, and $\mathcal{R}(\Omega)$ will be used below to fit the model case~(DWE)--(DE) into the abstract framework for~(HE)--(PE) presented in the previous section.
 In particular, if
$b\in\cB(\Omega)\cap\cB_D^\theta(\Omega)$,
then Hypothesis~\ref{hyp:main} holds, while
the condition $b\in\cR(\Omega)$ is relevant
to Theorem~\ref{thm:mild_DPE}(iii). 
For $\theta=1$, Proposition~\ref{prop:B_D_gives_DE}
shows that condition~\eqref{hyp:bound_2} in
Hypothesis~ \ref{hyp:main} implies condition~\eqref{DE-parab_est} in Hypothesis~\ref{hyp:main_2}. The converse is immediate.
Thus $b\in\mathcal B_D^1(\Omega)$ whenever
\eqref{hyp:bound_2} holds.
Also, if $b\in \cB_D^{\theta}(\Omega)\cap \cR(\Omega)$, then by Theorem~\ref{thm:mild_DPE}(iii), there exist slowly decaying orbits of the corresponding semigroup $\cT$.

In the case $\theta=1$, Theorem~\ref{thm: mild to DE}(iv) provides a stronger result: $b\in\mathcal{R}(\Omega)$ if and only if $0\in\sigma(L)$.
Combining this with a standard result from the theory of holomorphic semigroups (see Remark \ref{rem:a_statements}(a2)), we conclude that for every $b \in \cB_D^{1}(\Omega)\setminus \mathcal{R}(\Omega)$, any non-exponentially decaying solution $u$ to (DWE) does not admit a diffusive profile described by a solution to (DE).
 
Moreover, if for a given $\Omega$ the classes introduced in~\eqref{def:B_class}--\eqref{def:R_class} have non-empty intersections, then for $b\in \cB(\Omega)\cap \cB_D^\theta(\Omega)\cap \cR(\Omega)\cap W^{1,\infty}(\Omega)$ all assumptions of Theorems~\ref{thm:mild_DPE}--\ref{thm: diffusion} for a given $\theta$ are satisfied, and so these results can be reformulated in terms of classes~\eqref{def:B_class}--\eqref{def:R_class}. For instance in the toy case when $b\ge b_0>0$ almost everywhere on $\Omega$, one has
$b\in\mathcal B(\Omega)\cap\mathcal B_D^1(\Omega)$. If, in addition, $\Omega$ does not have the Poincar\'e property, then also $b\in\mathcal R(\Omega)$. Theorem~\ref{thm: diffusion intro} applies to the initial data specified in its statement. If moreover $b\in W^{1,\infty}(\Omega)$, the compatibility assumption of Theorem~\ref{thm: diffusion}  holds. The additional assertions of Theorem~\ref{thm:diffusion_theta=1}  apply under their respective hypotheses.

For the reader's convenience, we supplement the discussion by the following statements and provide a brief overview of the organization of Section~\ref{sec:examples}.

\smallskip

(s1) For any $b\in L^\infty_+(\Omega)\setminus\{0\}$, problem (DWE) 
corresponding to $A=-\Delta$ and $B=M_b$ is well posed, the associated first-order $C_0$-semigroup $\cU$ has at most linear growth, and the phase space coincides with $W^{1,2}_0(\Omega)\times L^2(\Omega)$; see Proposition \ref{prop:linear_growth}, which also implies that $\{\Re \lambda >0\}\subset \rho(P)$.
Proposition~\ref{prop:interm_freq} provides further conditions on $\Omega$ and $b$ guaranteeing that also $i\mathbb{R}\setminus \{0\} \subset \rho(P)$.

The study of the asymptotic behaviour of $R(is)$ as $|s|\to\infty$, often referred to as \emph{high-frequency estimates}, is a delicate matter and constitutes a research topic in its own right; cf., for instance, \cite{BuJo16} and the references therein. Both the geometry of $\Omega$ and the properties of the damping coefficients $b$ play a crucial role in such studies. Here, we do not aim at obtaining new high-frequency estimates for the pencil $P$, but merely use results available in the literature which provide conditions ensuring that $b$ belongs to the class 
\begin{equation}\label{def:B}
\cB_\infty(\Omega): = \{ b\in L^\infty_+(\Omega)\setminus\{0\}:  i\mathbb{R}\setminus \{0\} \subset \rho(P)  \text{ and } 
 \sup_{|s| \geq 1} \|\, i s R(is) \,\| < \infty \}.
\end{equation} 
In Lemma~\ref{lem:B_infty_class}, we adapt two results from \cite{AlIbKh15, BuJo16} ensuring that $b\in \cB_\infty(\Omega)$. This serves as the basis for Theorem~\ref{thm:diffusion_for_Delta} below, where we describe $\cB(\Omega)$ in two cases: when $\Omega=\mathbb{R}^n$ and when $\Omega$ is an \emph{exterior domain} in $\mathbb{R}^n$. 

\smallskip

(s2) Simple toy examples for $X=\ell^2(\N)$ illustrate that the location of
the resolvent sets and the asymptotic behaviour of the pencils $P$ and
$P_D$ associated with~(HE) and~(PE) may differ significantly in the
high-frequency regime. However, Lemma~\ref{lem near 0} shows that,
under suitable assumptions, the norms of $R(\la)$ and $R_D(\la)$ are
equivalent near zero, even in the abstract setting. In particular, $\cB_D^{\theta}(\Omega) \cap {\mathcal{B}}_\infty (\Omega)\subset \cB(\Omega)$. It remains unclear whether the classes $\mathcal B(\Omega)$ and $\mathcal B_D^\theta(\Omega)$ are comparable under inclusion.
 This motivates the study of the class $\cB_D^{\theta}(\Omega)$ independently of results concerning~$\cB(\Omega)$. 

Theorem~\ref{thm: pencil est} (see also Proposition~\ref{ex: of b} and Corollaries~\ref{cor:examples_Omega_b} and~\ref{cor:examples_b}) identifies a large class of functions $b$ that belong to $\cB_D^{1/2}(\Omega)$, while Propositions~\ref{prop:empty_resolvent_set_PD} and~\ref{prop:imaginary_exis_in_specturm} identify functions $b$ that do not belong to this class. For $\cB_D^\theta(\Omega)$ with $\theta \in (\frac12,1]$, see Proposition \ref{cond: 1/t at infty}. 

\smallskip

(s3) In Subsection~\ref{subs:ran_cond} we identify functions $b$ from $\cR (\Omega)$; see Proposition~\ref{prop:RninR}, where, in particular, we prove that $\cR(\R^n) =  L^\infty_+(\R^n) \setminus \{0\}$ for $n\leq 4$, and Proposition \ref{prop:correction_b} and Corollary~\ref{cor:correction_b} for $\Omega \neq \R^n$. 
We introduce the following class related to the Weyl sequences for the square root $(-\Delta)^{1/2}$ of $-\Delta$:
\begin{align}
\label{def:W_class}
\nonumber \cB_\cW(\Omega)
:=
\Bigl\{
b  \in & L^\infty_+(\Omega)
 :  \,
 \exists
(\varphi_j)\subset W^{1,2}_0(\Omega)
\text{ with pairwise disjoint supports, }
\\
& \|\varphi_j\|_{2}=1, \,
\lim_{j\to \infty}\|\nabla \varphi_j\|_{2} = 0, \text{ and }
\liminf_{j\to\infty} \operatorname*{ess\,inf}_{\xi\in\supp(\varphi_j)} b(\xi)>0 \Bigr\}.
\end{align}
Proposition~\ref{prop:correction_b} asserts that $\cB_\cW(\Omega) \subset \cR(\Omega)$.  If $\Omega$ contains a sequence of balls $B_j$, $j\in \N$, with radii tending to infinity and if $b(\xi) \geq b_0 > 0$ for all $\xi \in \bigcup_j B_j$, then $b \in \cB_\cW(\Omega)$. 
In the context of condition~\ref{def:W_class}, we also refer the reader to Propositions~\ref{prop:empty_resolvent_set_PD} and \ref{prop:imaginary_exis_in_specturm}, which indicate that the assumption that the coefficient $b$ is bounded away from zero on the supports of Weyl sequences cannot be completely omitted.
However, there exist damping coefficients $b\in\cB_D^\theta(\Omega)\cap\cB(\Omega)$ such that $b\notin\cB_\mathcal{W}(\Omega)$; see the example constructed in the proof of Proposition~\ref{prop:example_b}. Nevertheless, Proposition~\ref{prop:ncc_implies_R} shows that some such coefficients still belong to $\mathcal{R}(\Omega)$, a fact that we make use of in the following result.

Theorem \ref{thm:diffusion_for_Delta} presents the main result on the diffusion phenomena for solutions to (DWE) with Dirichlet Laplacian. We restrict our attention to two cases: $\Omega=\mathbb{R}^n$ and
$\Omega$ an exterior domain in $\mathbb{R}^n$.

Throughout the paper, an exterior domain in $\mathbb{R}^n$ is understood to be the unbounded connected component of $\mathbb{R}^n\setminus\overline{\cO}$, where
$\cO$ is a bounded, non-empty open subset of $\mathbb{R}^n$. Additional
regularity assumptions on the boundary are stated explicitly when
needed.
We say that a function $b\in L^\infty_+(\R^n)$ satisfies the \emph{averaged geometric control condition} (AGCC) (see, e.g. \cite{BuJo16}), if there exists a constant $\tau>0$ such that
\begin{equation*}\label{cond:agcc}
\inf_{\xi\in \R^n,\|\zeta\| = 1} \int_0^\tau b (\xi + t \zeta)\,\ud t >0. \tag{AGCC}
\end{equation*} 

We say that a function $b\in L^\infty_+(\Omega)$, where $\Omega$ is an exterior domain in $\R^n$, satisfies the \emph{exterior geometric control condition} (EGCC) (see, e.g. \cite{AlIbKh15}), if the following conditions hold: 
\begin{equation*}
\text{\normalfont (EGCC)}
\left\{
\begin{aligned}
 & \text{\it for every } r>0 \text{ \it there exists } t_r > 0 \text{ \it such that every generalized geodesic}\\ 
& \text{\it starting in } \Omega\cap B(0,r) \text{ \it either meets } \{\xi: b(\xi) > 0 \}\,  \text{\it at some time } t \leq t_r \\
& \text{\it or leaves } B(0,r)  \text{ \it at some time } t> t_r. 
\end{aligned}
\right. 
\end{equation*}

We refer the reader to, for instance, \cite{BLR92} for the notion of \emph{generalized geodesics} associated with (DWE); see also \cite{HorIII} for propagation of singularities and \cite{Gerard91} for the defect-measure framework.

In our next main result, we provide sufficient conditions on the damping coefficient $b$ in (DWE) under which the assumptions of Theorems~\ref{thm:mild_DPE}--\ref{thm: diffusion} are satisfied for $\theta  = \frac12$, and hence the diffusion phenomenon occurs. We recall that for $\theta =\frac12$ we have $Y^\theta_D = X^1_D$, $Z_D^\theta = X$, and $\wt Z_D^\theta = A^{-1}\big(BX\big)$; see \eqref{def_YD2} and \eqref{def_wtZD}. 

\begin{theorem}\label{thm:diffusion_for_Delta}
\noindent\emph{(i)} Assume that $b\in C^1(\R^n) \cap L^\infty_+(\R^n)$ satisfies \emph{(AGCC)} and $\|\nabla b\| \in   L^\infty(\R^n)$. 
Then the conclusions of Theorems \ref{thm:mild_DPE}, \ref{thm: diffusion intro}, and \ref{thm: diffusion} with $\theta = \frac12$ hold for $A = -\Delta$, $B = M_b$, $X = L^2(\R^n)$, $V = W^{1,2}(\R^n)$. 

 In particular, for every $u_0\in W^{2,2}(\R^n)$ and every $v_0 \in W^{2,2}(\R^n)$ such that \, $ \Delta v_0 = b  w_0$ for some $w_0 \in L^2(\R^n)$, setting $u_1: = b v_0 - b u_0$, we have
\[
\|u(t; u_0,u_1) - v(t; v_0)\| \lesssim 1/t, \quad t>0,
\]
where $u(\cdot; u_0,u_1)$ and $v(\cdot;v_0)$ denote the corresponding solutions to \emph{(DWE)} and \emph{(DE)}, respectively.

Furthermore, for each non-increasing  function $r: \R_+ \to \R_+$ such that $r(t)\rightarrow 0$ and $1/t = o(r(t))$ as $t\to \infty$ there exists $v_0 \in W^{2,2}(\R^n)$ such that for every $u_0 \in W^{2,2}(\R^n)$ and $u_1: = b v_0 - b u_0$, we have 
\[
\|u(t; u_0,u_1) - v(t; v_0)\| \lesssim 1/t, \quad t >0, \quad \textrm{and} \quad \|u(t;u_0,u_1)\|\gtrsim r(t)
\]
for all sufficiently large $t>0.$
\smallskip

\noindent\emph{(ii)} 
Let $\Omega\subset \R^n$ be an exterior domain with smooth boundary.  Let $b \in C^1({\Omega})\cap L^\infty_+(\Omega)$ satisfy \emph{(EGCC)}, $\|\nabla b\| \in   L^\infty(\Omega)$, and $\mathrm{ess\,inf}_{\|\xi\| > r_0} b(\xi) > 0$ for some  $r_0 > 0$.
 
Then the conclusions of Theorems \ref{thm:mild_DPE}, \ref{thm: diffusion intro}, and \ref{thm: diffusion} with $\theta =\frac12$ hold for $A = -\Delta$, $B = M_b$, $X = L^2(\Omega)$, $V = W^{1,2}_0(\Omega)$, and, in particular, statements analogous to those in \emph{(i)} (with $\Dom(A) = W^{2,2}(\Omega)\cap W^{1,2}_0(\Omega)$) also hold.
\end{theorem}

In the specific context considered in Theorem \ref{thm:diffusion_for_Delta}, the construction of $v_0$ given in Corollary \ref{cor:explicit_did} and based on Weyl sequences corresponding to $0 \in \sigma(-\Delta)$ leads to explicit \emph{diffusive} initial data for both (DE)  and (DWE). However, the resulting diffusion phenomenon described in Corollary \ref{cor:explicit_did} is somewhat weaker than that provided by parts (i) and (ii) of Theorem \ref{thm: diffusion}. It gives explicit initial data for which the parabolic component dominates the remainder along suitable sequences of times, rather than yielding the full asymptotic relation.

Although the conditions on the initial data in
Theorems \ref{thm 1} and \ref{thm: diffusion intro} may appear elaborate, they
admit useful concrete interpretations. In the
classical case $B=I$, with $A$ non-negative and
selfadjoint, one can take $\theta=1$, and all
three auxiliary spaces in \eqref{def_YD2} and \eqref{def_wtZD}
coincide with $X$. This includes, in particular, the
constant-damping framework of Chill and Haraux \cite{Chill03}.
These and other results go well beyond the
framework of commuting operator coefficients.
In particular, for the Dirichlet Laplacian
they apply to spatially varying, possibly
vanishing damping coefficients, for which
the multiplication operator $B=M_b$ generally
does not commute with $A=-\Delta$; see
Theorem \ref{thm:diffusion_for_Delta} below. In this setting, one can take
$\theta=1/2$, and the corresponding spaces
of initial data admit explicit descriptions.

\section{Diffusion phenomena for (DWE)  with Dirichlet Laplacian} \label{sec:examples}

In this section, we develop the program outlined in  Subsection~\ref{subs:applications} .
We begin with preliminaries on the Dirichlet Laplacian $\Delta = \Delta^\Omega_D$ in $L^2(\Omega)$.
Here and throughout, $\Omega$ denotes a non-empty open subset of $\mathbb{R}^n.$ If, in addition, $\Omega$ is connected, then we say that $\Omega$ is a domain. In general, we do not impose any regularity assumptions on the boundary of $\Omega$.
In the subsequent subsections we study classes $\cB_D^{\theta}(\Omega)$, $\cB(\Omega)$,  and $\cR(\Omega)$, following the description given in~(s1)--(s3) of Subsection~\ref{subs:applications}. To preserve the flow of exposition, we postpone the proofs of some technical results until Section \ref{sec:proofs_of_thm1.6}; for instance, Subsection~\ref{subs:diffusion_for_Delta} contains the proof of Theorem~\ref{thm:diffusion_for_Delta}.

\smallskip

\paragraph*{\bf Model example: the Dirichlet Laplacian on Euclidean domains} 
For a~given $\Omega$,  
recall that  $W^{1,2}_0(\Omega)$ is equipped with the standard $W^{1,2}$-norm given by 
\[
\|u\|^2_{W^{1,2}(\Omega)} : = \int_\Omega |u|^2 \ud \xi + \int_\Omega \|\nabla u\|^2 \ud \xi, \qquad u\in W^{1,2}(\Omega).
\]
Clearly, this completion can be identified with a subspace of the Sobolev space $W^{1,2}(\Omega)$. Let $\Delta = \Delta_D^\Omega$ denote the Dirichlet realization of the distributional Laplace operator $\Delta_{\cD'(\Omega)}$ in $L^2(\Omega)$, that is, 
\[
\Delta u = \Delta_{\cD'(\Omega)} u \qquad u\in \Dom(\Delta): = \{ u \in W^{1,2}_0(\Omega): \Delta_{\cD'(\Omega)} u \in L^2(\Omega) \}.
\]
Also,  $\Delta$ admits an equivalent definition via the associated sesquilinear form as the Friedrichs extension; see, e.g.\ \cite{Grubb09}.
The spectral properties of $\Delta$, which depend essentially on the geometry
of the underlying domain $\Omega$, are well studied; see, e.g., 
\cite{EdEv18}.  For instance, it is well known that $\Delta$ is the generator of a cosine function on $L^2(\Omega)$; see, e.g. \cite[Example 3.14.1]{ABHN11}. If $\Omega$ is of class $C^{1,1}$,  then 
\[
\Dom(\Delta) = W^{2,2}(\Omega)\cap W^{1,2}_0(\Omega),
\]
see, e.g., \cite[Section 9.6]{GiTr01}.

\subsection{Class ${\mathcal{B}_D^\theta}(\Omega)$}\label{subs:ex_b_for_0} 
We begin by presenting Proposition~\ref{prop:B_D_gives_DE}, which establishes the connection between condition~\eqref{hyp:bound_2} and estimate~\eqref{DE-parab_est} appearing in our abstract Theorem~\ref{thm:mild_DPE}. The proof is postponed to Section~\ref{sec:proofs_of_thm1.6}. The remaining results in this subsection provide conditions ensuring that
$b \in \cB_D^{\theta}(\Omega)$ or $b \notin \cB_D^{\theta}(\Omega)$.

For $\theta\in(0,1]$, we let 
\[
H_\theta := \left\{ \la\in \C\setminus\{0\}:  \Re \la \geq - \min(| \Im \la|, |\Im \la|^{\theta} )/(2[R_D B]_\theta) \right\}.
\]

\begin{proposition}\label{prop:B_D_gives_DE}
Let $\theta \in (0,1]$ and $b \in L^\infty_+(\Omega)\setminus \{0\}$. Suppose that $A = -\Delta$ and $B = M_b$ satisfy condition \eqref{hyp:bound_2} in Hypothesis~\ref{hyp:main}. Then 
\begin{equation}\label{cond:almost_bound_2}
\begin{gathered}
H_\theta \subset \rho(P_D) \quad \textrm{ and } \quad \sup_{\lambda\in H_\theta} \min(|\lambda|,|\lambda|^\theta) \|R_D(\lambda)M_b\|<\infty.
\end{gathered}
\end{equation}
In particular, when $\theta = 1$, then for each $b \in L^\infty(\Omega)$ condition \eqref{hyp:bound_2} in Hypothesis~\ref{hyp:main} implies  
condition \eqref{DE-parab_est} in Hypothesis \ref{hyp:main_2}, 
 and thus the conclusions of Theorems \ref{thm:mild_DPE}  and \ref{thm: mild to DE} hold for $A = -\Delta$, $B = M_b$, and $X = L^2(\Omega)$.
\end{proposition}
Using condition~\eqref{cond:almost_bound_2}, one can define a differentiable semigroup associated with the operators $A$ and $B$ in a way analogous to the construction of the semigroup $\cT$ in~\eqref{formula:semig}, by suitably modifying the contour $\Gamma$ lying within $H_\theta$; cf.~\cite[Chapter~3]{FaYa99}. Since the class of coefficients $b$ naturally arising in applications (e.g., in Theorem~\ref{thm: pencil est}) satisfies the stronger condition~\eqref{hyp:bound_2}, we restrict our attention to the holomorphic semigroups $\cT$.

The next Theorem \ref{thm: pencil est} should be read in connection with Proposition \ref{ex: of b}. The latter provides a large class of functions $b$ and open sets $\Omega$ for which the assumptions of Theorem \ref{thm: pencil est} are satisfied. 
We say that an open (non-empty) set $U\subset \R^n$ has the \emph{Poincar\'e property}, if there exists a constant $c_P = c_P(U)\geq 0$ such that
\[
\|f\|_{L^2(U)}\leq c_P \|\nabla f\|_{L^2(U)}\qquad \textrm{ for all } f\in W_0^{1,2}(U).
\]

\begin{theorem}\label{thm: pencil est} Let $\Omega$ be an open subset of $\mathbb{R}^n$. Let $U\subseteq \Omega$ be an open set having the Poincar\'e property. Let $U'\subset U$ be such that there exists a function $\chi \in W^{1,\infty}(\Omega)$ satisfying $\supp \chi \subset U$, $0\leq \chi \leq 1$, and $\chi(\xi) \equiv 1$ on $U'$. If $b \in L^\infty_+(\Omega)\setminus \{0\}$ is such that 
\begin{equation}\label{cond: zeros of b} 
\operatorname*{ess\,inf}_{\xi \in \Omega \setminus U'} b(\xi) > 0
\end{equation} then $b\in \cB_D^{1/2}(\Omega)$. More precisely, 
the following assertions hold: 
\begin{itemize}
\item [(i)]  $\rho (P_D) \supset \C\setminus (-\infty, 0]$ and for any $\omega\in (0, \pi)$ there exists $\varepsilon_\omega>0$ such that 
\[
\rho(P) \supset \{\la \in \C\setminus \{0\}: |\arg \la|\leq \omega, |\la|\leq \varepsilon_\omega\};
\]

\smallskip

\item [(ii)] for any $\omega \in (0, \pi)$, there exists $\varepsilon = \varepsilon_\omega > 0$ such that 
\begin{equation}\label{est: 1/t at 0}
\sup_{|\arg \la |\leq \omega \atop 0<|\la|\leq 1} \|\la R_D(\la)\|<\infty \qquad \textrm{and} \qquad \sup_{|\arg \la |\leq \omega \atop 0<|\la|\leq \varepsilon} \|\la R(\la)\|<\infty,
\end{equation}
\begin{equation}\label{est: 1/sqrt at infty}
\sup_{|\arg \la |\leq \omega \atop |\la|\geq 1} \left\|  R_D(\la)\right\|<\infty \qquad\textrm{ and }\qquad \sup_{|\arg \la |\leq \omega \atop |\la|\geq 1} \left\| \sqrt{|\la|} R_D(\la) M_b \right\|<\infty.
\end{equation}
\end{itemize}
\end{theorem}

The proof of Theorem~\ref{thm: pencil est}, preceded by the auxiliary Lemmas  \ref{lem near 0} and \ref{lem: est via form}, is postponed to Section \ref{sec:proofs_of_thm1.6}. Lemma \ref{lem near 0} also shows that the first bound in \eqref{est: 1/t at 0} yields the second one. 
Regarding \eqref{est: 1/sqrt at infty}, we point out that the first bound is not required in the definition of $\cB_D^{\theta}(\Omega)$. Nevertheless, this additional estimate for $R_D$ allows us to strengthen our conclusions; see the final statement of Theorem~\ref{thm 1}(i). The following remark shows that the first bound in \eqref{est: 1/sqrt at infty} is, in general, optimal, that is, we do not expect $\|R_D(\la)\|$ decays to zero as $|\la|\to\infty$.

\begin{remark} 
Assume that there exists a nonempty open set $G\subset\Omega$ such that $b=0$ a.e.\ on $G$. Pick $\varphi\in C_c^\infty(G)$, $\varphi \neq 0$. Then $\varphi\in \Dom(\Delta)$ and, since there are no non-zero compactly supported smooth harmonic functions, we have $\Delta \varphi \neq 0$.
Therefore, since $ M_b \varphi = 0$, for every $\lambda\in \rho(P_D)$, we have
\[
\|R_D(\la)\| \ \ge\ {\|\varphi\|_2}/{\|(-\Delta+\lambda M_b)\varphi\|_2}
= {\|\varphi\|_2}/{\|\Delta \varphi\|_2} >0.
\]
Thus $\| R_D(\la)\|$ cannot decay to zero as $|\lambda|\to\infty$ and $\la \in \rho(P_D)$ although for $b$ as in Theorem \ref{thm: pencil est} we have  $\|R_D(\la)M_b\|$ does decay to zero by the second estimate in \eqref{est: 1/sqrt at infty}.
$\hfill\Diamond$\end{remark}

It was observed by Favini and Yagi in \cite[Example 3.3, p.76]{FaYa99} for the case when $\Omega$ is bounded with the smooth boundary that assuming more smoothness and some order of vanishing of $b$, $\|R_D(is)M_b\|$ decays faster than $1/\sqrt{|s|}$ as $|s|\rightarrow \infty$. 
With the help of the first bound in \eqref{est: 1/sqrt at infty}, the argument from \cite{FaYa99} leading to such faster decay is adapted in the next proposition to our more general setting described in Theorem \ref{thm: pencil est}.

\begin{proposition}\label{cond: 1/t at infty} Let $\Omega$ be an arbitrary open set and assume that $b\in L^\infty_+(\Omega)\setminus \{0\}$ is such that for some $\omega\in (0,\pi)$
\begin{equation}\label{asum:RD_bound}
\{\la \in \C\setminus\{0\}: |\arg \la|\leq \omega\} \subset \rho(P_D)\quad \textrm{and}\quad  \sup_{|\arg \la|\leq \omega \atop |\la|\geq 1} \|R_D(\la)\|<\infty.
\end{equation}
 In addition, assume that for some $\theta \in [\frac12,1]$ and $c >0$
\begin{equation}\label{cond: b'<b}
b\in W^{1,\infty}(\Omega) \quad\textrm{ and }\quad \|\nabla b (\xi)\|\leq c\, b(\xi)^{2- 1/\theta} \quad \textrm{for a.e. } \xi\in \Omega.
\end{equation}
Then
\begin{equation}\label{est: RD theta}
\sup_{|\arg \la| \leq \omega \atop  |\la|\geq 1} |\la|^\theta \|R_D(\la)M_b\| < \infty.
\end{equation}
In particular, every function $b \in L^\infty_+(\Omega)$ satisfying the assumptions of Theorem~\ref{thm: pencil est} together with condition~\eqref{cond: b'<b} belongs to $\mathcal{B}_D^\theta(\Omega)$. 
 
\end{proposition}

The proof of Proposition \ref{cond: 1/t at infty} is postponed to Section \ref{sec:proofs_of_thm1.6}.

In the formulation of Theorem \ref{thm: pencil est} we do not exclude the possibility that $\Omega$ has the Poincar\'e property; then one can take $U=U'=\Omega$ and $\chi\equiv 1$. Consequently, in this case we view condition \eqref{cond: zeros of b} as automatically satisfied and hence the conclusion of the theorem follows  for every nonzero $b \in L^\infty_+(\Omega)$. A more interesting situation arises when $\Omega$ does \emph{not} possess the Poincar\'e property, and we treat this case in the remaining part of this subsection. First, we formulate a direct consequence of Theorem~\ref{thm: pencil est} and Proposition~\ref{cond: 1/t at infty} for later reference.

We recall that a domain $\Omega \subset \mathbb{R}^n$, is called a \emph{waveguide} if there exist an integer $n_0 \in \mathbb{N}$, $n_0 < n$, and an open, bounded, connected subset $\Omega_0 \subset \mathbb{R}^{n-n_0}$ such that $\Omega = \mathbb{R}^{n_0} \times \Omega_0.$

\begin{corollary}\label{cor:examples_Omega_b} Let $\Omega$ be a non-empty, open subset of $\R^n$. 
\begin{itemize} 
\item [(i)] (waveguides) If $\Omega$ has the Poincar\'e property, then every nonzero function $b\in L_+^\infty(\Omega)$ is in $\cB_D^{1/2}(\Omega)$. In particular, if $\Omega$ is a waveguide in $\R^n$, then this conclusion holds.

\item [(ii)] (exterior domains) If $\Omega = \R^n\setminus \overline{\cO}$ for an open and bounded $\cO \subset \R^n$, then each function $b\in L_+^\infty(\Omega)$ such that $\mathrm{ess\,inf}_{\|\xi\|\geq N} b(\xi) >0$ for some $N>0$ is in $\cB_D^{1/2}(\Omega)$. 
\end{itemize}
If $b$ in assertions~\emph{(i)} and \emph{(ii)} satisfies \eqref{cond: b'<b} with $\theta \in(\frac12, 1]$, then $b \in \cB_D^{\theta}(\Omega)$. 
\end{corollary}

\begin{proof}
(i) The first statement in (i) follows from Theorem \ref{thm: pencil est} by taking $U = U' = \Omega$ and $\chi \equiv 1$.
For the second statement, we observe that each waveguide is of finite width, and thus has the Poincaré property; see, e.g. \cite[Theorem 5.3.22]{EdEv18}. 

For (ii), every bounded open set has the Poincar\'e property. Moreover, one can easily construct  $\chi \in C^\infty(\R^n)$ as needed in Theorem \ref{thm: pencil est} for $U := \{\|\xi\| < 2N\} \cap \Omega$ and $U' := \{\|\xi\| < N\} \cap \Omega$. Hence, the claim follows directly from Theorem \ref{thm: pencil est} and the assumption on $b$.

The additional assertion in both cases (i) and (ii) follows from Theorem~\ref{thm: pencil est} and Proposition~\ref{cond: 1/t at infty}. 
\end{proof} 

We state the following simple complementary to Corollary~\ref{cor:examples_Omega_b} result.
\begin{proposition}\label{prop:b_positive_everywhere} Let $\Omega$ be a non-empty open set of $\R^n$. 
Then every function $b\in L^\infty_+(\Omega)$ that is bounded away from zero
belongs to $\mathcal{B}(\Omega)\cap\cB_D^{1}(\Omega)$. More precisely,  if $b\in L^\infty_+(\Omega)$ and $b\ge b_0$ a.e.\ on $\Omega$ for a constant $b_0>0$, then $i\R\setminus \{0\} \subset \rho(P) \cap \rho(P_D)$ and for every $s\neq 0$ we have
\begin{equation}\label{cond:b_positive}
\|is R_D(is) M_b\|/\|b\|_\infty \leq \|is R_D(is)\|\le 1/b_0 \quad \textrm{ and } \quad \|is R(is)\|\le 1/b_0.
\end{equation}
\end{proposition}
\begin{proof}
For $s\neq 0$ and $u\in \Dom(P_D(is)) =  \Dom(\Delta)$ we have
\[
|s|b_0\|u\|_2^2 \le |s|\langle bu,u\rangle = |\Im\langle P_D(is)u,u\rangle|
\le \|P_D (is)u\|_2\,\|u\|_2.
\]
Hence, $\|P_D(is)u\|_2\ge |s|b_0\|u\|_2$. The same estimate holds for the adjoint operator $(P_D(is))^* =-\Delta-is M_b $. 
Thus $\Ran (P_D(is))$ is closed, $P_D(is)$ and $P_D(is)^*$ are injective. Therefore $\Ran (P_D(is)) = L^2(\Omega)$ and $P_D(is)^{-1} \in \cL(L^2(\Omega))$, proving $i\R\setminus\{0\}\subset \rho(P_D)$.  Since $\|M_b\|_{\cL(L^2(\Omega))} = \|b\|_\infty$ we obtain the first estimate in \eqref{cond:b_positive} with the stated bound. For $P(is)$ we have $\Im\langle P(is)u,u\rangle = s \langle bu,u\rangle$ and $P(is)^* = P(-is)$, so the same argument applies yielding the second inequality in~\eqref{cond:b_positive}.
\end{proof}

Theorem~\ref{thm: pencil est} shows that whether a damping coefficient $b \in L^\infty_+(\Omega)$ is positive or not on a Poincar\'e subset $U$ of $\Omega$ 
is irrelevant to the fact that $b \in \cB_D^{1/2}(\Omega)$. 
Assumption \eqref{cond: zeros of b} on positivity of $b$ on a complement of $U'$ is somehow essential for the spectral inclusion in (i) of Theorem \ref{thm: pencil est} to hold. Indeed, as our next result, Proposition \ref{prop:empty_resolvent_set_PD}, shows the resolvent set $\rho(P_D)$ of $P_D$ may be even empty provided  the damping coefficient $b$ decays along a sequence of supports of functions $\varphi_j$ forming a \emph{Weyl sequence} associated with $0\in \sigma(-\Delta)$.

We recall that $\Delta$ is selfadjoint, and so $\sigma(-\Delta) = \sigma_{\mathrm{ap}}(-\Delta)$. 
For $\mu \in \sigma_{\mathrm{ap}}(-\Delta)$, a sequence $(\varphi_j) \subset \Dom(\Delta)$ satisfying 
\begin{equation}\label{Weyl_seq_Delta}
\|\varphi_j\|_{2} = 1,\, j\in \N,  \quad \textrm{ and } \quad \|(\mu + \Delta) \varphi_j\|_{2} \to 0 \quad \text{as } j \to \infty,
\end{equation}
is called the \emph{Weyl sequence for} $\mu \in \sigma_{\mathrm{ap}}(-\Delta)$.  The functions $\varphi_j$, $j\in\mathbb{N}$, are called \emph{approximate eigenfunctions} corresponding to $\mu$. 
Clearly, $\Omega$ does not have the Poincar\'e property if and only if 
\[
0 = \inf\left\{ \|\nabla \varphi\|_{2}:  \varphi \in W^{1,2}_0(\Omega) \textrm{ with } \|\varphi \|_2 = 1 \right\},
\] or, equivalently, $0\in \sigma_{\mathrm{ap}}(- \Delta)$.

\begin{proposition}\label{prop:empty_resolvent_set_PD} Let $\Omega$ be an open subset of $\R^n$ and let $b\in L_+^\infty(\Omega)$.
If for some $\mu \in \sigma_{\mathrm{ap}}(-\Delta)$ there exists a Weyl sequence $(\varphi_j)$ such that 
\begin{equation}\label{cond:inf_b=0}
\liminf_{j\to\infty}\operatorname*{ess\,sup}_{\xi\in\operatorname{supp}(\varphi_j)} b(\xi)=0.
\end{equation}
then $0 \in \sigma_{\mathrm{ap}}(-\mu + \la M_b - \Delta)$ for all $\la\in \C$. In particular,  $\pm \sqrt{\mu} i  \notin \rho(P)$. 

If, in addition, $\Omega$ does not have the Poincar\'e property and \eqref{cond:inf_b=0} holds for a Weyl sequence $(\varphi_j)$ for $\mu = 0\in \sigma_{\mathrm{ap}}(-\Delta)$, then $\rho(P_D) = \emptyset$ and $0\notin \rho(P)$. 
\end{proposition}

For the proof, it suffices to observe that $\varphi_j$, $j \in \mathbb{N}$, satisfying \eqref{cond:inf_b=0} are approximate eigenfunctions of $-\mu + \lambda M_b - \Delta$.

\begin{remark}\label{rem:2.8}  Weyl sequences can be constructed explicitly for many classes of unbounded domains $\Omega$.
For our purposes, we recall the construction of Weyl sequences for domains $\Omega$ that are not quasi-cylindrical.
Following \cite[p.~184]{AdFo03} we call 
an unbounded domain $\Omega \subset \mathbb{R}^n$ \emph{quasi-cylindrical} if
\begin{equation}\label{def:quasicylind}
\limsup_{\xi \in \Omega,\ \|\xi\| \to \infty} \mathrm{dist} (\xi, \partial\Omega) < \infty.
\end{equation}
Clearly, every \emph{quasi-bounded} domain $\Omega$ 
(i.e., such that $\lim_{\xi \in \Omega, \|\xi\| \to \infty} \mathrm{dist}(\xi, \partial\Omega) = 0$) is quasi-cylindrical, as is every (unbounded) domain of \emph{finite width} (e.g., waveguides). 
If $\Omega$ is not quasi-cylindrical, then for any $\mu \in [0,\infty) = \sigma_{\textrm{ess}}(-\Delta)$ and for any sequence of balls $B_j \subset \Omega$, $j\in \N$, with radii $r_j\rightarrow \infty,$ one can construct a corresponding Weyl sequence $(\varphi_j)$ for $\mu \in \sigma_{\mathrm{ess}}(- \Delta)$ such that $
\mathrm{supp}(\varphi_j) \subseteq B_j$, $j\in \N$. For instance, one can take
\[
\varphi_j (\xi) : = r_j^{ - n/2} \varphi((\xi-\xi_j)/r_j)e^{i \sqrt{\mu} \, \omega \xi}, 
\]
where $\omega\in \R^n$ with $\|\omega\|=1$,  $\varphi\in C^\infty(\R^n)$ with $\mathrm{supp}(\varphi) \subseteq B(0,1)$ and $\|\varphi\|_2=1$, 
$\xi_j \in \Omega$ and $0< r_j <\mathrm{dist}(\xi_j, \partial \Omega)$ with $r_j\rightarrow \infty$. Indeed, a direct computation gives 
\[
\|\varphi_j\|_2 = 1 \quad \textrm{ and } \quad  \|( \mu +\Delta)\varphi_j\|_2 \leq 2 \sqrt{\mu} \, r_j^{-1} \|\nabla \varphi\|_2 +  r^{-2}_j \|\Delta \varphi\|_2, \quad j\in \N. 
\]
In particular, taking $\mu =0$ yields that if $\Omega$ is not quasi-cylindrical then it does not have the Poincar\'e property.
For a broader account of the relation between Weyl sequences for Dirichlet Laplacians and the geometry of the underlying domains, see Edmunds and Evans \cite[Chapters~VIII--X]{EdEv18}.
$\hfill\Diamond$\end{remark}

By combining the construction described in Remark \ref{rem:2.8} with Proposition~\ref{prop:empty_resolvent_set_PD}, we obtain the following result relating positivity of the damping coefficient $b$ to the resolvent sets of the associated operator pencils $P$ and $P_D$.

\begin{proposition}\label{prop:imaginary_exis_in_specturm} 
Let $\Omega$ be a non-quasi-cylindrical domain in $\R^n$, that is, 
\eqref{def:quasicylind} fails to hold. Let $b \in L^\infty_+(\Omega)$ be such that for a sequence $(r_j)\subset (0,\infty)$ with $r_j\rightarrow\infty$ and a sequence of points $(\xi_j)$ in $\Omega$ such that $B_j: = B(\xi_j,r_j)\subset \Omega$ one has
\[
\sup_{\xi \in B_j}  b(\xi) \rightarrow 0 \qquad \textrm{ as } j\rightarrow \infty.
\] 
Then 
\[
i\R \cap \rho(P) = \emptyset \quad \textrm{ and } \quad \rho(P_D) = \emptyset.
\]
In particular, $b\notin \cB(\Omega)\cup  \cB_D^{1/2}(\Omega)$.
\end{proposition}

In order to ensure that condition $\Ran(M_b)\setminus \Ran(\Delta)\neq \emptyset$ in Theorem~\ref{thm:mild_DPE}(iii) holds, or equivalently $b\in \cR(\Omega)$, we are mainly interested in domains $\Omega$ that do not satisfy the Poincar\'e property and, in particular, are not quasi-cylindrical. In the next result, we discuss the existence of \emph{large}  subsets $U$ with the Poincar\'e property in such $\Omega$ allowing to consider damping coefficients $b$ that are positive on a small portion of $\Omega$, while still yielding $b\in \cB_D^{1/2}(\Omega)$. 
A~characterization of the open sets $U \subset \mathbb{R}^n$ for which the Poincar\'e inequality holds can be found, for instance, in \cite[Theorem~2.10, Chapter~8]{EdEv18}. 
In applying Theorem~\ref{thm: pencil est} to a given $b$ and a set $U \subset \Omega$ with the Poincar\'e property, one must additionally ensure that, for a sufficiently large subset $U' \subset U$, there exists $\chi\in W^{1,\infty}(\Omega)$ such that $0\le\chi\le1$, $\supp\chi\subset U$, and $\chi=1$ on $U'$, while \eqref{cond: zeros of b} holds. This we obtain by imposing an additional regularity condition on $b$; see Proposition \ref{ex: of b}(ii).

In the next result, relying on Theorem~\ref{thm: pencil est} and on a sufficient geometric condition ensuring the Poincar\'e property (see  \cite[Theorem 6.33]{AdFo03}), we specify classes of unbounded open sets $\Omega$ and functions $b$ for which $b \in \cB_D^{1/2}(\Omega)$. 
Recall that for every measurable set $E\subseteq\mathbb{R}^n$, we denote its Lebesgue measure by $|E|$.

\begin{proposition}\label{ex: of b} Let $\Omega$ be an unbounded, open subset of $\R^n$ and let $b \in L^\infty_+(\Omega)$.
Suppose that there exist  $r>0$, $l\in \N$, and $\varepsilon >0$  such that for any lattice point $p \in (r\ZZ)^n$ if $(p+[0,r]^n) \cap \Omega \cap \{ \|\xi\|\geq rl\} \neq \emptyset$ then 
\begin{equation}\label{meas}
\left| \{\xi\in \Omega\cap (p + [0,r]^n):  b(\xi)\geq \varepsilon\} \right|\geq 2\varepsilon r^n.
\end{equation}
Then the following assertions hold:
\begin{itemize}
\item [(i)] There exists an open set $U\subset \Omega$ which posses the Poincar\'e property and such that 
\[
\{\xi\in \Omega: b(\xi)<\varepsilon \} \subset U\quad \textrm{and} \quad \operatorname*{ess\,inf}_{\xi\in \Omega \setminus U} b(\xi)\geq \varepsilon.
\]

\item [(ii)] If, in addition, $b\in BUC(\Omega)$, then the conclusions of
Theorem~\ref{thm: pencil est} hold for $b$; in particular,
$b\in\cB_D^{1/2}(\Omega)$. 
\smallskip

\item [(iii)] If, in addition,  $b\in BUC(\Omega)$ and \eqref{cond: b'<b} holds for some $\theta \in (\frac12,1]$, then the conclusions of Theorem \ref{thm: pencil est} hold for $b$ and $b\in \cB_D^{\theta}(\Omega)$.
\end{itemize}
\end{proposition}

We postpone the proof of Proposition~\ref{ex: of b} to Section \ref{sec:proofs_of_thm1.6}. As an immediate consequence, we obtain the following result. 

\begin{corollary}[$\Omega = \R^n$]\label{cor:examples_b}
For every lattice periodic, continuous function $b_{\textrm{per}}\in L^\infty_+(\R^n)\setminus \{0\}$ and for every function $b_\infty\in L^\infty_+(\R^n)$ that vanishes at infinity, the conclusions of Theorem \ref{thm: pencil est} hold for $b: = b_{\textrm{per}} + b_\infty$, and so $b\in \cB_D^{1/2}(\R^n)$. If, in addition, $b$ satisfies \eqref{cond: b'<b} with $\theta \in(\frac12,1]$, then $b\in  \cB_D^{\theta}(\R^n)$.
\end{corollary}

\begin{proof} Since $b_{\mathrm{per}}$ is continuous,
non-zero and lattice-periodic, the condition \eqref{meas}
holds for $b_{\mathrm{per}}$ with suitable parameters.
Indeed, there are $a,\rho>0$ and $z\in\R^n$ such that
$b_{\mathrm{per}}\ge a$ on $B(z,\rho)$, hence on every lattice translate
of this ball. If $r_L$ is a covering radius of the period lattice and
$h>2(r_L+\rho)$, then each cube $p+[0,h]^n$ with $p\in(h\mathbb Z)^n$
contains such a ball. Thus \eqref{meas} holds with $r=h$ and any
$0<\varepsilon\le\min\{a,|B(0,\rho)|/(2h^n)\}$ (and any $l\ge1$).
Hence Proposition~\ref{ex: of b}(ii) shows that the hypotheses
of Theorem~\ref{thm: pencil est} are satisfied for
$b_{\mathrm{per}}$. They remain satisfied for
$b=b_{\mathrm{per}}+b_\infty$, because the only hypothesis
of Theorem~\ref{thm: pencil est} involving the damping
coefficient is the lower bound \eqref{cond: zeros of b}, and
$b\ge b_{\mathrm{per}}$. The final assertion follows from
Proposition~\ref{cond: 1/t at infty}. 
\end{proof}

\subsection{Class $\cB_\infty(\Omega)$} \label{subs:ex_b_for_infty} 
In this section, we provide conditions under which domains $\Omega \subset \mathbb{R}^n$ and damping coefficients $b \in L^\infty_+(\Omega)$ are such that $b\in \cB_\infty(\Omega)$; see \eqref{def:B}. We restrict our attention to two cases: the whole space $\Omega = \mathbb{R}^n$ and the exterior domain 
$\Omega = \mathbb{R}^n \setminus \overline{\cO}$, for $\cO$ open and bounded. Even in such cases, the study of sufficient conditions for $b$ to belong to $\cB_\infty (\Omega)$ is a delicate matter, and so we just adapt some results from the literature. In the case of the whole space we refer to a result of Burq and Joly \cite{BuJo16}, and in the case of exterior domains to a result of Aloui, Ibrahim, and Khenissi \cite{AlIbKh15}.
Our supplementary result, Proposition~\ref{prop:interm_freq}, which provides conditions on a general domain $\Omega$ and on $b$ guaranteeing that $i\mathbb{R}\setminus \{0\}\subset \rho(P)$, appears to be new in this context. 

We begin with a result that allows us to transfer the estimates~\eqref{pencils est} for the resolvent pencil from one damping coefficient to another. Its proof is postponed to Section~\ref{sec:proofs_of_thm1.6}.

\begin{proposition}\label{prop:a_la_ui} \emph{(i)} Let $\Omega$ be an open subset of $\R^n$. Let $b_*, b\in L^\infty_+(\Omega)\setminus \{0\}$ be such that for some constant $c_*>0$ we have $b_*(\xi)\leq c_* \, b(\xi)$ for a.e. $\xi\in \Omega$. Let $s\in \R$.

If there exists a constant $c_s>0$ such that
\begin{equation}\label{ineq:b_*}
|s| \|\varphi\|_2 \leq c_s \|(-s^2 + is b_* - \Delta) \varphi\|_2 \quad \textrm{ for every }\varphi \in \Dom(\Delta),
\end{equation}
 then 
\begin{equation}\label{ineq:b}
|s| \|\varphi\|_2 \leq  \wt c_s \| (-s^2 + is b -\Delta)\varphi\|_2\quad \textrm{ for every }\varphi \in \Dom(\Delta),
\end{equation} where  $\wt c_s : = 4 \max\big( c_s, \,c_s^2 |b\|_\infty(1+c_*)^2 \big)$. 
Consequently, if $b_*\in \cB(\Omega)$ (resp. $b_*\in \cB_\infty(\Omega)$), then $b\in \cB(\Omega)$ (resp. $b\in \cB_\infty(\Omega)$).
\smallskip

\emph{(ii)} 
Let $\Omega\subset\mathbb R^n$ be an exterior domain and let $b\in C(\Omega)\cap L^\infty_+(\Omega)$. Assume that there exist constant $b_0$ and $r_0$ such that $b(\xi)\geq b_0 >0$ for all $\xi \in \Omega$ with $\|\xi\|\geq r_0$.  

Then there exist $r_*>r_0$, $c_*>0$,  and a function $b_* \in C^\infty(\Omega)\cap L^\infty_+(\Omega)$ such that $b_*(\xi) = 1$ when $\|\xi\|\geq r_*$, $b_*\leq c_* b$ on $\Omega$, and  $\{\xi\in\Omega: b_*(\xi)>0\} = \{\xi\in\Omega:b(\xi)>0\}$.
Consequently, $b_*$ satisfies \emph{(EGCC)}
if and only if $b$ does.
\end{proposition}

The next result is essentially contained in \cite[Proposition 3.3]{BuJo16} (statement (i)) and in \cite[Propositions 7 and 8]{AlIbKh15} (statement (ii)). We refer the reader to Subsection~\ref{subs:applications} for the definitions of conditions (AGCC) and (EGCC).

\begin{lemma} \label{lem:B_infty_class}  
\noindent \emph{(i)}
 Every non-negative function $b \in \mathrm{BUC}(\mathbb{R}^n)$ satisfying condition \emph{(AGCC)} belongs to $\cB_\infty(\mathbb{R}^n)$.
\smallskip

\noindent \emph{(ii)} Let $\Omega$ be an exterior domain in $\R^n$ with smooth boundary. Every $b \in L_+^\infty(\Omega)\cap C(\Omega)$ satisfying condition \emph{(EGCC)} and such that $\mathrm{inf}_{\|\xi\| > r_0} b(\xi) > 0$ for some  $r_0 > 0$ belongs to $\cB_\infty(\Omega)$.
\end{lemma}

\begin{proof} The proof can be reduced to the arguments in~\cite{BuJo16}~and~\cite{AlIbKh15}.

\noindent (i) For $c\in L^\infty_+(\mathbb R^n)$ and $t\in\mathbb R$, put
\[
 \widetilde P_c(it)=-t^2+itM_c+I-\Delta.
\]
Suppose that $c$ is uniformly continuous and satisfies (AGCC).
By \cite[Theorem~1.1]{BuJo16}, the semigroup
generated on $H^1(\mathbb R^n)\times L^2(\mathbb R^n)$ by
\[
 \mathcal A_c=
 \begin{bmatrix}0&I\\ \Delta-I&-M_c\end{bmatrix},
 \qquad D(\mathcal A_c)=H^2(\mathbb R^n)\times H^1(\mathbb R^n)
\]
is exponentially stable. Hence its resolvent is uniformly
bounded on $i\mathbb R$. Solving
$(it-\mathcal A_c)(u,v)=(0,f)$ gives $v=itu$ and
$\widetilde P_c(it)u=f$. It follows that
$\widetilde P_c(it)$ is invertible on $L^2(\mathbb R^n)$
for every $t\in\mathbb R$, with
\[
 \sup_{t\in\mathbb R}
 \|it\widetilde P_c(it)^{-1}\|<\infty.
\]

Applying this to $c=b$, set
\[
 [\widetilde R]:=
 \sup_{s\in\mathbb R}\|is\widetilde P_b(is)^{-1}\|<\infty.
\]
Since $P(is)=\widetilde P_b(is)-I$, we have
$\|\widetilde P_b(is)^{-1}\|\le [\widetilde R]/|s|$.
For $|s|\ge 2[\widetilde R]$, the Neumann series gives
\[
 P(is)^{-1}=\bigl(I-\widetilde P_b(is)^{-1}\bigr)^{-1}
               \widetilde P_b(is)^{-1},
 \qquad \|isP(is)^{-1}\|\le 2[\widetilde R].
\]

To establish invertibility for every $s\ne0$, set
\[
 t_s=\operatorname{sgn}(s)\sqrt{1+s^2},
 \qquad b_s=\frac{|s|}{\sqrt{1+s^2}}\,b.
\]
The coefficient $b_s$ is a positive scalar multiple of $b$,
so it is uniformly continuous and satisfies (AGCC).
The preceding consequence of \cite[Theorem~1.1]{BuJo16},
applied to $c=b_s$, shows that
$\widetilde P_{b_s}(it_s)$ is invertible. Since
\[
 \widetilde P_{b_s}(it_s)
 =-t_s^2+it_sM_{b_s}+I-\Delta
 =-s^2+isM_b-\Delta=P(is),
\]
we obtain $i\mathbb R\setminus\{0\}\subset\rho(P)$.
The map $s\mapsto P(is)^{-1}$ is norm-continuous for
$s\ne0$, and hence bounded on every compact subset
of $\mathbb R\setminus\{0\}$. Together with the
high-frequency estimate above, this yields
\[
 \sup_{|s|\ge1}\|isR(is)\|<\infty,
\]
which proves (i).

(ii) Under the additional assumptions that the damping coefficient $b$ is smooth and satisfies $b(\xi)=1$ for all sufficiently large $\|\xi\|$, the proof of the inclusion $b\in\cB_\infty(\Omega)$ is given in~\cite{AlIbKh15}; see \cite[Theorem~1]{AlIbKh15} for the standing assumptions in that paper and \cite[Propositions~7 and 8]{AlIbKh15} for the \emph{a priori} estimates that readily imply $b\in\cB_\infty(\Omega)$.
For a general damping coefficient $b$ satisfying the weaker assumptions imposed in~(ii), Proposition~\ref{prop:a_la_ui}(ii) shows the existence of a function $b_*\in C^\infty(\Omega)\cap L^\infty_+(\Omega)$ satisfying condition~(EGCC), such that $b_*(\xi)=1$ for all sufficiently large $\|\xi\|$, and $b_*\leq c_*b$ on $\Omega$ for some constant $c_*>0$. Then,  $b_*\in\cB_\infty(\Omega)$ by \cite[Theorem~1]{AlIbKh15}, and   Proposition~\ref{prop:a_la_ui}(i) yields $b\in\cB_\infty(\Omega)$. 
\end{proof}

\begin{remark}\label{rem:ncc} 
Following \cite{BuJo16} (see Theorem 1.3 therein), we say that a function $b \in L_+^\infty(\mathbb{R}^n)$ 
 satisfies the \emph{network control condition} (NCC) 
if there exist constants $b_0, \nu, r > 0$ and a sequence $(\xi_j) \subset \mathbb{R}^n$ such that
\[
b(\xi) \ge b_0  \, \textrm{ for } \xi\in  \bigcup_j B(\xi_j; r)
\quad
\text{and}
\quad
\mathrm{dist}(\xi,\{\xi_j\}) < \nu \, \textrm{ for all } \, \xi \in \mathbb{R}^n \;\; . \tag{NCC}
\]
Under this condition, it was proved in \cite{BuJo16} that for any $s \in \mathbb{R}$ the operator $- s^2 + i s M_b - \Delta + I = P(is) + I$
is invertible on $L^2(\Omega)$; see the remark on the low-frequency estimates in \emph{Third step} of \cite[Section~2]{BuJo16}. The argument presented therein also shows that for all $b \in L^\infty_+(\mathbb{R}^n)$ satisfying (NCC), we have $i\mathbb{R} \setminus \{0\} \subset \rho(P)$. 
$\hfill\Diamond$\end{remark}

\begin{corollary}\label{cor:agcc_implies_thm21} Let $b\in BUC(\R^n)$ be a non-negative function satisfying \emph{(AGCC)}. Then the conclusions of
Theorem~\ref{thm: pencil est} hold for $b$, and consequently $b\in \cB(\R^n) \cap \cB_D^{1/2}(\R^n)$. If, in addition, $b$ satisfies \eqref{cond: b'<b} with some $\theta \in (\frac12,1]$, then $b\in \cB(\R^n) \cap \cB_D^{\theta}(\R^n)$.
\end{corollary}
\begin{proof}
 By Lemma~\ref{lem:B_infty_class}(i)
we obtain $b\in\mathcal B_\infty(\mathbb R^n)$.
It is readily seen that (AGCC) for a function
$b\in\mathrm{BUC}(\mathbb R^n)$ implies (NCC).
In turn (NCC) implies the assumptions of
Proposition~\ref{ex: of b}. Consequently, by part (ii)
of Proposition~\ref{ex: of b},
$b\in\mathcal B_D^{1/2}(\mathbb R^n)$ and
$b\in\mathcal B(\mathbb R^n)$. The additional statement
$b\in\mathcal B_D^\theta(\mathbb R^n)$ follows from
part (iii) of Proposition~\ref{ex: of b}. 
\end{proof}

We conclude this subsection with a result providing a condition on a function $b \in L^\infty_+(\Omega)$ that guarantees that $i\mathbb{R} \setminus \{0\} \subset \rho(P)$.  
We say that an open set $U \subset \mathbb{R}^n$ has the \emph{Rellich property} if the embedding $W_0^{1,2}(U) \hookrightarrow L^2(U)$ is compact. 
For example, every open set $U\subset\mathbb R^n$ of finite Lebesgue
measure has the Rellich property, even if $U$ is unbounded;
see \cite{Lair1976}. Another sufficient condition, due to Clark
\cite{Clark1966}, is the existence of
positive functions $d$,
$\delta$ on $\mathbb R_+$, and a constant $M>0$ such that
$d(R)+\delta(R)\to0$ as $R\to\infty$, $d(R)\leq M\delta(R)$,
and, for every $x\in U$ with $\|x\|>R$, there is $y\in\mathbb R^n$
with $\|x-y\|<d(R)$ and $B(y,\delta(R))\cap U=\varnothing$.
Quasi-boundedness is necessary for the Rellich property on unbounded
open sets, but is not sufficient when $n\geq2$; see \cite{Adams1970}.
For a characterization of the Rellich property in terms of capacity,
one may consult \cite[Theorem~3.1, Chapter~8]{EdEv18}.

\begin{proposition}\label{prop:interm_freq} Let $\Omega$ be a domain in $\mathbb{R}^n$.  Let $U\subset\Omega$ have the Rellich property and $|\Omega\setminus U|>0$.  Let $U'$ be any subset of $U$ such that there exists  a function $\chi \in W^{1,\infty}(\Omega)$ such that $\supp \chi \subset U$, $0\leq \chi \leq 1$,  and $\chi(\xi) \equiv 1$ on $U'$.
Let $b\in L^\infty_+(\Omega)$ be such that $\mathrm{ess\,inf}_{\xi \in \Omega\setminus U'} b(\xi) >0$. Then $i\R\setminus\{0\}\subset\rho(P)$.
\end{proposition}

The proof of Proposition~\ref{prop:interm_freq}, preceded by the auxiliary Lemma \ref{lem:interm_freq}, is postponed until Section \ref{sec:proofs_of_thm1.6}. We note that Proposition~\ref{prop:interm_freq} provides an alternative fixed frequency criterion under a Rellich localization of the damped region. Such a localization does not follow from (NCC) in general: for distributed or periodic network damping, the relevant weakly damped region need not have the Rellich property. Consequently, Proposition~\ref{prop:interm_freq} and Lemma~\ref{lem:B_infty_class} cover different, though overlapping, geometric classes.

\subsection{Class $\cR(\Omega)$} \label{subs:ran_cond} 
 In this subsection, we discuss conditions under which the inclusion $b \in \mathcal{R}(\Omega)$ holds, where $\cR(\Omega)$ is defined in \eqref{def:R_class}. The following result describes two types of domains $\Omega$ for which the relation between $\Ran(\Delta)$ and $\Ran(M_b)$ is partially controlled. In particular, parts (ii) and (iii) of our next result show that, even in the case of $\Omega = \mathbb{R}^n$, $n\in\mathbb{N}$, the structure of $\mathcal{R}(\mathbb{R}^n)$ depends on the dimension $n$. The case $n\geq 5$  serves as a motivation for the subsequent results in this subsection.
 
\begin{proposition}\label{prop:RninR} 
The following assertions hold. 
\begin{itemize}
\item [(i)] If $\Omega\subset \R^n$ has the Poincar\'e property, then $\cR(\Omega) = \emptyset$.
\smallskip

\item [(ii)] If $n\leq 4$, then $\cR(\R^n) = L^\infty_+(\R^n)\setminus \{0\}$. 
\smallskip

\item [(iii)] If $n\geq 5$, then $\cR(\R^n) \cap \big(L^\infty_+(\R^n)\cap L^q(\R^n)\big)  = \emptyset$ for any $1\leq q < n/2$, and so $\cR(\R^n) \neq  L^\infty_+(\R^n)\setminus \{0\} $. 
\end{itemize}
\end{proposition} 

\begin{proof} (i) If $\Omega$ has the Poincar'e property, then $0\in \rho(\Delta)$, and hence $\Ran(\Delta)=L^2(\Omega)\supset \Ran(M_b)$. Therefore, for every $b\in L^\infty_+(\Omega)$ we have $b\notin \mathcal{R}(\Omega)$. 
\smallskip

(ii) Let $b\in L^\infty_+(\mathbb{R}^n)$, $b\neq 0$. Then there exists a bounded measurable set $G\subset\mathbb{R}^n$ of positive measure such that $\mathrm{ess,inf}_{\xi\in G} b(\xi)>0$.
Consequently, $L^2(G)$, identified in the usual way with a subspace of $L^2(\mathbb{R}^n)$, is contained in $\Ran(M_b)$. 
We claim that every function $g\in L^2(G)\subset L^1(\mathbb R^n)$ satisfying
\[
\int_{\mathbb R^n} g(\xi)\,\ud \xi=(\cF g)(0)\neq 0
\]
does not belong to $\Ran(\Delta)$ whenever $n\leq 4$. Hence, $\Ran(M_b)\not\subset \Ran(\Delta)$, which shows that $b\in\mathcal{R}(\mathbb{R}^n)$.  To prove the claim, suppose, seeking a contradiction, that there exists $\varphi\in\Dom(\Delta)$ such that
$\Delta\varphi=g$. Then $-\|\xi\|^2(\cF\varphi)(\xi)=(\cF g)(\xi)$ for a.e. $\xi\in\mathbb R^n$. A direct calculation using the assumption
$n\leq 4$ shows that
\[
\left|\frac{1}{|B(0,r)|}\int_{B(0,r)}(\cF g)(\xi)\,\ud \xi\right|
\leq \frac{1}{|B(0,r)|} \left(\int_{B(0,r)}\|\xi\|^4\,\ud \xi\right)^{1/2}
\|\cF\varphi\|_{L^2(B(0,r))} \to 0
\]
as  $r\to 0+$. Since $\cF g$ is continuous at $0$, we obtain a contradiction.  

(iii) Since
\[
L^\infty(\R^n)\cap L^{\widetilde q}(\R^n)
\subset
L^\infty(\R^n)\cap L^{q}(\R^n)
\qquad
\text{whenever }1\leq \widetilde q\leq q\leq\infty,
\]
we may, without loss of generality, assume that $b\in L^\infty_+(\R^n)\cap L^q(\R^n)$ for some $2\leq q<n/2$. Let $\varphi = b\psi$ for some $\psi\in L^2(\R^n)$.
Then,  $b \psi \in L^p(\R^n)\cap L^2(\R^n)$ for any $\psi \in L^2(\R^n)$, where $p := 2q/(q+2)$. We fix  $\psi\in L^2(\R^n)$ and set $\varphi := b\psi$. We claim that $\varphi \in \Ran(\Delta)$. Then, $\Ran(M_b)\subset \Ran(\Delta)$, 
and consequently $b\notin\cR(\mathbb R^n)$.

To show the claim, we first note that the function $h := \|\cdot\|^{-2}\cF\varphi$ belongs to $L^2(\R^n)$. Indeed, since 
$|h(\xi)|\leq \max(\|\xi\|^{-2},1)|(\cF\varphi)(\xi)|$, $\xi \in \R^n$, 
it suffices to show that $\|\cdot\|^{-2} \cF\varphi \in L^2(B(0,1))$. Since $p   \in [1,  2n/(n+4))\subset [1,2]$, the Hausdorff--Young theorem yields   $\cF\varphi \in L^{p'}(\R^n)$. Moreover, since $1/p' + 1/q = 1/2$ and ${|\cdot|^{-2}}~\in L^q(B(0,1))$,  H\"older's inequality implies that $h\in L^2(\R^n)$, as desired. Let $f : = -\cF^{-1}h$. Since $f$, $\cF f=-h$, and ${\|\cdot\|^2}\cF f=-\cF\varphi$ all belong to $L^2(\R^n)$, the Fourier characterization of $W^{2,2}(\R^n)$ yields $f\in W^{2,2}(\R^n) = \Dom(\Delta)$.
Moreover, $\mathcal F(\Delta f)=\|\cdot\|^2 h=\mathcal F\varphi$ yields $\Delta f=\varphi$, as desired.
\end{proof}

Since we do not know whether, for a general domain $\Omega\subset\mathbb{R}^n$ without the Poincar\'e property, the inclusion $\cB_D^{1/2}(\Omega)\subset\cR(\Omega)$ holds, we provide a sufficient condition ensuring that a function $b$ belongs to $\cR(\Omega)$.  
More precisely, our main result in this direction, Proposition \ref{prop:correction_b}, shows that for domains $\Omega$ without the Poincar\'e property, countable-valued functions of the form $\sum_j c_j \mathbf{1}_{\mathrm{supp}(\varphi_j)}$,
built on the supports of functions from a Weyl sequence $(\varphi_j)$ for $0 \in \sigma\bigl((-\Delta)^{1/2}\bigr)$, belong to $\mathcal{R}(\Omega)$, and any damping coefficient $b$ that majorizes such a countable-valued function also belongs to $\mathcal{R}(\Omega)$. In particular, this majorization is valid for the coefficients $b$ which are bounded away from zero at infinity.

In what follows we consider (unbounded) domains $\Omega$ without the Poincar\'e property. This is equivalent to the condition
$0 = \inf \sigma (-\Delta)$. Moreover, we denote by 
$
\mathcal{W}(\Omega) := \mathcal{W}\bigl(\Omega; (-\Delta)^{1/2}\bigr)
$
the family of Weyl sequences associated with 
$0 \in \sigma \bigl((-\Delta)^{1/2}\bigr)$, 
this is, $(\varphi_j) \in \mathcal{W}(\Omega)$ if and only if
\[
(\varphi_j) \subset W^{1,2}_0(\Omega), \quad \|\varphi_j\|_{2} = 1  
\quad \text{and} \quad \|\nabla \varphi_j\|_{2} \to 0 \quad \text{as } j \to \infty.
\]
Clearly, if $0\in \sigma(-\Delta)$ then $\cW(\Omega)$ is nonempty and closed under taking subsequences; that is, 
if $(\varphi_j) \in \mathcal{W}(\Omega)$, then for any subsequence $(j_\ell) \subset \mathbb{N}$ we also have $(\varphi_{j_\ell}) \in \mathcal{W}(\Omega)$, and so there exists $(j_\ell) \subset \N$ such that $(\|\nabla \varphi_{j_\ell}\|_{2}) \in \ell^2(\mathbb{N})$.
For a set $E\subseteq\mathbb{R}^n$, we denote its characteristic
function by $\mathbf{1}_E$.
\begin{proposition}\label{prop:correction_b} Let $\Omega \subset \R^n$ be an open set without the Poincar\'e property. For any Weyl sequence $(\varphi_j)\in \cW(\Omega)$ such that supports of $\varphi_j$, $j\in \N$, are mutually disjoint, there exists a sequence 
$(c_j) \subset [0,\infty)$ with $c_j \to 0$ as $j\to \infty$ such that any function $b\in L^\infty_+(\Omega)$ with $b(\xi) \geq \sum_j c_j \mathbf{1}_{\mathrm{supp}(\varphi_j)}(\xi)$ for a.e. $\xi\in \Omega$ belongs to $\cR(\Omega)$.
\end{proposition}
It is worth noting that, for certain Weyl sequences, the sequence
$(c_j)$ appearing in Proposition~\ref{prop:correction_b} can be
constructed explicitly; see, e.g., the concluding discussion in the paragraph
{\emph{Cubes with ball-shaped holes}} in this subsection.
However, as Proposition \ref{prop:RninR}(iii) shows, in the case $\Omega=\mathbb{R}^n$ with $n\geq 5$, the function $\sum_j c_j \mathbf{1}_{\operatorname{supp}(\varphi_j)}$ does not belong to $L^q(\mathbb{R}^n)$ for any $q\in [1,n/2)$.

\begin{proof}
Passing to a subsequence if necessary, we may assume that 
$(\|\nabla \varphi_j\|_{2}) \in \ell^2(\mathbb{N})$. 
Then there exists a positive sequence $(\eta_j)\in \ell^2(\mathbb{N})$ such that $\eta_j/\|\nabla \varphi_j\|_2 \rightarrow \infty$ as $j\rightarrow \infty$. Indeed, let $(j_\ell)_{\ell=0}^\infty$ be an increasing sequence of integers such that
$
\sum_{j=j_\ell+1}^{j_{\ell+1}} \|\nabla \varphi_j\|_2^2 \leq 8^{-\ell}$, $\ell \in \mathbb{N}$. 
It is sufficient to define $\eta_j := 1$ for $j\in \{1,\dots,j_1\}$,
and $\eta_j := 2^\ell \|\nabla \varphi_j\|_2$ for $j=j_\ell+1,\dots,j_{\ell+1}, \ \ell\in\mathbb{N}$.

Let $\varphi: = \sum_j \eta_j \varphi_j$. Since the supports of $\varphi_j$ are mutually disjoint, $\varphi \in L^2(\Omega)$.  We claim that $\varphi\notin \Ran (\Delta)$. Assume, on the contrary, that there exists $\psi\in \Dom(\Delta)\subset W^{1,2}_0(\Omega)$ such that $\Delta \psi = \varphi$. Then,  for any $j\in \N$ we have
\[
\eta_j = \ip{\eta_j \varphi_j}{\varphi_j} = \ip{\varphi}{\varphi_j} = \ip{\Delta \psi}{\varphi_j} \leq \|\nabla \psi\|_2 \|\nabla \varphi_j\|_2.
\]
However, this contradicts the choice of $\eta_j$, proving the claim.

Furthermore, there exists a sequence $(c_j)\subset (0,1)$ such that $c_j \to 0$ as $j\to\infty$ and $(\eta_j/c_j)\in \ell^2(\mathbb{N})$. Indeed, let $c_j: = 1/2$ when $j=1, \dots, j_1$, and $c_j:=2^{-\ell/4}$ when $j\in \{j_{\ell}+1, \dots , j_{\ell +1}\}$, $\ell \in \N$.
Then, for every $\ell\in\mathbb N$,
\[
 \sum_{j=j_\ell+1}^{j_{\ell+1}}\left(\frac{\eta_j}{c_j}\right)^2
=2^{5\ell/2}\sum_{j=j_\ell+1}^{j_{\ell+1}}\|\nabla\varphi_j\|_2^2
 \le 2^{-\ell/2}.
\]
Hence $(\eta_j/c_j)\in\ell^2(\mathbb N)$ as desired. 
Let $b\in L^\infty_+(\Omega)$ with $b \geq \sum_j c_j \mathbf 1_j$, where $\mathbf 1_j$ stands for the characteristic function of $\mathrm{supp\,} (\varphi_j)$, $j\in \N$. 
Moreover, we define $\wt \psi(\xi): = \varphi(\xi)/b(\xi)$ for all $\xi$ with $b(\xi)\neq 0$ and $\wt \psi(\xi): = 0$ otherwise. Clearly, $\wt \psi \in L^2(\Omega)$ and $M_{b} \wt \psi = \varphi$. Therefore, $b\in \cR(\Omega)$. \end{proof}

Combining Persson's theorem~\cite[Theorem~2.1]{Persson60} with Proposition~\ref{prop:correction_b}, we readily obtain the following corollary. We omit a direct proof.

\begin{corollary}\label{cor:correction_b} Let $\Omega \subset \R^n$ be an open set without the Poincar\'e property. 
For any $\varepsilon>0$ there exists a countable-valued function $b\in L^\infty_+(\Omega)$ such that 
\begin{itemize}
\item [(i)]  $\mathrm{supp\,} b\subset \Omega_\varepsilon := \{\xi \in \Omega: \mathrm{dist\,}(\xi, \partial \Omega)<\varepsilon\, \textrm{ or }\, \|\xi\|>1/\varepsilon\}$,
\item [(ii)] $\|b\|_\infty\leq \varepsilon$ and  $|b(\xi)|\rightarrow 0$ as $\|\xi\|\rightarrow \infty$, and
\item [(iii)] for every $\wt b \in L^\infty_+(\Omega)$, if $\wt b \geq b$ then $\wt b \in \cR(\Omega)$.
\end{itemize}
\end{corollary}

\begin{remark}\label{rem:on_correction_of_b} (a)  In the context of Theorem~\ref{thm: pencil est}, Proposition~\ref{prop:correction_b} shows, in particular, that if there exists an open subset $U \subset \Omega$ such that $\Omega \setminus \overline{U}$ does not have the Poincar\'e property and
$\mathrm{ess\,inf}_{\xi \in \Omega \setminus U} b(\xi) > 0$,
then $b\in \cR(\Omega)$.
However, there exist domains $\Omega$ without the Poincar\'e property that admit a~decomposition into two subdomains, each of which have the Poincar\'e property. More precisely, there exists an open set $U \subset \Omega$ with the Poincar\'e property such that $\Omega \setminus \overline{U}$ also has the Poincar\'e property. We refer to the example coming from damping coefficients $b$ satisfying (AGCC) condition (and so $b\in \mathcal{B}(\mathbb{R}^n)\cap \cB_D^{1/2}(\mathbb{R}^n)$) discussed in Proposition~\ref{prop:example_b}. However, more direct examples can be given to illustrate such decompositions. 
\smallskip

(b) In the study of~(DWE), it is often important to control the regions where the damping coefficient $b$ is positive. The set $\{\xi \in \Omega : b(\xi)>0\}$ is called the \emph{damping region}. 
For a given $\Omega$ and $b\in L^\infty_+(\Omega)$, the explicit construction of Weyl sequences $(\varphi_j)\in \mathcal{W}(\Omega)$ allows one to identify regions of $\Omega$ such that, if $b$ is bounded below away from zero on these regions, then $b\in \mathcal{R}(\Omega)$. 
If this is not easy to determine, that $b$ is bounded  away from zero on the supports of any Weyl sequence in $\mathcal{W}(\Omega)$, this motivates the study of well-localized \emph{corrections} $b_c$ of $b$ such that $b_c + b \in \mathcal{R}(\Omega)$ while keeping the region  $\{ \xi \in \Omega: b_c(\xi) + b(\xi)>0\}$ under control. $\hfill\Diamond$\end{remark}

In view of Propositions~\ref{prop:RninR} and \ref{prop:correction_b}, 
for the damping coefficients $b$ of the type constructed in the proof of
Proposition~\ref{prop:example_b} below (their higher-dimensional analogues with $n\geq 5$), the set $\{\xi\in\mathbb{R}^n:b(\xi)>0\}$
does not contain $\bigcup_j\operatorname{supp}(\varphi_j)$,
where $(\varphi_j)$ is a Weyl sequence associated with 
$0\in\sigma((-\Delta)^{1/2})$. Therefore, Proposition~\ref{prop:correction_b} does not allow us to infer that such a coefficient $b$ belongs to $\cR(\mathbb{R}^n)$. The following result establishes, in particular, that for each $b\in L^\infty_+(\mathbb{R}^n)$ satisfying (NCC) (see Remark~\ref{rem:ncc}) belongs to $\mathcal{R}(\mathbb{R}^n)$.

\begin{proposition}\label{prop:ncc_implies_R}
Let $\Omega\subset\R^n$ be an open set. Let $(r_j)$ be an increasing
sequence of positive numbers such that $r_j\to\infty$. Suppose that
there exists a sequence $(\xi_j^*)\subset\R^n$ such that
$B(\xi_j^*;2r_j)\subset\Omega$ for every $j\in\N$, and the balls
$B(\xi_j^*;2r_j)$ are pairwise disjoint.

Let $b\in L^\infty_+(\Omega)$  satisfy the network control condition with respect to $\Omega_*: = \bigcup_j B(\xi_j^*;  r_j)$, that is, 
there exist constants $b_0, \nu, r > 0$ and a sequence $(\xi_j) \subset \Omega$ such that $\bigcup_j B(\xi_j; r) \subset \Omega$, 
\[
b(\xi) \ge b_0  \, \textrm{ for } \xi\in  \bigcup_j B(\xi_j; r) 
\quad
\text{and}
\quad
\mathrm{dist}(\xi,\{\xi_j\}) < \nu \, \textrm{ for all } \, \xi \in \Omega_*.
\]
Then $b\in \cR(\Omega)$. 
\end{proposition}

\begin{proof} 
Passing to a subsequence if necessary, we may assume without loss of generality that $r_j>j+\sqrt{n}\,2\nu$ for all $j\in\mathbb{N}$. For each $j\in\mathbb{N}$, we let $Q_j$ denote the cube centered at $\xi_j^*$ with side length $4\nu N_j$, where $N_j :=\left\lfloor\frac{r_j}{\sqrt{n} 2\nu}\right\rfloor$, that is, $Q_j$ is the maximal cube contained in $B(\xi_j^*,r_j)$, centered at $\xi_j^*$, whose side length is a multiple of $4\nu$. Each such cube $Q_j$ contains $N_j^n$ pairwise disjoint (open) cubes $Q_{j,\ell}$ of side length $4\nu$. We denote by $\widetilde{Q}_{j,\ell}$ the cube with the same center as $Q_{j,\ell}$ and side length $2\nu$.
Clearly, for all $\ell, \ell' \in \{1,... N_j^n\}$ we have
\[
\dist(\wt Q_{j,\ell}, \wt Q_{j,\ell'}) \geq 2\nu \quad \textrm{ and } \quad \wt Q_{j,\ell} \cap \{\xi_k\}_{k\in \N} \neq \emptyset.
\]

For each $j\in \N$ and $\ell\in\{1,\ldots,N_j^n\}$ we choose a point 
$\xi_{j,\ell}\in\widetilde Q_{j,\ell} \cap\{\xi_k:k\in\mathbb N\}.$
We choose $0\leq \varphi\in C^\infty(\mathbb{R}^n)\setminus\{0\}$ with $\supp(\varphi) \subset B(0;r_*)$, where $r_*: =\min(r,\nu)$. Let 
\[
a_j: = j^{-1} N_j^{- n/2} \quad \textrm{ and } \quad f(\xi) : = \sum_{j\in \N} a_j \sum_{\ell =1}^{N_j^n} \varphi(\xi - \xi_{j, \ell}), \quad \xi \in \Omega.  
\]

Since $B(\xi_{j,\ell};\nu)\subset Q_{j,\ell}$, the supports of $\varphi(\cdot-\xi_{j,\ell})$, with $j\in \N$ and $\ell\in\{1,\ldots,N_j^n\}$, are pairwise disjoint.
Hence, $
\|f\|_2^2=\|\varphi\|_2^2\sum_{j\in\mathbb{N}}j^{-2}<\infty$, so that $f\in L^2(\mathbb{R}^n)$. Moreover, since $b(\xi)\geq b_0$ for all $\xi\in\supp(f)$, it follows that $f\in\Ran(M_b)$.

We claim that $f\notin \Ran(\Delta)$. Suppose, for contradiction, that there exists $g\in \Dom(\Delta)$ such that $\Delta g = f$. Choose $0\leq \psi \in C^\infty(\R^n)$ such that $\psi(\xi) = 1$ when $\|\xi\|\leq 1$ and $\psi(\xi) = 0$ when $\|\xi\|>2$.  
We set
\begin{equation*}\label{def:packet_family}
\psi_j(\xi):=r_j^{-n/2}\psi\!\left({(\xi-\xi_j^*)}/{r_j}\right), \quad \xi \in \Omega, \, j\in \mathbb{N}.
\end{equation*}
It is readily seen that $\|\psi_j\|_2 = \|\psi\|_2$ and $\|\Delta\psi_j\|_2 =  r^{-2}_j \|\Delta \psi\|_2$, $j\in \N$. 
Therefore, for all $j\in \N$, 
\begin{equation}\label{eq:est_with_p}
|\ip{f}{\psi_{j}}| = |\ip{g}{\Delta \psi_{j} }| \leq r_j^{-2} \|g\|_2 \|\Delta \psi\|_2 \lesssim r_j^{-2} \lesssim j^{-2}. 
\end{equation}
On the other hand, since  $f(\xi)\psi_j(\xi) = r_j^{-n/2}a_j \sum_{\ell =1}^{N_j^n} \varphi(\xi - \xi_{j, \ell})$ for all $\xi \in \Omega$, 
\[
\ip{f}{\psi_j} = r_j^{- n/2} a_j   N_j^n \int \varphi \, \ud \xi \gtrsim j^{-1}.
\]
This contradicts \eqref{eq:est_with_p} for sufficiently large $j$. Hence,
$f\in\Ran(M_b)\setminus\Ran(\Delta)$, and therefore $b\in\cR(\Omega)$.
\end{proof}

Since, as already pointed out at the beginning of the proof of Corollary~\ref{cor:agcc_implies_thm21}, {(AGCC)} implies {(NCC)} for functions $b\in BUC(\mathbb{R}^n)$ (see Subsection~\ref{subs:applications} and Remark~\ref{rem:ncc} for the definitions), we obtain the following immediate consequence of Proposition~\ref{prop:ncc_implies_R}.

\begin{corollary}\label{cor:agcc_implies_r}
For every $n\in \N$ and every $b\in L^\infty_+(\R^n)\cap BUC(\R^n)$, if $b$ satisfies \emph{(AGCC)}, then $b\in \cR(\R^n)$. 
\end{corollary}

For the application of Proposition~\ref{prop:correction_b}, it is
sufficient to identify sequences of domains
$G_j\subset\Omega$, $j\in\mathbb N$, for which there exists a Weyl
sequence $(\varphi_j)$ satisfying
$\operatorname{supp}\varphi_j\subset G_j$ for every
$j\in\mathbb N$. In particular, no explicit representation of such 
Weyl sequences is required. Our next result, Lemma \ref{lem:steiner_formula}, shows that the construction of the Weyl sequences associated with $0 \in \sigma_{\mathrm{ap}}\bigl((-\Delta)^{1/2}\bigr)$ is related to the existence of a sequence of bounded open sets $G_j \subset \Omega$ such that 
${|\partial G_j|}/{|G_j|} \to 0$ as $j \to \infty$. Here, for any open set $G \subset \mathbb{R}^n$, we write $|\partial G|$ to denote the $(n-1)$-dimensional Hausdorff measure of its boundary $\partial G$.
For simplicity, we restrict ourselves to bounded Lipschitz subdomains of $\Omega$. 
For a sequence $(G_j)$ of bounded Lipschitz subdomains
of $\Omega$ with uniform Lipschitz character $(r,L)$
and $|\partial G_j|/|G_j|\to0$, we show how to
construct functions $\varphi_j\in W_0^{1,2}(G_j)$,
$j\in\mathbb N$, such that $\|\varphi_j\|_2=1$
 and $\|\nabla\varphi_j\|_2^2
\leq C {|\partial G_j|}/{|G_j|}$ for all $j$, where $C>0$ is independent of $j$; see Lemma~\ref{lem:steiner_formula}. Proposition~\ref{cor:diff_copy} below allows us to provide examples of sequences $(G_j)$ consisting of sets other than balls with radii tending to infinity, which can be used to determine whether 
$b\in\mathcal{R}(\Omega)$.

\smallskip

Below, we use a standard description of bounded Lipschitz
domains in terms of local Lipschitz graphs in coordinate
cylinders, which is more convenient for our purposes.
For related geometric formulations, see, e.g.,
\cite[Section~2.1]{MaVe11} and \cite[Section~2]{An24}.
We call a map $\Phi:\R^n \to \R^n$ a \emph{Euclidean isometry} (or a \emph{rigid motion}) if $\Phi(\xi) = \mathbb{O}\xi + \xi^*$.
for some orthogonal matrix $\mathbb{O}$ and a vector $\xi^*\in\R^n$. 
For a vector $\xi=(\xi_j)\in\R^n$ we write $\xi=(\xi',\xi_n)$, where 
$\xi'=(\xi_1,\dots,\xi_{n-1})\in\R^{n-1}$. 
For $r>0$ we denote by $B'(0;r)$ the open ball of radius $r$ centered at $0$ in $\R^{n-1}$. Moreover, we set $C(r):=B'(0;r)\times (-r,r)$.
We say that a bounded domain $G\subset\mathbb R^n$, $n\geq 2$, is a Lipschitz domain if there exist $r>0$ and $L>0$ such that, for every $\xi\in\partial G$, there exist a Euclidean isometry $\Phi_\xi$ with $\Phi_\xi(\xi)=0$ and a
Lipschitz function $\gamma_\xi:B'(0;r)\to\mathbb R$ with
$\mathrm{Lip}(\gamma_\xi)\le L$ such that
\begin{equation*}
\begin{split}
\Phi_\xi(G)\cap C(r)
&=
\{(x',x_n)\in C(r):x_n>\gamma_\xi(x')\},\\
\Phi_\xi(\partial G)\cap C(r)
&=
\{(x',x_n)\in C(r):x_n=\gamma_\xi(x')\}.
\end{split}
\end{equation*}
In this case, we say that $G$ has {\emph{Lipschitz character $(r,L)$}}. 
For a Lipschitz domain $G$, we define
\[
\mathrm{Lip}(\partial G) := \inf \left\{ L \ge 0 : G \text{ has Lipschitz character } (r,L) \text{ for some } r>0 \right\},
\] with the convention that $\inf \emptyset := \infty$. 
A sequence $(G_j)$ of Lipschitz domains in $\mathbb R^n$ is said to have \emph{uniform Lipschitz character} if there exist constants $r>0$ and $L\geq0$ such that each $G_j$ has Lipschitz character $(r,L)$.
For an individual bounded Lipschitz domain, the existence
of uniform parameters $r$ and $L$ follows from the usual
local graph description by compactness of the boundary.
For a family of bounded domains, however, a common choice
of these parameters is an additional assumption.

\begin{lemma}\label{lem:steiner_formula} Let $G\subset \R^n$ be a bounded domain with Lipschitz character $(r,L)$. Then there exists a function $\varphi_G\in W^{1,2}_0(G)$ with $\supp\varphi_G\Subset G$ such that 
\begin{equation}\label{eq:phi_G}
\|\varphi_G\|_2 = 1 \quad \textrm{ and } \quad \|\nabla \varphi_G\|^2_2 \leq \varepsilon^{-1} c(1+L)  \, |\partial G|/ |G|, 
\end{equation}
where 
\[\varepsilon := \min\big( r/(2\max(1,L)), |G|/(c (1 + L) |\partial G)|\big) \qquad \text{and}\qquad c : = 10(4\sqrt{2} +1)^n.
\]

In particular, suppose that $(G_j)$ is a sequence of
bounded domains in $\R^n$ with uniform Lipschitz character $(r,L)$ and
$|\partial G_j|/|G_j|\rightarrow 0$ as $j \to \infty.$ 
Then, for all sufficiently large $j\in\N$, there exists a function $\varphi_j\in W_0^{1,2}(G_j)$ with $\supp\varphi_j\Subset G_j$ satisfying
\begin{equation}\label{nabla_est}
\|\varphi_j\|_2=1,
\qquad
\|\nabla\varphi_j\|_2^2
\leq {2c(1+L)^2}{r^{-1}}\,{|\partial G_j|}/{|G_j|}.
\end{equation}
\end{lemma}

\begin{proof}
As we will see, the existence of the function $\varphi_G$ is a consequence of the following boundary layer estimate for $G$: \emph{ for any $0<\varepsilon \leq r/(2 \max(L,1))$}  
\begin{equation}\label{est:layer}
\left|\{\xi \in \R^n: \mathrm{dist} (\xi, \partial G)\leq \varepsilon \} \right|\leq  5^{-1}c \varepsilon  
 (1+L) \, |\partial G|.
\end{equation}
We will now prove \eqref{est:layer}. 
Since $G$ has Lipschitz character $(r,L)$, for every $\xi\in\partial G$ there exist a Euclidean isometry $\Phi_\xi$ and a Lipschitz function $\gamma_\xi:B'(0;r)\to\mathbb{R}$
with $\operatorname{Lip}(\gamma_\xi)\leq L$. We choose a maximal $r/2$-separated family $\{\xi_j\}_{j=1}^N$ in $\partial G$ and write $\Phi_j: = \Phi_{\xi_j}$, $U_j := U_{\xi_j} = \Phi_j^{-1}\big( C(r)\big)$, $\gamma_j:= \gamma_{\xi_j}$,   
Then, we have
\[
 \partial G\subset\bigcup_{j=1}^N B(\xi_j, r/2),
 \qquad B(\xi_j,  r)\subset  U_j\subset B(\xi_j,\sqrt2\,  r).
\]
Moreover, the family $\{U_j\}_{j=1}^N$ has the overlap constant bounded by a constant depending only on $n$. Indeed, since ${\xi_j}_{j=1}^N$ is $r/2$-separated, the balls $B(\xi_j,r/4)$ are pairwise disjoint. Let $x\in \R^n$. If $x\in U_j \subset B(\xi_j,\sqrt2\,r)$, then $B(\xi_j, r/4)\subset B(x,(\sqrt{2}+1/4)  r)$. Hence, by comparing volumes, 
\begin{equation}\label{overlaping}
\#\{j:x \in  U_j\}\leq \#\{j:x\in B(\xi_j,\sqrt{2} r)\} \leq  (4\sqrt{2}+1)^n.
\end{equation}
For $0<\varepsilon\le r/(2\max(1,L))$, we denote
\[
 S_j:=\partial G\cap \overline{B(\xi_j;r/2)},\qquad 
 T_{j,\varepsilon}:= S_j + \overline{B(0;\varepsilon)}.
\]

Since $\partial G$ is compact, we obtain 
\begin{equation}\label{cover}
 \{\xi \in \R^n: \operatorname{dist}(\xi,\partial G)< \varepsilon\}=\bigcup_{j=1}^N T_{j,\varepsilon}.
\end{equation}

Fix $\xi\in T_{j,\varepsilon}$ and choose $\zeta\in S_j$
such that  $\|\xi-\zeta\|\le \varepsilon$. Since
$|\zeta-\xi_j|<r/2$ and $\varepsilon\leq r/2$, we have
\[
 |\xi-\xi_j|\leq|\xi-\zeta|+|\zeta-\xi_j|<r.
\]
Thus both $\xi$ and $\zeta$ belong to $U_j$. Write
$\Phi_j(\xi)=(x',x_n)$ and $\Phi_j(\zeta)=(y',y_n)$.
Since $\Phi_j$ is an isometry with $\Phi_j(\xi_j)=0$,
we have $\|x'\|<r$ and $|y_n|<r/2$. Moreover,
$\zeta\in\partial G\cap U_j$ implies $y_n=\gamma_j(y')$,
and hence
\[
 \|x'-y'\|^2+|x_n-\gamma_j(y')|^2
 =|\xi-\zeta|^2\leq\varepsilon^2.
\]
In particular, $\|x'-y'\|\leq\varepsilon$. The Lipschitz
bound on $\gamma_j$, together with the choice of
$\varepsilon$, therefore gives
\[
 \begin{aligned}
 |\gamma_j(x')|
 &\leq |\gamma_j(y')|+L\|x'-y'\|\\
 &<\frac r2+L\varepsilon
 \leq\frac r2+\frac{Lr}{2\max\{1,L\}}
 \leq r.
 \end{aligned}
\]
Consequently, $x'\in D_j$, where
\[
 D_j:=\{u\in B'(0;r):|\gamma_j(u)|<r\}.
\]
The defining boundary representation yields
\[
 \Phi_j(\partial G\cap U_j)
 =\mathrm{graph}(\gamma_j|_{D_j}).
\]
Furthermore,
\[
 \begin{aligned}
 |x_n-\gamma_j(x')|   &\le |x_n- y_n|+L\|x'-y'\|\\
 &\le\sqrt{1+L^2}\, \bigl(|x_n-y_n|^2+\|x'-y'\|^2\bigr)^{1/2}\\
 &\le\varepsilon\sqrt{1+L^2}.
 \end{aligned}
\]
Consequently, we infer
\[
 \Phi_j(T_{j,\varepsilon})
 \subset\left\{(x',x_n):x'\in D_j,
                    |x_n-\gamma_j(x')|\le\varepsilon\sqrt{1+L^2}\right\}.
\]
Since $|\Phi_j(T_{j,\varepsilon})| = |T_{j,\varepsilon}|$, Fubini's theorem yields 
\[
 |T_{j,\varepsilon}| \leq 2 \varepsilon\sqrt{1+L^2}\,|D_j|.
\]
Finally, by the area formula for Lipschitz graphs, we obtain

\[
 |D_j|\le\int_{D_j}\sqrt{1+\|\nabla\gamma_j(u)\|^2}\,du
     =|\partial G\cap  U_j|.
\]
Therefore, since $\sqrt{1+L^2} \leq 1+L$,
\begin{equation}\label{volume}
 |T_{j,\varepsilon}|  \leq 2\varepsilon (1+L)\, |\partial G\cap U_j|.
\end{equation}
Summing \eqref{volume}, using
\eqref{cover} and \eqref{overlaping}, gives
\[
 \begin{aligned}
 |\{x:\dist(x,\partial G)\le\varepsilon\}|
 &\leq 2\varepsilon (1+L) \sum_{j=1}^N |\partial G\cap  U_j| \leq 2 (4\sqrt 2 +1)^n (1+L)\varepsilon |\partial G|.
 \end{aligned}
\]
This proves \eqref{est:layer}. 

We can now proceed with construction of the function $\varphi_G$. Let $d(\xi):=\operatorname{dist}(\xi,\partial G)$ and let $\varepsilon$ be as in the statement. With $\wt c:=5^{-1}c(1+L)|\partial G|$ we have $\wt c \varepsilon\leq 5^{-1}|G|$. We set
\[
 \psi(\xi) :=
 \begin{cases}
  1, & d(\xi)>\varepsilon,\\
  2d(\xi)/\varepsilon-1, & \varepsilon/2<d(\xi)\leq\varepsilon,\\
  0, & d(\xi)\leq\varepsilon/2,
 \end{cases}
 \qquad \xi\in G.
\]
The function $\psi$ is $2/\varepsilon$-Lipschitz, and its support is contained in the compact set $\{\xi\in\overline G:d(\xi)\geq\varepsilon/2\}\subset G$. Thus, $\psi$ belongs to $W^{1,2}_0(G)$, and Rademacher’s theorem gives $|\nabla\psi|\leq 2/\varepsilon$ almost everywhere in $G$. Moreover, the gradient vanishes almost everywhere on $\{\xi: d(\xi)>\varepsilon\}$.
By \eqref{est:layer},
\[
 \|\nabla \psi\|_2^2
 \leq \frac{4}{\varepsilon^2}|\{\xi\in G:d(\xi)\leq\varepsilon\}|
 \leq {4 \wt c}/{\varepsilon}.
\]

For $0<t<1$, we have
$\{\xi\in G:\psi(\xi)\leq t\}=\{\xi\in G:d(\xi)\leq\varepsilon(1+t)/2\}$.
Integrating the measures of these sublevel sets gives
\[
 \begin{aligned}
 |G|-\|\psi\|_2^2
 &=\int_0^1 2t\,|\{\xi\in G:\psi(\xi)\leq t\}|\,dt\\
 &\leq \wt c\varepsilon\int_0^1 t(1+t)\,dt = 5 \wt c \varepsilon/6 \leq  |G|/6.
 \end{aligned}
\]
Thus $\|\psi\|_2^2\geq 5 |G|/6$. The normalized function
$\varphi_G=\psi/\|\psi\|_2$ has compact support in $G$ and satisfies
\[
 \|\nabla\varphi_G\|_2^2
 \leq \frac{24\wt c }{ 5 \varepsilon|G|} = 
 \frac{24}{25}\frac {c(1+L)|\partial G|}{\varepsilon |G|}.
\]
Thus, the first assertion of the lemma follows. 

To prove the second assertion, we apply the first part to obtain, for each $j\in\N$, a function $\varphi_j$ with $\supp\varphi_j\Subset G_j$ such that
\begin{equation*}\label{eq:phi_j}
\|\varphi_{j}\|_2 = 1 \quad \textrm{ and } \quad \|\nabla \varphi_j\|^2_2 \leq \varepsilon_j^{-1} c(1+L)  \, |\partial G_j|/ |G_j|, 
\end{equation*}
where \[
\varepsilon_j := \min\big( r/(2\max(1,L)), |G_j|/(c (1 + L) |\partial G_j)|\big). 
\]
Since $|G_j|/|\partial G_j|\to\infty$ as $j\to\infty$, we have $\varepsilon_j=r/(2\max(1,L))\geq r/2(1+L)$ for all sufficiently large $j$, which yields the desired estimate \eqref{nabla_est}.
\end{proof}

\noindent The following result, Proposition~\ref{cor:diff_copy}, provides a more specific version of Lemma~\ref{lem:steiner_formula} and is better suited to our purposes described in Remark \ref{rem:on_correction_of_b}(b). 
To start for $N>0$ we set 
\[
Q(N) := (0,N)^n \quad \textrm{ and } \quad G(N;K) : = Q(N)\setminus K,
\]
where $K = K (N) \subset Q(N)$ is any compact set which is the closure of an open set $K^\circ$ with Lipschitz boundary. 

For a sequence $(N_j)\subset (0,\infty)$ we say that the family $K_j$, $j\in \N$, of compact subsets removed from the cubes $Q(N_j)$ is \emph{admissible} if the corresponding family $G_j: = G(N_j,K_j)$ has uniform Lipschitz character.

For an open set $\Omega\subset\mathbb R^n$ we say that $\Omega$ contains \emph{uniformly controlled (UC) copies}
$(\widetilde G_j)$ of $(G_j)$ if there exist pairwise disjoint open sets $\widetilde G_j\subset\Omega$, open neighbourhoods
$V_j\supset\overline{\widetilde G_j}$, and $C^1$-diffeomorphisms
$F_j:V_j\to F_j(V_j)$ satisfying $F_j(\widetilde G_j)=G_j$, for which
there is a constant $c\ge1$, independent of $j$, such that
\begin{equation}\tag{2.21}\label{bi-lip_constant}
  c^{-1}\le |\det DF_j(\xi)|\le c,
  \qquad \|DF_j(\xi)\|\le c,
  \qquad \xi\in V_j,\quad j\in\mathbb N.
\end{equation}
Here and below, the norm of a matrix is the operator norm induced by
the Euclidean norm. We call such a family $(F_j)$
\emph{uniformly controlled (UC)}. No disjointness assumption is imposed on the domains $G_j$.

\begin{remark}
The bounds \eqref{bi-lip_constant} also control the The bounds \eqref{bi-lip_constant} also control the derivatives of the inverse maps. Indeed, let $s_1\ge\cdots\ge s_n>0$ denote the
singular values of $DF_j(\xi)$. Then $s_k\le c$ and
$s_1\cdots s_n=|\det DF_j(\xi)|\ge c^{-1}$. It follows that
$s_n\ge c^{-n}$, and hence
\begin{equation*}
  \|D(F_j^{-1})(F_j(\xi))\|
  =\|(DF_j(\xi))^{-1}\|
  \le c^n,
  \qquad \xi\in V_j,\quad j\in\mathbb N.
\end{equation*}
These uniform derivative bounds give local bi-Lipschitz estimates,
but they need not give uniform global bi-Lipschitz estimates on
arbitrary domains. Conversely, a uniformly globally bi-Lipschitz
family of $C^1$-diffeomorphisms satisfies
\eqref{bi-lip_constant} for a suitable $M$.
$\hfill\Diamond$\end{remark}

An important special case of such maps $F_j$ is given by Euclidean
isometries. Namely, for a sequence of orthogonal $n\times n$-matrices $\mathbb{O}_j$ and $\xi^*\in \R^n$ we set
\[
F_j(\xi) = \mathbb{O}_j \xi +  \xi_j^*,\quad   \xi \in \mathbb{R}^n.
\]
Then $F_j$ is an isometry of $\mathbb{R}^n$, and in particular
\[
DF_j(\xi) = \mathbb{O}_j,  \quad |\det DF_j(\xi)| = 1, \quad 
\mathbb{O}_j^T \mathbb{O}_j = I .
\]

More generally, for a matrix $\mathbb{A} \in \mathbb{R}^{n\times n}$ and a vector $\xi^* \in \mathbb{R}^n$ we define the affine map
\[
F_\mathbb{A} (\xi) = \mathbb{A}  \xi + \xi^*, \quad \xi \in \mathbb{R}^n.
\]
It is readily seen that for a sequence of invertible matrices $(\mathbb{A}_j)_{j\in\mathbb{N}}$  in $\mathbb{R}^{n\times n}$ and $(\xi_j^*)\subset \R^n$, the associated affine maps $F_j(\xi)=A_j\xi+\xi_j^*$ form a UC-family if and only if
\[
\sup_j \|\mathbb{A}_j\| < \infty \quad\text{and}\quad \sup_j \|\mathbb{A}_j^{-1}\| < \infty.
\]
Our next result, Proposition~\ref{cor:diff_copy}, together with
Proposition~\ref{prop:correction_b}, shows that if a domain
$\Omega$ contains a UC-copy \(\widetilde G_j\)
of a sequence $(G_j)=(G_j(N_j,K_j))$ satisfying
$|\partial G_j|/|G_j|\to0$ as $j\to\infty$, and if $(K_j)$ is
admissible, then $\sum_{j\in\mathbb N}\mathbf{1}_{\widetilde G_j}\in\mathcal R(\Omega)$. Relying on this result, one can construct several particular sequences $(G_j)$ which, under UC-diffeomorphisms, lead to Weyl sequences and, consequently, to examples of functions $b\in\mathcal R(\Omega)$; see the last three items (i), (ii), (iii) below in this subsection.

\begin{proposition}\label{cor:diff_copy}
Let $\Omega$ be an open subset of $\R^n$ that contains UC-copies $(\widetilde G_j)$ of $G_j : = G(N_j; K_j)$, $j\in\mathbb{N}$, where $K_j$, $j\in \N$, is an admissible family of compact subsets removed from the cubes $Q(N_j)$, $j\in \N$. 

Suppose that 
\[
\frac{|\partial G_j|}{|G_j|} = 
\frac{2nN_j^{n-1} + |\partial K_j|}{N_j^n - |K_j|} \rightarrow 0 \quad \textrm{ as } j\rightarrow \infty.
\] 
Then there exists a sequence $(c_j)\subset (0,\infty)$ with $c_j\rightarrow 0$ as $j\rightarrow \infty$ such that the function $\sum_j c_j \mathbf{1}_{\wt G_j}$ belongs to $\cR(\Omega)$. 
\end{proposition}

\begin{proof}
By Lemma \ref{lem:steiner_formula} for each $j\in \N$ there exists $\varphi_j = \varphi_{G_j} \in W^{1,2}_0(G_j)$ with 
\[
\|\varphi_j\|_2 = 1 \quad \textrm{ and } \quad \|\nabla \varphi_j\|_2 \rightarrow 0 \textrm{ as } j\rightarrow \infty. 
\]
Let  $\psi_j(\xi): =\varphi_j \circ F_j (\xi)$, $\xi \in  \wt G_j$, $j\in \N$, where $F_j$, $j\in\mathbb N$ is the corresponding UC-family of diffeomorphisms. Since, for any $j\in \N$, 
\[
\|\psi_j\|_{L^2(\widetilde G_j)}^2
=
\int_{G_j} |\varphi_j(\xi)|^2 |\det DF^{-1}_j(\xi)|\,\ud \xi 
\]
and
\[
\nabla\psi(\xi)=DF_j(\xi)^T\nabla\varphi_j(F_j(\xi)),
\]
Combining this with \eqref{bi-lip_constant} we infer that
\begin{equation}\label{nabla_est_tilde}
c\leq \|\psi_j\|_{L^2(\widetilde G_j)}^2 \quad \textrm{ and } \quad 
\|\nabla\psi_j\|_{L^2(\widetilde G_j)}
\le c' \,\|\nabla\varphi_j\|_{L^2(G_j)}.
\end{equation}
where $c>0$ and $c'>0$ depend only on the uniform bounds in \eqref{bi-lip_constant}, and are independent of $j$. In particular,  
 $\|\nabla\psi_j\|_2 \rightarrow 0$ as $j\rightarrow \infty$, and so $\varphi_{\wt G_j}: = \psi_j/\|\psi_j\|_2$, $j\in \N$, forms a Weyl sequence in $\cW(\Omega)$. Therefore, the required assertion follows from Proposition~\ref{prop:correction_b}.
\end{proof}

We consider several special cases of Proposition \ref{cor:diff_copy}. 

\smallskip

(i) {\bf{Cubes with ball-shaped holes.}}
Here we specialize Proposition~\ref{cor:diff_copy} to the case where the removed sets $K_j$ are unions of closed balls. To start, let $N_j>0$, $j\in \N$, be given. Let $K_j = \bigcup_{\ell = 1}^{\ell_j} K_{j,\ell}$ with $K_{j, \ell}$ being pairwise disjoint closed balls with radii $r_{j,\ell}$ contained in $Q_j = Q(N_j)$,  $\ell = 1, ..., \ell_j$, for some $\ell_j\in \N$. 
We emphasize that the number $\ell_j$ and the sequence of radii 
$(r_{j,\ell})_{\ell=1}^{\ell_j}$ depend on the parameter $N_j$ and are chosen so that $K_{j,\ell}$ are pairwise disjoint and contained in $Q_j$.

Since for any ball $B(\xi, r)$ and a cube $Q$ we have $\mathrm{Lip}(\partial B) = 0$ and $\mathrm{Lip}(\partial Q) = 1$, we infer that $\mathrm{Lip}(\partial G_j) =1$ for $G_j: = G(N_j;K_j) := Q_j\setminus K_j$. 
To ensure admissibility of $(K_j)$ one must also impose a condition of \emph{uniform geometric separation} on $K_{j,\ell}$ and $\partial Q_j$. A sufficient condition is the existence of constants $r_*, r^* >0$ such that for any $j$ and $\ell, \ell'$ 
\begin{equation}\label{cond:admissible}
r_{j,\ell} \geq r_*  \quad \textrm{ and }\quad  \dist( K_{j,\ell}, K_{j,\ell'}), \, \dist (K_{j,\ell}, \partial Q_j) \geq r^*.
\end{equation}
\begin{lemma}\label{lem:admissible}
Let $(Q_j)$ and $(K_j)$ be defined as in the preceding paragraph. Suppose that condition \eqref{cond:admissible} is satisfied. Then the family $G_j=Q_j\setminus K_j$, $j\in \N$, has uniformly Lipschitz character, that is, is admissible.
\end{lemma}
A straightforward argument shows that every point of $\partial G_j$ admits a Lipschitz chart of radius $r = \frac{1}{10}\min\{r_*, r^*\}$. We omit the details of the proof. 

Furthermore, we have
\[
|\partial G_j| = 2n N_j^{n-1} + n\omega_n \sum_\ell r^{n-1}_{j, \ell}
\quad \textrm{ and } \quad
|G_j| = N_j^n - \omega_n \sum_{\ell} r^n_{j, \ell},
\]
where $\omega_n := \pi^{n/2}\Gamma\!\left(1+n/2\right)^{-1}$ is the volume of the unite ball in $\R^n$. 
We introduce a number $\kappa_j$ that represents the fraction of the cube $Q_j$ occupied by the removed balls $K_{j,\ell}$:
\[
\kappa_j : = \kappa(N_j, \ell_j, (r_{j,\ell}) ): = N_j^{-n} \omega_n \sum_{\ell=1}^{\ell_j} r_{j,\ell}^n.
\]
Clearly, $\kappa_j \in (0,1)$ and $|G_j| = (1- \kappa_j)N_j^n$. Using H\"older's inequality, we easily obtain the following upper bound:
\[
\sum_{\ell=1}^{\ell_j} r_{j, \ell}^{n-1}  \le \Big(\sum_{\ell=1}^{\ell_j} r_{j, \ell}^n \Big)^{\frac{n-1}{n}} \, \ell_j^{1/n} 
 = \Big(\kappa_j / \omega_n \Big)^{\frac{n-1}{n}}\, \ell_j^{1/n} \,N_j^{n-1}. 
\]

Hence 
\[
|\partial G_j|/|G_j| \leq (1-\kappa_j)^{-1} \Big( 2n + n \omega_n^{1/n}  \kappa_j^{1 - 1/n} \, \ell_j^{1/n} \Big) N_j^{-1}.
\]

Therefore, if there exists $\varepsilon\in (0,1)$ such that
\begin{equation}\label{cond:N_ell_decay}
\sup_j \kappa_j\leq 1-\varepsilon \quad \textrm{ and } \quad \ell^{1/n}_j/N_j \rightarrow 0 \textrm{ as } j\to \infty
\end{equation}
then $|\partial G_j|/|G_j| \lesssim \ell_j^{1/n}N^{-1}_j\rightarrow 0$ as $j\to \infty$.
To summarize, for any choice of $N_j>0$, $\ell_j \in \mathbb{N}$, and $(r_{j,\ell})_{\ell=1}^{\ell_j}$ satisfying \eqref{cond:admissible} and \eqref{cond:N_ell_decay}, 
 every family of UC-copies $(\widetilde G_j)$ of $(G_j)$ gives rise to a function $\sum_j c_j \mathbf{1}_{\wt G_j} \in \mathcal{R}(\Omega)$; see Proposition~\ref{cor:diff_copy}.

More specifically, suppose that  $N_j = j^2$, $j\in \N$, and $\vartheta, \varepsilon \in (0, 1/2)$. 
For any $j\in \N$ we chose a family of mutually disjoint closed balls $K_{j,\ell}$, $\ell = 1, ..., \ell_j$,  in the cube $Q_j$ satisfying \eqref{cond:admissible} and suppose that $\ell_j \leq j^{\vartheta n}$ and $\kappa_j \leq 2 - \varepsilon$.
Hence, the sequence $\ell_j^{1/n} N_j^{-1} = 1/j^{2 - \vartheta}$, $j\in \N$, is in $\ell^1(\N)$. The Weyl sequence $(\varphi_{\wt G_j})$ corresponding to $(\wt G_j)$, constructed as in the proof of Proposition \ref{cor:diff_copy}, satisfies $\|\nabla \varphi_{\wt G_j}\|_2^2 \lesssim {1}/{j^{2-\vartheta}}$,  $j\in\mathbb{N}$; see \eqref{nabla_est} and \eqref{nabla_est_tilde}.
 Then, for every $\nu \in (0,1 - (\vartheta +1)/2)$, by an argument as in the proof of Proposition \ref{prop:correction_b} with a sequence $\eta_j: = j^{-\kappa}$ for some $\kappa \in (1/2+ \nu, 1 -\vartheta/2)$, we can take $c_j: = 1/j^\nu$, i.e., $\sum_j j^{-\nu} \mathbf{1}_{\wt G_j} \in \mathcal{R}(\Omega)$.
\medskip

(ii) {\bf{Balls with equal radii.}} Assuming now that all $r_{j,\ell}$ are equal, i.e., $r_{j,\ell} = r > 0$ for all $j\in \N$ and  $\ell \in \{ 1,\dots,\ell_j\}$, we have
\[
|\partial G_j| = 2n N_j^{n-1} +  n \omega_n \ell_j r^{n-1}, \qquad
|G_j| = N_j^n -  \omega_n  \ell_j r^n,\qquad  \kappa_j = \omega_n \ell_j N_j^{-n} r^n.
\]
Therefore, 
\[
r = \left(\frac{\kappa_j N_j^n}{\omega_n \ell_j}\right)^{1/n} 
\quad \textrm{ and } \quad 
 \ell_j \, r^{n-1}  
= \left(\kappa_j/\omega_n\right)^{\frac{n-1}{n}} \, \ell^{1/n}  N^{\,n-1}.
\]
Hence, 
\begin{align*}
|\partial G_j|/|G_j| 
&  =   (1-\kappa_j)^{-1}\left( 2n/\ell_j^{1/n}  + n\omega_n \left(\kappa_j\omega_n\right)^{\frac{n-1}{n}}\right)\, \ell_j^{1/n} N_j^{-1}
\end{align*}
Let $N_j\geq 2$ be integer for all $j\in \N$. Then the number of lattice points in the cube $Q_j = (0,N_j)^n$ is equal to $( N_j-1 )^n$. 
We fix $r\in (0,1/2)$. Then a sequence of unions of all $(N_j-1)^n$ closed balls of radius $r$ centred at the lattice points $\ZZ^n \cap Q_j$ is admissible; see Lemma \ref{lem:admissible}.  Clearly, 
\[
\kappa_j = \omega_n ( N_j - 1)^n N_j^{-n} r^n < \omega_n r^n < 1,
\]
so $\sup_j \kappa_j<1 $. However, for the corresponding sequence $(G_j)$ the ratios $|\partial G_j|/|G_j|$ do not converge to $0$. 
  
If we take $\ell_j =  ( \lfloor (N_j-1)/h_j \rfloor)^n$, where $(h_j)\subset \N$ is an arbitrary sequence convergent to the infinity, and $G_j$ denote $Q_j$ with removed $( \lfloor (N_j-1)/h_j \rfloor)^n$ closed balls with radius $r$ and centered at lattice points of $Q_j$, then
\[
|\partial G_j|/|G_j| \lesssim \ell_j^{1/n} N_j^{-1} \leq  1/h_j \rightarrow 0 \textrm{ as } j\to \infty.
\]  
Therefore, the conclusion of Proposition~\ref{cor:diff_copy} applies to any domain $\Omega$ containing
 UC-copies $(\widetilde G_j)$ of $(G_j)$
a uniformly bi-Lipschitz copy $(\widetilde{G}_j)$ of $(G_j)$ and thus gives rise to functions in $\mathcal{R}(\Omega)$.
\medskip

(iii) {\textbf{Rectangles instead of balls.}}
Instead of balls as in items (i) and (ii), we may remove rectangular boxes from the cube. 
Let $\Pi = \mathbb{O}\bigl(\Pi_{k=1}^n[-d_k/2, d_k/2]\bigr) \subset \R^n$, where $\mathbb{O}\in O(n)$ an orthogonal matrix, be a fixed rectangle with side lengths $d_1, ..., d_n>0$. Clearly, 
\[
|\Pi| = d_1...d_n \quad \textrm{ and }\quad |\partial \Pi| = 2\sum_{k=1}^n 
\prod_{_{i\neq k}} d_i.
\]
Let $K_j = \bigcup_{\ell = 1}^{\ell_j} (\xi_{j,\ell} + \mathbb{O}_{j, \ell} \Pi)$, where $\xi_{j, \ell} \in \R^n$ and $\mathbb{O}_{j, \ell} \in O(n)$ are such that $K_{j, \ell}: = \xi_{j,\ell} + \mathbb{O}_{j, \ell} \Pi \subset Q_j = Q(N_j)$, $\ell = 1, ..., \ell_j \in \N$, are pairwise disjoint. We set again $G_j := Q_j \setminus K_j$, $j\in \N$.
Since orthogonal maps preserve Lebesgue measure and surface measure, 
we have that $|K_j| = \ell_j |\Pi|$ and $|\partial K_j| = \ell_j |\partial \Pi|$, and consequently
\[
|G_j| = N_j^n - \ell_j|\Pi| \quad \textrm{ and } \quad |\partial G_j| = 2n N_j^{n-1} + \ell_j |\partial \Pi|.
\]
Therefore, if $\ell_j/N_j^n \to 0$ as $j\to \infty$, then
 $|\partial G_j|/|G_j| \to 0$ as $j\rightarrow \infty$. As in the case of balls, the admissibility requires a uniform separation assumption: if the rectangles have uniformly bounded aspect ratios, fixed side lengths (or more generally side lengths bounded above and below away from $0$), and stay at a fixed positive distance from one another and from $\partial Q_j$, then the family $(G_j)$ has uniform Lipschitz character. 
 
A particularly convenient special case is obtained by removing from the cubes $Q_j := (0, N_j)^n$, with $N_j \in \mathbb{N}$, $\ell_j$ translates $\xi_{j,\ell} + \Pi$, where $\xi_{j,\ell} \in {Q_j} \cap \mathbb{Z}^n$ (or even $\xi_{j,\ell} \in \overline{Q_j} \cap \mathbb{Z}^n$), of a fixed (small) rectangle $\Pi$. The resulting sets
$G_j := Q_j \setminus \bigcup_{\ell = 1}^{\ell_j} (\xi_{j,\ell} + \Pi)$, $j \in \mathbb{N}$,  have uniform Lipschitz character. Hence Proposition~\ref{cor:diff_copy} applies provided ${|\partial G_j|}/{|G_j|} \to 0$ as $j \to \infty$.

\subsection{Examples: admissible damping coefficients}
In this subsection, we present several explicit examples of damping
coefficients belonging to the classes $\mathcal{B}(\Omega)$, $\cB_D^\theta(\Omega)$, and $\cR(\Omega)$, illustrating the theory developed in the previous sections. 
\smallskip

\paragraph*{\bf Damping coefficients bounded away from zero} 
Proposition \ref{prop:b_positive_everywhere} asserts that for any open (non-empty) set $\Omega\subset\R^n$, $n\in\N$, every function $b\in L^\infty_+(\Omega)$ bounded away from zero belongs to $\cB(\Omega)\cap\cB_D^{1}(\Omega)$. Consequently, if this is the case then our abstract results, namely Theorem~\ref{thm:mild_DPE} with $\theta=1$, its refined version, Theorem~\ref{thm: mild to DE}, and Theorem~\ref{thm: diffusion intro} with $\theta=1$, apply with $A=-\Delta$, $B=M_b$, $X = L^2(\Omega)$, and $V=W^{1,2}_0(\Omega)$. Moreover, $Y= L^2(\Omega)$, and $\cT = T$ is generated by the operator $L = M_b^{-1} \Delta$.

If $\Omega$ has the Poincar\'e property (equivalently, when $0 \in \rho(-\Delta)$), then $b\notin \cR(\Omega)$, i.e., $\Ran(M_b) = \Ran(-\Delta) = L^2(\Omega)$. Moreover, the semigroup $T$ decays exponentially, i.e., $\|T(t)\|\lesssim e^{-\varepsilon t}$ for some $\varepsilon>0$ and all $t\geq 0$. Consequently, since for any solution $u$ to (DWE) and any solution $v$ to (DE), $\|u(t)\|\leq 2 \max( \|v(t)\|, \|u(t) - v(t)\|)$, the particular statements  of Theorem~\ref{thm: diffusion intro}(under corresponding conditions on the initial values $u_0$, $u_1$, and $v_0$) assert respectively that: 
(i) $\|u(t; u_0, u_1)\| = O(1/t)$; 
(ii) $\|u(t; u_0, u_1)\| = O(1/t^2)$;
(iii) $\|u'(t,u_0,u_1)\| = O(1/t^2)$ (as $t\to \infty$).

If $\Omega$ does not have the Poincar\'e property (equivalently,
$0\in\sigma(-\Delta)$), then $\Ran(-\Delta)\neq L^2(\Omega)=\Ran(M_b)$,
so that the range condition in Theorem~\ref{thm: diffusion} is
satisfied. Moreover, since $\Delta$ is injective, and hence so is
$L=M_b^{-1}\Delta$, we have $0\in \sigma(L)\setminus\sigma_{\mathrm p}(L)$.
Consequently, if $b\in W^{1,\infty}(\Omega)$, then
Theorems~\ref{thm: diffusion} and
\ref{thm:diffusion_theta=1}(i),(ii) ensure the existence of solutions
to the corresponding (DWE) whose diffusive asymptotic profiles
are solutions of (DE). Moreover, if $b\in W^{2,\infty}(\Omega)$, then
Theorem~\ref{thm:diffusion_theta=1}(iii) applies as well. 

\smallskip

\paragraph*{\bf Damping coefficients on $\R^n$ and (AGCC) condition}
Since the (AGCC) condition has a simple geometric interpretation, it can
be readily verified in many concrete examples. In particular, if $E$ is
a closed subset of the cube $(-1/2,1/2)^n$ and $|E|>0$, then the lattice-periodic function $\sum_{p\in\mathbb{Z}^n}\mathbf{1}_{E+p}$ 
does not satisfy, whereas $1-\sum_{p\in\mathbb{Z}^n}\mathbf{1}_{E+p}$
does satisfy (AGCC). Clearly, the class of functions satisfying (AGCC) is
closed under passing to larger functions, in the sense that if $b_*$ satisfies (AGCC) and $b\geq b_*$, then $b$ satisfies (AGCC) as well. 
Corollary~\ref{cor:agcc_implies_thm21} shows that the conclusions of
Theorem~\ref{thm: pencil est} hold for every damping coefficient
$b\in BUC(\mathbb{R}^n)$ satisfying (AGCC). Consequently, $b\in\mathcal{B}(\mathbb{R}^n)\cap \mathcal{B}_D^{1/2}(\mathbb{R}^n)$.

For instance, for any $\ell\in\mathbb{N}$, the function
\[
b_1(\xi):=\sum_{j=1}^{n}\sin^{2\ell}(\xi_j),
\qquad \xi=(\xi_j)_{j=1}^{n}\in\mathbb{R}^n,
\]
belongs to $\cB(\R^n)\cap \cB^{1/2}_D(\R^n)$.
Moreover, $b_1\in C^\infty(\mathbb{R}^n)$ and all of its derivatives
are bounded. By Corollary \ref{cor:agcc_implies_r}, $b_1\in \cR(\R^n)$. However, the corresponding operator $M_{b_1}$ is injective, but not invertible. 
Therefore, the conclusions of Theorem~\ref{thm:diffusion_for_Delta}(i)
hold for this choice of $b_1$. 
Regarding Proposition~\ref{cond: 1/t at infty}, one can show that
condition~\eqref{cond: b'<b} is satisfied by $b_1$ with 
$\theta=\frac{2\ell}{2\ell+1}$, whereas it fails for
$\theta>\frac{2\ell}{2\ell+1}$. Therefore, Theorems \ref{thm:mild_DPE} and \ref{thm: diffusion intro} apply with $\theta = 2\ell/(2\ell +1)$,  $A = -\Delta$, $B = M_{b_1}$, $X = L^2(\R^n)$, $V = W^{1,2}(\R^n)$. Moreover, for $\ell \geq 2$, we have $Y_D^\theta = Z_D^\theta = \wt Z_D^\theta =L^2(\R^n)$; see \eqref{Y_D}.

To construct a more degenerate damping coefficient (with a nontrivial
kernel of $M_b$), we introduce the following auxiliary functions. Let
\begin{equation}\label{5.1}
\psi(t) := \begin{cases}
0,& t\leq0,\\
e^{-1/t},&t>0,\end{cases}
   \qquad h(t) :=\frac{\psi(t)}{\psi(t)+\psi(1-t)}, \quad t\in \R.
\end{equation}
Moreover, let $q: \R\to \R$ be the $1$-periodic function defined on $[-\frac12,\frac12]$ by $q(t): =h(4|t|-1)$. Then, $q\in C^\infty(\R)$, $q(t)=0$ when $\dist(t,\mathbb Z)\leq\frac14$, and $q>0$ near $\mathbb Z+\frac12$. 
A straightforward computation shows that for every $\eta\in[0,1)$,
there exists $c_\eta>0$ such that
\begin{equation}\label{5.3}
   |q'(t)|\leq c_\eta q(t)^\eta, \qquad t\in\R.
\end{equation} 
Using Gr\"onwall's inequality, one can show that condition~\eqref{5.3} cannot hold for $\eta=1$. Let 
\begin{equation}\label{5.4}
   b_2(\xi)
   =
   \sum_{j=1}^nq(\xi_j),
   \qquad \xi\in\R^n.
\end{equation}

The corresponding operator $M_{b_2}$ is not injective; indeed, $\Ker(M_{b_2})=L^2(\mathcal Q)$ where $\mathcal Q$ is the union of all translates $(-1/4,1/4)^n+p$, $p\in\mathbb Z^n$.  
Then $b_2$ satisfies the assumptions of Theorem~\ref{thm:diffusion_for_Delta}(i).  Consequently, Theorems~\ref{thm:mild_DPE}--\ref{thm: diffusion intro} apply for every $\theta<1$. 

\smallskip

In the following proposition, we provide an example of an admissible damping coefficient $b$ whose support is contained in a union of cubes and for which the measure of the support can be prescribed within certain limits.  Although a simpler construction of a function satisfying the properties stated in the proposition is possible, we consider a damping coefficient $b$ whose support has a more intricate structure, consisting of periodically translated copies of a unit cube in the plane. As this result is only illustrative, we restrict our attention to the case $n=2$. The support of $b$ is, in a sense, complementary to that of $b_2$ modulo lattice translations (see Figure~1 below), yet it is straightforward to verify that $b$ satisfies (AGCC).

\begin{proposition}\label{prop:example_b}  Let $\varepsilon\in(0,1)$ and $r>1$ be given. For each of the two alternatives in $(\mathrm{iii})$, there exists a function \(b\in L^\infty_+(\mathbb R^2)\cap C(\mathbb R^2)\) satisfying $(\mathrm{i})$, $(\mathrm{ii})$, and the corresponding alternative.
\begin{itemize}
\item [\textup{(i)}] $b$ satisfies \emph{(AGCC)} and $b\in \cB(\R^2) \cap \cB_D^{\theta}(\R^2) \cap \cR(\R^2)$ for every $\theta \in (0,1)$. 
\smallskip

\item [(ii)] The sets $\{\xi\in \R^2: b(\xi)>0\}$ and $\R^2\setminus \mathrm{supp}(b)$ have the Poincar\'e property.  In particular, the set $\{\xi\in \R^2: b(\xi) >0\}$ does not contain the union of the supports $\mathrm{supp}(\phi_j)$ for any Weyl sequence $(\phi_j)$ associated with $0\in \sigma((-\Delta)^{1/2})$.
\smallskip

\item [\textup{(iii)}] 
One of the following estimates holds:
\[
\sup_{\xi \in \mathbb R^2}|B(\xi;r)\cap\supp(b)|
\le\varepsilon \qquad \text{or} \qquad 
\sup_{x \in \mathbb R^2}|B(\xi;r)\setminus\supp(b)|
\le\varepsilon.
\]
\end{itemize}
Moreover, $b$ can be chosen from $C^\infty(\R^2)$ with all derivatives bounded. In particular, the assumptions of Theorem  \ref{thm:diffusion_for_Delta}\emph{(i)}  are satisfied. 
\end{proposition}

\begin{proof} We elaborate the geometric setting depicted in \cite[Fig.~2]{BuJo16}, to illustrate the (AGCC) condition. For any $\eta\in(0,1/2)$, we set
\[
Q_\eta:=Q\setminus(1-2\eta)\overline{Q},
\]
where $Q:=(-1/2,1/2)^2$ and $(1-2\eta)Q$ denotes the square with the same center as $Q$ and side length $1-2\eta$. For convenience, we denote $Q_{1/2}:=Q$.
Fix $\eta\in(0,1/2]$ and $\mu,\nu\in(0,1)$, and set
\[
 \alpha :=(1+\mu,\nu),\qquad \beta:=(-\nu,1+\mu),\qquad
 \Lambda:=\{\ell\alpha+k\beta:\ell,k\in\mathbb Z\}.
\]
We write
\[
 Q_{\ell,k}:=Q_\eta+\ell\alpha+k\beta,\qquad
 \mathcal Q=\mathcal Q_\eta:=\bigcup_{\lambda\in\Lambda}(Q_\eta+\lambda).
\]
The sets $Q_{\ell,k}$ are pairwise disjoint (see Figure~1 when $\eta=1/2$).
We fix a nonnegative function $b_\eta\in C(\R^2)$ supported in $\overline Q$ such that $\{\xi: b_\eta(\xi)>0\} = Q_\eta$ and set 
\begin{equation}\label{def:b_bQ}
 b(\xi):=\sum_{\lambda\in\Lambda} b_\eta(\xi-\lambda), \qquad\qquad \xi \in \R^2.
\end{equation} 
We assume that $\mu$ and $\nu$ satisfy $\nu(1-\nu)>\mu(1+\mu)$, or equivalently
\begin{equation}\label{mu_nu}
0<\mu<\frac{\sqrt2-1}{2},
\qquad
\frac{1-\sqrt{1-4\mu-4\mu^2}}{2}
<\nu<
\frac{1+\sqrt{1-4\mu-4\mu^2}}{2}.
\end{equation}
This condition has a clear geometric interpretation: the line passing through the upper-right vertices of the cubes $Q + k\beta$, $k\in\mathbb Z$, intersects the interior of every cube $Q+ \alpha + k\beta$, $k\in\mathbb Z$; see Figure~1. It follows that, for such choices of $\mu$ and $\nu$, the function $b$ satisfies condition (AGCC). More precisely, one can show that
\[
  \inf_{a\in\R^2,\ \zeta\in S^1}\int_0^6 b(a+t\zeta)\,dt>0,
  \qquad S^1:=\{\zeta\in\R^2\colon\norm{\zeta}_2=1\}.
\]
Indeed, we first show that every segment of length $6$ meets $\cQ_{1/2}$ at an interior point. 
For $\ell,k\in\mathbb Z$, put
$P_{\ell,k}:=Q+\ell\alpha+k\beta$. Fix $a\in\R^2$ and $\zeta\in S^1$, let $z=a+3\zeta$ be the midpoint of the segment $a+[0,6]\zeta$, and choose
$\omega\in S^1$ with $\omega\cdot\zeta=0$.
We consider the closed cube
\[
K:=\{s\alpha+t\beta\colon s,t\in[0,1]\}.
\]
Since $K+\Lambda$ covers $\R^2$, translation by an element of $\Lambda$
allows us to assume that $z\in K$. 
The vertices $0,\alpha,\beta,\alpha+\beta$ of $K$ are the centers of $P_{0,0},P_{1,0},P_{0,1}$,
and $P_{1,1}$, respectively.
 Let $\pi\colon\mathbb R^2\to\mathbb R$ denote the scalar projection onto the direction $\omega$, defined by $\pi(\xi):=\xi\cdot\omega$. We show that the projections of these four squares cover $\pi(K)$. We write
\[
  s_\alpha:=\alpha\cdot\omega,\qquad s_\beta:=\beta\cdot\omega,\qquad
  p:=\min(|s_\alpha|,|s_\beta|),\qquad q:=\max(|s_\alpha|,|s_\beta|).
\]
The minimum of the four projected vertex values is
$s_* :=\min(0,s_\alpha)+\min(0,s_\beta)$. After subtracting $s_*$, these values are the sums of one element of $\{0,|s_\alpha|\}$ and one
element of $\{0,|s_\beta|\}$. Thus, in nondecreasing order and counting multiplicities, the projected vertex values are 
\[
  s_*,\qquad s_* + p,\qquad s_* + q,\qquad s_* + p + q.
\]
It follows that $\pi(K)=[s_*,s_*+p+q]$, and that the largest gap between consecutive projected vertex values is $\max(p,q-p)$. 
For $\omega=(\omega_1,\omega_2)$, we have
$$
  \pi(Q)=(-\tau/2,\tau/2),\qquad
  \tau=|\omega_1|+|\omega_2|=: \|\omega\|_1.
$$

Let $d: = \|\alpha\|_2 = \|\beta\|_2$. The assumption on
$\mu$ and $\nu$ gives
$$
d^2=(1+\mu)^2+\nu^2<1+\mu+\nu\le\sqrt{2}\,d.
$$
In particular, $d<\sqrt{2}$. Since $\alpha/d$ and $\beta/d$ form an orthonormal basis of $\R^2$,
$p^2+q^2=d^2$. Consequently,
$$
  p\le\frac{d}{\sqrt{2}}<1=\norm{\omega}_2\le\tau.
$$
Moreover, the identity $d^2\omega=s_\alpha \alpha + s_\beta \beta$, together with
the reverse triangle inequality in the $\ell^1$ norm and
$\norm{\alpha}_1=\norm{\beta}_1=1+\mu+\nu$, yields
$$
  d^2\tau
  =\|s_\alpha \alpha + s_\beta \beta\|_1
  \ge\bigl| \|s_\alpha \alpha\|_1- \|s_\beta \beta\|_1\bigr| =(q-p)(1+\mu+\nu).
$$
Hence
$$
  q-p\le\frac{d^2}{1+\mu+\nu}\,\tau<\tau.
$$
Thus every gap between consecutive projected vertex values is strictly smaller than~$\tau$. The projections of the four squares $P_{\ell,k}$,
$\ell,k\in\{0,1\}$, are open intervals of length $\tau$ centered at the projected vertex values.
The subsequent intervals overlap, and their union contains both extreme vertex values. Therefore,
$$
  \pi(K)\subset \bigcup_{\ell,k\in\{0,1\}}\pi(P_{\ell,k}).
$$
Since $z\in K$, there exist a vertex $\lambda$ of $K$ and a point $\xi\in Q + \lambda$ such that $\pi(\xi)=\pi(z)$. As
$\omega\cdot\zeta=0$, this equality implies that
$\xi\in z+\R\zeta$. Furthermore,
$$
  \|\xi - z\|_2
  \leq \|\xi-\lambda\|_2+\|\lambda-z\|_2
  < {1}/{\sqrt{2}} +\sqrt{2}\,d<3,
$$
where we used $\operatorname{diam}K=\sqrt{2}\,d$ and
$d<\sqrt{2}$. Hence $\xi= a + t_0\zeta$ for some $t_0\in(0,6)$. This proves the required intersection property for $\cQ_{1/2}$. Since every line segment of length $6$ intersecting $\cQ_{1/2}$ also intersects $\cQ_\eta$, the same property holds for every $\eta\in(0,1/2)$. We now fix $\eta\in(0,1/2]$. We set
$$
 F(a,\zeta):=\int_0^6 b(a+t\zeta)\,\ud t,
  \qquad (a,\zeta)\in\R^2\times S^1.
$$
The intersection property gives $b(a+t_0\zeta)>0$ for some
$t_0\in(0,6)$. By continuity and nonnegativity of $b$, it follows that $F(a,\zeta)>0$ for every $(a,\zeta)\in\R^2\times S^1$.
Moreover, $F$ is continuous  and $\Lambda$-periodic in its first
variable. Compactness of $K\times S^1$ therefore gives
$$
  \inf_{a\in\R^2,\ \zeta\in S^1}\int_0^6 b(a+t\zeta)\,dt
  =\min_{(a,\zeta)\in K\times S^1}F(a,\zeta)>0.
$$ 
Therefore, $b$ satisfies the (AGCC) condition, as claimed.
By Corollary~\ref{cor:agcc_implies_thm21}, we have 
$b\in \cB(\mathbb{R}^2)\cap \cB_D^{1/2}(\R^2)$. Thus, by Proposition \ref{prop:RninR}(ii) assertion~(i) is satisfied with $\theta  = 1/2$.  
A suitable choice of $b_\eta$ in \eqref{def:b_bQ} yields a function
$b\in\cB_D^\theta(\R^2)$ for every $\theta\in(0,1)$. Indeed, for each $\eta\in(0,1/2)$, we use \eqref{5.1} to define
\[
g_\eta(t) :=\psi(t^2-(\tfrac12-\eta)^2), \qquad f(t) :=\psi(\tfrac14-t^2),
\]
and 
\[
 b_\eta(\xi):=f(\xi_1)f(\xi_2) \big(g_\eta(\xi_1)+g_\eta(\xi_2)\big), \qquad
 b_{1/2}(\xi) :=f(\xi_1)f(\xi_2).
\]
Then $b_\eta\in C^\infty(\R^2)$ and $\{\xi: b_\eta(\xi) >0\} = Q_\eta$ for all $\eta \in (0,1/2]$.  If $b$ is defined by \eqref{def:b_bQ}, then straightforward calculations show that $b$ satisfies condition~\eqref{cond: b'<b} for every $\theta\in[1/2,1)$. Proposition \ref{cond: 1/t at infty} yields $b\in \cB_D^\theta(\R^2)$ for any $\theta \in (0,1)$. 

To prove (ii), we first observe that $\{b>0\}:=\{\xi: b(\xi)>0\}=\mathcal Q_\eta$. Each component $Q_\eta+\lambda$ is contained in the unit square
$Q+\lambda$, $\lambda\in\Lambda$.  The Dirichlet Poincar\'e inequality on
these squares, followed by summation, gives the Poincar\'e property of 
$\{b>0\}$. For $\R^2\setminus\supp b$, consider the translate $K_0:=K-\frac{\alpha+\beta}{2}$ of the square $K$ introduced above.
The squares $K_0+\lambda$, $\lambda\in\Lambda$, cover $\R^2$ with disjoint interiors.  
Put $E:=K_0\cap Q_\eta$. Since \(B(0,1/2)\subset K_0\) and \(|B(0,1/2)\cap Q_\eta|>0\), we have \(|E|>0\). 
Let $\varphi\in C_c^\infty(\R^2)\setminus\supp (b)$ and let
$\varphi_{K_0+\lambda}$ denote its mean over $K_0+\lambda$.  As $\varphi=0$ on
$E+\lambda$, the Poincar\'e--Wirtinger inequality gives
\[
 |E|^{1/2}|\varphi_{K_0+\lambda}|
 \le\|\varphi - \varphi_{K_0+\lambda}\|_{L^2(K_0+\lambda)}
 \le c\|\nabla \varphi\|_{L^2(K_0+\lambda)},
\]
where the constant $c$ is independent of $\lambda$. 
Consequently,
\[
 \|\varphi\|_{L^2(K_0+\lambda)}
 \le c\Bigl(1+\big({|K_0|}/{|E|}\big)^{1/2}\Bigr)
      \|\nabla \varphi\|_{L^2(K_0+\lambda)},
\]
hence summation proves the Poincar\'e property of $\R^2\setminus\supp b$.  Finally, the existence of a Weyl sequence supported in $\{b>0\}$ would contradict the Poincar\'e inequality on $\{b>0\}$.

To show (iii), we use $b$ given by \eqref{def:b_bQ}. For the first alternative, we fix admissible \(\mu,\nu\) and let \(\eta\downarrow0\), keeping the lattice \(\Lambda\) fixed. For the second, we take \(\eta=1/2\), set \(\mu=\nu^2\), and let \(\nu\downarrow0\), so that \(\Lambda\) varies. 
We first bound the number of lattice translates meeting a fixed ball.
Both $Q$ and $K_0$ are contained in $B(0,1)$, since
$d<\sqrt2$.  The balls $B(\lambda,1/2)$, $\lambda\in\Lambda$, are
pairwise disjoint.  If $D=Q$ or $D=K_0$ and
$(D+\lambda)\cap B(\xi;r)\ne\varnothing$, then
$B(\lambda,1/2)\subset B(\xi;r+3/2)$.  Comparing areas, at most
\(
 N_r:=4(r+3/2)^2
\)
translates $D+\lambda$ meet $B(\xi;r)$, for either choice of $D$.

For the first alternative, fix admissible $\mu,\nu$ and let $\eta\downarrow0$.
Since $\supp b=\bigcup_{\lambda\in\Lambda}
(\overline{Q_\eta}+\lambda)$,
\[
\begin{aligned}
 |B(\xi;r)\cap\supp b|
 &\le \sum_{\substack{\lambda\in\Lambda\\
             (Q+\lambda)\cap B(\xi;r)\ne\varnothing}}
       |\overline{Q_\eta}+\lambda| \\
 &\le N_r\bigl(1-(1-2\eta)^2\bigr)
   =4N_r\eta(1-\eta)\rightarrow0 \qquad \text{as}\qquad  \eta \downarrow 0.
\end{aligned}
\]
For the second alternative, take $\eta=1/2$ and $\mu=\nu^2$,
with $\nu>0$ sufficiently small to satisfy (2.29).  The translates
$K_0+\lambda$ cover $\R^2$, whereas $\supp b$ contains
$\overline Q+\lambda$ for each $\lambda\in\Lambda$.  Hence
\[
 B(\xi;r)\setminus\supp b
 \subset \bigcup_{\substack{\lambda\in\Lambda\\
        (K_0+\lambda)\cap B(\xi;r)\ne\varnothing}}
       ((K_0\setminus\overline Q)+\lambda).
\]
Since $K_0\subset(1+\mu+\nu)\overline Q$, we obtain
\[
\begin{aligned}
 |B(\xi;r)\setminus\supp b|
 &\le \sum_{\substack{\lambda\in\Lambda\\
        (K_0+\lambda)\cap B(\xi;r)\ne\varnothing}}
             |(K_0\setminus\overline Q)+\lambda| \\
 &\le N_r|K_0\setminus\overline Q| \\
 &\le N_r\bigl((1+\mu+\nu)^2-1\bigr)
   \longrightarrow0, \qquad \nu\downarrow0.
\end{aligned}
\]
Both estimates are uniform in $\xi$, proving (iii). \end{proof}

\begin{figure}[ht]
\centering
\begin{subfigure}{0.48\textwidth}
\centering

\begin{tikzpicture}[scale=0.85]
\def\muv{0.10}
\def\nuv{0.60}

\coordinate (a) at ({1+\muv},{\nuv});
\coordinate (b) at ({-\nuv},{1+\muv});

\foreach \i in {-1,0,1}{
  \foreach \j in {-1,0,1}{

    \pgfmathsetmacro{\x}{\i*(1+\muv)-\j*\nuv}
    \pgfmathsetmacro{\y}{\i*\nuv+\j*(1+\muv)}

    \draw[black,thin,fill=gray!10]
      ({\x-0.5},{\y-0.5})
      rectangle
      ({\x+0.5},{\y+0.5});
  }
}

\node[anchor=north east] at (0.5,0.5) {$Q_{0,0}$};

\pgfmathsetmacro{\xA}{1+\muv}
\pgfmathsetmacro{\yA}{\nuv}
\node[anchor=north east] at ({\xA+0.5},{\yA+0.5}) {$Q_{1,0}$};

\pgfmathsetmacro{\xB}{-\nuv}
\pgfmathsetmacro{\yB}{1+\muv}
\node[anchor=north east] at ({\xB+0.5},{\yB+0.5}) {$Q_{0,1}$};

\pgfmathsetmacro{\xC}{1*(1+\muv)-1*\nuv}
\pgfmathsetmacro{\yC}{1*\nuv+1*(1+\muv)}
\node[anchor=north east] at ({\xC+0.5},{\yC+0.5}) {$Q_{1,1}$};

\coordinate (P00) at (0.5,0.5);
\coordinate (P01) at ({0.5-\nuv},{1.5+\muv});

\draw[black,thick,shorten <=-18mm,shorten >=-6mm]
  (P00) -- (P01);

\end{tikzpicture}
\end{subfigure}

\caption{Geometry of the periodic damping region ($\mu = 3/20$, $\nu = 3/5$, and $\eta = 1/2$) } 
\label{fig:periodic}
\end{figure}
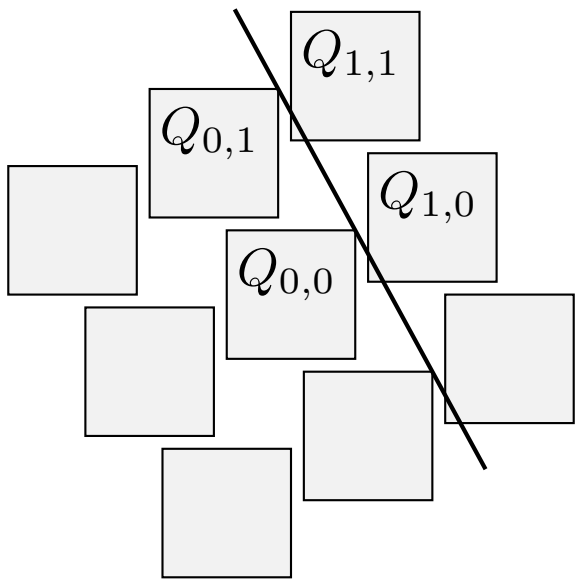

\paragraph{\bf Damping coefficients on exterior domains and (EGCC)}

Let $\Omega :=\R^n\setminus\overline{B(0;1)}$. We define
\begin{equation}\label{6.1}
 b(\xi): =h(\|\xi\|-2),
\end{equation}
where $h$ is given by \eqref{5.1}. Note that $b(\xi)=0$ when $1<|\xi|\leq2$, and $b(\xi)=1$ when $|\xi|\geq 3$. Such damping coefficient satisfies all assumptions of Theorem~\ref{thm:diffusion_for_Delta}(ii). In particular, Theorems~\ref{thm:mild_DPE}--\ref{thm: diffusion intro} apply with any $\theta < 1$.

\subsection{Explicit diffusive initial data for (DWE)}\label{subs:explict}
We now provide an explicit construction of \emph{diffusive initial data} for (DWE). We begin with an abstract auxiliary result; its proof is postponed to Section \ref{sec:proofs_of_thm1.6}.

\begin{proposition}\label{cor:fixed-scale} Let $A$ and $B$ satisfy Hypothesis~\ref{hyp:main_2}, and let $\cT$ denote the associated semigroup defined by~\eqref{formula:semig}. Let $\eta\in(0,1)$ and $y\in X$. Then the following conditions are equivalent:
\begin{itemize}
\item [(i)] $\sup_{t\geq 1} t^\eta \norm{\cT(t)y}<\infty$,
\item [(ii)] For some, equivalently for all $\omega'\in (0,\omega)$, we have 
\[
\sup_{\la \in \Sigma_{\omega' + \pi/2} \atop |\la|\leq 1}|\la|^{1-\eta}\norm{(\la B + A)^{-1}By}<\infty.
\]
\end{itemize}  
\end{proposition}

We now return to our model example, namely, $X=L^2(\Omega)$, $A=-\Delta$, and $B=M_b$.

\begin{proposition}\label{prop:Weyl_seq_Delta} Let $\Omega \subset \R^n$ be a domain such that $0\in \sigma(-\Delta)$ and let $b\in L^\infty_+(\Omega)$. Suppose that $(\varphi_j)\subset \Dom(\Delta)$ is a Weyl sequence corresponding to $0\in \sigma(-\Delta)$, i.e., 
\[
\|\varphi_j\|_2 = 1 \quad \textrm{ and }\quad s_j: = \|\Delta \varphi_j\|_2 \to 0 \quad \textrm{ as } j\rightarrow \infty.
\]
Suppose that the sets $G_j := \mathrm{supp}(\varphi_j)$ are pairwise disjoint and  
$\operatorname*{ess\,inf}_{\xi \in \bigcup_j G_j} b(\xi) = b_0$ for some $b_0 >0$. Let $(a_j)\in \ell^2(\N)$ and put $v_a : = \sum_{j\geq 1} a_j \varphi_j$.

Then $v_a \in \Delta^{-1}(\Ran(M_b))$ and if $(i s_j)_j \subset \rho(P_D)$ then
\[
\|(is_j M_b - \Delta)^{-1} M_b\, v_a\|_2\geq c |a_j|/s_j, \quad j\in \N, 
\]
where $c: = b_0/(1 + \|b\|_\infty)$. 
\end{proposition}

\begin{proof} By the closedness of $\Delta$, we infer that $v_a\in\Dom(\Delta)$, and hence $\Delta v_a=\sum_j a_j\Delta\varphi_j\in L^2(\Omega)$. Since $\supp(\Delta\varphi_j)\subset G_j$ and $b(\xi)\geq b_0$ for all $\xi\in G_j$, we have $\Delta\varphi_j\in\Ran(M_b)$ for every $j$, and consequently $\Delta v_a\in\Ran(M_b)$.
Let $x_j:=R_D(i s_j)M_b v_a$, $j\in \mathbb{N}$. Since $G_j$, $j\in \N$,  are pairwise disjoint and $M_b v_a = (i s_j M_b - \Delta) x_j$ we obtain
\begin{align*}
a_j \ip{ M_b \varphi_j}{\varphi_j}
= \ip{ M_b v_a}{\varphi_j} =  \ip{(i s_j M_b - \Delta)x_j}{\varphi_j}
=\ip{x_j}{(-i s_j M_b - \Delta)\varphi_j}.
\end{align*}
Since $|\ip{ M_b \varphi_j}{\varphi_j}|\geq b_0$ and $\|b\varphi_j\|_2 \le \|b\|_{\infty}$, 
\begin{equation}\label{eq:beta_a_j}
b_0 \, |a_j| 
\le  s_j \, (\|b\varphi_j\|_2+ 1) \|x_j\|_2 \leq s_j \, (\|b\|_\infty + 1) \|x_j\|_2,
\end{equation}
thereby showing the desired estimate. 
\end{proof}

Proposition \ref{prop:Weyl_seq_Delta} can be specialized to a given Weyl sequence for $0 \in \sigma(-\Delta)$. For instance, let $\varphi \in C_c^\infty(\R^n)$ with $\|\varphi\|_2 = 1$. Let $\Omega \subset \R^n$ be a domain such that for any $j\in \mathbb{N}$ there exist $\xi_j\in \Omega$ and $r_j>0$ such that $r_j\to \infty$ as $j\to \infty$, and 
\begin{equation}\label{def:G_j}
\xi_j+r_j \, \mathrm{supp}(\varphi)   \subset \Omega
\end{equation}
are pairwise disjoint sets. We set
\begin{equation}\label{def:packet_family}
\varphi_j(\xi):=r_j^{-n/2}\varphi \!\left({(\xi-\xi_j)}/{r_j}\right), \quad \xi \in \Omega, \, j\in \mathbb{N}.
\end{equation}
It is readily seen that $\|\varphi_j\|_2 = 1$ and $\|\Delta\varphi_j\|_2 =  r^{-2}_j \|\Delta \varphi\|_2 \to 0$ as $j\to \infty$.

\begin{proposition}\label{explict_parabolic_data}
Let $\Omega \subset \R^n$ be a domain and let $b\in L^\infty_+(\Omega)$. Assume that condition \eqref{DE-parab_est} holds with $A = -\Delta$ and $B = M_b$ on $X=L^2(\Omega)$. Suppose that $0\in \sigma(-\Delta)$ and $(\varphi_j)\subset \Dom(\Delta)$ is a Weyl sequence with pairwise disjoint supports of $\varphi_j$ and such that $\operatorname*{ess\,inf}_{\xi \in \bigcup \mathrm{supp}(\varphi_j)} b(\xi) >0$. 
Then the following assertions hold.
\begin{itemize}
\item [(i)] For every $\eta \in (0,1)$ and every $(a_j)\in \ell^2(\N)$ such that $\sup_{j\in \N} |a_j| /\|\Delta \varphi_j\|_2^{\eta} = \infty$ one has
\[
\limsup_{t\to \infty} t^\eta \|\cT(t)v_a\|_2 = \infty, 
\]
where $v_a = \sum_j a_j \varphi_j$ and $\cT$ denotes the semigroup associated with $A$ and $B$ (see \eqref{formula:semig}).

\item [(ii)] There exists a vector $v_0\in \Delta^{-1}(\Ran(M_b))$ such that for any $\eta\in (0,1)$ we have 
\[
\limsup_{t\to \infty} t^\eta \|v(t; v_0)\|_2 = \infty. 
\]
\end{itemize}
\end{proposition}

\begin{proof}
Statement (i) follows directly from Propositions \ref{prop:Weyl_seq_Delta} and \ref{cor:fixed-scale}. For (ii), since $\|\Delta \varphi_j\|_2 \to 0$ as $j\to \infty$, there exists an increasing subsequence $(j_k)\subset \N$ such that for any $\eta >0$ we have $ k\|\Delta \varphi_{j_k}\|^\eta_2 \to 0$ as $k\to \infty$. Let $a_j = 1/k$ if $j = j_k$ for some $k\in \N$, and otherwise $a_j=0$. Clearly, $(a_j)\in \ell^2(\N)$ and, for every $\eta \in (0,1)$, we have $\sup_j a_j/ \|\Delta \varphi_j\|^\eta_2 =\infty$. Setting $v_0 = \sum_j a_j \varphi_j$ and applying (i), we obtain~(ii).
\end{proof}

Based on Proposition \ref{explict_parabolic_data} and Theorem \ref{thm: diffusion intro}, one can formulate a weaker version of Theorem \ref{thm: diffusion}, and consequently of Theorem \ref{thm:diffusion_for_Delta}, with an explicit description of the initial data leading to the diffusion profiles. The following result gives such a variant. 
\begin{corollary}\label{cor:explicit_did}
\emph{(i)} Let $b \in L^\infty_+(\R^n)$ satisfy the assumptions of Theorem \ref{thm:diffusion_for_Delta}(i). In addition, suppose that there exists a Weyl sequence $(\varphi_j)$ for $0 \in \sigma(-\Delta)$ with pairwise disjoint supports such that $\operatorname{ess\,inf}_{\xi \in \bigcup_j \supp(\varphi_j)} b(\xi) > 0$.

Then for every $a = (a_j)\in \ell^2(\N)$ such that $\sup_{j\geq 1} |a_j|/\|\Delta \varphi_j\|_2^{\eta} = \infty$ for some $\eta \in (0,1)$ and for any $u_0 \in W^{2,2}(\R^n)$, denoting $u_1: = b v_0 -  b u_0$ and $v_0 : = \sum_{j\geq 1} a_j \varphi_j$,  we have
\[
\|u(t; u_0,u_1) - v(t; v_0)\| \lesssim 1/t, \, t>0, \quad \textrm{ and } \quad \limsup_{t\to \infty} t^\eta \|u(t; u_0,u_1)\| = \infty, 
\] where $u(\cdot;u_0, u_1)$ and $v(\cdot; v_0)$ denote the corresponding solutions to \emph{(DWE)} and \emph{(DE)}, respectively. 

\smallskip

\emph{(ii)} Let $\Omega$ and $b\in L^\infty_+(\Omega)$ satisfy the assumption of Theorem \ref{thm:diffusion_for_Delta}(ii). 
Let $\varphi \in C^\infty_c(\Omega)$ with $\|\varphi\|_2 = 1$.
Let $\varphi_j$, $j\in \N$, be given by \eqref{def:packet_family} with $r_j : = 2^j$ and $\xi_j$ such that \eqref{def:G_j} holds,
and  $\operatorname{ess\,inf}_{\xi \in \bigcup_j \supp(\varphi_j)} b(\xi) > 0$.
Let $v_0 : = \sum j^{-1}\varphi_j$. 

Then for any $u_0 \in W^{1,2}_0(\Omega)\cap W^{2,2}(\Omega)$, if $u_1: = b v_0 -  b u_0$, then 
\[
\|u(t; u_0,u_1) - v(t; v_0)\| \lesssim 1/t \, \textrm{ and } \, \limsup_{t\to \infty} t^\eta \|u(t; u_0, u_1)\| = \infty \, \textrm{ for all $\eta \in (0,1)$}, 
\] where $u(\cdot;u_0, u_1)$ and $v(\cdot; v_0)$ denote the corresponding solutions to \emph{(DWE)} and \emph{(DE)}, respectively.

Furthermore, in both cases \emph{(i)} and \emph{(ii)}, the corresponding solution $v(\cdot;v_0)$ to \emph{(DE)} admits a bounded holomorphic extension to every sector $\Sigma_\omega$, where $\omega\in(0,\pi/2)$.

\end{corollary}
We only briefly comment on the proof. First, Proposition~\ref{prop:Weyl_seq_Delta},
Theorem~\ref{thm:mild_DPE}(ii), and definition~\eqref{def_YD2}
yield $v_0\in \Delta^{-1}\big(\operatorname{Ran}(M_b)\big)
\subset X^1_D\subset Y^{1/2}_D$
for every $v_0$ specified in either statement (i) or (ii) above. 
Since each Weyl sequence associated with $0\in\sigma(-\Delta)$ is also a Weyl 
sequence for $0\in\sigma((-\Delta)^{1/2})$, the assumptions of Corollary 
\ref{cor:explicit_did} imply those of Proposition~\ref{prop:correction_b}. 
Consequently, $b\in\mathcal{R}(\mathbb{R}^n)$ in case (i) and $b\in\mathcal{R}(\Omega)$ in case (ii).

\section{Proof of the abstract multiplier result}\label{sect3}

\subsection{Auxiliary lemmas}
The proof of Theorem \ref{thm 1} uses the following vector-valued version of \cite[Proposition 4.4.2]{St93}. It converts Mihlin-type estimates for a symbol into regularity and decay of its distributional inverse Fourier transform away from zero. For a Banach space $Y$, we write $\mathcal S'(Y):=\mathcal S'(\R;Y)$ and regard $L^1(\R;Y)$ and $L^\infty(\R;Y)$ as subspaces of $\mathcal S'(Y)$. The Fourier transforms are understood with the convention given in the Notation paragraph. For background on vector-valued distributions and their Fourier transforms, see \cite[Section III.4]{Am95}.

Throughout, for functions $a, k \in C^n(\R\setminus \{0\}; Y)$, $n\in \N_0$,  we write $a\in \m^n(Y)$ and $k\in \cK^n(Y)$, if respectively 
\[
 [a]_{\m^n}:=\max_{\ell =0,...,n}\sup_{s\neq 0} \|s^\ell a^{(\ell)}(s)\|_{Y}<\infty \qquad \textrm{and} \tag{$\m^n$}
\]
\[
[k]_{\cK^{n}} := \max_{\ell =0,...,n} \sup_{t\neq 0} \|t^{\ell + 1} k^{(\ell)}(t)\|_{Y}<\infty.
\tag{$\cK^{n}$}
\]

\begin{lemma}\label{CZ cond}
Let $Y$ be a Banach space and $n\ge2$ an integer.
For every $a\in\m_n(Y)$ there exists a unique
$k\in\cK_{n-2}(Y)$ such that
\begin{equation}\label{k}
  \left.\cF^{-1}a\right|_{\R\setminus\{0\}}=k
  \quad\text{in }\mathcal D'(\R\setminus\{0\};Y),
\end{equation}
and
\[
 [k]_{\cK_{n-2}}\le C_n[a]_{\m_n},
\]
the constant $C_n$ depends only on $n$ (and on the fixed dyadic partition).
\end{lemma}
In what follows, $\cF^{-1}a(t)$ for $t\ne0$ denotes the representative $k$ given by \eqref{k}.

This lemma is a vector-valued variant of \cite[Proposition~4.4.2(a)]{St93}. An operator-valued version was stated in \cite[Lemma~4.1]{Krol23}, where the proof was omitted. Since the result plays a fundamental role in our approach, we give a detailed proof.

\begin{proof}[{Proof of Lemma \ref{CZ cond}}]
Assume $a\in \m^n(Y)$. Let $\varphi\in C^\infty(\mathbb{R})$ be a function supported in $\{1/2\le |s|\le 2\}$ such that the following partition of unity holds:
\[
\sum_{j\in\mathbb{Z}}\varphi(2^{-j}s)=1, \qquad s\neq 0.
\]
Clearly, for any $n\in \N$ we have $\varphi \in \m^n(\R)$, that is, $[\varphi]_{\m^n} < \infty$.  
For each $j\in \ZZ$ we set 
\[
\varphi_j(s) : =\varphi(2^{-j}s) \quad \textrm{ and }\quad
a_j(s):=a(s)\varphi_j(s),\qquad s\neq 0.
\]
Then $a_j\in C^n(\mathbb{R};Y)$, $\operatorname{supp} a_j\subset\{2^{j-1}\le |s|\le 2^{j+1}\}$, and
$\sum_j a_j(s)=a(s)$ for $s\neq 0$ (locally finite sum on $\mathbb{R}\setminus\{0\}$). 

By the dominated convergence theorem the sequence $S_N: = \sum_{|j|\leq N} a_j$, $N\in \N$, converges to $a$ in $S'(Y)$ as $N\rightarrow \infty$.
Therefore, setting $k_j:=\mathcal{F}^{-1}a_j\in \mathcal{S}'(Y)$, we have
\[
K_N: = \sum_{|j|\leq N} k_j = \cF^{-1}S_N  \rightarrow \cF^{-1}a \qquad \textrm{ in } S'(Y) \textrm{ as } N\rightarrow \infty.
\]
We claim that for any $\ell \in \{0, ..., n-2\}$ the sequence of the derivatives $K_N^{(\ell)}$, $N\in \N$, 
converges locally uniformly on $\R\setminus\{0\}$ to a function
$\widetilde k_\ell\in C(\R\setminus\{0\};Y)$ and for any  $t\neq 0$ and $N\in \N$, 
\begin{equation}\label{K_N est}
\|K_N^{(\ell)}(t)\|_Y\leq C  [a]_{\m^n}|t|^{-(\ell+1)}
\end{equation} with  a constant $C$ independent on $N$ and $t$. 
Then, setting $k := \widetilde k_0$, we obtain the desired conclusion:  
Indeed, for $n = 2$ this follows directly. 
If $n \ge 3$, then for every $\ell \in \{0, \dots, n-3\}$ we have
\[
K_N^{(\ell)}(t) - K_N^{(\ell)}(\sgn(t))
= \int_{\sgn(t)}^{t} K_N^{(\ell+1)}(\tau)\, \mathrm d\tau,
\qquad t \neq 0.
\]
Letting $N \to \infty$, we obtain
\[
\widetilde k_{\ell}(t) - \widetilde k_{\ell}(\sgn(t))
= \int_{\sgn(t)}^{t} \widetilde k_{\ell+1}(\tau)\, \mathrm d\tau,
\qquad t \neq 0.
\]
This shows that $k := \widetilde k_0$ is the desired function. 
Thus, it remains to justify the claim. Fix $\ell \in \{0,...,n-2\}$ and $j\in \ZZ$. It is readily seen that
\begin{equation}\label{aj}
\left(\cF^{-1}a_j\right)^{(\ell)}(t)  = \cF^{-1}(g_{\ell,j})(t), 
\end{equation}
where 
\[
g_{\ell,j}(s):= (is)^{\ell} a_j(s), \qquad s\neq 0.
\]
Hence, for any $t\in \R$ we have 
\begin{align}\label{n=0}
\nonumber \left\|\left(\cF^{-1}a_j\right)^{(\ell)}(t)\right\|_Y &  \leq \frac{1}{2\pi}\int_{\{ 2^{j-1} \leq |s| \leq 2^{j+1}\}} |s|^\ell |\varphi (2^{-j}s)| \|a(s)\|_Y  \ud s\\
& \leq 2^{\ell + 1} \sup_{s\neq 0} \|a(s)\|_Y \sup_{s\neq 0} |\varphi(s)| \left(2^{\ell +1}\right)^j.
\end{align}

Moreover, we note that 
\begin{equation}\label{sn aj}
t^n \left(\cF^{-1}a_j\right)^{(\ell)}(t)  = i^n \cF^{-1}\left(g_{\ell,j}^{(n)}\right)(t), \qquad t\in \R.
\end{equation}

Further, by Leibniz' rule, 
\begin{align*}
g_{\ell,j}^{(n)}(s) &= \sum_{q=0}^n \binom nq \left( (is)^\ell\right)^{(q)}
\left( a_j(s) \right)^{(n-q)}\\
& = i^\ell\sum_{q=0}^\ell \binom nq \frac{\ell !}{(\ell - q)!} s^{\ell - q}
\left( a_j(s) \right)^{(n-q)},
\end{align*}
\[
a_j^{(n-q)}(s)=\sum_{r=0}^{n-q} \binom {n-q}r\, a^{(r)}(s)\, \varphi_j^{(n-q-r)}(s), \qquad s\neq 0.
\]
Since for any $r\in \{0,..., n\}$ we have $[a]_{\m^{r}}\leq [a]_{\m^n} < \infty$, 
\begin{align*}
\left\| s^{n-q} (a_j(s))^{(n-q)} \right\|_Y & = \left\| \sum_{r=0}^{n-q} \binom {n-q}r\, s^r a^{(r)}(s)\, (s2^{-j})^{n-q-r}\varphi^{(n-q-r)} (2^{-j}s) \right\|_Y\\
& \leq \sum_{r=0}^{n-q} \binom {n-q}r [a]_{\m^r} [\varphi]_{\m^{n-q-r}} \chi_{j}(s)\\
& \leq   C_{\varphi, n} [a]_{\m^n} \, \chi_j(s),
\end{align*}
where $\chi_j$ stands for the characteristic function of the shell $\{2^{j-1}\leq |s| \leq 2^{j+1}\}\supset \supp \varphi_j$, and, consequently, 
\begin{align*}
\left\|g_{\ell,j}^{(n)}(s)\right\|_Y & = \left\| \sum_{q=0}^\ell \binom nq \frac{\ell !}{(\ell - q)!} s^{\ell - n} s^{n - q} 
\left( a_j(s) \right)^{(n-q)}\right\|_Y\\
& \leq (\ell +1) C_{\varphi, n} [a]_{\m^n} \, |s|^{\ell-n} \chi_j(s), \qquad s\neq 0,
\end{align*}
for an appropriate constant $C_{\varphi, n}$ depending only on $n$ and $\varphi$.

Returning to \eqref{sn aj} we can write for every $t\in \R$
\begin{align*}\label{n=n}
\left\| t^n \left(\cF^{-1}a_j\right)^{(\ell)}(t) \right\|_Y & \leq \frac{ C_{\varphi, n}[a]_{\m^n}}{2\pi} \int_\R |s|^{\ell - n} \chi_j(s) \ud s\\
& \leq   C_{\varphi, n}[a]_{\m^n} 2^{(j-1)(-n+ \ell +1)}.
\end{align*}
Combining this estimate with \eqref{n=0}, for an appropriate constant $\wt C$ and all $t\neq 0 $, $j\in \ZZ$ and $\ell \in \{0,...,n-2\}$ we infer
\[
\left\|k_j^{(\ell)}(t) \right\|_Y = \left\| \left(\cF^{-1}a_j\right)^{(\ell)}(t) \right\|_Y  \leq \wt C [a]_{\m^n}  \min\left( |t|^{-n} 2^{(-n + \ell +1)j} , \, 2^{(\ell +1)j} \right).
\]
This estimate readily yields the uniform convergence of 
$(K_N^{(\ell)})_{N\in\mathbb N}$, for $\ell\in\{0,\dots,n-2\}$, 
on compact subsets of $\mathbb R\setminus\{0\}$ 
to a function in $C(\mathbb R\setminus\{0\};Y)$. 
 Indeed, let $N, N'\in \N$ with $N'<\N$ and $\varepsilon>0$. Then
for any $|t|\geq \varepsilon$ 
\begin{align*}
\left\|K^{(\ell)}_N(t) - K^{(\ell)}_{N'}(t)\right\|_Y & \leq \sum_{N'<|j|\leq N} \left\|k^{(\ell)}_j(t)\right\|_Y\\
& \leq \wt C [a]_{\m^n} \left( \varepsilon^{-n} \sum_{N'<j\leq N} 2^{(-n +\ell +1)j}  +  \sum_{-N\leq j < -N'} 2^{(\ell +1)j} \right).
\end{align*}
Since $n \ge 2$ and $n-\ell-1 \ge 1$, the term in the bracket 
tends to zero as $N, N' \to \infty$.  
For the pointwise estimate \eqref{K_N est}, we fix $t \neq 0$ and choose 
$N_0 \in \mathbb{Z}$ such that 
$
2^{N_0}|t| \le 1 < 2^{N_0+1}|t|$. 
Then  
\begin{align*}
\left\| K^{(\ell)}_N(t) \right\|_Y & \leq  \sum_{|j|\leq N} \left\|k^{(\ell)}_j(t)\right\|_Y \\ 
& \leq \wt C [a]_{\m^n}  \sum_{j\leq N_0} 2^{(\ell +1)j} + C [a]_{\m^n}  \sum_{j > N_0} |t|^{-n} 2^{(-n + \ell +1)j} \\
& \leq 2 \wt C [a]_{\m^n} \left(  2^{N_0(\ell +1)} +    |t|^{-n}   2^{(N_0 +1 )(-n + \ell +1)} \right)\\ 
& \leq   6  \wt C  [a]_{\m^n}  {|t|^{-(\ell + 1)}}
\end{align*}
as claimed in \eqref{K_N est}.  

To identify the limit with the distributional Fourier transform, let
$\psi\in C_c^\infty(\R\setminus\{0\})$.
Since $K_N\to\cF^{-1}a$ in $\mathcal S'(\R;Y)$ and
$K_N\to k$ uniformly on $\supp\psi$, we have
\[
 \langle\cF^{-1}a,\psi\rangle
 =\lim_{N\to\infty}\int_\R K_N(t)\psi(t)\,dt
 =\int_\R k(t)\psi(t)\,dt.
\]
Thus $k$ represents the restriction of $\cF^{-1}a$ to
$\R\setminus\{0\}.$
\end{proof}

For later purposes, we introduce the following auxiliary notation.  We say that a function $a\in C^n(\R\setminus \{0\}; Y)$, $n \in \N$,  satisfies condition $(\m^n)$ {\it near zero}, and write $a \in \m^n_{0}(Y)$, if 
\[
[a]_{\m^n, 0}:=\max_{\ell =0,...,n}\sup_{0<|s|\leq 1} \|s^\ell a^{(\ell)}(s)\|_{Y}<\infty.  \tag{$\m^n_{0}$}
\]
We say that $a$ satisfies condition $(\m^n)$ at infinity and write $a \in \mathcal{M}_\infty^n(Y)$~if
\[
[a]_{\m^n, \infty}: = \max_{\ell = 0, \dots, n} \sup_{|s| \geq 1} \| s^\ell a^{(\ell)}(s) \|_{Y} < \infty. \tag{$\mathcal{M}_\infty^n$}
\]

Below, we denote by $\arg \la$ the principal argument of $\la\in \C$ taking values in $(-\pi,\pi]$. We use the notation introduced in Subsection~\ref{subs:main_results}. Moreover, if $i\R \setminus \{0\} \subset \rho(P)\cap \rho(P_D)$, then for any $s\neq 0$ we write 
\begin{equation}\label{a i d}
a(s): = is \alpha(s) \quad \textrm{ and } \quad d(s): = is \delta(s) B,
\end{equation}
where 
\begin{equation}\label{a i d2}
\alpha(s): = R(is) = (-s^2 + is B +A)^{-1} \quad \textrm{ and } \quad \delta(s): = R_{D}(is) = (is  B +A)^{-1}. 
\end{equation}
The following lemma provides an extension of estimates \eqref{pencils est} and \eqref{hyp:bound_2} for $R$ and $R_DB$, which allows us to show that $m = m(\cdot; x_0, y_0)$ is in $\m_{0}^n(X)$ for any $n\in \N$ and all $x_0$ and $y_0$ in $X$. 

\begin{lemma}\label{lem on pencil est} Assume Hypothesis \ref{hyp:main}.
Then the following assertions hold.
\begin{itemize}
\item [(i)] 
There exists $\omega\in (0,\frac{\pi}{2})$ such that 
\[
G_{\omega}^R:= \left\{\la \in \C\setminus\{0\} :  \left| |\arg \la| - \frac{\pi}{2}\right|\leq \omega \textrm{ and }  |\Re \la| \leq  1/(12[R])  \right\} \subset \rho(P)
\]
and
\begin{equation}\label{pencil est ex}
\sup_{\la \in G_\omega^R} \| \la R(\la) \| \leq 4[R].
\end{equation} 
Moreover, the function $a$ from \eqref{a i d} is in the class $\m^n_{0}( \cL(X))$ for any $n \in \N$.

\item [(ii)] We have the inclusion
\[
G_\theta : =  \left\{ \la\in \C\setminus\{0\}:  |\Re \la|\leq  \min(| \Im \la|, |\Im \la|^{\theta} )/(2[R_D B]_\theta) \right\} \subset \rho(P_D)
\]
and
\begin{equation}\label{pencil ext small}
\sup_{\la \in G_\theta} \min(|\la|, |\la|^\theta) \|  R_D(\la)B \| \leq 1+ 2[R_D B]_\theta.
\end{equation} 
Moreover, the function $d$ from \eqref{a i d} is in the class $\cM^n_0( \cL(X))$ for any $n \in \N$. If $\theta = 1$, then $d$ is in $\cM^n(\cL(X))$ for any $n\in \N$. 
\end{itemize}
\end{lemma}

\begin{proof} (i) To prove $G^R_\omega \subset \rho(P)$, let $s \neq 0$ and 
$r \in \mathbb{R}$ be arbitrary and set \[
T: = -s^2 + is B + A  \quad \textrm{ and } \quad S := r^2 + 2i r s + r B.
\]
In particular, by our assumption \eqref{pencils est}, $T = P(is)$ is invertible and
\[
\|S T^{-1}\| \leq {|r|}{|s|^{-1}} \| r + 2is + B\| [R] < {2^{-1}}
\]
provided
\begin{equation}\label{region G}
|r|\leq |s|^\frac{1}{2}/(6[R])^{1/2} \quad \textrm{and} \quad    |r| \leq  {1}/{(12 [R])}\quad \textrm{and} \quad \|B\| |r| \leq   |s|/ (6 [R]). 
\end{equation}
Consequently, for any such $r$ and $s$ the operator  
\[
(ST^{-1} + I)T = T+S = (r+is)^2 +  (r+is) B + A = P(r+is)\]
is invertible; see, e.g., Kato \cite[Section IV, Theorem 1.16, p.196]{Ka95}.
Moreover, it is readily seen that there exists $\omega>0$ such that the set $G^R_\omega$ given in (i) is contained in the region of $\la$'s described by \eqref{region G}. Thus, we have that $\la^2 + \la B + A$ is invertible for all $\la \in G^R_\omega$ and, using the Neumann series expansion, $
\|R(\la)\| \leq 2 \|R(i \Im \la)\|$. 
Taking $\omega$ small, if necessary, we have $\sup_{r+is\in G^R_\omega} |r+is|/|s|\leq 2$, and the desired estimate \eqref{pencil est ex} holds.

To prove that $a\in \m^n_0(\cL(X))$ for every $n\in \N$, note that for any $N>0$ there exists $0<\omega'<\omega$ such that 
\[
\left\{ \la\in \C\setminus\{0\}:  \left| |\arg \la| - \frac{\pi}{2}\right|\leq \omega' \right\}  \cap \{|\la|\leq 2N\} \subset G^R_\omega. 
\]
Since $R$ is analytic in $G^R_\omega$,  for any $0<|s|\leq N$, by Cauchy's formula,
\[
R^{(\ell)}(is)  = \frac{\ell!}{2\pi i} \int_{|\la - is|= |s|\sin \omega' }\frac{R(\la)}{(\la-is)^{\ell +1}} \ud \la. 
\]
Consequently, for any $\ell \in \N_0$ we have, by \eqref{pencil est ex},
\begin{equation}\label{Cauchy est}
\sup_{0<|s|\leq N} \| s^{\ell+1} R^{(\ell)} (is)\| \leq  \frac{2\pi\cdot \ell!}{(1-\sin\omega')\sin^\ell\omega'} 4[R]. 
\end{equation}
Recall that $\alpha(s) = R(is)$, $s\neq 0$. Clearly, for each $\ell\in \N_0$, $  \alpha^{(\ell)}(s) = i^\ell R^{(\ell)}(is) $, $s\neq 0$. By \eqref{Cauchy est}, for each $\ell\in \N_0$
\begin{equation}\label{loc est alpha}
\sup_{0<|s|\leq  N} \| s^{\ell +1} \alpha^{(\ell)} ( s)\| < \infty.
\end{equation}
Since for each $\ell \in \N$ 
\begin{align*}
a^{(\ell)}(s) & = is \alpha^{(\ell)}(s) + i\ell \alpha^{(\ell -1)}(s),
\end{align*}
we get that $a\in \m^n_{0}(\cL(X))$ for all $n \in \N$.
This completes the proof of (i).

(ii) To show that $G_\theta \subset \rho(P_D)$ one can proceed similarly as in the proof of (i). For any $r\in \R$ and $s\neq 0$, we set 
\[
T = is B + A \qquad \textrm{ and } \qquad S = r B.
\]
Let $\phi(\la) : = \min(|\la|, |\la|^\theta)$, $\la \in \C$. First, we note that 
\[
\|T^{-1}S\| = {|r|}{\phi(s)^{-1}}\|\phi(s) (is B +A )^{-1} B\| \leq {|r|}\phi(s)^{-1} [R_D B]_\theta \leq {2^{-1}}
\]
for any $|r|\leq \phi(s)/(2[R_D B]_\theta)$. Hence, using the Neumann series expansion, $(I+T^{-1}S)^{-1}T^{-1}$ is a bounded, injective operator on $X$ whose range is contained in the domain $D(T)=D(T+S)$. 
Moreover, 
\[
(T+S) (I + T^{-1}S)^{-1} T^{-1}  =   I,
\]
\[
(I+T^{-1}S)^{-1} T^{-1}(T+S) = (I+T^{-1}S)^{-1}(I + T^{-1}S) \subset I,
\]
where $I$ stands for the identity operator on $X$.
Therefore, for any $s\neq 0$ and $|r|\leq {\phi(s)}/{2[R_D B]_\theta}$, the operator $T+S = (r+is)B +A = P_D(r+is)$ is invertible, i.e., $G_\theta\subset \rho(P_D)$. 

In particular, for any  $|r|\leq {\phi(s)}/{2[R_D B]_\theta}$ and $s\neq 0$,
we obtain
\begin{align*}
\phi(r+is) \|R_D(r+is) B\| & \leq  {2 \phi(r+is)}{\phi(s)^{-1}}  \| \phi(s) R_D(is) B \| \leq 1+ 2 [R_DB]_\theta. 
\end{align*}
This proves \eqref{pencil ext small}. 
The same argument based on Cauchy's formula shows that
for every $\ell\in \N_0$ 
\begin{equation}\label{est on beta}
\|s^{\ell +1} R_D^{(\ell)}(is)B\| \leq \frac{2\pi  \cdot \ell!  }{(1-\sin \omega) \sin^\ell\omega} (1+ 2 [R_DB]_\theta).
\end{equation}
holds for all $0<|s|\le1$ if $\theta\neq1$, and for all $s\neq0$ if $\theta=1$.
Recalling $\delta(s) = R_D(is)$, $s\neq 0$, since
\begin{align*}
d^{(\ell)}(s) & = (is \delta(s) B)^{(\ell)} =  is \delta^{(\ell)}(s)B  + i\ell \delta^{(\ell -1)}(s)B,
\end{align*}
the desired assertion concerning $d$ follows.
\end{proof}

\subsection{Proof of Theorem \ref{thm 1}} 
We follow the notation introduced in \eqref{a i d} and \eqref{a i d2}, and recall that $\mathcal{M}_0^n(X)$ and $\mathcal{M}_\infty^n(X)$ denote the classes of functions in $C^n(\R\setminus\{0\};X)$ satisfying the Mihlin estimates of order $n \in \mathbb{N}$ on the regions $\{0 < |s| \leq 1\}$ and $\{|s| \geq 1\}$, respectively. Given $x_0, y_0\in X$, the function $m$ defined in \eqref{def of m} could be written as
\begin{align}\label{def of m2}
m(s)&=is\alpha(s) x_0+\alpha(s)By_0-\delta(s)By_0,\, \text{ or as}\\
m(s)&=a(s)x_0 - a(s) d(s)y_0. \label{def of m3}
\end{align}
The latter identity follows from the formula
\[\alpha(s)-\delta(s) = \alpha(s)\big((isB+A)-(-s^2+isB+A)\big)\delta(s) =s^2\alpha(s)\delta(s), \quad s\neq0,\]
which also shows that $\alpha(s)$ and $\delta(s)$ commute.
We also define 
\begin{equation}\label{def: m_123}
m_1(s) = is \alpha(s) x_0, \quad m_2(s) = \alpha(s)By_0, \quad m_3(s) = \delta(s) By_0, \quad s\neq 0.
\end{equation}
Thus, by \eqref{def of m2}, $m(s) = m_1(s) + m_2(s) - m_3(s)$, $s \neq 0$.

To prove Theorem~\ref{thm 1} it suffices to show, under the corresponding assumptions on $x_0$ and $y_0$, that in Theorem~\ref{thm 1}(i) one has $m \in \mathcal{M}^2(X)$, in Theorem~\ref{thm 1}(iii) one has $m \in \mathcal{M}^3(X)$, and in Theorem~\ref{thm 1}(ii) one has  $M\in \mathcal{M}^3(X)$, where 
\begin{equation}\label{def:M}
M(s) := (is)^{-1}m(s), \qquad s \neq 0.
\end{equation}
Once this is established, the desired bounds in (i) and (iii) follow directly  from Lemma \ref{CZ cond}. For (ii), Lemma \ref{CZ cond} implies  $\cF^{-1}((is)^{-1}m)\in\cK^1(X)$, and thus the function 
$t\mapsto t^2\big(\cF^{-1}((is)^{-1}m)\big)'(t)$ is uniformly bounded on $\mathbb{R} \setminus \{0\}$. Since $(\cF^{-1}m)(t)=\big(\cF^{-1}((is)^{-1}m)\big)'(t)$, the desired bound in (ii) follows.  

We begin by collecting some immediate consequences of Lemma~\ref{lem on pencil est}, together with several observations that facilitate the verification of the above conditions on~$m$.

\begin{lemma}\label{lem:reduction_step} Assume Hypothesis~\ref{hyp:main}. Let $n \in \mathbb{N}$ and $x_0, y_0 \in X$. Let $m = m(\cdot; x_0, y_0)$, and let $m_j$, $j=1,2,3$, and $M$ be defined as in \eqref{def: m_123} and \eqref{def:M}. Then the following assertions hold.
\begin{itemize}
\item [(1)] $m\in\cM_0^n(X)$. Hence, $m\in\cM^n(X)$ whenever $m_1,m_2,m_3\in\cM_\infty^n(X)$.
\item [(2)] If $m\in \cM_\infty^n(X)$, then $M \in \cM_\infty^n(X)$.
\item[(3)] $\Dom(A)\subset X_2$ and $A^{-1}\big(BX\big) \subset X_D^1$. Moreover, if $\sup_{|s|\ge1}\|\delta(s)\|<\infty$, then $\Dom(A)\subset X_D^1$.
\item [(4)] If $n(1-\theta)\leq 1$ and $y_0\in X_D^{n(1-\theta)}$, then $m_3 \in \cM^n_\infty(X)$. In particular,  if $\theta  = 1$, then $m_3\in \cM_\infty^n(X)$.
 
\item [(5)] If $y_0 \in \Dom(A)$, $\theta \in [1/2, 1]$, and $\sup_{|s|\geq 1} \|\delta(s)\|<\infty$ (or $Ay_0 \in \Ran(B)$), then $m_3\in \cM^2_\infty(X)$.
\end{itemize}

\end{lemma}
\begin{proof}
(1) Combining Lemma~\ref{lem on pencil est} with \eqref{def of m3} and applying the product rule to $(ad)^{(n)}$, we obtain that $m\in\cM_0^n(X)$. The second assertion follows directly from formula \eqref{def of m2}; see also \eqref{def: m_123}.
 
(2) follows directly by Leibniz' rule. 
(3) For any $x \in \Dom(A)$ we may write
\[
s^2 \alpha(s)x = is\,\alpha(s)Bx + \alpha(s)Ax - x,
\]
and thus $\Dom(A) \subset X_2$ follows from \eqref{pencils est}. Moreover, if $x\in \Dom(A)$ with $Ax = By$ for some $y\in X$, then 
\begin{equation}\label{eq:basic_delta}
is \delta(s)Bx =  x - \delta(s)Ax = x - \delta(s)By 
\end{equation} 
showing that $A^{-1}\big(BX\big) \subset X_D^1$ and $\Dom(A) \subset X_D^1$ whenever \eqref{hyp:bound_2} and $\sup_{|s|\geq 1}\|\delta(s)\|<\infty$, respectively.
 
(4) A direct verification shows that for all  $s\neq 0$ and $n\in \N$ we have \begin{equation}\label{eq:deriv_of_delta}
\delta^{(n)} (s) = (-i)^n n! (\delta(s)B)^{n} \delta(s),
\end{equation}
which easily yields (4). 

(5) From (4) with $n=1$ we infer that $m_3\in \cM_\infty^1(X)$. Moreover, by \eqref{eq:deriv_of_delta}, for $y_0\in \Dom(A)$ we have
\begin{align*}
s^2 m''_3(s) & = - 2 (|s|^\theta \delta(s) B)^2 |s|^{2-2\theta} \delta(s) By_0\\
& = - 2 (|s|^\theta \delta(s) B)^2 |s|^{2-2\theta} (is)^{-1} (y_0 - \delta(s) Ay_0).
\end{align*}
Since by our assumption $\theta \geq 1/2$ and  $\sup_{|s|\geq 1} \|\delta (s)Ay_0\|\leq \sup_{|s|\geq 1} \|\delta (s)\| \|Ay_0\|<\infty$ (or $ \sup_{|s|\geq 1} \|\delta (s)Ay_0\| = \sup_{|s|\geq 1} \|\delta (s)B x\|<\infty$ for some $x\in X$), (5) holds.
\end{proof}

We follow the notation introduced in \eqref{def_YD2} and \eqref{def_wtZD}.

\begin{lemma}\label{lem:m3_in_Mikhlin} Let $\theta \in [0,1]$ and  let $A$ and $B$ satisfy \eqref{hyp:bound_2}. Then, $m_3 \in \cM^2_\infty(X)$ provided $y_0\in Y_D^{\theta}$ and $m_3 \in \cM^3_\infty(X)$ provided $y_0\in \wt Z_D^{\theta}$.

Moreover, if $y_0\in Z_D^{\theta}$ then, for every $\ell \in \{1,2,3\}$, we have 
\begin{equation}\label{cM3_for_M}
\sup_{|s|\geq 1} \|s^{\ell - 1} m_3^{(\ell)}(s) \| <\infty. 
\end{equation}
\end{lemma}
\begin{proof} 
Since $A^{-1}\big(BX\big)\subset X_D^1$ by Lemma \ref{lem:reduction_step}(3), a direct verification shows that for any  $\theta \in [0,1]$ one has $\wt Z_D^{\theta}\subset Y_D^{\theta}\subset X_D^{1-\theta}$.

Let $y_0\in X_D^{1-\theta}$. By \eqref{hyp:bound_2} we have 
the uniform boundedness of $m_3(s)$ on $\{|s|\geq 1\}$.  Moreover, by \eqref{eq:deriv_of_delta}, 
\[
|s| m_3'(s) = -i |s|^\theta \delta(s)B |s|^{1-\theta}\delta(s)B.
\]
Thus, the uniform boundedness of $sm_3'(s)$ on $\{|s|\geq 1\}$ follows directly from $y_0\in X_D^{1-\theta}$. Assume now that $y_0 \in Y_D^{\theta}$. By \eqref{eq:deriv_of_delta} we have
\[
s^2 m_3''(s) = -2 (|s|^\theta\delta(s)B)^2 |s|^{2-2\theta}\delta(s) By_0.
\]  
Since $2-2\theta\leq \theta$ for $\theta \in [\frac23 ,1]$, $[R_DB]_\theta<\infty$ gives the desired bound for such $\theta$. Since $2-2\theta \in [\frac23, 1]$ for $\theta \in [\frac12, \frac23]$, condition $y_0\in X_D^{2-2\theta}$ suffices to obtain the uniform boundedness of $s^2m''_3(s)$ on $\{|s|\geq 1\}$. For $\theta \in [0,\frac12)$ we have $2-2\theta>1$, however, taking $y_0 \in A^{-1}(BX_D^{1-2\theta})$, i.e., $y_0\in D(A)$ with $Ay_0 = By_1$ for some $y_1 \in X_D^{1-2\theta}$, we can write (see \eqref{eq:basic_delta})
\begin{align*}
s^2 m_3''(s) & = -2 (|s|^\theta\delta(s)B)^2 |s|^{2-2\theta}(is)^{-1} ( y_0 -  \delta(s) By_1)\\
& = 2i \sgn(s) (|s|^\theta\delta(s)B) (|s|^{1-\theta} \delta(s) B y_0)\\
&\qquad - 2i \sgn(s) (|s|^\theta\delta(s)B)^2 (|s|^{1-2\theta}\delta(s) By_1),
\end{align*}
and thereby showing the uniform boundedness of $s^2m_3''(s)$ on $\{|s|\geq 1\}$.

Finally, assume that $y_0\in \wt Z_D^{\theta}$. By \eqref{eq:deriv_of_delta}, we have
\[
|s|^3m_3'''(s)=i6\bigl(|s|^\theta\delta(s)B\bigr)^3|s|^{3-3\theta}\delta(s)By_0.
\]
Hence, the case $\theta\in[\frac23,1]$ follows immediately from the assumption on $y_0$.
 For $\theta \leq \frac23$, if $y_0\in \Dom(A)$ with $Ay_0 = By_1$ and $y_1\in X$, then by \eqref{eq:basic_delta} we have
\begin{align}\label{eq:delta_y1}
\nonumber |s|^3m_3'''(s) & = 6\sgn(s) \bigl(|s|^\theta\delta(s)B\bigr)^3|s|^{2-3\theta}(y_0 - \delta(s)By_1)\\
\nonumber & = 6\sgn(s) \bigl(|s|^\theta\delta(s)B\bigr)^2 |s|^{2-2\theta}\delta(s)By_0 \\
& \qquad - 6\sgn(s) \bigl(|s|^\theta\delta(s)B\bigr)^3|s|^{2-3\theta} \delta(s)By_1.
\end{align}
Therefore, since $2-2\theta\leq 1$ and $2-3\theta\leq \theta$ if $\theta\in [\frac12,\frac23]$, for every $y_0\in A^{-1}(BX)$ and $\theta \in [\frac12,\frac23]$, we obtain the uniform boundedness of $s^2m_3'''(s)$ on $\{|s|\geq 1\}$. Since 
\begin{align}\label{eq:delta_y0}
\nonumber\bigl(|s|^\theta\delta(s)  B\bigr)^2|s|^{2-2\theta}\delta(s)By_0 & = -i\sgn(s)\bigl(|s|^\theta\delta(s)B\bigr)|s|^{1-\theta}\delta(s)By_0\\
& \qquad + i\sgn(s)\bigl(|s|^\theta\delta(s)B\bigr)^2|s|^{1-2\theta}\delta(s)By_1
\end{align}
and $y_0\in X_D^1\subset X_D^{1-\theta}$, $y_1\in X_D^{2-3\theta}\subset X_D^{1-\theta}$ for every $\theta \in [\frac13,\frac12)$, we infer the desired uniform bound for $s^3m_3'''(s)$ on $\{|s|\geq 1\}$ if $y_0\in A^{-1}(BX_D^{2-3\theta})$ and $\theta\in [\frac13,\frac12)$. 
Finally, if $y_1\in A^{-1}\big(BX_D^{1-3\theta}\big)$ with $\theta\in\left[0,\frac13\right)$, similar arguments, applying \eqref{eq:basic_delta} to \eqref{eq:delta_y1} and \eqref{eq:delta_y0}, yield the desired bound and complete the proof of the first part of the lemma. 

For the second statement, since
$m_3'(s)=-i(\delta(s)B)^2y_0$, its uniform boundedness on
$\{|s|\geq1\}$ follows directly from \eqref{hyp:bound_2}, yielding
\eqref{cM3_for_M} with $\ell=1$ for every $\theta\in[0,1]$ and
$y_0\in Z_D^\theta$. Relying on \eqref{eq:deriv_of_delta}, we can write
\begin{equation}\label{cM3_M>1/3}
\begin{split}
|s|m_3''(s) & = -2|s|^{1-3\theta}
\bigl(|s|^\theta\delta(s)B\bigr)^3 y_0,\\
|s|^2m_3'''(s) & = 3i
\bigl(|s|^\theta\delta(s)B\bigr)^3
|s|^{2-3\theta}\delta(s)B y_0.
\end{split}
\end{equation}This shows that \eqref{cM3_for_M} holds for $\ell=2,3$, every
$\theta\in[1/3,1]$, and $y_0\in Z_D^\theta=X_D^{2-3\theta}$.

Let $\theta \in [0,1/3)$ and $y_0 \in Z^\theta_D = A^{-1}\big(BX^{1-3\theta}_D)$, i.e., $y_0\in \Dom(A)$ and $Ay_0 = B y_1$ for some $y_1 \in X^{1-3\theta}_D$. Combining \eqref{eq:basic_delta} and \eqref{cM3_M>1/3}, we obtain 
\begin{equation}
\begin{split}
|s|m_3''(s) & = -2 
\bigl(|s|^\theta\delta(s)B\bigr)^2 |s|^{1-2\theta} \delta(s)B y_0,\\
|s|^2m_3'''(s) & = 3 \sgn(s) \big(|s|^\theta \delta(s) B)^2 |s|^{1-2\theta} \delta(s)By_0\\
& \qquad 
 - 3 \sgn(s)  \big(|s|^\theta \delta(s)B)^3 |s|^{1-3\theta} \delta(s)By_1.
\end{split}
\end{equation}
Therefore, since $\Dom(A)\subset X^1_D$,  \eqref{cM3_for_M} also holds for $\ell=2,3$ when
$\theta\in[0,1/3)$ and $y_0\in Z^\theta_D$. \end{proof}

Below, we repeatedly use the following formulas; their verification is straightforward. For any $s \neq 0$, we have
\begin{equation}\label{formula:diriv_alpha}\begin{split}
\alpha'(s) & = \alpha(s)(2s - iB)\alpha(s)\\
\alpha''(s) & = 2\alpha'(s)(2s -iB)\alpha(s) + 2\alpha(s)^2\\
\alpha'''(s) & =
2\alpha''(s)(2s - iB) \alpha(s) + 6 \alpha'(s) \alpha(s)\\ 
& \qquad + 
2 \alpha'(s) (2s - iB) \alpha'(s) + 2 \alpha(s)\alpha'(s).\\ 
\end{split}
\end{equation}
The following auxiliary result is useful for the study of condition $(\cM_\infty^\ell)$.
Let $k \geq 0$ and $\ell \geq 1$ be integers, and consider an expression of the form $\pi(s)\alpha(s)$, $s \neq 0$, where $\pi(s) = \pi(s; k,\ell)$ is a product consisting of $s^k$, $\ell - 1$ factors of $\alpha(s)$, and finitely many factors of $B$, arranged in arbitrary order (we recall that $\alpha(s)$ and $B$ need not commute).
We say that $\pi(s)\alpha(s)$ is of \emph{type 0} if $k - \ell \leq 0$. Otherwise, that is, when $k - \ell > 0$, we say that $\pi(s)\alpha(s)$ is of \emph{type $k - \ell$}; see, for example, \eqref{def:J} below.

The following result provides conditions on $x \in X$ ensuring the uniform boundedness of $\pi(s)\alpha(s)x$ on $\{|s|\geq 1\}$, depending on the type of the expression $\pi(s)\alpha(s)$.

\begin{lemma}\label{lem:types}
Assume Hypothesis~\ref{hyp:main}. Then the following statements hold.
\begin{itemize}
\item[(i)] Any expression $\pi(s)\alpha(s)$ of type~0 is uniformly bounded on $\{|s|\geq 1\}$.

\item[(ii)] Let $\pi(s)\alpha(s)$ be of type~1. Then, for any $x \in X_2$, the function $s \mapsto \pi(s)\alpha(s)x$ is uniformly bounded on $\{|s|\geq 1\}$.

\item[(iii)] 
Let $\pi(s)\alpha(s)$ be of type~$2$ and suppose that $\pi(s)=\widetilde{\pi}(s)\alpha(s)B^j$ for some $\widetilde{\pi}(s)$ and $j\in\N_0$. Then $s\mapsto \pi(s)\alpha(s)x$ is uniformly bounded on $\{|s|\geq 1\}$ for all $x\in\Dom(A)$ such that $B x\in X_2$ and $B^j x\in X_2$.
\end{itemize}
\end{lemma}
\begin{proof}
(i) The proof follows directly from \eqref{pencils est} and the boundedness of $B$. 

(ii) If the expression $s^{-2}\pi(s)$ occurs $s^{\ell-1}$ and $\ell-1$ factors of $\alpha(s)$ then by \eqref{pencils est} it is bounded on $\R\setminus\{0\}$. Since  
\[
\pi(s)\alpha(s)x = s^{-2}\pi(s)\,[s^2 \alpha(s)x],
\]
taking $x\in X_2$ yields the desired bound. 

(iii) We first note that 
\begin{equation}\label{eq: t2 alpha incl}
s^2 \alpha(s) \supset \alpha(s)(isB + A) - I,
\end{equation}
which for $x \in \Dom(A)$ yields
\[
\pi(s)\alpha(s)x = i s^{-1} \pi(s)\alpha(s)Bx + s^{-2} \pi(s)\alpha(s)Ax - s^{-2}\pi(s)x.
\]
Clearly, $s^{-1}\pi(s)\alpha(s)$ is of type~1. Hence, by part~(ii), the term $s^{-1}\pi(s)\alpha(s)Bx$ is uniformly bounded on $\{|s|\geq 1\}$, provided $Bx \in X_2$. 
Moreover, $s^{-2}\pi(s)\alpha(s)$ is of type~0 and thus uniformly bounded on $\R\setminus\{0\}$ by part~(i). It remains to consider $s^{-2}\pi(s)x$. Writing
\[
s^{-2}\pi(s) = s^{-2}\widetilde{\pi}(s)\alpha(s)B^j,
\]
we observe that $s^{-2}\widetilde{\pi}(s)\alpha(s)$ is of type~1. Therefore, by part~(ii), $s^{-2}\pi(s)x$ is uniformly bounded on $\{|s|\geq 1\}$ provided $B^j x \in X_2$.
\end{proof}
Expressions with $k - \ell = 3$ are handled by iterating \eqref{eq: t2 alpha incl}, thereby reducing them to the lower-order types. Since this situation occurs only once (see~\eqref{theta1}), we do not treat it right away.

\medskip

{\it Proof of Theorem \ref{thm 1}\emph{(i)}.} By Lemma~\ref{lem:reduction_step}(1) and Lemma \ref{lem:m3_in_Mikhlin} it suffices to show that $m_1, m_2 \in \cM_\infty^2(X)$ provided $x_0 \in \Dom(A)$ and $Bx_0,\, By_0\in X_2$ (the assumption $y_0\in Y_D^\theta$ is needed only to guarantee that $m_3\in\mathcal{M}^2_\infty(X))$). 
Using the first identity in \eqref{formula:diriv_alpha}, we easily arrive at
\begin{equation}\label{eq:m_in_M1}
\begin{split}
s m_1'(s) & = is \alpha(s)x_0 +  2i s\alpha(s) [s^2\alpha(s)x_0] +  s^2 \alpha(s) B \alpha(s) x_0,\\
s m_2'(s) & =  2 (s\alpha(s))^2 B y_0 - is (\alpha(s)B)^2 y_0.
\end{split}
\end{equation}
Thus, since $x_0\in \Dom(A)\subset X_2$, see Lemma \ref{lem:reduction_step}(3), by parts (i) and (ii) of Lemma \ref{lem:types}, we obtain $m_1, m_2\in \cM_\infty^1(X)$. Further, we note that
\[
s^2 m_1''(s) = 2 s m_1'(s) - 2 is \alpha(s)x_0 + is^3 \alpha''(s)x_0.
\]
By \eqref{formula:diriv_alpha}, we can write 
\begin{equation}\label{def:gamma}
\gamma(s): = s^3 \alpha''(s) =  2  \gamma_1(s) - 2 \gamma_2(s) - 4i \gamma_3(s) - 4i \gamma_4(s) +  8 \gamma_5(s),
\end{equation}
where 
\begin{equation}\label{def:J}
\begin{tabular}{lll}
$\gamma_1(s) := s^3\alpha(s)^2$, 
& $\gamma_2(s) := s^3(\alpha(s)B)^2 \alpha(s)$, 
& $\gamma_3(s) := s^4 \alpha(s) B \alpha(s)^2$, \\[0.3em]
$\gamma_4(s) := s^4 \alpha(s)^2 B \alpha(s)$, 
& $\gamma_5(s) := s^5 \alpha(s)^3$, 
& $s \neq 0$.
\end{tabular}
\end{equation}
Clearly, $\gamma_2(s)$ is of type~0, while $\gamma_1(s)$, $\gamma_3(s)$, and $\gamma_4(s)$ are of type~1. Hence, since $x_0 \in X_2$, parts~(i) and~(ii) of Lemma~\ref{lem:types} imply that $\gamma_j(s)x_0$ is uniformly bounded on $\{|s|\geq 1\}$ for all $j=1,\dots,4$.
Since $\gamma_5(s)$ is of type 2 and $x_0\in \Dom(A)$ with $Bx_0 \in X_2$, by Lemma \ref{lem:types} (iii) with $j=0$,  we infer that 
$\gamma_5(s)x_0$ is uniformly bounded on $\{|s|\geq 1\}$. Thus,  $m_1\in \cM_\infty^2(X)$. 
Finally, for $m_2$ we have that
\begin{align*}
s^2 m_2''(s) &  =  s^2 \alpha''(s)By_0  = -i s^{-1} \gamma(s)By_0. 
\end{align*}
Clearly, $s^{-1}\gamma_j(s)$, $j=1,\dots,4$, are of type~0. Hence, $s^{-1}\gamma_j(s)By_0$ is uniformly bounded on $\{|s|\geq 1\}$ for any $y_0 \in X$. 
Since $s^{-1}\gamma_5(s)$ is of type~1, by Lemma \ref{lem:types} the term $s^{-1}\gamma_5(s)By_0$ is uniformly bounded on $\{|s|\geq 1\}$ provided $By_0 \in X_2$. 
This proves that $m_2 \in \cM_\infty^2(X)$ and completes the proof of Theorem \ref{thm 1}(i).
\medskip

{\it Proof of Theorem \ref{thm 1}\emph{(iii)}.}
Again, by Lemma~\ref{lem:reduction_step}(1) and Lemma \ref{lem:m3_in_Mikhlin}, it suffices to show that $m_1, m_2 \in \mathcal{M}_\infty^3(X)$, provided that $x_0$ and $y_0$ satisfy the conditions in \eqref{cond_kinetic}. From Theorem \ref{thm 1}(i) we already know that $m_1,m_2$ are in $\mathcal{M}_\infty^2(X)$. Thus it suffices to prove that $s^3m_j'''(s)$, $j=1,2$, are uniformly bounded on $\{|s|\geq 1\}$. 
For $m_1$, we note that
\begin{equation}\label{s3m1'''}
s^3 m_1'''(s) = 3i s^3 \alpha''(s) x_0 + i s^4 \alpha'''(s)x_0 = 3i \gamma(s)x_0  + i s^4 \alpha'''(s)x_0, \qquad s\neq 0,
\end{equation}
where $\gamma$ is defined in \eqref{def:gamma}, and the uniform boundedness of $\gamma(s) x_0$ on $\{|s| \ge 1\}$ was already established in the proof of part~(i). 
For the second term, by \eqref{formula:diriv_alpha},  we observe that 
\begin{equation}\label{beta}\begin{split}
\beta(s) & := s^4 \alpha'''(s)\\
 &  =  6s^4 \alpha'(s)\alpha(s) + 2 s^4 \alpha(s)\alpha'(s)\\
& \quad  + 2 s^4 \alpha'(s)(2s -i B)\alpha'(s) + 2s^4 \alpha''(s) (2s - iB)\alpha(s)\\
& =: 6\beta_1(s) + 2\beta_2(s) + 2\beta_3(s) +  2\beta_4(s). 
\end{split}
\end{equation}
Using \eqref{formula:diriv_alpha}, we have
\begin{equation}\label{beta12}\begin{split}
\beta_1(s) & = s^4 \alpha'(s)\alpha(s) = 2 s^5 \alpha(s)^3 - i s^4 \alpha(s)B \alpha(s)^2 = 2 \gamma_5(s) -i \gamma_3(s),\\
\beta_2(s) & = s^4 \alpha(s)\alpha'(s) = 2 s^5 \alpha(s)^3 -  i s^4 \alpha(s)^2 B \alpha(s) = 2 \gamma_5(s) +i \gamma_4(s).
\end{split}
\end{equation}
Thus the uniform boundedness of $\beta_1(s)x_0$ and $\beta_2(s)x_0$ on $\{|s|\geq 1\}$, follows from part (i) provided $x_0\in \Dom(A)$ with $Bx_0\in X_2$. For $\beta_3(s)x_0$, we note that
\begin{align}\label{beta3}
\beta_3(s) & = 2 s^5 \alpha'(s)\alpha'(s) - is^4 \alpha'(s)B\alpha'(s) =: 2 \vartheta_1(s) - i \vartheta_2(s),
\end{align}
where
\begin{align}\label{theta1}
\nonumber \vartheta_1(s) & = 4 s^7\alpha(s)^4 - 2i s^6 \alpha(s)^3B\alpha(s) - 2i s^6\alpha(s)B\alpha(s)^3 - s^5\alpha(s)B\alpha(s)^2 B \alpha(s)\\
& =: 4\vartheta_{1,1}(s) -2i \vartheta_{1,2}(s) -2i \vartheta_{1,3}(s) - \vartheta_{1,4}(s). 
\end{align}
The expression $s^7\alpha(s)^4$ is of type~3. As $x_0, Bx_0 \in \Dom(A)$, we iterate \eqref{eq: t2 alpha incl} to obtain  
\begin{align}\label{theta11}
\vartheta_{1,1}(s) x_0 & = s^5 \alpha(s)^4 (isB  +A)x_0 -  s^5\alpha(s)^3x_0\\
\nonumber & = is^4 \alpha(s)^3\big(\alpha(s)(isB  +A) - I\big)Bx_0 + s^5 \alpha(s)^4 Ax_0 -   \gamma_5(s)x_0\\
\nonumber & = - (s^3 \alpha(s))^3 [s^2\alpha(s)B^2 x_0] - is^4 \alpha(s)^3Bx_0 +   i(s\alpha(s))^4 ABx_0\\
& \quad +  (s \alpha(s))^3 [s^2\alpha(s)Ax_0] - \gamma_5(s)x_0.
\end{align}
Hence, since $x_0, Bx_0 \in \Dom(A)$ and $Ax_0, B^2x_0\in X_2$, it follows from that $\vartheta_{1,1}(s) x_0$ is uniformly bounded on $\{|s|\geq 1\}$. 
Clearly, all $\vartheta_{1,j}x_0$, $j=2,3,4$, are of type at most $2$. Since $x_0 \in \Dom(A)$ with $Bx_0\in X_2$, again by Lemma \ref{lem:types} we infer that $\vartheta_{1,j}(s)x_0$, $j=2,3,4$, are uniformly bounded on $\{|s|\geq 1\}$.   

Further, applying \eqref{formula:diriv_alpha}, we obtain 
\begin{align*}
\vartheta_2(s) &
 = 4 s^6 \alpha(s)^2 B \alpha(s)^2 - 2i s^5 \alpha(s)^2 (B \alpha(s))^2 -
2i s^5 (\alpha(s)B)^2 \alpha(s)^2 - s^4 (\alpha(s)B)^3\alpha(s)\\
& =: 4\vartheta_{2,1}(s) - 2i \vartheta_{2,2}(s) - 2i \vartheta_{2,3}(s) - \vartheta_{2,4}(s). 
\end{align*}
The term $\vartheta_{2,1}(s)$ is of type 2. Since $x_0\in \Dom(A)$ with $Bx_0\in X_2$, Lemma \ref{lem:types}(iii) (with $j=0$) yields the uniform boundedness of $\vartheta_{2,1}(s)x_0$ on $\{|s|\geq 1\}$. 
Similarly, by Lemma \ref{lem:types}(i) and (ii) we obtain the uniform boundedness of $\vartheta_{2,j}(s)x_0$ on $\{|s|\geq 1\}$ for any $j=2,3,4$.

Ultimately, we are left with showing boundedness of $\beta_4(s)x_0$. We note that 
\begin{equation}\label{beta4}
\beta_4(s) = s^4 \alpha''(s) (2s - iB)\alpha(s) = - 2i s^2 \gamma(s) \alpha(s) - s \gamma(s)B\alpha(s).   
\end{equation}
From the proof of part~(i) we recall that $\gamma(s)$ is a linear combination of $\gamma_j(s)$, $j=1,\dots,5$; see \eqref{def:J}. Hence, we may write 
\begin{equation}\label{beta41}\begin{split}
s^2\gamma(s)\alpha(s) & =  2i  s^5 \alpha(s)^3 - 2i  s^2 (\alpha(s)B)^2 \alpha(s)^2 + 4 s^6 \alpha(s) B \alpha(s)^3 \\
& \quad +  4 s^6 \alpha(s)^2 B \alpha(s)^2 +  8 i s^7 \alpha(s)^4,\\
& = 2i  \gamma_5(s) - 2i  s^5 (\alpha(s)B)^2 \alpha(s)^2 + 4 \vartheta_{1,3}(s) \\
& \quad +  4 s^6 \alpha(s)^2 B \alpha(s)^2 +  8 i \vartheta_{1,1}(s),
\end{split}
\end{equation}
Note that the boundedness of the new expressions, i.e., $s^2 (\alpha(s)B)^2 \alpha(s)^2x_0$ and $s^6 \alpha(s)^2 B \alpha(s)^2x_0$ follows from part (i) and (ii)(with $j=0$) of Lemma \ref{lem:types}, respectively. 
Further, for the second term in \eqref{beta4}, we have that
\begin{equation}\label{beta42}\begin{split}
s\gamma(s) B\alpha(s) & =     2i  s^4 \alpha(s)^2 B \alpha(s)  - 2i s (s\alpha(s)B)^3 \alpha(s) + 4 s^5 \alpha(s) B \alpha(s)^2 B \alpha(s)\\
& \quad +  4 s^5 \alpha(s) (\alpha(s) B)^2 \alpha(s) 
+  8 i s^6 \alpha(s)^3 B\alpha(s).
\end{split}
\end{equation}
The first term is of type~0, and the next three are of type~1. Hence, since $x_0 \in X_2$, parts~(i) and~(ii) of Lemma~\ref{lem:types} apply. 
For the fifth term, which is of type~2, we apply part~(iii) of Lemma~\ref{lem:types} with $j=1$ to obtain the desired bound for $s^6 \alpha(s)^3 B\alpha(s)x_0$. This shows that $m_1 \in \cM_\infty^3(X)$.

Finally, by \eqref{beta}, we have
\begin{equation}\label{s3alpha3}
s^3 m_2'''(s) = s^3 \alpha'''(s) By_0 = s^{-1}\beta(s) By_0.
\end{equation}
It is readily seen from the above representation of $\beta(s)$ (see \eqref{beta12}-\eqref{beta42}) that the only expression of type~3 is $\vartheta_{1,1}(s)$, which appears in the representations of $\beta_3(s)$ and $\beta_4(s)$; see \eqref{theta1} and \eqref{beta41}. Consequently, all remaining terms are of type at most~2 and, upon multiplication by $s^{-1}$, become expressions of type at most~1. Hence, by Lemma~\ref{lem:types}(i) and~(ii), we obtain the desired bound, provided that $By_0 \in X_2$.

Making use of $By_0 \in \Dom(A)$ and applying \eqref{eq: t2 alpha incl}, we can write (cf.~\eqref{theta11})
\[
s^{-1} \vartheta_{1,1}(s) By_0 = i s^5 \alpha(s)^4 B^2 y_0 + (s \alpha(s))^4 A B y_0 - s^4 \alpha(s)^3 B y_0.
\]
Therefore, by Lemma~\ref{lem:types}(i) and~(ii), and since $\Dom(A) \subset X_2$ and $B^2 y_0 \in X_2$, it follows that $s^{-1} \vartheta_{1,1}(s) By_0$ is uniformly bounded on $\{|s|\geq 1\}$. Hence, $m_2 \in \mathcal{M}_\infty^3(X)$, which completes the proof of Theorem \ref{thm 1}(iii). 
\medskip

{\it Proof of Theorem \ref{thm 1}\emph{(ii)}.} Since the conditions on $x_0$ and $y_0$ in (ii) imply those in (i), it follows from part~(i) and Lemma \ref{lem:reduction_step}(2) that $M \in \mathcal{M}_\infty^2(X)$. It remains to verify that $M \in \mathcal{M}_0^3(X)$ and $\sup_{|s| \geq 1} \| s^3 M'''(s) \| < \infty$. 

We begin with the proof of $M \in \mathcal{M}_0^3(X)$ and make use of formula \eqref{def of m3}, cf. \eqref{a i d}. For all $z\in \Dom(A)$ we have 
\begin{align}\label{cond:a_range}
\alpha(s)Az & =  z -  a(s)(is+B)z\quad  \textrm{ and } \quad
\delta(s)Az  = z- d(s)z.
\end{align}
Since $a$ and $d$ are uniformly bounded by \eqref{pencils est}, it follows from \eqref{cond:a_range} that for any $y \in \Ran(A)$ the functions
\begin{equation}\label{range A}
s\mapsto \alpha(s)y, \, \delta(s)y \quad \textrm{ are uniformly bounded on   } \{0<|s|\leq 1\}.
\end{equation}

In particular, the assumptions in~(ii) ensure that \eqref{range A} holds for $y = x_0$ and $y = By_0$.

Thus, since $
M(s) = (is)^{-1} m(s) = \alpha(s)x_0 - is\, \alpha(s)\delta(s) By_0$ (see \eqref{def of m3}), $M$ is uniformly bounded on $\{0 < |s| \leq 1\}$.
Moreover, it is readily seen that
\begin{align}\label{M'''}
\nonumber sM'(s) & = - M(s) - i m'(s),\\
\nonumber s^2M''(s) & = - 2 s M'(s) - i s\,m''(s),\\
s^3M'''(s) & = - 3 s^2 M''(s) - is^2 m'''(s), \quad s\neq 0.
\end{align}
Therefore, to show $M\in \cM^3_0(X)$ it suffices to show that for each $\ell\in \{1,2,3\}$,
\begin{equation}\label{eq:ell_123}
\sup_{0<|s|\leq 1} \|s^{\ell-1}m^{(\ell)}(s)\|<\infty.
\end{equation}
For $\ell=1$ we note that, by \eqref{def of m2}, \[
m'(s) = a'(s)x_0 - a'(s) d(s) y_0- a(s) d'(s) y_0 . 
\]
Moreover, for any $s\neq 0$, using \eqref{formula:diriv_alpha},
\begin{equation}\label{a'}
a'(s)x_0  = i\alpha (s)x_0 + is \alpha'(s)x_0 = i [\alpha(s) x_0]  - a (s) (-2s + iB) [\alpha(s) x_0]
\end{equation}
\begin{align*}
d'(s)y_0 & = i\delta(s)By_0 + is \delta'(s)By_0 = i[\delta(s) B y_0] -i d(s)[\delta(s)By_0]\\
a'(s)d(s)y_0 & = [is a'(s)] [\delta(s) By_0].   
\end{align*}
By Lemma~\ref{lem on pencil est} we, in particular, know that $a \in \mathcal{M}_0^1(\mathcal{L}(X))$, thus we infer from \eqref{range A} that \eqref{eq:ell_123} holds for $\ell=1$, i.e.,  $M \in \mathcal{M}_0^1(X)$.
For $\ell =2$, we first observe that
\[
s m''(s) = s a''(s) x_0 - s \left(a''(s) d(s) + 2a'(s) d'(s) + a(s) d''(s)\right)y_0.
\]
Since 
\begin{equation}\label{a''}
s a''(s) = is^2  \alpha''(s) + 2 i s \alpha'(s),
\end{equation}
combining \eqref{formula:diriv_alpha} with \eqref{cond:a_range}, we obtain that $s a''(s) x_0$ is bounded uniformly on $\{0<|s|\leq 1\}$. 
 Moreover, the same argument as in the case $\ell = 1$ yields the uniform boundedness of  $s a'(s)\big[d'(s)y_0\big]$ and $s a''(s) d(s) y_0 = i s^2 a''(s) [\delta(s)By_0]$ on $\{0<|s|\leq 1\}$. For the last term, using \eqref{formula:diriv_alpha}, we have
\[
s d''(s)y_0 = is^2 \delta''(s)By_0 + 2i s \delta'(s)By_0 =  - 2i (s \delta(s)B)^2 [\delta(s)By_0] + 2 s\delta(s)B [\delta(s)By_0].
\]
Thus $s a(s) d''(s)y_0$ is uniformly bounded on $\{0<|s|\leq 1\}$ as well. 
This proves \eqref{eq:ell_123} for $\ell=2$. 
For $\ell =3$, we observe that 
\begin{equation}\label{s2m'''}
s^2 m'''(s) = s^2a'''(s)x_0 - s^2\left( a'''(s) d(s) + 3 a''(s) d'(s) + 3 a'(s) d''(s) + a(s) d'''(s) \right) y_0,
\end{equation}
For $s^2a'''(s)x_0$, using \eqref{a'} and \eqref{a''}, we have
\[
s^2 a'''(s) = i s^3 \alpha'''(s) + 3s a''(s) - 6 a'(s) + 6i \alpha(s).
\]
Hence, to bound $s^2 a'''(s)x_0$ on $\{0 < |s| \leq 1\}$, it suffices to bound $s^3 \alpha'''(s)x_0$ on the same region. Using \eqref{formula:diriv_alpha}, we set
\begin{align*}
s^3 \alpha'''(s) & =  s^3\alpha''(s)(2s - iB) \alpha(s) + 6 s^3\alpha'(s) \alpha(s)\\
&\qquad + 2 s^3\alpha'(s) (2s - iB) \alpha'(s) + s^3 \alpha(s) \alpha'(s)\\
& =:  \varrho_1(s)+ 6 \varrho_2(s) + 2 \varrho_3(s) +2 \varrho_4(s).
\end{align*}
Since $s^3 \alpha''(s)$ is uniformly bounded on $\{0 < |s| \leq 1\}$ (see \eqref{loc est alpha}), we infer from \eqref{cond:a_range} the desired bound for $\varrho_1(s)x_0$. Using \eqref{formula:diriv_alpha}, we have
\[
\varrho_3(s)x_0 = \big(s\alpha(s)(2s - iB)\big)^3 [\alpha(s)x_0],
\]
so \eqref{cond:a_range} and \eqref{pencils est} yield uniform boundedness on $\{0 < |s| \leq 1\}$. Finally, the terms $\varrho_2(s)$ and $\varrho_4(s)$ are uniformly bounded on $\{0 < |s| \leq 1\}$ by \eqref{pencils est}.

 For the remaining terms in the expression for $s^2 m'''(s)$ in \eqref{s2m'''}, it suffices to observe that
\[
s^2 \left(a'''(s) d(s) + 3 a''(s) d'(s) + 3 a'(s) d''(s)\right) y_0
\]
is uniformly bounded on $\{0<|s|\leq 1\}$, which follows from $a \in \mathcal{M}_0^3(X)$ (see Lemma~\ref{lem on pencil est}(i)) and the estimates established above.  For the term $s^2a(s)d'''(s)y_0$, 
it is readily seen that
\[
\delta'''(s) =  6 i (\delta(s)B)^{3}\delta(s)\quad \textrm{ and } \quad
d'''(s) =  6i (\delta(s) B)^{2}\left(d(s) - I\right)\delta(s)B. 
\]
Combining this with \eqref{cond:a_range}, we obtain the desired boundedness of this term. Consequently, this completes the proof of $M \in \mathcal{M}_0^3(X)$. 
Therefore, it remains to show that $\sup_{|s| \geq 1} \| s^3 M'''(s) \| < \infty$, which, by \eqref{M'''}, reduces to $\sup_{|s| \geq 1} \| s^2 m'''(s) \| < \infty$, and by \eqref{def of m2} $\sup_{|s| \geq 1} \| s^2 m_j'''(s) \| < \infty$ for each $j=1,2,3$. Since we assume that $y_0\in  Z_D^\theta$, the case $j=3$ follows from the second claim in Lemma \ref{lem:m3_in_Mikhlin}. Therefore, it remains to consider $j=1,2$. 

To this end, we appeal to the formulas for $m'''(s)$ provided in the proof of Theorem \ref{thm 1}(iii) and verify that, under the assumptions on $x_0$ and $y_0$ in Theorem \ref{thm 1}(i), those formulas yield $\sup_{|s| \geq 1} \| s^2 m_j'''(s) \| < \infty$ for each $j=1,2$.
We again make use of the formula \eqref{def of m2}; see also \eqref{def: m_123}. By \eqref{s3m1'''}, we have 
\[
s^2m'''_1(s) = 3is^{-1} \gamma(s)x_0 + i s^3\alpha'''(s)x_0.
\]
As in the proof of the boundedness of $s^2 m_2''(s)$ in Theorem \ref{thm 1}(i), we note that $s^{-1}\gamma(s)$ is a linear combination of expressions of type at most~1. Hence, $x_0 \in X_2$ implies that $s^{-1}\gamma(s)x_0$ is uniformly bounded on $\{|s|\geq 1\}$.

Furthermore, since $s^3 \alpha'''(s) = s^{-1} \beta(s)$, the proof of the boundedness of $s^3 \alpha'''(s)x_0$ follows exactly along the same lines as the proof of the boundedness of \eqref{s3alpha3}, with $By_0$ replaced by $x_0$. This shows that the condition $x_0 \in \Dom(A)$ with $Bx_0 \in X_2$ suffices for the uniform boundedness of $s^2m_1'''(s)x_0$ on $\{|s|\geq 1\}$.

Finally, we first note that 
\[
s^2 m_2'''(s) = s^2 \alpha'''(s)By_0 = s^{-2}\beta(s)By_0
\]
Using the observation on the order type of the individual terms appearing in the proof of the boundedness of \eqref{s3alpha3} in part~(iii), we infer that $s^{-2}\beta(s)$ can be written as a linear combination of expressions of type at most~1 (in particular, the term $\vartheta_{1,1}(s)$ is of this type). Hence, the condition $Bx_0 \in X_2$ suffices to ensure the uniform boundedness of $s^2 m_2'''(s)$ on $\{|s|\geq 1\}$. This completes the proof of Theorem \ref{thm 1}(ii).

\smallskip

We conclude with the following direct consequence of the proof of
Theorem~\ref{thm 1}.
\begin{corollary}\label{cor:decay_of_mH}
Let $A$ and $B$ satisfy Hypothesis~\ref{hyp:main}. For any $x_0,y_0\in X$, define 
\[
m_H(s):=R(is)(isx_0+By_0)
\quad\text{and}\quad
m_P(s):=R_D(is)By_0,
\quad s\neq0.
\]
If $x_0\in\Dom(A)$, $Bx_0,By_0\in X_2$, and $By_0\in\Ran(A)$, then
$m_H\in\mathcal{M}^2(X)$. Consequently,
\[
\mathcal{F}^{-1}m_H\in C(\mathbb{R}\setminus\{0\};X) \qquad \textrm{ and } \quad
\|(\mathcal{F}^{-1}m_H)(t)\|\lesssim t^{-1},
\quad t>0.
\] 

In particular, under the assumptions of Theorem~\ref{thm 1}(ii), we have
$m_H,m_P\in\mathcal M^2(X)$ and 
\[
\|\cF^{-1}m_H (t) - \cF^{-1}m_P(t)\|\lesssim 1/t^2\,  \textrm{ and } \, \|\cF^{-1}m_H(t)\|,\,  \|\cF^{-1}m_P(t)\| \lesssim 1/t, \quad t>0. 
\]
\end{corollary}

\begin{proof}
Using the notation introduced in \eqref{def: m_123}, we have
$m_H(s)=m_1(s)+m_2(s)$ for $s\neq0$. If $x_0\in\Dom(A)$ and $Bx_0,By_0\in X_2$, then the proof of Theorem~\ref{thm 1}(i) shows that
$m_1,m_2\in\mathcal{M}_\infty^2(X)$. By Lemma~\ref{lem on pencil est}(i),
$m_1=a(\cdot)x_0\in\mathcal{M}_0^2(X)$, since $a\in\mathcal{M}_0^n(\mathcal{L}(X))$ for every $n\in\mathbb{N}$. Moreover, the assumption $By_0\in\Ran(A)$ implies that there exists $z\in\Dom(A)$ such that
\[
m_2(s)=z-a(s)(B+is)z,
\qquad s\neq0,
\]
which implies $m_2\in\mathcal{M}_0^2(X)$. Therefore, $m_H\in\mathcal{M}^2(X)$, and the conclusion follows from Lemma~\ref{CZ cond}. 

For the last statement, the proof of Theorem~\ref{thm 1}(i) shows
that $m\in\cM^2(X)$, and hence $m_P=m_H - m \in\cM^2(X)$. Therefore, Theorem~\ref{thm 1}(ii) and Lemma~\ref{CZ cond} yield the desired estimates.
\end{proof}

\subsection{On the decay rate $O(1/t^{\ell})$ with $\ell \geq 3$}\label{subs:limit-of-decay}
One could naturally try to apply an idea similar to that used in the proof of Theorem~\ref{thm 1}(ii) by considering the $\ell$-fold antiderivative of the function $\mathcal{F}^{-1} m = k$.
That is, assuming that $k \in C(\mathbb{R}_+;X)$ and defining
\[
K_\ell(t) := \int_0^t k(s)\,\frac{(t-s)^{\ell-1}}{(\ell-1)!}\,\ud s,
\quad t\geq 0, \quad \textrm{ and } \quad K_\ell(t) : = 0, \quad t<0,
\]
Lemma~\ref{CZ cond} implies that if $M_\ell := \mathcal{F} K_\ell$ belongs to the class $\mathcal{M}^{\ell+2}(X)$, then 
\[
\|K^{(\ell)}(t)\| = \|k(t)\| \lesssim  {[M_\ell]_{\m_{\ell+2}}}/{t^{\ell+1}},  \quad t >0.  
\]
 If the distribution $\mathcal F^{-1}m$ is precisely the extension by zero
of its positive time representative $k$, then
$M_\ell(s)=m(s)/(is)^\ell$ for $s\ne0$.
In this case a necessary first step is to determine whether, for
$\ell\ge2$, the quotient $m(s)/s^\ell$ is bounded on near $s=0$.

The following considerations show that, under the assumptions of
Theorem~\ref{thm 1}, an estimate of the form
\begin{equation}\label{faster decay}
\sup_{0<|s|\leq 1} \left\|s^{-2} m(s;x_0,y_0)\right\| < \infty
\end{equation}
is very restrictive with respect to the choice of vectors
$x_0, y_0 \in X$. Note that if, in addition to \eqref{faster decay}, $M_2=m/(is)^2$ satisfies the remaining
assumptions of Lemma~\ref{CZ cond} for $\mathcal M^4(X)$, then for
$k=\mathcal F^{-1}m=(\mathcal F^{-1}M_2)''$ one has
$\|k(t)\|\lesssim t^{-3}$. 

\begin{proposition}\label{prop: sharp}
Assume Hypothesis \ref{hyp:main} with $\theta =1$. If, in addition,  $X$ is reflexive, and the operators $A$ and $B$ are injective, then the following assertions hold.
\begin{itemize}
\item [(i)] For any $x_0\in \overline{\Ran(B)}$ and $y_0 \in X$ 
\[
m(s; x_0,y_0) \rightarrow 0 \in X \qquad \textrm{ as } |s|\rightarrow 0.
\]

\item [(ii)] For any $x_0 \in \Ran (A)$ and $y_0\in X$ with $By_0 \in \Ran (A)$ and $A^{-1} By_0 \in \overline{\Ran(B)}$ 
\[
(is)^{-1} m(s; x_0, y_0) \rightarrow A^{-1} x_0 \qquad \textrm{ as } |s|\rightarrow 0.
\]

Moreover, if the range of $B$ is dense in $X$, then for any $x_0 \in \Ran (A) \setminus \{0\}$ and $y_0 \in X$ such that $By_0\in \Ran (A)$ we have 

\begin{equation}\label{asymp of d}
\left\||s|^{-\eta} m (s; x_0,y_0)\right\| \rightarrow \left\{ 
\begin{array}{ll} 
\|A^{-1} x_0\|, & \textrm{ if } \eta = 1, \\
\infty, & \textrm{ if } \eta > 1,\\
0, & \textrm{ if } \eta< 1,
\end{array}
\right.\qquad \textrm{ as } |s|\rightarrow 0.
\end{equation}

\item [(iii)] For any $y_0\in X$ with $By_0 \in \Ran (A^2)$ such that  $BA^{-1}By_0 \in \Ran (A)$ and $A^{-1} BA^{-1} B y_0 \in \overline{\Ran(B)}$ we have
\[
s^{-2} m(s; 0, y_0) \rightarrow A^{-2} B y_0 \qquad \textrm{ as } |s|\rightarrow 0.
\]
\end{itemize}
\end{proposition}

\begin{proof} We start with introducing some auxiliary notation. 
For $x_0, x_1\in X$ we set for $\la \in \rho(P)\cap \rho(P_D)$
\[
\beta (\la): = \beta(\la; x_0, x_1) := R(\la)(\la x_0 + Bx_0 + x_1) - R_D(\la)(Bx_0 + x_1). 
\]
Clearly, if  $i\R\setminus \{0\} \subset \rho(P)\cap \rho(P_D)$ and $x_1 = By_0 - Bx_0$ for some $y_0\in X$, then 
\[
m(s ; x_0, y_0) = \beta (is; x_0,x_1), \quad s\neq 0.
\]
Since $[R]<\infty$, by Lemma \ref{lem on pencil est}(i) and  Lemma \ref{lem near 0}, there exists $\varepsilon>0$ and $\omega \in (0,\frac{\pi}{2})$ such that $U:= G_\omega \cap \{|\la|<\varepsilon\}\subset \rho(P) \cap \rho(P_D)$  and 
\[
\sup_{\la \in U} \|\la^2 R_D(\la)\| \leq {1}/{2}.
\] 
Moreover, note that for any $\la\in U$ 
\begin{equation}\label{3.49}
R(\la) = R_D(\la)\left( I + \la^2R_D(\la)\right)^{-1}.
\end{equation}

Thus, if $x_1 = By_0 - Bx_0$ with $y_0\in X$, we can write (see \eqref{def of m3}) 
\begin{align*}
\beta(\la) & = \left( I + \la^2 R_D(\la) \right)^{-1} \left[ \la R_D(\la) x_0 - \la^2 R_D(\la)^2 B y_0 \right], \qquad \la \in U.
\end{align*}
Clearly, 
\begin{equation}\label{conv uniform}
\left( I + \la^2 R_D(\la) \right)^{-1} \rightarrow I \quad  \textrm{ as } |\la|\rightarrow 0, \, \la \in U,
\end{equation} with convergence in the operator norm.
We claim that
\begin{equation}\label{conv on cl BX}
 \lambda R_D(\lambda)x\longrightarrow0
 \qquad (x\in\operatorname{Ran}(B))
\end{equation}
as $|\lambda|\to0$, $\lambda\in U$.  For $y\in X$, set
$z_\lambda:=\lambda R_D(\lambda)By$.  By Lemma~\ref{lem on pencil est}(ii),
the family $\{z_\lambda\}_{\lambda\in U}$ is bounded.  Moreover,
\[
 Az_\lambda=\lambda B(y-z_\lambda)\rightarrow0,
\]
as $\lambda \to 0$ in $U$. Since $X$ is reflexive and $A$ is closed and injective,
$z_\lambda \to 0$ weakly as $\lambda \to 0$ in $U$.  Consequently,
\[
 AR_D(\lambda)By=By-Bz_\lambda \to By \quad \text{weakly as}\quad\lambda \to0, \lambda \in U.
\]
As $AR_D(\lambda)By\in\operatorname{Ran}(A)$, we obtain
$\operatorname{Ran}(B)\subset\overline{\operatorname{Ran}(A)}$.

On the other hand, \eqref{3.49} gives
\[
 \lambda R_D(\lambda)
 =\lambda R(\lambda)(I+\lambda^2R_D(\lambda)),
\]
so $\{\lambda R_D(\lambda)\}_{\lambda\in U}$ is uniformly
bounded.  For $x=Aw\in\operatorname{Ran}(A)$,
\[
 \lim_{\lambda \to 0, \lambda \in U}\lambda R_D(\lambda)x
 =\lambda w-\lambda^2R_D(\lambda)Bw=0.
\]
By uniform boundedness, the same limit holds on
$\overline{\operatorname{Ran}(A)}$, and hence on
$\overline{\operatorname{Ran}(B)}$.  This proves \eqref{conv on cl BX}.

For the proof of (i), we note that for any $y_0\in X$ 
\[
\la^2 R_D(\la)^2 By_0 = (\la R_D(\la))( \la R_D(\la) By_0) \rightarrow 0 \quad \textrm{ as } |\la|\rightarrow 0, \la \in U.
\]
Consequently, by \eqref{conv uniform} and \eqref{conv on cl BX}, for any $x_0\in \overline{\Ran(B)}$ and $y_0\in X$, 
\begin{equation}\label{conv delta in U}
\beta(\la)  \rightarrow 0 \qquad \textrm{ as } |\la|\rightarrow 0, \, \la \in U.
\end{equation}
This proves (i).
For (ii), we first note that for $x_0 \in \Ran (A)$, if $w\in \Dom(A)$ and $Aw = x_0$, then 
\[
R_D(\la)x_0 = R_D(\la) A w = w - \la R_D(\la)Bw \rightarrow w\quad \textrm{ as } |\la| \rightarrow 0, \, \la \in U.
\]
Moreover, if $y_0 \in X$ is such that $By_0 = A y$ with $y \in \Dom(A)\cap \overline{\Ran(B)}$, then
\[
\la R_D(\la)^2 B y_0 =  [\la R_D(\la)]( y -  \la R_D(\la)B y) \rightarrow 0 \quad \textrm{ as } |\la| \rightarrow 0, \, \la\in U.
\]
Consequently, 
\begin{equation}\label{conv 1/la delta}
\beta(\la)/\la \rightarrow w \quad \textrm{ as } |\la|\rightarrow 0, \lambda \in U.
\end{equation}
This completes the proof of (ii).
For (iii), if $By_0 = A^2 w$ for $w\in \Dom(A^2)$ and $BAw = A z$ for $z\in \Dom(A)$, then  
\begin{align*}
R_D(\la)^2 By_0 & = R_D(\la) \left( A w - \la R_D(\la) BAw\right)\\
& = w - \la R_D(\la) Bw - \la R_D(\la) (z - \la R_D(\la) Bz).
\end{align*}
Therefore, 
\[
 \beta(\la)/{\la^2} \rightarrow -w \quad \textrm{ as } |\la|\rightarrow 0, \, \la\in U.
\]
\end{proof}

\section{Well-posedness of (HE) and (PE)}
In this section, we present preparatory results on $(\mathrm{HE})$ and $(\mathrm{PE})$ needed for the proofs of Theorems~\ref{thm:ind-solu}--\ref{thm: diffusion}. Some of these results are of independent interest.
We begin by introducing the notion of \emph{well-posedness} of $(\mathrm{HE})$ used in our abstract framework, which can be roughly interpreted as the possibility of reducing $(\mathrm{HE})$ to the first-order system $(\mathrm{CP})$, see Subsection \ref{subsec well posed}. We prove that the phase space of a well-posed $(\mathrm{HE})$ is independent of the bounded operator $B$; see Theorem~\ref{thm: DWE well-posed}. This result also shows that $(\mathrm{HE})$ is well-posed if and only if $-A$ generates an operator cosine function; see Corollary~\ref{cor:DWE_well_posed}.
We then estimate the norm of $R$ in terms of the resolvent of the generator associated with $(\mathrm{CP})$; see Proposition~\ref{lem: pencil vs resolvent}. In particular, Proposition~\ref{prop on equiv cond} provides a significant extension of a result known in the Hilbert space setting, obtained by a different argument; cf. also Theorem~\ref{thm on equiv cond}.

\subsection{Well-posedness of (HE)}\label{subsec well posed} 
As a primary reference for the theory of abstract second-order problems, we use the monograph by Fattorini~\cite{Fa85}. By well-posedness of $(\mathrm{HE})$ we mean the notion introduced in \cite[p.~270-271]{Fa85}, with the additional requirement that \cite[Assumption VIII.3.1]{Fa85} is satisfied, which guarantees the existence of a \emph{phase space} for $(\mathrm{HE})$. For the convenience of the reader, we restate this definition in an equivalent form in the setting of this paper, that is, under the assumption that $B$ is bounded. 
The formulations in \cite{Fa85, Xiao98} impose an a priori estimate on every classical solution, whereas we require uniqueness explicitly. For bounded \(B\), these conditions are equivalent.
For related formulations of well-posedness of higher-order abstract Cauchy problems we refer the reader to \cite[Section 3]{Neub86} and \cite[Chapter~2]{Xiao98}. 

Let $A$ be a closed, densely defined operator on a Banach space $X$ and let $B \in \mathcal{L}(X)$. 
We consider the following conditions~(a) and~(b), and say that $(\mathrm{HE})$ is \emph{well-posed} whenever both of them are satisfied.

\smallskip
    (a) \emph{For any $y, x \in \Dom(A)$, there exists a unique {classical solution} $u = u(\cdot; y, x)$ to $(\mathrm{HE})$ satisfying
\[
u(0) = y, \quad u'(0) = x, \quad \text{and} \quad \|u(t)\| \leq c_t \bigl(\|y\| + \|x\|\bigr), \quad t \geq 0,
\]
for some constant $c_t > 0$ independent of $y$ and $x$.}
\smallskip

Here, we call $u \in C^2(\mathbb{R}_+; X)$ a \emph{classical solution} if  $u(t) \in \Dom(A)$ for all $t \geq 0$ and $(\mathrm{HE})$ is satisfied. 
The estimate in (a) allows us to define (uniquely) the \emph{propagators} (or \emph{solution families}) $(C(t))_{t\geq 0} \subset \mathcal{L}(X)$ and $(S(t))_{t\geq 0} \subset \mathcal{L}(X)$ associated with $(\mathrm{HE})$~by
\[
C(t)y = u(t; y, 0) \qquad \text{and} \qquad S(t)x = u(t; 0, x), \quad y, x \in \Dom(A).
\]

(b) \emph{$S$ is continuously differentiable in strong operator topology, i.e., $S(\cdot)x\in C^1(\R_+; X)$ for any $x\in X$.}
\smallskip

It  turns out that condition (b) guarantees exponential growth bounds for the solutions to $(\mathrm{HE})$ and the existence of the (maximal) \emph{phase space} $V \times X$; see~\cite[Theorem~4.2, Chapter VIII]{Fa85}.
The \emph{state space} $V$ is defined as the space of all $y \in X$ such that $C(\cdot)y \in C^1(\mathbb{R}_+; X)$, and is endowed with the norm 
\begin{equation}\label{norm_V}
\|y\|_V : = \|y\| + \sup_{0\leq t \leq 1} \|C'(t) y\|.
\end{equation}
In particular, \cite[Theorem~4.3, Chapter~VIII]{Fa85} asserts that $\mathrm{Dom}(A)$ is dense in $V$ and 
\begin{equation}\label{matrix_cU}
\mathcal{U}(t) := 
\begin{bmatrix}
C(t) & S(t)\\
C'(t) & S'(t)
\end{bmatrix}
\end{equation}
is a $C_0$-semigroup on $V \times X$ generated by the operator $\mathcal{A}$  given by 
\[
\mathcal{A} =
\begin{bmatrix}
0 & I\\
-A & -B
\end{bmatrix} \qquad \textrm{ with the domain } \quad \mathrm{Dom}(\mathcal{A}) = \mathrm{Dom}(A)\times V.
\]
Further, there exists a constant $\omega \geq 0$ such that $\{\Re \la > \omega\}\subset \rho(P)\cap \rho(\cA)$ and for any 
$\Re\la >\omega$
\begin{equation}\label{pencil vs resolvent}
(\la- \cA)^{-1} =  \left[\begin{array}{cc}
R(\la) (\la + B)  & R(\la) \\
\la R(\la)(\la + B) - I  & \la R(\la)
\end{array}\right].
\end{equation}
Equivalently, for all $(y,x)\in V\times X$ this means
\begin{equation}\label{eq: resolvent}
(\la - \cA)^{-1}(y,x) = \bigl( R(\la)(\la + B)y + R(\la)x, \la R(\la)(\la + B)y - y + \la R(\la) x \bigr).
\end{equation}
A straightforward verification shows that, for any $\lambda \in \mathbb{C}$, one has $\lambda \in \rho(P)$ if and only if $\lambda \in \rho(\mathcal{A})$. 
Indeed, the closedness of $A$ implies the closedness of both $P(\lambda)$ and $\lambda - \mathcal{A}$, and a direct algebraic manipulation shows that $P(\lambda)$ is bijective if and only if $\lambda - \mathcal{A}$ is bijective.

We denote by (CP) the Cauchy problem associated with the operator $\cA$ on $V\times X$, i.e., 
\begin{equation*}
\vec{u}'(t) = \cA \vec{u}(t),\,\, t\geq 0, \qquad \vec{u}(0) = (x_0,x_1) \in V\times X.\tag{CP}
\end{equation*}

In the case $B = 0 \in \mathcal{L}(X)$, we refer to $(\mathrm{HE})$ as $(\mathrm{HE}_0)$, that is,
\[
u''(t) + A u(t) = 0, \quad t \geq 0, \quad \text{with } u(0) = y,\, u'(0) = x. \tag{HE$_0$}
\]
If $(\mathrm{HE}_0)$ is well-posed, we denote by $V_0 \times X$, $\mathcal{U}_0$, $\mathcal{A}_0$, and $(C_0(t))_{t\geq 0}$ and $(S_0(t))_{t\geq 0}$ the corresponding phase space, semigroup, generator, and solution families, respectively.
We recall the well-known result (see~\cite[Theorem~3.14.11]{ABHN11}) that $(\mathrm{HE}_0)$ is well-posed if and only if the operator $-A$ generates a cosine function $\mathrm{Cos}(t) = C_0(|t|)$, $t \in \mathbb{R}$, on $X$. We refer the reader to~\cite{ABHN11,Fa85} for background on the cosine function theory. 

\smallskip

The main point of the following result is the identification of the maximal phase space of (HE) with that of the undamped equation (HE$_0$), with equivalent norms. This identification will be used in the proofs of Theorems~\ref{thm: diffusion intro}--\ref{thm:diffusion_for_Delta}.

\begin{theorem}\label{thm: DWE well-posed} Let $-A$ be the generator of a cosine operator function on a Banach space $X$, and let $V_0 \times X$ denote the phase space of $(\mathrm{HE}_0)$. Then for any $B \in \mathcal{L}(X)$ with $B \neq 0$ the problem $(\mathrm{HE})$ is well-posed, and its phase space $V \times X$ coincides with $V_0 \times X$, with equivalent norms.
\end{theorem}

\begin{proof}
By Kisyński's theorem (see~\cite[Theorem~3.14.11]{ABHN11}), the problem $(\mathrm{HE}_0)$ is well-posed. 
Since $\cB: = \left[\begin{smallmatrix}0&0\\0&-B\end{smallmatrix}\right]$ is a bounded operator on $V_0\times X$, the operator $\cA: = \cA_0 + \cB$ with the domain $\Dom(\cA) = \Dom(\cA_0) = \Dom(A)\times V_0$ generates a $C_0$-semigroup,  $\cU$, on $V_0\times X$.
Since $\cU(t)$, $t \geq 0$, is a bounded operator on $V_0 \times X$, it admits a block operator representation with components
\[
U_1(t) \in \cL(V_0), \quad
U_2(t) \in \cL(X, V_0), \quad
U_3(t) \in \cL(V_0,X), \quad
U_4(t) \in \cL(X).
\]
that is,  for all $t\geq 0$,  
\begin{equation}\label{cU_forumula}
\cU(t)(y,x) = ( U_1(t)y + U_2(t)x,\, U_3(t)y + U_4(t)x), \quad (y,x)\in V_0\times X.
\end{equation}
We first collect properties of the operator-valued functions $U_j$, $j \in \{1,\dots,4\}$, which will be used in the proof of the well-posedness of $(\mathrm{HE})$. Since $\mathcal{U}$ is a $C_0$-semigroup on $V_0 \times X$, each $U_j$ is strongly continuous on $\mathbb{R}_+$, and
\begin{equation}\label{cU(0)=I}
U_1(0) = I_{V_0}, \, U_4(0) = I_X, \, U_2(0) = 0 \in \mathcal{L}(X), \, U_3(0) = 0\in \cL(V_0,X).
\end{equation}
Moreover, for any $t\geq 0$,
\begin{equation}\label{pgs est 1}
\|U_1(t)\|_{\cL(V_0)}, 
\|U_2(t)\|_{\cL(X,V_0)},  
\|U_3(t)\|_{\cL(V_0,X)}, 
\|U_4(t)\| \leq  \|\cU(t)\|_{\cL(V_0 \times X)}. 
\end{equation}
Since for each $(y,x)\in \Dom(\cA) = \Dom(A)\times V_0$ the semigroup orbit $\cU(\cdot) (y,x)$ is the classical solution to the problem 
$\vec u\, ' (t) = \cA \vec u (t)$,  $t\geq 0$,  with  $\vec u(0) = (y,x)$,
we infer that $\cU(t) (y,x) \in \Dom(A)\times V_0$ for all $t \geq 0$, $\cU(\cdot)(y,x) \in C^1(\R_+ ; V_0\times X)$, and $\cA \cU(\cdot)(y,x) \in C(\R_+ ; V_0 \times X)$; see, e.g., \cite[Proposition 3.1.9(h)]{ABHN11}.

This implies that for each $y\in \Dom(A)$ and $x\in V_0$ the functions
\begin{equation}\label{diff_of_U12}
\begin{split}
& \R_+ \ni t\mapsto U_1(t)y + U_2(t)x \in V_0,\\ 
& \R_+ \ni t\mapsto U_3(t)y + U_4(t)x \in X \quad \textrm{are continuously differentiable with }
\end{split}
\end{equation}
\begin{equation}\label{derivatives of Uj}
\begin{split}
\left(U_1(t)y + U_2(t)x\right)'  & = U_3(t)y + U_4(t)x,\\
\left( U_3(t)y + U_4(t)x \right)' & =  -A \left(U_1(t)y + U_2(t)x\right) -B\left(U_3(t)y + U_4(t)x\right), 
\end{split}
\end{equation}
and 
\begin{equation}\label{dom_cond_A}
U_1(t)y, \, U_2(t)x \in \Dom(A), \qquad t\geq 0. 
\end{equation}
Combining \eqref{cU(0)=I}, \eqref{diff_of_U12}, \eqref{dom_cond_A}, and recalling that $V_0\hookrightarrow X$, we conclude that for any $(y,x) \in \mathrm{Dom}(A) \times V_0$, the function
\begin{equation}\label{classical_sol}
u(t) = u(t;y,x) := U_1(t)y + U_2(t)x, \quad t \geq 0,
\end{equation}
is a classical solution to $(\mathrm{HE})$ with $u(0)=y$ and $u'(0)=x$. 

In particular, 
for any $x\in V_0$ the function $\R_+\ni t\to U_2(t)x \in X$ is continuously differentiable and $(U_2(t)x)' = U_4(t)x$, $t\geq 0$.
Since, by \eqref{pgs est 1}, the functions $\R_+ \ni t\to U_2(t),\,U_4(t)\in \cL(X)$ are uniformly bounded on compact subsets of $\R_+$ and $V_0$ is dense in $X$, a standard argument shows that $U_2(\cdot)x$ is in $C^1(\R_+; X)$ for any $x\in X$ and $(U_2(t)x)' = U_4(t)x$, $t\geq 0$. This implies condition~(b) for the well-posedness of~{(HE)}. 

To verify condition~(a), we first show that for every $(y,x) \in \Dom(A)\times V_0$ the function $u$ given by \eqref{classical_sol} is the unique classical solution to (HE) with initial data $y$ and $x$. Indeed, let $v$ be any classical solution to (HE) with the same initial data. Since $Av=-v''-Bv'$ is continuous as an $X$-valued function on $\R_+$, the function $v$ is continuous on $\R_+$ with values in $\Dom(A)$ equipped with its graph norm, and hence also with values in $V_0$. Thus, $\vec v:=(v,v')$ belongs to $C(\R_+;V_0\times X)$ and satisfies
$$
\int_0^t\vec v(s)\,ds\in\Dom(\cA),\qquad
\cA\int_0^t\vec v(s)\,ds=\vec v(t)-\vec v(0),\quad t\ge0.
$$
Consequently, $\vec v$ is a mild solution of (CP), and uniqueness of mild solutions yields $\vec v(t)=\cU(t)\vec v(0)$ for all $t\ge0$, as desired. Therefore, since $U_2(t)\in\cL(X)$ for every $t>0$, it remains to show that
$$
\|u(t,y,0)\|=\|U_1(t)y\|\leq c\|y\|,\qquad y\in\Dom(A),
$$
with a constant $c>0$ independent of $y$. Equivalently, since $\Dom(A)$ is dense in $X$, we have to show that $U_1(t)$, $t>0$, extends to a bounded operator on $X$. Since $B \in \mathcal{L}(X)$ and $\mathrm{Dom}(A)$ is dense in $X$, it suffices to show that
\begin{equation}\label{crucial_formula}
U_1(t)y = U_4(t)y + U_2(t)By, \quad t>0,\; y \in \mathrm{Dom}(A).
\end{equation}
Since $U_4(t)y = (U_2(t)y)'$, $t\geq 0$, and $U_2(0)y=0$, this is equivalent to show that 
\begin{equation}\label{identity:u1}
U_2(t)y = \int_0^t (U_1(s)y - U_2(s)By) \ud s, \quad t>0,\; y \in \mathrm{Dom}(A).
\end{equation}
We fix $y\in \Dom(A)$ and let 
\[
u(t): =  \int_0^t \left(U_1(s)y - U_2(s) B y\right) \, \ud s, \qquad t\geq 0.
\]
To show \eqref{identity:u1}, it suffices to prove that the function $(u,u')$ is a mild solution to $\vec{u}'(t) = \cA \vec{u}(t)$, $t\geq 0$, with $\vec{u}(0)= (0,y)$, that is, 
\begin{equation}\label{orbit_of_cU}
\vec{u} \in C(\R_+; V_0 \times X), \, \int_0^t \vec{u}(s)\ud s \in \Dom(\cA) \, \textrm{ and }  \, \cA\int_0^t \vec{u}(s)\ud s = \vec{u}(t) - (0,y), \, t\geq 0.
\end{equation}
Then, by \cite[Theorem 3.1.11]{ABHN11}, $\vec{u}(t) = \cU(t) (0,y)$, $t\geq 0$. On the other hand, by \eqref{cU_forumula}, $\cU(t)(0,y) = (U_2(t)y, U_4(t)y)$, so we obtain \eqref{identity:u1}.

Therefore, we have to establish~\eqref{orbit_of_cU} for $\vec{u} := (u,u')$. As noted above, $U_1$ and $U_2$ are strongly continuous on $\mathbb{R}_+$, and hence $\vec{u} \in C(\mathbb{R}_+; V_0 \times X)$. 
Since $(U_2(t)By)' = U_4(t)By$, $t\geq 0$  and \eqref{cU(0)=I}, we readily see that $u\in C^2(\R_+; X)$ and 
\begin{equation}\label{deriv_of_u}
\begin{split}
u'(t) & = \int_0^t \left(U_3(s)y - U_4(s)By \right) \ud s + y, \\
u''(t) & =  U_3(t)y - U_4(t)By, \quad t\geq 0. 
\end{split}
\end{equation}
Combining \eqref{deriv_of_u} with the following identity (see, e.g., \cite[Proposition 3.1.9(e)]{ABHN11}),  
\begin{equation}\label{semig_formula}
\cU(t)(y, -By) - (y,-By) = \cA \int_0^t \cU(s)(y,- By) \ud s, \quad t\geq 0,
\end{equation}
we infer that $u$ is a classical solution to (HE) with $u(0) = 0$ and $u'(0) = y$. 
This proves \eqref{orbit_of_cU}, hence
\eqref{identity:u1} and \eqref{crucial_formula}. This establishes the remaining assertion of condition~(a), thereby completing the proof of well-posedness of (HE).

For $t\geq 0$, we set $S(t) := U_2(t)$ and let $C(t)$ denote the extension of $U_1(t)$ to the whole space $X$. By~\eqref{crucial_formula}, we then have
\begin{equation}\label{C_vs_S'}
C(t) :=  S'(t) + S(t)B, \quad t \geq 0,
\end{equation}
where $S'$ denotes the derivative taken with respect to the strong operator topology on $\mathcal{L}(X)$. It remains to verify that $V_0 \times X$ is the maximal phase space for $(\mathrm{HE})$, that is,
\begin{equation}\label{eq:V0}
V_0 = \{ x \in X : C(\cdot)x \in C^1(\mathbb{R}_+; X)\} =: V,
\end{equation}
and that the norm on $V$ (see~\eqref{norm_V}) is equivalent to the norm on $V_0$. 
Since $S(\cdot)$ is strongly differentiable on $\R_+$, it follows from \eqref{C_vs_S'} that, for every $x\in X$, $S'(\cdot)x\in C^1(\R_+;X)$ if and only if $C(\cdot)x\in C^1(\R_+;X)$.
By \cite[Proposition 3.1.9(h)]{ABHN11}, for every $x\in X$, we have $(0,x)\in\Dom(\cA)=\Dom(\cA_0)=\Dom(A)\times V_0$ if and only if the orbit $t\mapsto\cU(t)(0,x)=(S(t)x,S'(t)x)$ belongs to $C^1(\R_+;V_0\times X)$. 
Thus, if $x\in V_0$ then 
$S'(\cdot)x\in C^1(\R_+;X)$, and hence $V_0\subset V$. 
Conversely, suppose that $S'(\cdot)x\in C^1(\R_+;X)$ for some $x\in X$. The perturbation formula
\[
 \cU (t) \vec z=\cU_0 (t) \vec z+\int_0^t \cU_0(t-s)\cB \cU(s) \vec z\,ds, \quad \vec z \in V_0\times X,
\]
applied to $\vec z=(0,x)$ and projected onto the second component yields
\[
C_0(t)x =  S'(t)x + \int_0^tC_0(t-s)BS'(s)x\,\ud s.
\]
Since $S_0'(t)=C_0(t)$, $S_0(0)=0$, and $S'(0)=I$,
integration by parts gives
\[
\int_0^t C_0(t-s)BS'(s)x\,ds
=S_0(t)Bx+\int_0^t S_0(t-s)BS''(s)x\,ds.
\]
The right-hand side belongs to $C^1(\R_+;X)$, as
$S_0$ is strongly $C^1$ and $BS''(\cdot)x$ is continuous.
Thus $C_0(\cdot)x\in C^1(\R_+;X)$, so $x\in V_0$.

Finally, since both spaces $V_0$ and $V$ are continuously embedded in $X$, the equivalence of their norms follows from the closed graph theorem applied to the identity map.
\end{proof}

 The argument of Theorem~\ref{thm: DWE well-posed} also yields the following characterization of well-posedness of (HE); see also \cite[ Theorem~2.5.1]{Xiao98}.

\begin{corollary}\label{cor:DWE_well_posed}
Let $A$ be a linear densely defined, closed operator on a Banach space $X$. Then the following assertions are equivalent:
\begin{itemize}
\item[(i)]The problem \emph{(HE)} associated with $A$ and $B$ is well-posed for some (equivalently, for all) $B \in \mathcal{L}(X)$.
\item[(ii)] $-A$ is the generator of a cosine function on $X$.
\end{itemize}
\end{corollary}

\subsection{Linear growth of the semigroup $\cU$}\label{subs:linear_grothw} 
In this subsection, we present a corollary of Theorem~\ref{thm: DWE well-posed} that allows us to verify the assumptions of Theorems~\ref{thm: diffusion intro} and~\ref{thm: diffusion} concerning the well-posedness of $(\mathrm{HE})$ and the polynomial growth of the semigroup $\mathcal{U}$ associated with $(\mathrm{HE})$.
We restrict our attention to the case where $X$ is a Hilbert space,
$A$ is a positive selfadjoint operator, and $B\in\mathcal L(X)$ is non-negative.

\begin{proposition}\label{prop:linear_growth}
Let $A$ be a positive selfadjoint operator on a Hilbert space $X$ and  $B$ be a bounded non-negative operator on $X$.
Then the corresponding problem \emph{$(\mathrm{HE})$} is well posed and its phase space is $\Dom(A^{1/2}) \times X$.
Moreover, the $C_0$-semigroup $\mathcal{U}$ associated with \emph{$(\mathrm{HE})$} satisfies the estimate
\begin{equation}\label{linear_growth_U}
\|\cU(t)\|_{\cL(V\times X)} 
\leq  \sqrt{3}\max(1,t),\qquad t\geq 0.
\end{equation}
\end{proposition}

\begin{proof} Since $A$ is positive, $-A$ generates an operator cosine function on $X$ (see, e.g., \cite[Example~3.14.16]{ABHN11}), and hence the well-posedness of $(\mathrm{HE})$ follows from Theorem~\ref{thm: DWE well-posed}. By \cite[Proposition~3.16.3 and Theorem~3.16.7]{ABHN11}, the phase space for $(\mathrm{HE}_0)$ is $\Dom(A^{1/2}) \times X$. We equip $\Dom(A^{1/2})$ with the graph norm of $A^{1/2}$, making it a Hilbert space with the associated inner product.
Consequently, Theorem~\ref{thm: DWE well-posed} implies that the phase space $V \times X$ for $(\mathrm{HE})$ coincides with $\Dom(A^{1/2}) \times X$, with equivalent norms.
Let $\cA$ denote the generator of a $C_0$-semigroup $\cU$ on $V\times X$, that is,
\[
\cA(x,y): = (y, -A x - By) \, \textrm{ with } \Dom(\cA): = \Dom(A)\times \Dom(A^{1/2}). 
\]
Moreover, let $C$ and $S$ stand for the \emph{solution families} for the (HE). 

Since $A$ is injective, so is $A^{1/2}$, and $\|A^{1/2}\cdot\|$ defines a norm on $\Dom(A^{1/2})$. We denote by $\dot V$ the space $\Dom(A^{1/2})$ equipped with the norm $\|A^{1/2}\cdot\|$ and the corresponding inner product $\langle\cdot, \cdot \rangle_{\dot V}$. Furthermore, we let $\widetilde V$ to denote the completion of~$\dot V$. That is, $\wt V$ is the Hilbert space of all (equivalence classes of) Cauchy sequences $(x_j)_{j}$ in $\dot V$ with the inner product defined by $\langle (x_j)_j, (y_j)_j \rangle_{\wt V} : = \lim_{j\to \infty} \langle x_j, y_j \rangle_{\dot V}$.
We denote by $J \colon \dot V \ni x \mapsto (x)_{j} \in \widetilde V$
the corresponding isometric isomorphism from $\dot V$ onto a dense subspace of $\widetilde V$.

We first define a densely defined operator $\wt\cA_0$ on $\wt V\times X$ by
\begin{equation}\label{def tA}
\wt\cA_0(Jx,y):=(Jy,-Ax-By),\qquad
\Dom(\wt\cA_0):=J\Dom(A)\times\Dom(A^{1/2}).
\end{equation}

The space $\wt V \times X$ is equipped with the 
(hilbertian) product norm $\|\left((x_j), x\right)\|_{\wt V \times X} : = \big( \|(x_j)\|^2_{\wt V} + \|x\|^2\big)^{1/2}$, $(x_j)\in \wt V$, $x\in X$, and with the corresponding inner product.

We claim that $\wt \cA_0$ is densely defined and dissipative 
and that its closure is maximal dissipative. To begin the proof of the claim we first note that $D(\wt \cA)_0$ is a dense subset of $\wt V \times X$. 
Indeed, $\Dom(A)$ is a core of $A^{1/2}$, i.e., for any $x\in \Dom(A^{1/2})$ there exists $(x_j) \subset \Dom(A)$ such that $x_j \rightarrow x$ and $A^{1/2} x_j \rightarrow A^{1/2} x$ in $X$ as $j\rightarrow \infty$.
In particular, $\Dom(A)$ is dense in $\dot V$, and so $J \Dom(A)$ is dense in the completion $\wt V$ of $\dot V$ while $\Dom(A^{1/2})$ is dense in $X$ as claimed.

To show that  $\wt \cA_0$ is dissipative note that, for any $(Jx , y) \in \Dom(\wt\cA_0)$, we have
\begin{align*}
 \lan \widetilde{\mathcal A}_0(Jx,y),(Jx,y)\ran_{\wt V \times X}
 &=\langle A^{1/2}y,A^{1/2}x\rangle
   -\langle Ax,y\rangle-\langle By,y\rangle\\
 &=\langle A^{1/2}y,A^{1/2}x\rangle
   -\langle A^{1/2}x,A^{1/2}y\rangle-\langle By,y\rangle.
\end{align*}
Hence
$\operatorname{Re}\langle\widetilde{\mathcal A}_0(Jx,y),(Jx,y)\rangle
=-\langle By,y\rangle\le0$.
Since $\wt {\mathcal A}_0$ is densely defined and dissipative, it is closable. Put $\wt {\mathcal A} := \overline {\wt {\mathcal A}_0}.$ Then
$\wt {\mathcal A}$ is closed and dissipative. Consequently,
\[
\|(I - \wt {\mathcal A})z\| \ge \|z\|, \qquad 
z \in \Dom( \wt {\mathcal A}),
\]
so $I - \wt {\mathcal A}$ is injective and has closed range.

Next, we show that $\Ran(I-\wt\cA)$ is dense in
$\wt V\times X$. Since $\wt\cA$ extends $\wt\cA_0$,
it suffices to prove that
\begin{equation}\label{inc}
\cD(\wt\cA_0)\subset \Ran(I-\wt\cA_0).
\end{equation}
Indeed, by the definition of $\wt\cA_0$ in (4.18), for any $x\in\Dom(A)$
and $y\in\Dom(A^{1/2})$, we have
\[
(I-\wt\cA_0)(Jx,y) =(Jx-Jy,\;y+Ax+By).
\]
Fix $(Jz,w)\in\cD(\wt\cA_0)$.
Since $A+B$ is selfadjoint and non-negative,
the operator $I+A+B$ has a bounded inverse on $X$.
Consequently, there exists $y\in\Dom(A)$ such that
\[
(I+A+B)y=Az+w.
\]
Putting $x=z+y$, we have $x\in\Dom(A)$ and
\[
(I-\wt\cA_0)(Jx,y)=(Jz,w).
\]
This proves \eqref{inc}. Since
$\cD(\wt\cA_0)$ is dense in $\wt V\times X$, so is
$\Ran(I-\wt\cA)$. As the latter range is also closed,
we conclude that
\[
\Ran(I-\wt\cA)=\wt V\times X.
\]
Thus $\wt\cA$ is maximal dissipative. By the Lumer-Phillips theorem, see, e.g., \cite[Theorem 3.4.5]{ABHN11}, $\wt\cA$  generates  a contraction $C_0$-semigroup on $\wt V \times X$ that we will denote by $\wt \cU$.

We now provide the relation between $\cA$, $\cU$ on $\Dom(A^{1/2})\times X$ and $\wt\cA, \wt \cU$ on $\wt V \times X$. Let 
\[
\jmath :  \Dom(A^{1/2}) \times X \ni (x,y) \mapsto (Jx, y) \in \t V\times X. 
\]
Clearly, $\jmath$ is a contraction from $\Dom (A^{1/2}) \times X$ to $\wt V \times X$, and an isometry from $\dot V \times X$ to $\wt V \times X$. Moreover, for all $t\geq 0$, we have
\[
\jmath \, \cU(t) = {\t \cU}(t) \jmath \, \textrm{ on } \Dom(A^{1/2}) \times X \quad \textrm{and} \quad \jmath \cA = \wt \cA \jmath \,  \textrm{ on } \Dom(\cA). 
\]

In particular, for any $x\in \Dom(A^{1/2})$ and $y\in X$ we have
\begin{align*}
\| \cU(t)(x,y)\|_{\dot V \times X}^2 & = \| \t\cU (Jx, y) \|^2_{\wt V \times X}  \leq \|x\|_{\dot V}^2 + \|y\|_X^2 = \|(x,y)\|^2_{\dot V \times X}, 
\end{align*}
that is, $\cU(t)$ is a contraction on $\dot V\times X$.
Therefore, by \eqref{matrix_cU}, for all $x\in \Dom(A^{1/2})$ and $y\in X$
\begin{equation}\label{pgs est 22}
\begin{split}
\|C(t)x\|_{\dot V} & \leq \| x\|_{\dot V}, \qquad \|S(t) y\|_{\dot V} \leq \|y\|_X,\\
\|C'(t) x\|_X & \leq \|x\|_{\dot V},  \qquad \|S'(t) y\|_X\leq \|y\|_X.
\end{split}
\end{equation}
Using this we obtain
\begin{align}\label{est:cU_in_VxX}
\nonumber\| \cU(t)(x,y)\|_{V\times X}^2 & =  \| \cU(t)(x,y)\|_{\dot V \times X}^2 +  \|C(t)x + S(t) y\|_{X}^2\\
& \leq \|x\|_{\dot V}^2 + \|y\|_X^2 + \|C(t)x + S(t) y\|_{X}^2.
\end{align}
It remains to estimate $\|C(t)x + S(t) y\|_{X}^2$.
For $x\in \Dom(A^{1/2})$ and $y \in X$ the functions $C(\cdot)x$ and $S(\cdot)y $ are continuously differentiable on $\R_+$ by Theorem \ref{thm: DWE well-posed}. Therefore, since
\[
C(t)x = \int_0^t C'(s)x \ud s + x \quad \textrm{ and } \quad S(t)y = \int_0^t S'(s)y \ud s,  \quad t \geq 0,
\]
estimates \eqref{pgs est 22} yield
\begin{equation*}
\begin{split}
\|C(t)x\|_X  & \leq t \sup_{s\in [0,t]} \|C'(s)x\|_X + \|x\|_X \leq t \|x\|_{\dot V} + \|x\|_X \leq \max\{1,t\} \|x\|_V,\\
\|S(t)y\|_X & \leq t \sup_{s\in [0,t]} \|S'(s)y\|_X \leq t \|y\|_{X}, \quad t\geq 0.
\end{split}
\end{equation*}
Consequently,   
\[
\|C(t)x + S(t) y\|_{X}^2\leq  2 \max\{1,t^2 \} \|x\|^2_V + 2t^2 \|y\|_X^2.
\]
Combining the last estimate with \eqref{est:cU_in_VxX} yields the required estimate \eqref{linear_growth_U}.
\end{proof}

In the context of the assumptions of Theorem \ref{thm: diffusion}, the following remark is in order.

\begin{remark}\label{Goldstein example} Under the assumption that (HE) is well posed, the associated semigroup $\cU$, in general, is not uniformly bounded on the phase space $V\times X$. For instance, in the case when $B=0$ and $A$ is the negative Laplacian on $X = L^2(\R^n)$ it is known that $\|\cU (t)\|_{\cL(V\times X)} \simeq t$ for $t$ large; see \cite{GoWa03}. 
On the other hand, for any generator $-A$ of a uniformly bounded cosine function on an arbitrary Banach space $X$, due to a result by Fattorini  \cite[Theorem 1.3]{Fa85}, the linear growth holds, that is,   $\|\cU_0(t)\|_{\cL(V\times X)}\leq 4(1+t)$, $t\geq 0$. $\hfill\Diamond$\end{remark}

\subsection{Pencil estimates and the inclusion $V\subset X_2$}\label{subsec43}
We recall that Theorem~\ref{thm 1}(i) asserts that
$\Dom(A)\subset X_2$. This follows from \eqref{pencils est} in
Hypothesis~\ref{hyp:main} by a simple algebraic argument; see the proof
of Lemma~\ref{lem:reduction_step}(iii).  
Among other things, the next three results show that, under additional
assumptions on (HE) and when the underlying space~$X$ is UMD, the inclusion
$\Dom(A) \subset X_2$ can be strengthened to $V\subset X_2$. Combining this with Theorem~\ref{thm: DWE well-posed} and Fattorini's theorem on the square root reduction of cosine function generators (see, e.g., \cite[Theorem 3.16.4]{ABHN11}), we obtain $ \Dom\big((\omega+A)^{1/2}\big) = V \subset X_2$ for some $\omega\geq 0$. Consequently, our results, Theorems~\ref{thm: diffusion intro}--\ref{thm:diffusion_for_Delta},
can be formulated in terms of the spaces associated with $A$,
rather than the auxiliary space $X_2$.
Moreover, working with $V$ instead of $X_2$ allows us to
formulate more natural conditions on $B$ in concrete applications,
ensuring that $Bx_0, By_0\in X_2$ for every
$x_0,y_0\in \Dom(A)$. This is important in the adaptation of
Theorem~\ref{thm 1} to concrete situations. For instance, in the
model example studied in Section~\ref{sec:examples}, it is sufficient
to assume that $b\in W^{1,\infty}(\Omega)$, to obtain that $B(\Dom(A)) \subset V \subset X_2$ for the corresponding $A$ and $B$.

We start with a result which, in particular, identifies a sufficient condition (see \eqref{eq: pencil vs resolvent} below) ensuring that $V\subset X_2$ in the case when $A$ and $B$ satisfy \eqref{pencils est} in Hypothesis \ref{hyp:main}. 
Then, in Theorem~\ref{thm on equiv cond}, we show that the assumption \eqref{eq: pencil vs resolvent} is automatically satisfied when $X$ is a UMD space. 
The UMD property of a Banach space \(X\) is characterized by the boundedness of the Hilbert transform on \(L^p(\mathbb R;X)\) for some (equivalently, every) \(p\in(1,\infty)\). The class of UMD spaces includes Hilbert spaces and \(L^p\)-spaces, \(1<p<\infty\), as well as their closed subspaces and quotient spaces. For background on UMD spaces and their applications to evolution equations, see e.g. \cite[Sections~1.8 and~3]{KuWe04} and 
 \cite[Chapter I]{ABHN11}.

\begin{proposition}\label{lem: pencil vs resolvent} Let $A$ and $B$ be such that \emph{(HE)} is well-posed with the phase space $V\times X$ as described in Subsection \ref{subsec well posed}. Let us assume that $is \in \rho(P) \cap \rho(\cA)$ for all $|s|> s_0$,  where $s_0\geq 0$ is a constant that for some integer $k\geq -1$ one has
\begin{equation*} 
\sup_{|s|> s_0} \| s^{-k} R(is)\|<\infty.
\end{equation*}
If there exists a constant $c>0$ such that
\begin{equation}\label{eq: pencil vs resolvent}
\|(is - \cA)^{-1}\|_{\cL(V\times X)} \leq c \|is R(is) \|, \qquad |s|> s_0,
\end{equation}
then for any $x\in X$ with $(0,x)\in \Dom(\cA^{\ell})$, where $\ell\in \N$, and any $|s|> s_0$, we have
\begin{equation}\label{est: pencil via resolvent}
\| is R(is) x\|_{X} \leq {|s|^{-1}}\left( \sum_{j=0}^{\ell-1} {|s|^{-j}} {\|\cA^j(0,x)\|_{V\times X} } + c [R]_{-k, \infty} \|\cA^\ell(0,x)\|_{V\times X} |s|^{k+ 2 - \ell} \right).
\end{equation}
Moreover, for $j=1,2$, the following assertions hold:
\begin{itemize}
\item [(i)] if $(0,x) \in \Dom(\cA^{k+j})$ then 
\[
\sup_{|s|> s_0} \| s^j R(is) x\|_{X} <\infty.
\]
\item [(ii)] For any $x\in X$, 
\begin{equation*}\begin{split}
(0,x) \in \Dom(\cA) \quad  & \Leftrightarrow \quad x\in V,\\
(0,x) \in \Dom(\cA^2)  \quad & \Leftrightarrow \quad x\in \Dom(A), \quad Bx \in V,\\
(0,x) \in \Dom(\cA^3) \quad & \Leftrightarrow \quad x,\, Bx \in \Dom(A),\, -Ax + B^2x \in V.
\end{split}
\end{equation*}
\end{itemize}
\end{proposition}

\begin{proof} Setting $\la = is$ with $|s|> s_0$
by \eqref{eq: resolvent} we obtain for any $x\in X$, 
\[
(\la - \cA)^{-1}(0,x) = \bigl(R(\la) x,  \la R(\la) x \bigr)\quad \textrm{and}\quad 
\|\la R(\la)x\| \leq \|(\la -\cA)^{-1}(0,x)\|_{V\times X}.
\]
Moreover, for any $\ell\in \N$ and $\bf x \in \Dom(\cA^\ell)$, 
\[
(\la - \cA)^{-1} {\bf x} = \sum_{j=1}^\ell {\la^{-j}} \cA^{j-1}{\bf x} + {\la^{-\ell}} (\la - \cA)^{-1} \cA^\ell {\bf x}.
\]
Combining this with \eqref{eq: pencil vs resolvent} yields
\eqref{est: pencil via resolvent}. Statement~(i) now follows readily from \eqref{est: pencil via resolvent}. For~(ii), we note that  
\begin{equation*}\begin{split}
\cA(0,x) & = (x,\, -Bx), \qquad \cA^2(0,x) = (-B x,\, -A x - B(-Bx)),\\
\cA^3(0,x) & = (-Ax + B^2x, \, ABx + BAx - B^3 x). 
\end{split}
\end{equation*}
Since $\Dom(\cA) = \Dom(A)\times V$ and $\Dom(A)\subset V$, we get the desired equivalences in~(ii).
\end{proof}

To establish condition~~\eqref{eq: pencil vs resolvent} for UMD spaces, we shall use a resolvent comparison which is also of independent interest. In the self-adjoint Hilbert-space setting, comparisons between the resolvent of the quadratic pencil and that of the associated first-order operator have been used in the study of decay of damped waves; see, for example, \cite[Lemma~4.6]{AnLa14} and the references therein. Proposition~\ref{prop on equiv cond} extends such comparisons to sectorial operators $A$ on reflexive Banach spaces and arbitrary $B\in\cL(X)$, and gives estimates at every nonzero $\lambda\in\rho(P)$.

We recall that an operator $A$ on a Banach space $X$ is called \emph{sectorial} if  $(-\infty,0)\subset \rho(A)$ and $\sup_{\lambda>0} \|\lambda(\lambda + A)^{-1}\| < \infty$. Following the notation of~\cite{KuWe04}, for a sectorial operator $A$ and $\vartheta \in \mathbb{R}$, we set 
\[
(X_{\vartheta, A}, \|\cdot\|_{X_{\vartheta, A}}) : = \left\{ \begin{array}{ll}
\Dom(A^\vartheta), \|(I + A)^\vartheta \cdot\|, & \vartheta\geq 0,\\[0.5ex]
(X, \|(I + A)^{\vartheta}\cdot\|)^{\widetilde{  }}, & \vartheta <0.
\end{array}\right.
\]
Let $X^*$ denote the dual space to $X$ and $A^*$ the adjoint operator to $A$ with respect to the duality $\langle\cdot, \cdot \rangle_{X^*,X}$. 
We recall that for any $\vartheta>0$, $\|\cdot\|_{X_{\vartheta,A}}$ is equivalent to the graph norm of $A^\vartheta$; see, e.g., \cite[Lemma 15.22]{KuWe04}. 

We first prove an auxiliary result.
Its proof is based on extrapolation fractional scales of sectorial operators. We refer the reader to~\cite{Am95} and to~\cite{KuWe04} for background. We shall use the next partial case of the moment inequality for sectorial operators
(see, e.g.,~[44, Theorem~15.14]).
Note that since $X^*$ is reflexive, $A^*$ is a densely defined sectorial operator on $X^*$.
Since $A$ and $A^*$ are sectorial,
we fix $c_A\geq 1$ such that
\begin{equation}\label{eq:moment-fix}
 \|T^{1/2}z\|^2\leq c_A\|Tz\|\,\|z\|,\qquad z\in\Dom(T),
\end{equation}
for each of $T=A$, $A^*$, $I+A$, and $I+A^*$, with the norms taken
in the respective spaces $X$ or $X^*.$ If $X$ is a Hilbert space, and $A$ is non-negative, then $c_A=1.$

\begin{lemma}\label{prop on equiv cond_lemma}
Let $A$ be a sectorial operator on a reflexive Banach space $X$. Let $B\in \cL(X)$. Then the following assertions hold.

\begin{itemize}
\item [(i)] For any $0\neq \la \in \rho(P)$ and for every $x^* \in X^*$ and $x\in X$ we have
\begin{equation}\label{bd X into V}
\|R(\la)^*x^*\|_{(X^*)_{\frac{1}{2},A^*}} \leq  C(\la) \|x^*\|_{X^*}\quad \textrm{and} \quad \|R(\la)x\|_{X_{\frac{1}{2},A}} \leq  C(\la) \|x\|_{X},
\end{equation} 
where 
\begin{equation}\label{eq:C-replacement}
 C(\lambda)=\sqrt{c_A}\,
 \left(1+|\lambda|^{-2}+|\lambda|^{-1}\|B\|
       +|\lambda|^{-2}\|R(\lambda)\|^{-1}\right)^{1/2}
 \|\lambda R(\lambda)\|.
\end{equation}

\item [(ii)] For any $0\neq \la \in \rho(P)$ the operator $R(\la) \in \cL(X)$  (respectively, $A$) has an extension to the operator in $\cL(X_{-\frac{1}{2}, A},X)$ (respectively, $\cL(X_{\frac{1}{2},A},X_{-\frac{1}{2},A})$) with the norm less or equal to $C(\la)$ (respectively, $\|A(I+A)^{-1}\|$). 
\end{itemize}
\end{lemma}

\begin{proof}
For any $y^*\in \Dom(A^*)$, using the defining fractional-domain norm and \eqref{eq:moment-fix} with $T=I+A^*$, we have
\begin{align*}
 (c_A)^{-1}\|y^*\|_{(X^*)_{1/2,A^*}}^2
 &=(c_A)^{-1}\|(I+A^*)^{1/2}y^*\|_{X^*}^2\\
 &\leq \|(I+A^*)y^*\|_{X^*}\,\|y^*\|_{X^*}\\
& \leq  \|y^*\|_{X^*}^2 + \|A^*y^*\|_{X^*} \|y^*\|_{X^*}\\
& \leq  \|y^* \|_{X^*}^2 + (\|P(\la)^* y^*\|_{X^*} + \|\left(\la B^* + \la^2)\|_{\cL(X^*)} \|y^*\|_{X^*}\right) \|y^*\|_{X^*}\\
& \leq  \left(1 +  \|\la B^* + \la^2\|_{\cL(X^*)} \right) \|y^*\|^2 + \|P(\la)^* y^*\|_{X^*} \|y^*\|_{X^*}\\
& \leq     \left(1 +  \|\la B^* + \la^2\|_{\cL(X^*)} \right) \|R(\la)^*\|^2_{\cL(X^*)} \| P(\la)^* y^* \|_{X^*}^2\\
& \quad  + 
\|R(\la)^*\|_{\cL(X^*)} \| P(\la)^* y^*\|_{X^*}^2.
\end{align*}
Since $\|B^*\|_{\cL(X^*)} = \|B\|$ and $\|R(\la)^*\|_{\cL(X^*)} = \|R(\la)\|$ , setting $x^* = P(\la)^* y^*$ and  
\begin{align*}
C(\la) & : =  \sqrt{c_A}\left(\left(1 +  |\la|  \|B\| + |\la|^2 \right) \|R(\la)\|^2  + \|R(\la)\|\right)^{1/2}
\\
& =  \sqrt{c_A} \big( 1 + |\la|^{-2} +  |\la|^{-1} \|B\|   + |\la|^{-2} \|R(\la)\|^{-1}\big)^{1/2}  \| \la R(\la)\|
\end{align*}
completes the proof of the first estimate in \eqref{bd X into V}.  The second estimate follows by applying the same moment inequality to $A$ in place of $A^*$.

(ii)  Let $T$ denote the operator $R(\la)^*$ considered as an operator in $\cL(X^*, (X^*)_{\frac{1}{2}, A^*})$. Recall that $((X^*)_{\frac{1}{2}, A^*})^* = X_{-\frac{1}{2},A}$ with respect to the duality pairing $\langle \cdot, \cdot \rangle_{X,X^*}$; see, e.g., \cite[Proposition 15.23]{KuWe04}. Taking the adjoint $T^*$ of $T$ (again with respect to the  pairing $\langle \cdot, \cdot \rangle_{X,X^*}$), we infer that $T^*$ is an extension of $R(\la)$ to an operator in $\cL(X_{-\frac{1}{2},A}, X)$. Using the first bound in \eqref{bd X into V}, $
\|T^*\|_{\cL(X_{-\frac{1}{2},A},X)
} = \|T\|_{\cL(X^*, (X^*)_{\frac{1}{2},A^*})}\leq C(\la)$, as desired. 

Since $\Dom(A) \subset X_{\frac{1}{2},A}$ is dense and, for every $y\in\Dom(A)$, we have $\|Ay\|_{X_{-\frac{1}{2},A}}
\leq \|A(I+A)^{-1}\|_{\cL(X)}\|y\|_{X_{\frac{1}{2},A}}$, 
the operator $A$ has a unique extension to an operator in
$\cL(X_{\frac{1}{2},A},X_{-\frac{1}{2},A})$ with the norm less or equal to  $\|A(I+A)^{-1}\|$.
\end{proof}

We are now in a position to prove the announced resolvent estimate.

\begin{proposition}\label{prop on equiv cond}
Let $A$ be a sectorial operator on a reflexive Banach space $X$. Let $B\in \cL(X)$.  Let $\cX: =  X_{\frac{1}{2},A}\times X$ be equipped with the norm $\|(y,x)\|_\cX: = \|(I+A)^{1/2}y\| + \|x\|$, $y\in X_{\frac{1}{2},A}$, $x\in X$. Let $\cA$ be an operator on $\cX$ defined by  $\cA: = \left[\begin{smallmatrix}0&I\\-A&-B\end{smallmatrix}\right]$,  $\Dom(\cA): = \Dom(A)\times X_{\frac{1}{2},A}$. 

Then $\rho(P) \subset \rho(\cA)$, and there exists a constant $c>0$ such that, for every $0\neq\la\in\rho(P)$, we have
\begin{equation}\label{norm of A}
\|(\la - \cA)^{-1}\|_{\cL(\cX)} \leq c \max\big\{\|\la R(\la)\|, \, |\la|^{-1}\| R(\la)\|,\, |\la|^{-1} \big\}.
\end{equation} 
In particular, there exists a constant $C\ge 1$ such that, for every $\lambda\in\rho(P)$ with $|\lambda|\ge1$,
\begin{equation}\label{2sided}
 \|\lambda R(\lambda)\|_{\cL(X)}
 \le \|(\lambda-\cA)^{-1}\|_{\cL(\cX)}
 \le C\|\lambda R(\lambda)\|_{\cL(X)}.
\end{equation}
\end{proposition}

\begin{proof}
 Let $0\neq \la \in \rho(P)$ and let $\cG(\la)$ denote the matrix operator on the right-hand side of \eqref{pencil vs resolvent}. It is readily seen that $
 \cG(\la) (\la - \cA) \subset I_{\cX} =  (\la - \cA) \cG(\la)$,  
where $I_\cX$ stands for the identity operator on $\cX$. Therefore, it suffices to show that $\mathcal{G}(\lambda)$ is bounded on $\cX$ and satisfies the desired norm bound if $|\la|\geq 1$.

By \eqref{eq: resolvent}, for any $y\in \Dom(A^{1/2})$ and $x\in X$, we have 
\begin{align}\label{ineq:for_G}
\begin{split}
\|\cG(\la)(y,x)\|_{ \cX} &  \leq \| R(\la)(\la + B)y\|_{X_{\frac{1}{2},A}}
+ \|R(\la)x\|_{X_{\frac{1}{2}, A}}\\
&\qquad  + \| R(\la)(\la^2 + \la B)y - y\| + \|\la R(\la) x\|.
\end{split}
\end{align}  
To estimate the first term on the right-hand side of \eqref{ineq:for_G}, we begin with auxiliary bounds.  Let $t: = |\la|^2.$
By the moment inequality, we have
\begin{align*}
\|A^{1/2}(t+A)^{-1}z\|^2 & \leq c_A \|A(t+A)^{-1}z\| \|(t+A)^{-1}z\|\\
&\leq c_A \big(\|(t+A)^{-1}\| + t\|(t+A)^{-1}\|^2 \big)\|z\|^2, \quad z\in X.
\end{align*}
Thus,
\begin{equation}\label{ineq:for_A_1}
\|A^{1/2}(t+A)^{-1}\|\leq \sqrt{c_A}(\|(t+A)^{-1}\| + t \|(t+A)^{-1}\|^2)^{1/2}\leq \sqrt{2c_A} \mu t^{-1/2},
\end{equation}
where $\mu: = \sup_{r>0} \|r(r+A)^{-1}\|$.
Furthermore,  
\begin{equation}\label{ineq:for_A_2}
\|A^{1/2} y\| \leq \|A^{1/2}(I+A)^{-1/2}\| \|y\|_{X_{\frac{1}{2},A}} =: \nu \|y\|_{X_{\frac{1}{2},A}},\quad  y\in \Dom(A^{1/2})
\end{equation}

We fix $y\in \Dom(A)\subset \Dom(A^{1/2})$ and $x\in X$. Let $y_0 := t (t+A)^{-1} y$ and  $y_1 : = A(t + A)^{-1}y$, and so $y = y_0 +y_1$.  
Since $R(\la)(\la + B)y_0 = \la^{-1}(y_0 -  R(\la) Ay_0)$, by the second bound in \eqref{bd X into V}, $\|y_0\|_{X_{\frac{1}{2},A}}\leq \mu \|y\|_{X_{\frac{1}{2},A}}$, \eqref{ineq:for_A_1}, and \eqref{ineq:for_A_2}, we obtain
\begin{align*}
\|R(\la)(\la+B)y_0\|_{X_{\frac{1}{2},A}} & = |\la|^{-1}\|y_0 - R(\la) Ay_0\|_{X_{\frac{1}{2},A}}\\
& \leq |\la|^{-1}\big( \mu \|y\|_{X_{\frac{1}{2},A}} + C(\la) \|A y_0\|\big)\\
& \leq  \big(|\la|^{-1} \mu \|y\|_{X_{\frac{1}{2},A}} +  t^{1/2} C(\la) \|A^{1/2}(t+A)^{-1}\| \|A^{1/2} y\| \big)\\
& \leq \big(|\la|^{-1} \mu  +  \sqrt{2c_A} \mu \nu  C(\la) \big)\|y\|_{X_{\frac{1}{2},A}}.
\end{align*} 
Moreover, combining \eqref{ineq:for_A_1}, and \eqref{ineq:for_A_2}, we infer that
\begin{align}\label{ineq:y1_est}
\nonumber \|R(\la)(\la+B)y_1\|_{X_{\frac{1}{2},A}} & \leq  C(\la) \|(\la +B)\| \|y_1\| \\
& \leq  C(\la) (|\la| +\|B\|)  \|A^{1/2}(t+A)^{-1} A^{1/2} y\|\\
\nonumber & \leq \sqrt{2c_A}\mu \nu C(\la) \big(1 + |\la|^{-1}\|B\|\big)\|y\|_{X_{\frac{1}{2},A}}.
\end{align}

 Consequently, 
\begin{equation}\label{ineq:for_A_3}
\|R(\lambda)(\lambda+B)y\|_{X_{\frac12,A}}
\lesssim \big( |\la|^{-1} + C(\la) + \|B\| C(\la) |\la|^{-1}  \big)\|y\|_{X_{\frac12,A}}.
\end{equation}
where the implicit constant is independent of $\la$. Since $\Dom(A)$ is dense in $\Dom(A^{1/2})$, \eqref{ineq:for_A_3} holds for every $y\in \Dom(A^{1/2})$.

For the third term in \eqref{ineq:for_G},  by Lemma \ref{prop on equiv cond_lemma}(ii), the operator $R(\la)A$ extends to a~bounded operator in $\cL(X_{\frac{1}{2},A},X)$ with norm at most $C(\la)\|A(I+A)^{-1}\|\leq C(\la)~{(1+\mu)}$.
Therefore, for every $y \in \Dom(A)$, we have
\begin{equation}\label{ineq:for_A_4}
\|R(\la)(\la^2 + \la B)y - y\| = \|R(\la)Ay\| \lesssim  C(\la) \|y\|_{X_{\frac{1}{2},A}}. 
\end{equation}
Therefore, by combining \eqref{ineq:for_A_3} and \eqref{ineq:for_A_4} and using again the second bound in \eqref{bd X into V} for the second term in \eqref{ineq:for_G}, we infer that
\begin{align*}
\|\cG(\lambda)(y,x)\|_{\cX} & \lesssim \big( |\la|^{-1} + C(\la) + \|B\| C(\la) |\la|^{-1} + \|\la R(\la)\| \big)  \|(y,x)\|_{\cX}.
\end{align*}

Therefore, to establish \eqref{norm of A} for $|\la| \geq 1$ it suffices to show that $|\la|^{-1} + C(\la) \lesssim \|\la R(\la)\|$. Let $0\neq z\in \Dom(A)$ with $\|z\| = 1$. Then 
\begin{equation}\label{eq:for_C}
1 =\| R(\la) P(\la) z\| \leq \|R(\la)\| (|\la|^2 + |\la| \|B\| + \|Az\|)
\end{equation}
Thus, if $|\la|\geq 1$ , then 
\[
|\la|^{-1} \leq \big(1 + |\la|^{-1}\|B\| + |\la|^{-2} \|Az\|\big) \|\la R(\la)\| \lesssim \|\la R(\la)\|.
\] 

Moreover, \eqref{eq:for_C} shows that
$|\la|^{-2}\|R(\la)\|^{-1}\leq 1+\|B\|+\|Az\|$ if $|\la|\geq1$.
Hence, by the definition of $C(\la)$ in part (i), we obtain
$C(\la)\lesssim\|\la R(\la)\|$ if $|\la|\geq1$, as desired. 
Finally, if $|\lambda|\leq 1$, then $C(\lambda)\lesssim \|R(\lambda)\|$,
with the implicit constant independent of $\lambda$. This implies
\eqref{norm of A} for $|\lambda|\leq 1$.

To show \eqref{2sided} note that for $|\lambda|\ge1$, we have
$|\lambda|^{-1}\|R(\lambda)\|=|\lambda|^{-2}\|\lambda R(\lambda)\|\le\|\lambda R(\lambda)\|$.
Moreover, \eqref{eq:for_C} gives $|\lambda|^{-1}\lesssim\|\lambda R(\lambda)\|$. Thus \eqref{norm of A} yields the asserted upper bound in \eqref{2sided}. Conversely, by \eqref{eq: resolvent},
\[
 (\lambda-\cA)^{-1}(0,x)=\bigl(R(\lambda)x,\lambda R(\lambda)x\bigr),
 \qquad x\in X,
\]
so that $\|(\lambda-\cA)^{-1}\|_{\cL(\cX)}\ge\|\lambda R(\lambda)\|_{\cL(X)}$. This completes the proof of \eqref{2sided}.
\end{proof}

Combining Proposition \ref{prop on equiv cond} with Fattorini's theorem on the square root reduction of cosine function generators and Theorem \ref{thm: DWE well-posed}, in the following result we provide the desired inclusion $V \subset X_2$. 

\begin{theorem}\label{thm on equiv cond}
Let $-A$ be the generator of a cosine operator function on a UMD space
$X$, and let $B\in\mathcal L(X)$. If
\begin{equation}\label{eq:-1,infty}
\{is: |s|\geq1\}\subset \rho(P)\quad \textrm{ and } \quad \sup_{|s|\geq 1}|s| \| R(is)\|<\infty,
\end{equation} 
then, for every $\omega\geq 0$ such that $A + \omega$ is sectorial, we have
\[
V = \Dom\big((A + \omega)^{1/2})\big) \subset X_2.
\]
\end{theorem}

\begin{proof} By Theorem~\ref{thm: DWE well-posed} and Fattorini's square root reduction theorem for generators of cosine functions (see
\cite[Theorem~5.4]{Fa85}), we have $V=V_0=\Dom((A+\eta^2)^{1/2})$, with equivalent norms, for any $\eta\geq 0$ such that $\|C_0(t)\|\lesssim e^{\eta |t|}$, $t\in \R$. Here $C_0$ denotes the cosine function generated by $-A$. The fact that $\Dom((A + \eta^2)^{1/2}) = \Dom((A +  \omega)^{1/2})$ for any $\omega\geq 0$ such that $A+\omega$ is a sectorial operator follows readily, for instance, from \cite[Proposition 3.8.2]{ABHN11}. Thus, $V = \Dom((A +  \omega)^{1/2})$ with equivalent norms. 

Fix $\omega\geq0$ such that $A_\omega:=A+\omega$ is sectorial.
Applying Proposition~\ref{prop on equiv cond} to $A_\omega$, we
obtain a constant $c>0$ such that, for every $\lambda\in\C\setminus\{0\}$ with $|\lambda|\geq1$ for which $\lambda^2+\lambda B+A_\omega$ is invertible, 
\[
\|(\la - \cA_\omega)^{-1}\|_{\cL(V\times X)} \leq c \|\la (\la^2 + \la B + A_\omega)^{-1}\|,
\]
where $\cA_\omega = \cA - \left[\begin{smallmatrix}0& 0\\ \omega & 0\end{smallmatrix}\right]$. 
By \eqref{eq:-1,infty}, a perturbation argument shows that, for every $|s|>s_0$, where $s_0: = \max(1, 2\omega \sup_{|s|\geq 1}|s|\|R(is)\|)$,  the operator $-s^2+isB+A_\omega$ is invertible and 
\[
\|(-s^2+isB+A_\omega)^{-1}\|\le 2\|R(is)\|.
\] As observed after \eqref{pencil vs resolvent} (see also
\eqref{eq: resolvent}), we also have $\{is:|s|>s_0\} \subset \rho(\cA_\omega)$, and
\[
\|is(-s^2+isB+A_\omega)^{-1}\| \le \|(is-\cA_\omega)^{-1}\|_{V\times X}, \quad |s|>s_0.
\] Moreover, $A_\omega$ is a generator of a cosine function; see, e.g.,
\cite[Corollary~3.14.13]{ABHN11}. Hence, by
Corollary~\ref{cor:DWE_well_posed}, the problem (HE) associated with
$A_\omega$ and $B$ is well posed.
Therefore, the operators $A_\omega$ and $B$ satisfy the assumptions of Proposition~\ref{lem: pencil vs resolvent} with $k=-1$. Hence, by Proposition~\ref{lem: pencil vs resolvent}(i) with
$j=2$, for every $(0,x)\in\Dom(\cA_\omega)
=\Dom(\cA)=\Dom(A)\times V$, or equivalently for every
$x\in V$, we have 
\[
\sup_{|s|> s_0} s^2 \|(-s^2 + is B + A_\omega)^{-1}x\| <\infty.
\]
Since, for all $|s|>s_0$, we have
\[
R(is)=\big(I + \omega R(is)\big)(-s^2+ isB+A_\omega)^{-1},
\]
and $I + \omega R(is)$, $|s|>s_0$, is uniformly bounded, we conclude
that $V\subset X_2$, as desired.
\end{proof}

\subsection{Decay of the total energy of solutions to (HE)}
 The ideas employed in the proofs of Propositions~\ref{prop:linear_growth} and~\ref{prop on equiv cond}, together with
stability results for \(C_0\)-semigroups, have a noteworthy implication for the
decay of the total energy of solutions to (DWE), that is, in the
Dirichlet Laplacian case, for the decay of
\[
 \|(-\Delta)^{1/2}u(t;u_0,u_1)\|_2+\|u'(t;u_0,u_1)\|_2
 \simeq
 \|\nabla u(t;u_0,u_1)\|_2+\|u'(t;u_0,u_1)\|_2 .
\]
We refer the reader to \cite{BaChTo16} for a seminal paper on the stability
theory of $C_0$-semigroups in the presence of the two possible singularities
of the resolvent of the generator on the imaginary axis, at zero and at infinity. We point out that the existing results provide a precise correspondence between the decay of semigroup orbits and the asymptotic behaviour of the resolvent for bounded $C_0$-semigroups. Since the semigroup $\cU$ on the phase space $V\times X$ associated with (DWE) can exhibit linear growth, see Remark~\ref{Goldstein example} and Proposition \ref{prop:linear_growth}, we have to apply such stability results to the (abstract) semigroup $\wt\cU$ on the energy space $\wt V\times X$, and then transfer the resulting conclusions to $\cU$ on $V\times X$. 
The next result gives a homogeneous energy-space counterpart of Proposition~\ref{prop on equiv cond} and also establishes contractivity of the associated semigroup.

In the next result, we use the notation from the proof of Proposition \ref{prop:linear_growth}.
Thus, we let $A$ be a positive selfadjoint operator on a Hilbert space $X$
and let $B$ be a bounded non-negative operator on $X$. We denote by
$\mathcal A$ the generator of the associated $C_0$-semigroup $\mathcal U$
on the phase space $V\times X$ for (HE). Moreover,
\[
 \dot V:=\operatorname{Dom}(A^{1/2}),\qquad
 \|x\|_{\dot V}:=\|A^{1/2}x\|,
\]
$\widetilde V$ denotes the Hilbertian completion of $\dot V$,
$J:\dot V\to\widetilde V$ is the canonical isometric embedding,
and $\widetilde{\mathcal A}$ is the closure in $\widetilde V\times X$
of the operator $\widetilde{\mathcal A}_0$ defined by
\[
 \widetilde{\mathcal A}_0(Jx,y):=(Jy,-Ax-By),\qquad
 \operatorname{Dom}(\widetilde{\mathcal A}_0)
 :=J\operatorname{Dom}(A)\times\operatorname{Dom}(A^{1/2}).
\]

\begin{lemma}\label{prop:disipative_A}
Let $A$ be a positive selfadjoint operator on a Hilbert space $X$ and  $B$ be a bounded non-negative operator on $X$.

Then the following conditions hold: 
\begin{itemize}
\item [(i)] $\wt\cA$ is the generator of a  contraction $C_0$-semigroup $\wt\cU$ on $\widetilde V \times X$. 

\item [(ii)] For any $y\in \Dom(A^{1/2})$ and $x\in X$ we have
\begin{align*}
\|A^{1/2}u(t)\|_X^2 + \|u'(t)\|_X^2 = \| \cU(t)(y,x)\|_{\dot V \times X}^2 & = \| \t\cU (Jy, x) \|^2_{\wt V \times X}  
\end{align*}
where $u$ stands for the mild solution to the corresponding \emph{(HE)} with initial data $y$ and $x$. In particular, $\cU(t)$ is a contraction on a normed space $\dot V\times X$. 
\item [(iii)] There exists a constant $c>0$ such that, for any $\la \in \rho(P)\setminus\{0\}$, we have $\la\in \rho(\wt \cA)$ and 
\begin{equation}\label{eq:est_resolvent_tcA}
\|(\la - \wt \cA)^{-1}\|_{\cL(\wt V \times X)}\leq c \max\big(|\la|^{-1}, |\la|^{-2}\|\la R(\la)\|,  \|\la R(\la)\|\big).
\end{equation}
\end{itemize}
\end{lemma}

\smallskip

\begin{proof} A straightforward inspection of the proof of Proposition~\ref{prop:linear_growth} yields (i) and (ii). The proof of (iii) follows the same lines as the proof of \eqref{norm of A}. We provide only the main additional observation needed in this case.

For \(\lambda\in\rho(P)\setminus\{0\}\), we define $\wt G(\lambda): J\Dom(A^{1/2})\times X\longrightarrow \wt V\times X$
by 
\[
 \wt G(\lambda)(Jy,x):=  \left(JR(\lambda)(\lambda+B)y+JR(\lambda)x,\,
 \lambda R(\lambda)(\lambda+B)y-y+\lambda R(\lambda)x \right).
\]
We first note that it suffices to show that there exists a constant $c>0$ such that for every $0\neq \la \in \rho(P)$, $y\in \Dom(A)$, and $x\in X$, we have 
\begin{equation}\label{eq:tildeG}
\|\wt \cG(\la)(Jy, x)\|_{\wt V\times X} \leq c \max\big(|\la|^{-1}, |\la|^{-2}\|\la R(\la)\|,  \|\la R(\la)\|\big) \|(Jy,x)\|_{\wt V\times X}.
\end{equation}
Indeed, since $\Dom(A)$ is a core for $A^{1/2}$, this estimate shows that $\wt\cG(\lambda)$ extends to a bounded operator on $\wt V\times X$, which we continue to denote by $\wt\cG(\lambda)$.
Moreover, it is readily seen that for any $0\neq \la \in \rho(P)$, we have
\[
(\la - \wt \cA)\wt \cG (\la) (Jy,x) = (Jy,x) \quad \textrm{and } \quad \wt \cG (\la) (\la - \wt \cA) (Jz,y) = (Jz,y) \
\]
for all $z\in \Dom(A)$, $y\in \Dom(A^{1/2})$,  and $x\in X$. Since $\la - \wt\cA$ is closed on $\wt V\times X$ and $\wt \cG(\la)$ is bounded on $\wt V\times X$, we infer that $\la \in \rho(\wt \cA)$ and $(\la - \wt \cA)^{-1} = \wt \cG(\la)$. Therefore, since $J\Dom(A)\times X$ is dense in $\wt V\times X$, \eqref{eq:tildeG} yields \eqref{eq:est_resolvent_tcA}.

To show \eqref{eq:tildeG} we fix $0 \neq \la \in \rho(P)$, $y\in \Dom(A)$ and $x\in X$.  As in \eqref{ineq:for_G}, 
\begin{equation}\label{ineq:tildeG}
\begin{split}
\|\wt \cG(\la)(Jy,x)\|_{\wt V\times X} &  \lesssim \| A^{1/2} R(\la)(\la + B)y\| 
+ \|A^{1/2} R(\la)x\|\\
&\qquad  + \|  R(\la)(\la^2 + \la B)y - y\| + \|\la R(\la) x\|.
\end{split}
\end{equation}
For the first term in \eqref{ineq:tildeG}, as in the corresponding part of the proof of Proposition~\ref{prop on equiv cond}, we set $y_0: = t(t+A)^{-1}y$ and $y_1: = A(t+A)^{-1}y$, where $t=|\la|^2$. Let $\mu := \sup_{r>0}\|r(r+A)^{-1}\|$. We recall that 
\begin{equation}\label{ineq:z}
\|A^{1/2}z\|\lesssim \|z\|_{X_{\frac{1}{2},A}},\qquad z\in \Dom(A^{1/2}).
\end{equation}
Hence, since $R(\la)(\la + B) y_0 = \la^{-1} ( y_0 - R(\la)Ay_0)$,  combining \eqref{ineq:z}, the second estimate in \eqref{bd X into V}, and \eqref{ineq:for_A_1}, we obtain
\begin{align*}
\| A^{1/2} R(\la)(\la + B)y_0\| & \leq |\la|^{-1}\big(\mu \|A^{1/2} y\| + \|R(\la)Ay_0\|_{X_{\frac{1}{2},A}}\big)\\
& \leq \big( |\la|^{-1}\mu  + \sqrt{2}\mu C(\la)\big) \|A^{1/2} y\|
\end{align*}
Moreover, combining (4.36), the second inequality in (4.28), and (4.26), we obtain 
\begin{align*}
\|A^{1/2}R(\la)(\la + B)y_1\|\lesssim \|R(\la)(\la + B)y_1\|_{X_{\frac12,A}} \leq \sqrt{2} \mu C(\la) (1+ |\la|^{-1}\|B\|)\|A^{1/2}y\|,
\end{align*}
Consequently, since $y = y_0 + y_1$, we have 
\[
\|A^{1/2}R(\la)(\la+B)y\|\lesssim \big(|\la|^{-1} + C(\la) + |\la|^{-1} \|B\| C(\la)\big)\|A^{1/2}y\|.
\]
For the second term in \eqref{ineq:tildeG}, by \eqref{ineq:z} and the second bound in \eqref{bd X into V}, we obtain
\[
\|A^{1/2}R(\la)x\| \lesssim C(\la)\|x\|.
\]
For the third term, since $R(\la)(\la^2 +\la B)y - y = - R(\la)A^{1/2}(A^{1/2}y)$,  it suffices to show that $R(\la)A^{1/2}$ has the extension to an operator in $\cL(X)$. By Lemma \ref{prop on equiv cond_lemma}(ii), $R(\la)$ extends to an operator in $\cL(X_{-\frac12,A},  X)$ with the norm less than or equal to $C(\la)$. For any $z\in \Dom(A^{1/2})$, we have 
\[
\|A^{1/2} z\|_{X_{-\frac{1}{2},A}} = \|A^{1/2}(I+A)^{-\frac{1}{2}} z\|\leq  
 \|A^{1/2}(I+A)^{-\frac{1}{2}}\|  \|z\|.
\]
Since $\Dom(A^{1/2})$ is a dense subset of $X$, the operator $A^{1/2}$ admits the extension to an operator in $\cL(X,X_{-\frac{1}{2}, A})$, which yields the desired extension of $R(\la)A^{1/2}$. Thus, 
\[
\|R(\la)A y\| \leq C(\la)\|A^{1/2}(I+A)^{-1/2}\| \|A^{1/2} y\|.
\]
Applying to the fourth term in \eqref{ineq:tildeG} the bound
$\|\la R(\la)x\|\leq \|\la R(\la)\| \|x\|$, we have
\[
\|\wt \cG(\la)(Jy, x)\|_{\wt V\times X} \lesssim \big(|\la|^{-1} + C(\la) + |\la|^{-1} \|B\| C(\la) + \|\la R(\la)\| \big) \|(Jy, x)\|_{\wt V \times X}.
\]
The same arguments as those used at the end of the proof of Proposition~\ref{prop on equiv cond} yield the desired constant in \eqref{eq:tildeG}, thereby completing the proof.
\end{proof}
\begin{remark}
Under the assumptions of Lemma~\ref{prop:disipative_A}, the map
$W(Jx)=A^{1/2}x$, $x\in\dot V$, extends to a unitary
operator $W:\widetilde V\to X$. Under this unitary
transformation,
\[
(W\oplus I)\widetilde{\mathcal A}(W\oplus I)^{-1}
=
\begin{pmatrix}
0 & A^{1/2}\\
-A^{1/2} & -B
\end{pmatrix},
\]
with domain
$\Dom(A^{1/2})\times\Dom(A^{1/2})$.
Consequently,
\[
\Dom(\widetilde{\mathcal A})
=
\bigl\{(v,x)\in\widetilde V
\times\Dom(A^{1/2}):
Wv\in\Dom(A^{1/2})\bigr\}.
\]
This provides an explicit description of $\Dom(\widetilde{\mathcal A}).$
$\hfill\Diamond$\end{remark}

We recall the stability for semigroup orbits under the singularity at zero; see \cite[Theorem 7.6]{BaChTo16}.
Let $(\wt\cU(t))_{t\geq 0}$ be a bounded $C_0$-semigroup on a Hilbert space $\cX$
with generator $\wt\cA$. Assume that
$\sigma(\wt\cA)\cap i\mathbb R\subset\{0\}$, and let $\alpha>1$. Then 
\begin{equation}\label{equi:decay_U}
\|(is - \wt\cA)^{-1}\|
=
\begin{cases}
O\left(|s|^{-\alpha}\right),
& s\to0,\\
O(1), & |s|\to\infty,
\end{cases}\,\, \Longleftrightarrow \,\,
\|\wt\cU (t) \wt\cA(I-\wt\cA)^{-1}\|
=O\left(t^{-\frac{1}{\alpha}}\right),\,\, t\to\infty.
\end{equation}

Combining Lemma~\ref{prop:disipative_A} with the equivalence in \eqref{equi:decay_U}, we obtain the following result.

\begin{proposition}\label{prop:total_energy}
Suppose that $A$ and $B$ satisfy the assumptions of Lemma~\ref{prop:disipative_A}, $i\R\setminus\{0\} \subset \rho(P)$, and $\sup_{s\neq 0} \min(|s|^{1+\beta},|s|) \|R(is)\|_{\cL(X)}<\infty$ for some $\beta\geq 0$. 

Then, for every 
$y\in \Dom(A^{1/2})$ and $x\in \Ran(A) - By$, there exists a constant $c>0$ such that   
\[
\|A^{1/2}u(t)\|_X + \|u'(t)\|_X \leq c t^{-1/(2+\beta)}, \qquad t\geq 1,
\]
where $u$ denotes the mild solution to \emph{(HE)} with initial data $y$ and $x$.
\end{proposition}

\begin{proof}
Let $y\in\Dom(A^{1/2})$ and $x\in\Ran(A)-By$.  Choose $z\in\Dom(A)$ such that 
$x=Az-By$. Then $(Jy,x)=\wt\cA(-Jz,y)\in\Ran(\wt\cA)$. Since  $1\in\rho(\wt\cA)$, we have $\Ran(\wt\cA)
 =
 \Ran\bigl(\wt\cA(I-\wt\cA)^{-1}\bigr).
$
The resolvent hypothesis of Proposition~\ref{prop:total_energy} and
Lemma~\ref{prop:disipative_A}(iii) give 
\[
 \|(is-\wt\cA)^{-1}\| =
 \begin{cases}
 O(|s|^{-2-\beta}),&s\to0,\\
 O(1),&|s|\to\infty.
 \end{cases}
\]
Hence, by \cite[Theorem 7.6]{BaChTo16},
\[
 \|\wt\cU(t)\wt\cA(I-\wt\cA)^{-1}\| =  O(t^{-1/(2+\beta)}).
\]
Applying this to $(Jy,x)\in\Ran(\wt\cA)$ and using Lemma~4.9(ii),
\[
 \|\wt\cU(t)(Jy,x)\|_{\wt V\times X}^2
 =
 \|A^{1/2}u(t)\|^2+\|u'(t)\|^2,
\]
and the assertion follows.
\end{proof}

\subsection{Mild solutions to (PE) and (HE) via Laplace transform}\label{subsec:well-posed-DE}

We present auxiliary results concerning the Laplace transform representation of mild solutions to (HE) and (PE).

For $f\in L^1_{\mathrm{loc}}(\R_+;X)$, let $\operatorname{abs}(f)$ denote the abscissa of convergence of its Laplace transform:
\[
\operatorname{abs}(f):=
\inf\{\operatorname{Re}\lambda:(\cL f)(\lambda)\text{ exists}\},
\]
In the half-plane of convergence, the integral in the definition of $\cL f$ is understood as
\[
(\cL f)(\lambda)
=\lim_{\tau\to\infty}\int_0^\tau e^{-\lambda t}f(t)\, \ud t.
\]

\paragraph{\bf{The mild solutions to (PE)}}

Let $A$ be a linear, closed operator on a Banach space $X$ and $B\in \cL(X)$. 
The solvability of (PE) can be also studied within the framework of multivalued generators of semigroups; see Favini and Yagi~\cite{FaYa99}. As we are primarily interested in the case where $A$ and $B$ satisfy Hypothesis \ref{hyp:main}, we present an alternative approach based on the Laplace transform techniques. 

Following the standard semigroup theory (corresponding to $B=I$), by a {\it classical} solution of (PE) we mean a function $v\in C^1(\R_+; X)$ such that $v(t) \in \Dom(A)$ for any $t\geq 0$ and $v$ satisfies (PE). A function $v\in C(\R_+;X)$ is called a {\it mild} solution to (PE) if 
\[
\int_0^t v(s) \ud s \in \Dom(A) \qquad \textrm{and} 
\qquad  Bv(t) - By_0 = - A \int_0^t v(s) \ud s, \quad t\geq 0.
\]

The following result is a direct counterpart of \cite[Theorem~3.1.3]{ABHN11}. Since the proof carries over almost verbatim, we omit it.

\begin{proposition}\label{mild-DE} 
Let $A$ be a closed linear operator on a Banach space and $B \in \cL(X)$. 
Let $y_0\in X$ and $v\in C(\R_+;X)$ be Laplace transformable, that is, $\abs(v)<\infty$. Let $\omega \geq \max(0,\abs(v))$.
Then the following conditions are equivalent:
\begin{itemize}
\item [(i)] $v$ is a mild solution to \emph{(PE)} with the initial value $y_0$;
\item [(ii)] For all $\la >\omega$ 
\begin{equation} 
(\cL v)(\la) \in \Dom(A) \quad \textrm{ and  } \quad \la B (\cL v)(\la) + A (\cL v)(\la) = B y_0.
\end{equation} 
\end{itemize}
If $(\omega, \infty)\subset \rho(P_D)$, then for the mild solution $v$ to \emph{(PE)} we have 
\[
(\cL v)(\la) = (\la B +A)^{-1} B y_0,  \qquad \la > \omega.
\]
\end{proposition}

\paragraph{\bf{The mild solutions to (HE)}}
Let $A$ be a closed linear operator on a Banach space and $B \in \cL(X)$. 
Let $x_0,x_1\in X$.  
A {\it mild solution} to (HE) with initial data $x_0$ and $x_1$ is a function $u \in C(\mathbb{R}_+,X)$ such that
\[
\int_0^t \int_0^s u(r)\,dr\,ds
= \int_0^t (t-s)u(s)\,ds \in \Dom(A)
\]
and
\begin{equation}\label{mild sol formula}
u(t) = x_0 + t x_1 - B\int_0^t u(s) \ud s + t Bx_0  - A \int_0^t (t-s)u(s)\,ds \end{equation}
for all $t \ge 0$. If $u$ is a classical solution, then integrating (HE) twice shows that
$u$ is a mild solution. Conversely, if $u$ is a mild solution and
$u \in C^2(\mathbb{R}_+,X)$, then $u$ is a classical solution. This follows from \eqref{mild sol formula} and the fact that $A$ is closed.

\begin{proposition}\label{mild_DWE}
Let $u \in C(\mathbb{R}_+,X)$ with $\operatorname{abs}(u) < \infty$.
Let $\omega \geq \max\big(\mathrm{abs}(u),0 \big)$.
Then the following assertions are equivalent: 
\begin{itemize}
\item [(i)] $u$ is a mild solution to \emph{(HE)} with initial data $x_0$ and $x_1$.
\item [(ii)] For all $\la >\omega$
\begin{equation}
(\cL u)(\lambda) \in \Dom(A)
\quad \text{and} \quad
\lambda x_0 + B x_0 + x_1 = (\lambda^2  +  \la B +  A)(\cL u)(\lambda).
\end{equation}
\end{itemize}
\end{proposition}

\begin{proof}The proof is a consequence of the uniqueness theorem for the Laplace transform.  
Indeed, since for every $\la > \omega$
\[
(\cL u) (\la)  = \la^2 \int_0^\infty e^{-t\la} \int_0^t \int_0^s u(r) \ud r \ud s \ud t,
\]
and 
\begin{align*}
\la^2 \cL \Big( -u - B\int_0^t u(s) \ud s & + x_0 + t x_1 + t B x_0 \Big)(\la)\\
& = -\la^2 (\cL u)(\la) - \la B (\cL u)(\la) + \la x_0 + Bx_0 + x_1,
\end{align*}
the equivalence follows easily from \cite[Proposition 1.7.6]{ABHN11}.
\end{proof}

 The statements and proofs of Theorems \ref{thm: diffusion intro} - \ref{thm:diffusion_for_Delta} rely on the following result, which guarantees the existence of mild solutions to (HE).

\begin{lemma}\label{lem: mild s to DWE} If \emph{(HE)} is well posed and the corresponding $C_0$-semigroup $(\cU(t))_{t\geq 0}$ is polynomially bounded in $\cL(V\times X)$, then for any pair $(x_0, x_1) \in V\times X$ there exists a polynomially  bounded mild solution $u$ to \emph{(HE)} with the initial data $x_0$ and $x_1$.
\end{lemma}

\begin{proof} According to the definition of the well-posedness of {(HE)} given in Subsection~\ref{subsec well posed}, for any $(x_0,x_1)\in V\times X$ the function
\[
\mathbb{R}_+ \ni t \longmapsto \mathcal U(t)(x_0,x_1) \in V\times X
\]
is a mild solution to the Cauchy problem
\[
\vec u'(t)=\mathcal A\vec u(t),\quad  t\geq 0,
\qquad \vec u(0)=(x_0,x_1).
\]

We show that the first component of $\vec{u}$ that we denote by $u$ is a mild solution to {(HE)} with the initial data $x_0$ and $x_1$.
Indeed, since $V$ is continuously embedded into $X$, we have $u \in C(\mathbb{R}_+;V) \subset C(\mathbb{R}_+;X)$.
Moreover, since $\mathcal{U}$ is polynomially bounded in $\mathcal{L}(V\times X)$ and $\|u(t)\|_V \leq \|\mathcal{U}(t)(x_0,x_1)\|_{V\times X}$,
and using again the continuous embedding of $V$ into $X$, it follows that $u$ is polynomially bounded. In particular, the abscissa of convergence of the Laplace transform of $u$ in $X$ is less than or equal to $0$.
According to the observations made at the beginning of Subsection~\ref{subsec well posed}, we have $\{\Re\lambda>0\}\subset\rho(P)$. Since
\[
\mathcal L\big(\mathcal U(\cdot)(x_0,x_1)\big)(\lambda)
=(\lambda-\mathcal A)^{-1}(x_0,x_1),
\qquad \Re\lambda>0,
\]
in $V\times X$, it follows from
\eqref{pencil vs resolvent} that
\[
\mathcal L u(\lambda)
= R(\lambda)\big(\lambda x_0 + Bx_0 + x_1\big),
\qquad \Re\lambda>0,
\]
where the convergence of $\mathcal L u(\lambda)$ holds in $V$, and hence in $X$. By Proposition~\ref{mild_DWE}, the proof is complete.
\end{proof}

\subsection{A maximum principle}\label{sec:max-principle}

In this section, we establish the Phragmén--Lindelöf principle for holomorphic functions on the half-plane $\C_+ = \{z \in \C : \Re z > 0\}$ with a singularity at zero. The following proposition is a polynomial growth variant of \cite[Lemma~5.3]{BaChTo16}. We include the proof since the slightly more general form used here is convenient in several places below.

\begin{proposition}\label{thm:MP-fudge-JEMS-correct}
Let $X$ be a Banach space and let $f:\overline{\C_+}\setminus \{0\} \to X$ be a continuous function, holomorphic on $\C_+$, and such that
\[
M := \sup_{s \neq 0} \|f(is)\| < \infty.
\]

Assume that there exist constant $c > 0$ and integers $N, K \ge 0$ such that
\begin{equation}\label{eq:MP-growth-JEMS-correct}
\|f(z)\| \le c\, {(1 + |z|)^N}{(\Re z)^{-K}}, \qquad z\in \C_+.
\end{equation}

Then
\[
\|f(z)\| \le M, \qquad z\in \C_+.
\]
\end{proposition}

\begin{proof}
We set  $g(z):={f(z)}/{(1+z)^N}$,  $z\in \overline{\C_+}\setminus \{0\}$. Clearly, 
$\|g(is)\|\le M$, $s\neq 0$,  and, since $|1+z|\ge (1+|z|)/2$ for any $z\in \C_+$, we have
\begin{equation}\label{eq:N0-growth-JEMS}
\|g(z)\|
\le 2^N c\,/(\Re z)^{K}, \qquad z\in \C_+.
\end{equation}
For any $0<r<1<\rho$ we put $A_{r,\rho}:=\{z\in\C_+:\ r\le |z|\le \rho\}$ 
and define
\[
g_{r,\rho}(z):=
\bigl({z}{(1+z)^{-1}}\bigr)^{\!K}
\bigl(1+{z^2}{\rho^{-2}}\bigr)^{\!K}
\bigl(1+{r^2}{z^{-2}}\bigr)^{\!K} g(z),
\qquad z\in A_{r,\rho}.
\]
The function $g_{r,\rho}$ is holomorphic in the interior of $A_{r,\rho}$ and continuous up to the boundary $\partial A_{r,\rho}$. 
We estimate $g_{r,\rho}$ on the three components of $\partial A_{r,\rho}$. 
On the inner half-circle $\{|z|=r,\, \Re z\ge 0\}$, writing $z=re^{\ii\theta}$, $|\theta|\le\pi/2$, we have
\[
        \left|1+ {r^2}{z^{-2}}\right|
        =|1+e^{-2\ii\theta}|
        =2\cos\theta
        =2 {\Re z}/{r}.
\]
Moreover $|z/(1+z)|\le r/(1-r)$ and $|1+z^2/\rho^2|\le2$.  Hence, by \eqref{eq:N0-growth-JEMS},
\[
        \|g_{r,\rho}(z)\|
        \le
        \left({r}{(1-r)^{-1}}\right)^K
        2^K
        \left(2{\Re z}/{r}\right)^K
        2^N c \, (\Re z)^{-K}
        \le {2^N 4^K c}{(1-r)^{-K}}.
\]
On the outer half-circle $\{|z|=\rho,\ \Re z\ge0\}$, similarly,
\[
        \left|1+{z^2}{\rho^{-2}}\right|
        =|1+e^{2\ii\theta}|
        =2 {\Re z}/{\rho},
\]
while $|z(1+z)^{-1}|\le \rho/(\rho-1)$ and $|1+r^2z^{-2}|\le2$.  Therefore
\[
        \|g_{r,\rho}(z)\|
        \le
        \left({\rho}{(\rho-1)^{-1}}\right)^K
        \left(2{\Re z}/{\rho}\right)^K
        2^K 2^N c\, (\Re z)^{-K}
        \le {2^N 4^K c }{(\rho-1)^{-K}}.
\]
Finally, on the vertical part $\{\ii s: r\le |s|\le \rho\}$, we have
\[
        \left|{\ii s}{(1+\ii s)^{-1}}\right|\le1,
        \qquad
        \left|1-{s^2}{\rho^{-2}}\right|\le1,
        \qquad
        \left|1-{r^2}{s^{-2}}\right|\le1,
\]
and so
\[
        \sup_{r\le |s|\le \rho}\|g_{r,\rho}(\ii s)\|
        \le \sup_{s\ne0}\|g(\ii s)\|.
\]
The maximum principle for Banach-valued holomorphic functions implies that
\[
\sup_{z \in A_{r,\rho}} \|g_{r,\rho}(z)\| \le c_1,
\]
where $c_1$ is independent of $r$ and $\rho$, provided that $0 < r \le \tfrac{1}{2}$ and $\rho \ge 2$. Letting $r\rightarrow0$ and $\rho\rightarrow\infty$ with $z$ fixed yields
\begin{equation}\label{eq:weighted-G-bound}
        \left\|\left({z}{(1+z)^{-1}}\right)^K g(z)\right\|
        \le c_1,
        \qquad z\in\C_+.
\end{equation}
In particular, $g$ is bounded on $\{z\in\C_+:|z|\ge1\}$, while
\begin{equation}\label{eq:G-small-poly}
        \|g(z)\|\le 2^K c_1\, |z|^{-K}
        \quad \textrm{ for all } z \in \C_+ \textrm{ with } 0<|z|\le1.
\end{equation}
We now show that $g$ is uniformly bounded for all $z \in \C_+ \cap \{0<|z|\leq 1\}$. Let
\[
\wt g(z):=g(1/z), \qquad z\in \overline{\C_+}\setminus \{0\}.
\]
From the boundedness of $g$ we readily see that $\wt g$ is uniformly bounded on $\ii \R \setminus \{0\}$ and $\C_+ \cap \{|z| \leq 1\}$, and, by \eqref{eq:G-small-poly}, $\wt g$ has at most polynomial growth as  $|z|\rightarrow \infty$ and $z\in \mathbb{C}_+$.  Thus, by the standard Phragmen--Lindelof principle, $\wt g$ is uniformly bounded on $\overline{\C_+}\setminus\{0\}$. Equivalently, $g$ is bounded uniformly on $\overline{\C_+}\setminus\{0\}$. 
Consequently,   $f(z)=(1+z)^N g(z), z \in \mathbb C_+,$ is
bounded near zero and has at most polynomial growth as $|z|\to \infty$.  Another application of the Phragm\'en--Lindel\"of principle, together with the boundary bound $\sup_{s\ne0}\|f(is)\|\le M$,
yields $\|f(z)\|\le M$ for all $z\in\mathbb C_+$.
\end{proof}

\section{Proofs of the main abstract results}\label{sec:proofs}

\subsection{Proofs of Theorems \ref{thm:ind-solu} -- \ref{thm: diffusion} }

\begin{proof}[{Proof of Theorem \ref{thm:ind-solu}}]
The assumption on the polynomial growth of $u$ and $v$ implies that $\abs(u)$, $\abs(v)\leq 0$. Since $u$ is a mild solution to (HE), using Proposition \ref{mild_DWE}, for all $\la > 0$ we have
\begin{equation}\label{rep:u}
(\cL u)(\la) = R(\la)(\la x_0 + B x_0 + x_1).  
\end{equation}
Clearly, $\cL(u)$ is holomorphic on $\{\Re \la >0\}$ and since, by our assumption, the right-hand side of \eqref{rep:u} is also holomorphic on $\{\Re \la >0\}$, equality \eqref{rep:u} holds for all $\la$ with $\Re \la>0$.
Similarly, since $v$ is a mild solution to (PE), by Proposition \ref{mild-DE}, we have 
\begin{equation}\label{ref:v}
(\cL v)(\la) = R_D(\la)B y_0, \qquad \Re \la > 0.
\end{equation}

Let $\widetilde u$ and $\widetilde v$ denote the extensions of $u$ and $v$ from $\R_+$ to $\R$  by $0 \in X$. We set $\widetilde w : = \widetilde u - \widetilde  v$  and $\widetilde w_\varepsilon (t) : = e^{-\varepsilon t} \widetilde w (t) $, $t\in \R$ and $\varepsilon >0$. Since $\wt w$ has at most polynomial growth, $\wt w \in \cS'(\R;X)$. Moreover, $\|\widetilde w_\varepsilon(t)\| \leq \|\widetilde w(t)\|$ for all $t\in \R$ and $\varepsilon>0$, and $\widetilde w_\varepsilon \rightarrow \widetilde w$ uniformly on compact subset of $\R$. Hence, we have $\cF \widetilde w_\varepsilon \rightarrow \cF \widetilde w$ as $\varepsilon \rightarrow 0+$ in $\cS'(\R; X)$.  By the representations of $\cL u$ and $\cL v$ in \eqref{rep:u} and \eqref{ref:v}, for all $\varepsilon > 0$ and $s\in \R$ we infer
\begin{align*}
\cF \widetilde w_\varepsilon (s) & = (\cL u)(is +\varepsilon) - (\cL v)(is + \varepsilon).
\end{align*}
Further, if $N\in \N$ is a constant such that 
\[
\|u(t)\|, \|v(t)\| \lesssim 1+t^N, \quad t\geq 0,
\] 
then by a direct estimation of the Laplace transforms of $u$ and $v$ we obtain 
\[
\|(\cL u)(\la)\|, \|(\cL v)(\la)\| \lesssim {(\Re \la)^{-1}} +  {N!} {(\Re \la)^{-N-1}}, \quad \Re \la >0.
\]

Combining these estimates with \eqref{pencils est} and \eqref{hyp:bound_2} (cf. \eqref{def of m2} and \eqref{def of m3}), and using Proposition \ref{thm:MP-fudge-JEMS-correct}, we obtain 
\[
\sup_{\Re \la >0} \| (\cL u) (\la) - (\cL v)(\la) \|_{X}\lesssim [R]\|x_0\| + \max\big([R][R_D B]_\theta, [R]+ [R_D B]_\theta \big) \|y_0\|.
\]
Consequently,  $\cF \widetilde w_\varepsilon \in L^\infty(\R;X)$ are uniformly bounded for all $\varepsilon>0$.
Since, 
$\cF \widetilde w_\varepsilon (s) \rightarrow m(s)$ as $\varepsilon \rightarrow 0+$ in $X$ for all $s\neq 0$, by Lebesgue's domination convergence theorem, $\cF \widetilde w_\varepsilon \rightarrow m$ as $\varepsilon \rightarrow 0+$ in $\cS'(\R;X)$. Hence, $m = \cF \widetilde w$ or 
$\cF^{-1} m = \widetilde w$ in $\cS'(\R;X)$. 
\end{proof}

\begin{proof}[Proof of Theorem \ref{thm:mild_DPE}]  (i) For any $\omega'\in (0,\omega)$, the integral \eqref{formula:semig} is absolutely convergent in $\cL(X)$, uniformly for $z$ in compact subsets of $\Sigma_{\omega'}$, and thereby $z\to \cT(z)$ is holomorphic on $\Sigma_{\omega'}$. By Cauchy's theorem, the definition of $\cT(z)$ is independent on the choice of $\rho>0$ and $\omega' \in (0,\omega)$ (so long as $|\arg z|<\omega'$). Thus, $z \to \cT(z)$ is well-defined and holomorphic on $\Sigma_\omega$. 
 
The first bound in (i) follows directly from \eqref{formula:semig} by choosing $\rho:=1/|z|$ in $\Gamma(\rho, \omega')$ with any $\omega' \in (\wt \omega, \omega)$. For $|\arg z|\leq\widetilde\omega$ and $0<|z|<1$,
the same choice $\rho=|z|^{-1}$ yields
\[
\|\cT(z)\|_{\mathcal L(X)}=O(|z|^{\theta-1}),
\qquad |z|\to0,\ z\in\Sigma_{\wt \omega}.
\]
Combining Cauchy's theorem with the identity
\begin{equation}\label{resolvent_identity}
(\lambda B +A)^{-1}B-(\mu B +A)^{-1}B
=(\mu-\lambda)(\lambda B+A)^{-1}B(\mu B+A)^{-1}B,
\end{equation}
for all $\lambda, \, \mu\in\Sigma_{\omega+\pi/2},$ one concludes that the family $\cT(z)$, $z\in\Sigma_{\omega}$, satisfies the semigroup property. 

For the remaining assertion, the estimates on $\cT$ just established imply that the Laplace integral $(\mathcal{L}\cT)(\lambda)$ of $\cT$ is absolutely convergent for every $\lambda$ with $\Re\lambda>0$. Applying \eqref{formula:semig} with $\rho<|\lambda|$, followed by Fubini's theorem and Cauchy's theorem, yields the desired identity \eqref{eq:laplace_int}.

(ii) Let $q:\Sigma_{\omega+ \pi/2} \to X$ be given by $q(\la): = R_D(\la)Bx$, where $x\in \Dom(A)$ with $Ax = By$ for some $y\in X$. 
Applying condition \eqref{DE-parab_est}, we first note that 
\[
\la q(\la)  = x - R_D(\la)Ax  = x - R_D(\la) By \to x  \quad \textrm{ as } |\la| \to \infty, \la \in \Sigma_{\omega+\pi/2}.
\]
Moreover, by \eqref{DE-parab_est},
$\lambda q(\lambda)$ is uniformly bounded for
$0<|\lambda|\leq1$. Together with the preceding
identity, this yields
\[
\sup_{\lambda\in\Sigma_{\omega+\pi/2}}
\|\lambda q(\lambda)\|<\infty.
\]
Thus, combining \eqref{eq:laplace_int} with \cite[Theorems~2.6.1 and~2.6.4(a)]{ABHN11}, we conclude that $\cT(t)x\to x$ as $t\to0+$. Consequently, by Proposition~\ref{mild-DE}, the function $v$ is a mild solution of (PE) with initial value $x$. By (i), the mapping $\Sigma_{\omega}\ni z\longmapsto \cT(z)x$ is a holomorphic extension of $v$, and its boundedness on  $\Sigma_{\wt \omega}$ follows from the boundedness of $\lambda q(\lambda)$ on $\Sigma_\omega$; cf. \cite[Theorem~2.6.1]{ABHN11}.

\smallskip 
 
(ii') By \eqref{DE-parab_est}, the integral defining $J_\beta$ converges
absolutely in $\cL(X)$. 
Using \eqref{eq:laplace_int}, Fubini's theorem, and Euler's identity $\Gamma(\beta)\Gamma(1-\beta) = \frac{\pi}{\sin(\pi\beta)}$,
we obtain the following alternative representation of \(J_\beta\), also with absolute convergence in \(\mathcal L(X)\): 
\begin{equation}\label{eq:J-semigroup}
J_\beta = \frac{1}{\Gamma(\beta)} \int_0^\infty t^{\beta-1}e^{-t}\cT (t)\,\ud t.
\end{equation} 
By the semigroup law, for any $h>0$, we have
\[
\cT (h)J_\beta = \frac1{\Gamma(\beta)} \int_0^\infty t^{\beta-1} e^{-t} 
\cT(t+h)\,\ud t.
\]
Since $\cT$ is continuous on $(0,\infty)$ in operator norm and $\|\cT(t+h)\| \lesssim \max (1, t^{\theta -1})$ for every $t>0$ and $h\in (0,1)$, by the Lebesgue dominated convergence theorem, we infer that $\cT(h) J_\beta \to J_\beta$ in $\cL(X)$ as $h\to 0+$. That $v$ is a mild solution follows from
Proposition~\ref{mild-DE} and~\eqref{eq:laplace_int}, and the existence of its bounded holomorphic extension to $\Sigma_{\wt\omega}$ follows from part~(i).
The remaining assertion can be deduced either from \cite{ABHN11}
together with \eqref{eq:laplace_int}, or by a direct computation based
on \eqref{eq:J-semigroup} and \eqref{resolvent_identity}.

\smallskip

(iii) By Gelfand's formula for the spectral radius it suffices to show that 
$\Ran(B) \setminus \Ran(A) \neq \emptyset$ implies 
\begin{equation}\label{formula:spectral_radius}
\lim_{t\to \infty} t^{-1} \ln \|\cT(t)\| \geq 0.
\end{equation} 
Indeed, for a given $t>0$, by the semigroup property we have
\begin{align*}
r_{sp}(\cT(t)) = \lim_{n\to \infty} \|\cT(nt)\|^{1/n} = \exp\big(t \lim_{n\to \infty} (nt)^{-1} \ln \|\cT(nt)\|  \big).
\end{align*} 
To prove~\eqref{formula:spectral_radius}, we argue by contradiction. First, since the function $t\mapsto \ln\|\cT(t)\|$ is subadditive, by a standard argument we have
\[
\lim_{t\to \infty} t^{-1}\ln\|\cT(t)\|= \inf_{t\geq 1} t^{-1} \ln\|\cT(t)\|.
\] 
Therefore, suppose that $\Ran(B) \setminus \Ran(A) \neq \emptyset$ and  
$\lim_{t\to \infty} t^{-1}\ln\|\cT(t)\| < - \varepsilon$ for some $\varepsilon >0$. Thus there exists a constant $c>0$ such that 
$\|\cT(t)\| \leq c e^{-\varepsilon t}$ for all $t\geq 1$. 
 Combining this with the second condition in (i) shows
that the abscissa of convergence of the Laplace transform
$\mathcal LT$ of the operator-valued function
$(0,\infty)\ni t\mapsto T(t)$ is at most $-\varepsilon$,
and hence $\mathcal LT$ is holomorphic on the half-plane
$\{\operatorname{Re}\lambda>-\varepsilon\}$.
By \eqref{eq:laplace_int} and the uniqueness theorem
for the Laplace transform,
\(
(\mathcal LT)(\lambda)=(\lambda B+A)^{-1}B,\)
for \(\operatorname{Re}\lambda>0.
\)
In particular, $\sup_{\Re\lambda>0,\ |\lambda|\le1}
|(\lambda B+A)^{-1}B\|<\infty$. 
Therefore, there exists $\lambda$ with $\Re \la >0$ and sufficiently small $|\la|$ such that $\|\la (\la B +A)^{-1} B \|<1$. For each $y\in X$ we define 
\[
 w : = (I- \la (\la B+A)^{-1}B)^{-1} y \quad \textrm{and} \quad x : = (\la B + A)^{-1}B w. 
\] 
Clearly, $w - \la x = y$ and $x\in \Dom(A)$.  Moreover, 
\[
Ax  = (A + \la B - \la B) (\la B + A)^{-1}B w = B w - \la B x = By, 
\]
Thus, $\Ran(B)\subset \Ran(A)$, which contradicts the assumption and thereby completes the proof of the first statement in~(iii). 

Now the remaining assertion follows from an application of M\"uller's theorem, along the lines of the proof of \cite[Lemma 3.1.7]{JvN96}.   
A key auxiliary observation that guarantees the possibility of choosing $y\in \bigcap_{\ell \in \N} D_\ell$ is the following \emph{shifting} argument, which exploits the fact that 
\begin{equation}\label{inc:S_in_core}
\cT(t)X\subset D_\ell, \quad t>0, \, \ell \in \N.
\end{equation}
Namely, let the function $r$ be given. Let $\wt r:\R_+ \to \R_+$ be defined by $\wt r(t) : = r(0)$ if $t\in [0,1]$ and $\wt r(t) := r(t-1)$ for $t\geq 1$. Following the proof of \cite[Lemma 3.1.7]{JvN96} there exists $\wt y\in X$ such that $\wt r(t) \lesssim \|\cT(t)\wt y\|$ for all $t\geq 1$.  
Let $y := \cT(1)\wt y$. Then, for every $t>0$ we have
\[
r(t) = \wt r(t+1)  \lesssim \|\cT(t+1) \wt y\| \leq \|\cT(t) y\|.
\]
Therefore, it suffices to show $\cT(1)x \in D_\ell$ for every $\ell \in \N$. 
For \eqref{inc:S_in_core}, by the semigroup property and an inductive argument, the proof reduces to showing \eqref{inc:S_in_core} with $\ell =1$, i.e., 
\begin{equation}\label{inc:core}
\cT(t)X\subset A^{-1} B X = D_1, \quad t>0.
\end{equation}
To justify \eqref{inc:core}, we first note that for every $\lambda\in\rho(P_D)$ we have
\[
A(\lambda B+A)^{-1}B
=
B\bigl(I-\lambda(\lambda B+A)^{-1} B\bigr).
\]
Since the functions
\[
\Gamma\ni\lambda \mapsto e^{\lambda t}R_D(\la)B
\quad\text{and}\quad
\Gamma\ni\lambda \mapsto e^{\lambda t}\bigl(I-\lambda R_D(\la) B\bigr)
\]
are integrable on $\Gamma$, where $\Gamma$ is the contour appearing in the definition of the operator $\cT(t)$ (see~\eqref{formula:semig}), 
a standard argument based on the closedness of $A$ and Hille's theorem (see, e.g., \cite[Proposition~1.1.7]{ABHN11}) shows that for every $x\in X$ we have $\cT(t)x \in \Dom(A)$ and 
\[
A\cT(t)x = \frac{1}{2\pi i}\int_{\Gamma} e^{\la t} A R_D(\la) B x \,\ud \la = \frac{1}{2\pi i} B \int_{\Gamma} e^{\la t} \big(I - \la R_D(\la) B\big) x \, \ud \la \in BX.
\]
Thus, inclusion~\eqref{inc:core} holds. 
\end{proof}

\begin{remark}\label{rem:range_F_kappa} 
For $\theta\in(1/2,1]$ and $\beta\in(1-\theta,\theta)$, one has the inclusion  $A^{-1}\big(BX\big)\subset\Ran(J_\beta)$. 
Indeed, for any $\gamma >1 - \theta$ we set $J_\gamma: =  \frac{1}{\Gamma(\gamma)} \int_0^\infty t^{\gamma-1}e^{-t} \cT(t)\,\ud t$. Then, using Fubini's theorem and the change of variables $s = t+r$ and $t = \tau$, 
\[
\begin{aligned}
J_\beta J_\gamma
&= \frac1{\Gamma(\beta)\Gamma(\gamma)} \int_0^\infty\int_0^\infty
t^{\beta-1} r^{\gamma-1} e^{-(t+r)}
\cT (t+r)\,\ud t\,\ud r\\
&=
\frac1{\Gamma(\beta)\Gamma(\gamma)}
\int_0^\infty
e^{-s}\cT(s)
\int_0^s \tau^{\beta-1}(s-\tau)^{\gamma-1}\,\ud \tau\,\ud s\\
&= \frac1{\Gamma(\beta+\gamma)}
\int_0^\infty s^{\beta+\gamma-1} e^{-s }\cT(s)\,\ud s, 
\end{aligned}
\]
where we have used the identity
\[
\int_0^s
\tau^{\beta-1}(s-\tau)^{\gamma-1}\,\ud \tau
= s^{\beta+\gamma-1}\frac{\Gamma(\beta)\Gamma(\gamma)}
{\Gamma(\beta+\gamma)}.
\]
Thus, $J_\beta J_\gamma = J_{\beta+\gamma}$. Moreover, one can show that $J_1 = R_D(1) B = (B+A)^{-1}B$ and  $\Ran (J_1)  = A^{-1}(BX)$.
(Indeed, $y = (B+A)^{-1}B x$ implies $Ay = B(x-y)$, and $Ay=Bx$ implies $y=(B+A)^{-1}B(x+y)$.) Therefore, if $\theta>\frac12$, $1-\theta<\beta<\theta$, then $J_\beta$ and $J_{1-\beta}$ are bounded and $J_1 = J_\beta J_{1-\beta}$. 
Thus, $A^{-1}(BX)=\Ran J_1 \subset\Ran J_\beta$ as desired. 

Finally, if $B=I$ and $A$ is a sectorial operator on $X$, then $
J_\beta=(I+A)^{-\beta}$, and
\[
A^{-1}(BX)=\Dom(A)\subset\Dom(A^\beta)=\Ran(J_\beta).
\]
For every $0<\beta<1$, this inclusion is proper unless $A$ is bounded.
$\hfill\Diamond$\end{remark}

As announced in Section~\ref{subs:main_results}, the case $\theta=1$ leads to the following refinement of Theorem~\ref{thm:mild_DPE}; cf. the discussion preceding Theorem~\ref{thm:mild_DPE}. In particular, 
Theorem \ref{thm: mild to DE}(iii) identifies the generator $L$ of the solution propagator to (PE), that is, the $C_0$-semigroup $T$; see \eqref{eq:decomp_of_cT}. The spectral behaviour of the operator $L$ plays a crucial role in revealing the diffusion effect for certain solutions to (HE), cf. Remark~\ref{rem:a_statements} and Proposition \ref{thm:decying-orbits-of-T}  for further discussion. 
In particular, if $X$ is reflexive and $A$ is injective, then parts~(i) and (iv) of Theorem~\ref{thm: mild to DE} state that the condition
$0\in\rho(L)$ is equivalent to
$\Ran(B)\setminus\Ran(A) = \emptyset$, and is further equivalent to the
exponential decay of the $C_0$-semigroup $T$ (and hence also of $\cT = T\oplus 0$).

Let $Y = Y_{A,B}$ denote the subspace of $X$ defined by 
\begin{equation}\label{def:Y}
Y: = \overline{A^{-1}\big(BX\big)} = \overline{\{x\in \Dom(A): Ax \in \Ran( B)  \}}.
\end{equation}

\begin{theorem}\label{thm: mild to DE} Assuming Hypothesis~\ref{hyp:main_2} with $\theta =1$, the following assertions hold.
\begin{itemize}
\item [(i)] $Y\cap \Ker( B) = \{0\}$. If $X$ is reflexive, then  $X = Y \oplus \Ker (B)$, a direct sum.
\smallskip

\item [(ii)] For every $y_0 \in Y$ there exists a mild solution $v: = v(\cdot, y_0)$ to \emph{(PE)} with the initial data $y_0$. Furthermore, the operator-valued function $T: \R_+ \rightarrow \cL(Y)$ defined by 
\begin{equation}\label{def: of T}
 T(t)y_0 := v(t; y_0)  \quad \textrm{for } t\geq 0 \textrm{ and } y_0\in Y, 
\end{equation}
gives a sectorially bounded holomorphic $C_0$-semigroup on $Y$.
\smallskip

\item [(iii)] The generator $L$ of the $C_0$-semigroup $T$ is the composition
 \begin{equation}\label{def:L}
L = - (B|_Y)^{-1} \circ  A|_{\Dom(A)\cap Y}
\end{equation}
 defined on its maximal domain $\Dom(L) = \{x\in \Dom(A)\cap Y: Ax \in  BY \}$. 
\smallskip 
 
\item [(iv)]  $\Ker(L) = \Ker(A)$, and  $\Ran(L) \neq  Y$ if and only if $B Y \setminus A(\Dom(A)\cap Y) \neq \emptyset$.
\end{itemize}
\end{theorem}
\medskip

Condition \eqref{DE-parab_est} with $\theta =1$ is equivalent to the following one:
 $\{ \Re \la \geq 0\}\setminus \{0\} \subset \rho(P_D)$, $[R_D B]_1 < \infty$  and there exist constants $N,K\geq 0$ such that  
\begin{equation}\label{DE-parab_est_2}
 \|R_D(\la)  B\|\lesssim {(1+ |\la|)^N}{(\Re \la)^{-K}}, \quad \Re \la >0.
\end{equation}
This follows from a maximum principle for functions that may have singularities at the boundary; see Proposition~\ref{thm:MP-fudge-JEMS-correct}, which will also be used below in other contexts.
The estimate \eqref{DE-parab_est_2} is often easier to verify in concrete examples. In particular, for the model example considered in Section~\ref{sec:examples} (where $A$ is the negative Dirichlet Laplacian and $B$ is a multiplication operator), the most convenient way to verify Hypothesis~\ref{hyp:main_2} with $\theta =1$, under the assumptions
$i\R\setminus\{0\}\subset \rho(P_D)$ and $[R_DB]_1<\infty$, is to check that
\eqref{DE-parab_est_2} holds with $N=K=1$; see Proposition~\ref{prop:B_D_gives_DE}.

\begin{proof}[Proof of Theorem \ref{thm: mild to DE}]
(i) For any $y\in \Dom(A)$ such that $Ay = Bx$ for some $x \in X$, we have 
\[
is R_D(is) By = y - R_D(is) Ay = y -  R_D(is) B x \rightarrow y \quad \textrm{ as } |s| \rightarrow \infty. 
\]
Since $|s|\| R_D(is) B\|$, $s\neq 0$, are uniformly bounded, see \eqref{DE-parab_est}, for any $y \in Y$ we also have
\[
isR_D(is)B y \rightarrow y \quad \textrm{ as } |s|\rightarrow \infty.
\] 
This yields $Y \cap \Ker B = \{0\}$. 

For the second assertion in (i), let $X$ be reflexive and fix $x\in X$. Let $y$ be a weak limit point of the bounded sequence $x_n: = i n R_D(in) Bx$, i.e., $x_{n_j} \rightarrow y$ in the weak topology of $X$ for a subsequence $(n_j)\subset \N$. Since $Ax_n = in Bx - (in)^2 B R_D(in)Bx$, $(Ax_n)\subset \Ran(B)$, i.e. $(x_n)\subset Y$. Therefore, since $Y$ is weakly closed, $y\in Y$. To show that $x-y \in \Ker B$, it suffices to check that 
$R_D(\mu)B (x-y) = 0$ for some $\mu \in \rho(P_D)$. 
But for any $\mu, \nu \in \rho(P_D)$, $R_D(\mu)BR_D(\nu) = R_D(\nu)B R_D(\mu)$ and $R_D(\mu)Bx \in Y$, and so we have
\begin{align*}
R_D(\mu)B (x-y) & = R_D(\mu)B x - \operatorname*{\mathrm{w}\!-\!lim}_{j \to \infty} R_D(\mu)B \left( in_j R_D(i n_j) Bx \right)\\
&  = R_D(\mu)B x - \operatorname*{\mathrm{w}\!-\!lim}_{j \to \infty} 
 in_j R_D(i n_j) B R_D(\mu)Bx = 0.
\end{align*}

(ii) We fix $y_0 \in  Y$ and set $q(\la) := R_D(\la)By_0$, $\Re \la >0$. Combining Lemma \ref{lem on pencil est}(ii) and Proposition \ref{thm:MP-fudge-JEMS-correct}, the function $q$ satisfies the assumptions of \cite[Corollary 2.5.2]{ABHN11} (or even \cite[Theorem 2.6.1]{ABHN11}). Consequently, there exists a bounded continuous function $v: (0, \infty) \rightarrow X$ such that 
\[
q(\la) = (\cL v)(\la) \qquad \textrm{ for all } \Re \la >0.
\]
By Proposition \ref{mild-DE}, in order to show that $v$ is a mild solution to (PE) with the initial value $y_0$, it suffices to show that $\lim_{t\to 0^+} v(t) = y_0$.
This follows directly form the Tauberian Theorem; see e.g. \cite[Theorem 2.6.4]{ABHN11}.
Indeed, in the case when $y_0 \in \Dom(A)$ and $Ay_0 = Bx_0$ for some $x_0\in X$ we have 
\[
\la q(\la) = \la (\la B +A)^{-1} By_0 = y_0 - (\la B +A)^{-1}A  y_0 =  y_0 - (\la B +A)^{-1} B x_0 \rightarrow y_0
\]
as $\la \rightarrow \infty$ by using \eqref{hyp:bound_2} (with $\theta =1$).
For a general $y_0 \in Y$, the claim follows by a simple approximation argument: if $y_n \rightarrow y_0$ with $y_n\in \Dom(A)$ such that $A y_n \in \Ran(B)$ then
\[
\la q(\la) = \la(\la B +A)^{-1} (By_0 - By_{n}) + [\la (\la B+ A )^{-1}By_n - y_n ] + y_n  \rightarrow y_0
\]
as $\la \rightarrow \infty$. Therefore, $v(t) \rightarrow y_0$ as $t\rightarrow 0+$. 

To see that $T$ defined in \eqref{def: of T} is a $C_0$-semigroup on $Y$, we first note that its Laplace transform (in the strong operator topology) is a pseudo-resolvent, see \cite[Appendix B]{ABHN11}, that is, for all $\Re \la, \Re \mu >0$ we have
\begin{equation}\label{eq:pseudo_res}
(\la B + A)^{-1}B - (\mu B + A)^{-1}B = (\mu - \la) (\la B + A)^{-1} B (\mu B + A)^{-1}B.
\end{equation}
Since for any $x \in X$ and $\Re\la >0$  
\[
(\la B +A)^{-1} B x \in \Dom(A) \quad \textrm{ and } \quad A(\la B +A)^{-1} B x = Bx - \la B (\la B +A)^{-1} B x, 
\]
the subspace $Y$ is invariant whit respect to $(\la B  + A)^{-1}B$ for $\Re \la > 0$. Indeed, $AR_D(\la)Bx = B(x- \la R_D(\la)Bx)$, so $R_D(\la)B \in A^{-1}\big(BX\big) \subset Y$. Since $\Ker (\la B +A)^{-1}B$ is independent on $\la$ by \eqref{eq:pseudo_res}, and, as we have shown above, for all $y_0\in Y$ we have $\la q(\la) \rightarrow y_0$ as $\la \rightarrow \infty$, it follows that the family of the operators $(\la B +A)^{-1} B$, $\Re \la >0$, is in fact the resolvent of a closed linear operator $L$ on $Y$; see, e.g., \cite[Proposition B.6]{ABHN11}, that is, 
\begin{equation}\label{resolvent_L}
(\la - L)^{-1} = (\la B + A)^{-1}B\big|_{Y},\quad \Re \la >0.
\end{equation} 
By \cite[Theorem 3.1.7]{ABHN11}, $T$ is a $C_0$-semigroup on $Y$ with the generator $L$. The bounded holomorphic extension of $T$ to a sector $\{z\in \C\setminus: 0< |\arg z |<\omega' \}$ for some $\omega'>0$ follows now from the characterization of this class of semigroups; see e.g. \cite[Theorem 3.7.11]{ABHN11}.

(iii) By \eqref{resolvent_L} with $\la =1$, we have
\begin{equation}\label{FL}
L = I - \left((A+B)^{-1}B\big|_{Y}\right)^{-1}
\end{equation} and so $\Dom(L)$ is equal to $\Ran \big( (B+A)^{-1}B\big|_{Y}\big) \subset \Dom(A) \cap Y$. We take $x\in \Dom(L)$ and set $y: = \bigl((A+B)^{-1}B\big|_{Y}\bigr)^{-1} x$, so that, $y\in Y$ and $(A+B)^{-1}By = x$, or $By = Ax + Bx$.
Therefore, 
\[
B Lx = Bx - B \left((A+B)^{-1}B\big|_{Y}\right)^{-1}x = -Ax. 
\]
Since $Lx\in Y$, one has $x\in \Dom(A)\cap Y$ with $Ax = -BLx \in BY$. Thus, we have proved the inclusion $L \subset -(B\big|_{Y})^{-1} A_{|\Dom(A)\cap Y}$. 

For the converse inclusion, let $x\in \Dom(A)\cap Y$ be such that $Ax \in BY$, that is, $Ax = Bz$ for some $z\in Y$. Setting $y: = x+z$ we can write
$Ax + Bx = B y$. Since $y\in Y$, we have $x\in \Ran \big( (B+A)^{-1} B\big|_{Y}\big) = \Dom(L)$. Moreover, 
\[
- Ax = - Bz = Bx - By = B\left(x - [(A+B)^{-1}B\big|_{Y}]^{-1}x\right) = B Lx.
\]
This yields the desired inclusion and hence the equality of $L$ and $-(B\big|_{Y})^{-1} A$.

Since $\Ker (A) \subset Y$, statement (iv) is a direct consequence of  formula~\eqref{FL}. \end{proof}

\begin{proof}[Proof of Theorem \ref{thm: diffusion intro}] 
First, by Corollary~\ref{cor:DWE_well_posed}, $-A$ generates a cosine
function on $X$. Consequently, Theorem~\ref{thm on equiv cond} yields $
V=\Dom((A+\omega)^{1/2})\subset X_2$.

It remains to show assertions (i) through (iii) in Theorem \ref{thm: diffusion intro}. By Lemma \ref{lem: mild s to DWE} for any $x_0 \in V$ and $x_1\in X$ there exists a polynomially bounded mild solution $u$ to (HE) with the initial data $x_0$ and $x_1$. By Theorem \ref{thm:mild_DPE} (resp. Theorem~\ref{thm: mild to DE} when $\theta =1$) for any $y_0 \in A^{-1}\bigl(BX\bigr)$ or $y_0 \in \Ran(J_\beta)$ (resp. $\overline{A^{-1}(BX)}$) there exists a mild solution to (PE) with the initial value $y_0$ which is uniformly bounded on $\R_+$. Consequently, by Theorem \ref{thm:ind-solu},  we have $ u(t; x_0,x_1) - v(t; y_0) = \bigl(\cF^{-1}m(\cdot; x_0, y_0)\bigr)(t)$, $t>0$.

 Therefore, \eqref{decay of u-v} in part~(i) follows from Theorem \ref{thm 1}(i). 
The estimate~\eqref{decay of u-v t2} (resp. \eqref{est:thm4(ii)_u}) in part~(ii) follows from
Theorem~\ref{thm 1}(ii) (resp. Corollary \ref{cor:decay_of_mH}).  

Finally, the first estimate in part~(iii) follows from
Theorem~\ref{thm 1}(iii). To prove the second estimate we first recall that 
$v(\cdot; y_0) = \cT(\cdot)y_0$ has a bounded holomorphic extension on a sector $\Sigma_\omega$. By Cauchy's formula, for any $\omega'\in (0,\omega)$ and for every $t>0$ we have 
\[
v'(t;y_0) = \frac{1}{2\pi i} \int_{|\la - t| = t\sin \omega'} \frac{\cT(\la)y_0}{(\la - t)^{2}} \ud \la.
\]
Thus, $\|v'(t,y_0)\|\lesssim 1/t$, $t>0$. Together with the
estimate for $u'-v'$, this gives the desired estimate for $u'$.
\end{proof}

\begin{proof}[{Proof of Theorem \ref{thm: diffusion}}] 
By Theorem \ref{thm:mild_DPE}(iii) there exists $y_0\in D_2$  such that 
\begin{equation}\label{est_of_v}
\|\cT(t)y_0\| = \|v(t; y_0)\|\gtrsim r(t),\qquad t\geq 0.
\end{equation}
By Theorem~\ref{thm: diffusion intro}, we already know that
$\Dom((A+\omega)^{1/2})=V$ whenever $\omega\geq0$ is such that
$A+\omega$ is sectorial.  
We fix $x_0\in \Dom(A)$ and set $x_1 := By_0 - Bx_0$. By our assumption on $B$, since  $D_2 \subset \Dom(A)$ and  $B(\Dom(A)) \subset V$, we have $Bx_0, By_0 \in V$, and so also $x_1\in V$. Since $D_2\subset A^{-1}(BX)\cap Y_D^\theta$ for every $\theta \in (0,1]$, the vectors $x_0$, $x_1$ and $y_0$ satisfy the assumptions of Theorem \ref{thm: diffusion intro}(i). Consequently, $\|u(t;x_0;x_1) - v(t,y_0)\|\lesssim 1/t$, $t>0$.
Combining this with \eqref{est_of_v} and using the fact that $v$ is holomorphic, and hence $v(t;y_0)\neq 0$ for all $t\ge0$, we may write
\begin{align*}
\Big| \|u(t)\|/\|v(t)\| - 1 \Big| & \leq \|u(t) - v(t)\| /\|v(t)\|  \lesssim {1}/{t r(t)}, \quad t>0.
\end{align*} Combining this with the assumption on $r$, we obtain
\eqref{decay u/v} and the third estimate in \eqref{dp}.
\end{proof}

\subsection{Proof of Theorem \ref{thm:diffusion_theta=1}}\label{subsec:weak-orbits}
Theorem~\ref{thm: diffusion} corresponds to part~(i) of
Theorem~\ref{thm: diffusion intro}. It asserts that, for every
$\theta\in(0,1]$, under additional assumptions on the operator $B$, it is possible to choose initial data for which the norms of both
solutions $u$ and $v$ decay more slowly than the norm of their difference.

We do not know whether analogous results corresponding to parts~(ii)
and~(iii) of Theorem~\ref{thm: diffusion intro} hold for general
$\theta\in(0,1]$. We establish such counterparts only in the case
$\theta=1$, where the associated semigroup $\cT$ (see
Theorem~\ref{thm:mild_DPE}) admits a more detailed description,
which is carried out in Theorem~\ref{thm: mild to DE}. 
Before proceeding to the details, we emphasize the following points.

\begin{remark}\label{rem:a_statements}
{(a1)} Since $L$ is the generator of the sectorially bounded holomorphic $C_0$-semigroup $T$ on a Banach space $Y$,  we have $\sigma(L)\cap \{\Re \lambda \geq 0\} \subset {0}$ and for any $x\in \Ran(L^k)$, that is, such that $x = L^k y$ for some $y\in \Dom(L^k)$, we have the estimate
\[
\|T(t)x\| = \|L^k T(t) y\| \lesssim {\|y\|}/{t^k},\qquad t> 0.
\]
Therefore, regardless of any additional assumptions imposed on $X$ and $L$, the orbits $(T(t)x)_{t \ge 0}$ do not decay as $t \to \infty$ slower than $1/t$ (respectively, $1/t^2$) for vectors $x \in \Ran(L)$ (respectively, $x \in \Ran(L^2)$).
For a remainder of order $t^{-2}$, one must therefore find $y\in\Ran(L)$ whose orbit dominates $t^{-2}$, despite the upper bound $\|T(t)y\|=O(t^{-1})$.
\medskip

\noindent{(a2)} If $X$ is an arbitrary Banach space and $0 \in \rho(L)$, then the semigroup $T$ decays exponentially, that is, for some $\varepsilon >0$
\[
\|T(t)\|_{\mathcal{L}(Y)} \lesssim e^{-\varepsilon t},
\qquad t\ge 0.
\]
Therefore, for each $y_0 \in Y$ we have $v(t;y_0) = O(e^{- \varepsilon t})$ as $t\rightarrow \infty$. 
Consequently, in such case Theorem~\ref{thm: diffusion intro} provides just the decay rate for solutions $u = u(\cdot; x_0,x_1)$ to (HE). More precisely, for any initial data $x_0$ and $x_1$ satisfying the conditions of part~(i) (resp. part~(ii)), we have $\|u(t)\| = O(1/t)$ (resp. $\|u(t)\| = O(1/t^2)$) as $t\to \infty$. Moreover, for any $y_0 \in Y$, we have $v'(t;y_0) = LT(t)y_0 = T(t-1) LT(1)y_0$, $t\geq 1$, so that $\|v'(t;y_0)\| \lesssim e^{-\varepsilon t}$ for $t \ge 0$. Hence, if $x_0$ and $x_1$ satisfy the conditions of part~(iii) in Theorem~\ref{thm: diffusion intro}, we obtain the decay of the kinetic energy of (HE), i.e., $\|u'(t)\| = O(1/t^2)$ as $t\to \infty$.

\medskip

\noindent{(a3)} If $X$ is reflexive and $0 \in \sigma(L) \setminus \sigma_p(L)$, then $Y$ is also reflexive, and for any $y \in Y$ one has $\|T(t)y\| \to 0$ as $t \to \infty$, which is referred to as the \emph{strong stability} of $T$ (the Arendt--Batty--Lyubich--V\=u theorem; see, e.g., \cite[Theorem 5.5.5]{ABHN11}). In other words, each mild solution $v(\cdot; y_0) = T(\cdot)y_0$, with $y_0 \in Y$, to (PE) decays. Moreover, for $y_0 \in \Ran(L^k)$, we have $\|v(t; y_0)\| = O(t^{-k})$ as $t\to \infty$.
We recall that if a solution $u$ to (HE) admits the asymptotic profile $v = T(\cdot)y_0$ for some $y_0 \in Y$, then $u(t) \to 0$ as $t \to \infty$ with the same decay rate as $v(t)$; cf. the proof of Theorem~\ref{thm: diffusion}. 

In applications (see, e.g., Subsection~\ref{subs:applications}), the operator $A$ often has trivial kernel; hence, by Theorem~\ref{thm: mild to DE}(iv), so does $L$, i.e., we indeed have $0 \in \sigma(L)\setminus \sigma_p(L)$. 

\medskip

\noindent{(a4)} If $0\in \sigma_p(L)$ (and $X$ an arbitrary Banach space), then there exists $y\in Y$ such that $\|y\| = 1$ and $Ly = 0$. Then for any $n\in \N$ we have $y \in \Dom(L^n)$ and $T(t)y = y$ for all $t \ge 0$. Under the assumptions of Theorem \ref{thm: diffusion}, let $x_0\in \Dom(A)$ and  $0\neq y_0\in \Ker L$, and set $x_1 := By_0 - Bx_0$. Since $v(t;y_0)=y_0$ for all $t \ge 0$,  Theorem \ref{thm: diffusion intro}(i) yields 
\[
\|u(t; x_0, x_1) - y_0\| \lesssim 1/t, \qquad t>0. 
\]
In contrast to (a3), $u$ does not decay, i.e., $\|u(t)\| \to \|y_0\|$ as $t \to \infty$, and the constant function $v(t; y_0) = y_0$, $t \geq 0$, can be regarded as the asymptotic profile of $u$. 
$\hfill\Diamond$\end{remark}

For the proof of part~(i) of Theorem~\ref{thm:diffusion_theta=1}, we use
results from \cite{MuTo13}. For the convenience of the reader, we provide
below a reformulation of selected results from \cite{MuTo13}, adapted to our
particular setting of the generator $L$.

Roughly speaking, the following result asserts that if $0\in\sigma(L)$, then
every bounded function tending to zero at infinity can be dominated on a
\emph{large} subset $J$ of $\mathbb R_+:=[0,\infty)$
($J=\mathbb R_+$ when $X$ is a Hilbert space) by a certain
\emph{smooth weak orbit} $t\to \langle T(t)y,y^*\rangle$ of the corresponding semigroup $T$.
Regarding Theorem~\ref{thm: diffusion}, clearly, domination by weak orbits
$t\mapsto \langle T(t)y,y^*\rangle$ is a stronger property than domination by strong orbits $t\mapsto T(t)y$.

For a measurable set $J\subset \R_+$ we denote by ${\mathrm{Dens}}\,J$ and $\overline{\mathrm{Dens}}\,J$   the density of $J$ and the upper density of $J$, defined by
\[
\mathrm{Dens}\,J := \lim_{t \to \infty} {|J \cap [0, t]|}/{t}\quad \textrm{ and } \quad
\overline{\mathrm{Dens}}\,J := \limsup_{t \to \infty} {|J \cap [0, t]|}/{t}
\]
(assuming the existence of the first limit).

\begin{lemma}[{\cite[Theorem~6.4]{MuTo13}}]\label{MT-theorem} Let $L$ be the generator of a bounded $C_0$-semigroup $T$ on a Banach space $Y$. Assume that $0 \in \sigma(L)$. 
Then for each bounded function $r \colon \mathbb{R}_+ \to \mathbb{R}_+$ such that $\lim_{t\to\infty} r (t)=0$
there exist $y\in \bigcap_{\ell\in \N}\Dom(L^\ell)$, $y^\ast\in Y^\ast$, and a measurable set
$J\subset \mathbb{R}_+$ such that $\|y^\ast\|\le 1$,
$\overline{\mathrm{Dens}}\,J = 1$, and
\[
|\langle T(t)y, y^\ast\rangle| > r(t)
\qquad \text{for all } t\in J.
\]
Moreover, if $Y$ is a reflexive Banach space, one can take $J$ such that $\mathrm{Dens}\,J = 1$, and if $Y$ is a Hilbert space, one can take $J = \mathbb{R}_+$.
\end{lemma}

Theorem~\ref{thm 1}(ii) (a result leading to the $1/t^2$-decay of the difference between solutions to (HE) and (PE)), in particular, requires that $By_0 \in \Ran (A)$, so $y_0 \in \Ran(L)$, cf. Theorem \ref{thm: mild to DE}(iii). Consequently, the corresponding orbit $(T(t)y_0)_{t \ge 0}$ decays no slower than $1/t$; see Remark~\ref{rem:a_statements}(a1). 
The existence of vectors $y$ in $\Ran(L)$ for which the decay rate of $(T(t)y)_{t\geq 0}$ is arbitrarily close to $1/t$ does not appear to have been studied in the literature (in fact, $1/t^2 = o(\|T(t)y\|)$ as $t\to \infty$ would suffice to conclude the occurrence of the diffusion phenomenon). 
At present, it remains unclear whether direct analogues of Lemma~\ref{MT-theorem} hold for the orbits of $T$ corresponding to vectors in $\Ran(L)$. A systematic investigation of these counterparts is postponed to future work.

In the next two results, Proposition~\ref{thm:decying-orbits-of-T} and Corollary~\ref{cor:decay-orbits-L},
we apply a result of Deutsch and Hundal \cite{DeHu10} to show that,
under the stated spectral assumptions, for every positive
function $r:\mathbb R_+\to\mathbb R_+$ satisfying
\[
 r(t)=o(t^{-1}),\qquad t^{-2}=o(r(t)),
 \qquad t\to\infty,
\]
there exist $y\in\operatorname{Ran}(L)$ and a sequence
$t_j\to\infty$ such that
\[
 \|T(t_j)y\|\gtrsim r(t_j).
\]
Consequently, in Theorem~\ref{thm:diffusion_theta=1}(ii), we obtain a solution
$v$ to (PE) that provides an asymptotic diffusive profile
for the corresponding solution $u$ to (HE), but only
along a sequence of times.

First, to motivate the assumptions imposed in Proposition \ref{thm:decying-orbits-of-T}, we note that there may be no vectors $y \in \Ran(L)$ for which $T(t)y$ decays at a rate slower than exponential. Indeed, let $L$ be the generator of a sectorially bounded holomorphic semigroup $T$ on $Y$. Suppose that $0$ is an isolated point of $\sigma(L)$ and let $Y_0 \oplus Y_1$ be the spectral decomposition of $Y$ corresponding to the parts $\{0\}$ and $\sigma\setminus\{0\}$ of the spectrum of $L$; see, e.g., \cite[Theorem 6.17]{Ka95}. 
Let $L_0$ and $L_1$ denote the parts of $L$ in $Y_0$ and $Y_1$, so that $\sigma(L_0) = \{0\}$ and $\sigma(L_1) = \sigma(L)\setminus \{0\}$. It is easily seen that the restriction $T_j$ of $T$ to $Y_j$ is a semigroup on $Y_j$ generated by $L_j$, $j=0,1$.  Therefore, $T_1$ is exponentially decaying. If $L_0 = 0\in \cL(Y_0)$, then $T_0(t) = I_{Y_0}$, $t\geq 0$. Moreover, for any $y\in \Ran(L)$, since $y = L_0 x_0 + L_1 x_1$ for some $x_0 \in Y_0$ and $x_1 \in \Dom(L)\cap Y_1$, we have
\begin{align*}
\|T(t)y\|_Y & = \|T_0(t)L_0 x_0 +  T_1(t)L_1 x_1\|_Y\\ 
& \leq  \|L_1 T_1(t/2)\|_{Y_1} \| T(t/2)\|_{Y_1}  \|x_1\|_Y. 
\end{align*}
Hence, for any $y \in \Ran(L)$, the orbit $(T(t)y)_{t \ge 0}$ decays exponentially as $t \to \infty$. 
The following result shows that the situation described in this paragraph is the only possible one in which $0 \in \sigma(L)$ and no vectors $y \in \Ran(L)$ have orbits that decay slower than exponentially.

\begin{proposition}\label{thm:decying-orbits-of-T} Let $L$ be the generator of a sectorially bounded, holomorphic $C_0$-semigroup $T$ on a Banach space $Y$. Suppose that one of the following two conditions holds: 
\begin{itemize}
\item [(i)] $0$ is a limit point of $\sigma(L)$,
\item [(ii)] $0$ is an isolated point of $\sigma(L)$ and $L_0\neq 0 \in \cL(Y_0)$.
\end{itemize}
Then there exists a sequence $(t_j)_{j\in \N}\subset (0,\infty)$ with $t_j\rightarrow \infty$ such that for any sequence $(\gamma_j)_{j\in \N}$ of positive real numbers with $\gamma_j\rightarrow 0$ as $j\rightarrow \infty$, there exists a vector $y\in \Ran(L)\cap \bigcap_{\ell\in \N} \Dom(L^\ell)$ such that 
\[
\|T(t_j)y\|_Y \gtrsim {\gamma_j}/{t_j}\quad \textrm{ for infinitely many } j\in \N.
\]
\end{proposition}

\begin{proof} For (i), as in the proof of \cite[Theorem~6.9]{BaChTo16}, by an argument based on the spectral mapping theorem and the spectral radius formula, we infer that there exist $(\wt t_j)_{j\in\mathbb{N}} \subset (1,\infty)$ with $\wt t_j\rightarrow \infty$ and a constant $c>0$ such that  for all $j \in \mathbb N$ 
\begin{equation}\label{eq:decay_orbit_1}
\|\wt t_j T( \wt t_j) L(I-L)^{-1} \|_{\cL(Y)} \ge c.
\end{equation}
This implies the desired conclusion. Indeed, if there exists a vector $\tilde{y}\in Y$ such that 
\begin{equation}\label{eq:decay_orbit}
\limsup_{j\rightarrow \infty} \| \wt t_j  T(\wt t_j) L(I-L)^{-1}\tilde{y}\|_Y >0
\end{equation}
then, since $L(1-L)^{-1}T(1)\wt y \in \bigcap_{\ell\in \N} \Dom(L^\ell)$, we can take $y := y(\tilde{y}): = L(I-L)^{-1}T(1)\tilde{y}$ and $t_j := \wt t_j -1$ to get the desired conclusion (even for $(\gamma_j)_{j\in \N}$ being constant). If there is no $\wt y$ such that \eqref{eq:decay_orbit} holds, then the sequence of operators 
\[
(\wt t_j -1) T(\wt t_j) L(I-L)^{-1},\qquad j\in \N,
\] converges strongly to $0\in \cL(Y)$, but not in operator norm. Then \cite[Theorem 3.1]{DeHu10} yields a vector $\tilde{y}$ such that $y (\tilde{y})$ has the desired properties.

For (ii), the space $Y$ admits a spectral decomposition $Y_0\oplus Y_1$ corresponding to the separated parts $\{0\}$ and $\sigma(L)\setminus\{0\}$ of the spectrum of $L$. Let $L_j$ and $T_j$, $j=0,1$, be the respective parts of $L$ and $T$. If 
\[
\limsup_{t\rightarrow \infty}\| t T_0(t) L_0(I-L_0)^{-1}\|_{\cL(Y_0)} > 0, 
\] 
then the same argument as in (i) gives the desired vector $y$. The case 
\[
\lim_{t\rightarrow \infty} \|t L_0 T_0(t)(1-L_0)^{-1}\|_{\cL(Y_0)} = 0 
\]
is not possible since then 
\begin{align*}
\|t L_0 T_0(t)\|_{\cL(Y_0)} & \leq  \|t L_0 T_0(t)(1-L_0)^{-1}\|_{\cL(Y_0)} \| 1 - L_0\|_{\cL(Y_0)} \rightarrow 0 \textrm{ as } t\rightarrow\infty, \textrm{ and }
\end{align*}
by \cite[Theorem~2.1]{KM-SOT04}, we have $L_0 = 0$ in $\mathcal{L}(Y_0)$, which contradicts  assumption~(ii). 
\end{proof}

For every sectorial operator (here, $-L$), there exists a decomposition of the underlying Banach space that allows one to restrict the operator to a subspace on which its part is injective and has dense range; see, e.g., \cite[Proposition 15.2]{KuWe04}. This yields Corollary \ref{cor:decay-orbits-L}, which can be viewed as a variant of Proposition \ref{thm:decying-orbits-of-T}. We recall that, for a linear operator $L$ on a Banach space $Y$ and a subspace $\widetilde Y \subset Y$, the part $\widetilde L$ of $L$ in $\widetilde Y$ is defined as
the restriction of $L$ to the domain
\[
D(\widetilde L) := \{ y \in \widetilde Y \cap \Dom(L) : Ly \in \widetilde Y \}.
\]

\begin{corollary}\label{cor:decay-orbits-L}
Let $L \neq 0$ be the generator of a sectorially bounded holomorphic $C_0$-semigroup $T$ on a Banach space $Y$. Suppose that $0\in \sigma(\tilde{L})$, where  $\tilde{L}$ denotes the part of ${L}$ in $\overline{\Ran(L)}$.  
Then the conclusion of Proposition \ref{thm:decying-orbits-of-T} holds.

In particular, if $Y$ is reflexive and $L$ is injective with $0\in \sigma(L)$, then the conclusion of Proposition \ref{thm:decying-orbits-of-T} holds.
\end{corollary}
\begin{proof} The part $\tilde{L}$ of ${L}$ in $\widetilde Y : = \overline{ \Ran(L)} \neq \{0\}$  is an injective generator of sectorially bounded  holomorphic $C_0$-semigroup $\tilde{T}$ on $\tilde{Y}$. Moreover, for any $t\geq 0$ $\tilde{T}(t)$ is the restriction of $T(t)$ to $\tilde{Y}$, and $\tilde{L}$ has dense range; see, e.g., \cite[Proposition 15.2]{KuWe04}. Clearly, $\tilde{L} \neq 0\in \cL(\widetilde Y)$. 
If $0$ is a limit point of $\sigma(\widetilde{L})$, then the desired conclusion follows from Proposition~\ref{thm:decying-orbits-of-T}.
If $0$ is an isolated point in $\sigma(\wt L)$, then $\widetilde Y$ admits a spectral decomposition $\widetilde Y_0$ and $\widetilde Y_1$ corresponding to $\{0\}$ and $\sigma(\widetilde{L})\setminus \{0\}$. Let $\widetilde L_j$ and $\widetilde T_j$, $j=0,1$, be the corresponding decompositions of $\widetilde L$ and $\widetilde T$. The semigroup $\widetilde T_1$ decays exponentially in the norm of $\mathcal{L}(\widetilde Y_1)$, whereas $\widetilde T_0$ does not. For $\wt T_0$, either 
\[
\limsup_{t\rightarrow \infty}\| t \wt T_0(t) \wt L_0(1- \wt L_0)^{-1}\|_{\cL(\wt Y_0)} > 0,
\]
which leads to the desired conclusion as in the proof of Proposition~\ref{thm:decying-orbits-of-T}, 
or 
\begin{equation}\label{eq:decay_orbits_2}
\limsup_{t\rightarrow \infty}\|t \wt T_0(t) \wt L_0(1- \wt L_0)^{-1}\|_{\cL(\wt Y_0)} = 0,
\end{equation}
This gives that $\wt L_0 = 0$ in $\cL(\wt Y_0)$ as in the proof of Proposition \ref{thm:decying-orbits-of-T}. Since $\wt L$ is injective, we get $\wt Y_0 = \{0\}$ and, consequently, $\wt Y_1 = \wt Y$ and $\wt L_1 = \wt L$. Since $0\notin \sigma(\wt L_1)$ we see that \eqref{eq:decay_orbits_2} contradicts  our assumption.
For the last statement we recall that on a reflexive Banach spaces any injective sectorial operator (here, $-L$) has dense range; see, e.g., \cite[Proposition 15.2]{KuWe04}.
\end{proof}

We are now in a position to prove Theorem~\ref{thm:diffusion_theta=1}.

\begin{proof}[{Proof of Theorem \ref{thm:diffusion_theta=1}}] We begin with a general observation; cf. \eqref{eq:decomp_of_cT}. Since $X$ is reflexive, by combining  Theorem~\ref{thm:mild_DPE}(i) with $\theta=1$ and Theorem~\ref{thm: mild to DE}, we infer that the semigroup $\cT$ associated with $A$ and $B$ in Theorem~\ref{thm:mild_DPE} admits the decomposition $\cT(t)=T(t)\oplus0$, $t>0$, with respect to the decomposition $X=Y\oplus\Ker(B)$, where $Y:=\overline{A^{-1}(BX)}$. Moreover, $T$ is a $C_0$-semigroup on $Y$ with generator $L$ and $\Dom(L)\subset A^{-1}\big(BX\big)$; see Theorem~\ref{thm: mild to DE}.

For (i), since $Y$ is a Hilbert space, by relying on Lemma~\ref{MT-theorem} instead of \eqref{est_of_v}, the proof of Theorem~\ref{thm:diffusion_theta=1}(i) follows by essentially the same argument as that of Theorem \ref{thm: diffusion}, and we omit the details.

(ii) Since $X$ is reflexive, and hence so is $Y$, it follows that $Y$ is the direct sum of $\Ker(L)$ and $\overline{\Ran(L)}$ (see, e.g.,  \cite[Proposition 15.2]{KuWe04}). Since by our assumption $0\in \sigma(L)\setminus \sigma_p(L)$, we have $0\in \sigma(\wt L)$, where $\wt L$ stands for the part of $L$ in $\overline{\Ran(L)}$. Therefore, $L$ satisfies the assumptions of  Corollary \ref{cor:decay-orbits-L}. Let $(t_j)$ be a sequence from the conclusion of Proposition~\ref{thm:decying-orbits-of-T}. We put $\gamma_j: = r(t_j)$, $j\in \N$. Hence, there exists a vector $y_0\in \Dom(L) \cap \Ran(L)$ such that 
\begin{equation}\label{est:D_H_est}
\limsup_{j\rightarrow \infty} (t_j/\gamma_j) \|v(t_j; y_0)\| >0. 
\end{equation}
Since $y_0\in \Ran(L)$, we have that  $By_0\in \Ran (A)$ and $\|v(t; y_0)\|\lesssim 1/t$ for all $t > 0$ (see Remark \ref{rem:a_statements}(a1)).  
We fix a vector $x_0 \in \Dom(A)\cap \Ran (A)$ and set $x_1: = By_0 - Bx_0$. 
Since $Y_D^1 = Z_D^1=X$, applying  $B(\Dom(L)) \subset B(\Dom(A)) \subset V$, we infer that the vectors $x_0$, $x_1$ and $y_0$ satisfy the assumptions of Theorem \ref{thm: diffusion intro}(ii). Consequently,  
$\|u(t;x_0;x_1) - v(t,y_0)\|\lesssim 1/t^2$, $t>0$. Similarly, as in the proof of Theorem \ref{thm: diffusion}, relying on \eqref{est:D_H_est}, we obtain 
\[
\limsup_{j\rightarrow \infty}\, (t_j/\gamma_j) \|u(t_j; x_0, x_1)\| >0 
\quad \textrm{ and } \quad
\|u(t;x_0,x_1)\|\lesssim 1/t, \qquad t>0.
\] 
(iii) Let $r$ be given. By Corollary \ref{cor:decay-orbits-L} there exists a vector $y\in \Ran(L)\cap \Dom(L)$ such that 
\[
\limsup_{t\rightarrow \infty}\, t r(t)^{-1} \left\| T(t)y\right\|>0.
\]
Let $y_0 \in \Dom(L)$ be such that $y = Ly_0$. Since $(T(t)y_0)' = T(t)Ly_0 = T(t) y$, we infer that 
\begin{equation}\label{est:bound_v'}
\limsup_{t\rightarrow \infty}\, t r(t)^{-1} \left\| v'(t,y_0)\right\|>0,
\end{equation}
where $v(t,y_0) = T(t)y_0$, $t\geq 0$. 
Since $D(L) \subset A^{-1}\big( BX \big)\subset \Dom(A)$, $\wt Z_D^1 = X$, and $B(\Dom(L)) \subset B(\Dom(A)) \subset \Dom(A) \subset V = \Dom((A+\omega)^{1/2})$, a direct verification shows that if the vectors $x_0\in X$ and $y_0$ satisfy the assumptions of part~(iii), and $x_1:=By_0-Bx_0$, then $x_0$, $x_1$, and $y_0$ satisfy the assumptions of part~(iii) of Theorem~\ref{thm: diffusion intro}.
Consequently,  
$\|u'(t;x_0;x_1) - v'(t,y_0)\|\lesssim 1/t^2$, $t>0$. Moreover, 
$\|v'(t,y_0)\| = \|LT(t) y_0\| \lesssim 1/t$, $t>0$; see Remark \ref{rem:a_statements}(a1). Therefore, arguing as in the proof of Theorem~\ref{thm: diffusion} and relying on~\eqref{est:bound_v'}, we obtain the remaining estimates for $u$ in~(iii).
\end{proof}

\section{Proofs of main results on (DWE) with Dirichlet Laplacian}\label{sec:proofs_of_thm1.6}

\subsection{Proofs of Proposition \ref{prop:B_D_gives_DE},  Theorem \ref{thm: pencil est} and Proposition \ref{cond: 1/t at infty}}

\begin{proof}[Proof of Proposition \ref{prop:B_D_gives_DE}] Recall that, by Lemma~\ref{lem on pencil est}, condition~\eqref{hyp:bound_2} implies that $G_\theta\subset\rho(P_D)$ and that $\min(|\lambda|,|\lambda|^\theta)\|R_D(\lambda)M_b\|$ is uniformly bounded for $\lambda\in G_\theta$; see~\eqref{pencil ext small}.  We first show that 
\begin{equation}\label{cond:Re_est_PD}
\{\Re \la >0 \} \subset \rho(P_D) \quad \textrm{ and } \quad \|R_D(\la)M_b\| \lesssim \max\left( 1, 1/\Re \la\right), \quad \Re \la> 0, 
\end{equation}
To this end, we prove that $1\in \rho(P_D)$, that is, that $M_b - \Delta$ is invertible on $L^2(\Omega)$. Suppose, to the contrary, that $0\in \sigma(M_b - \Delta)$.
Since $M_b - \Delta$ is selfadjoint and non-negative, it follows that 
$0\in \sigma_{\mathrm{ap}}(M_b - \Delta)$,
that is, there exists a sequence $(\varphi_j)\subset \Dom(M_b - \Delta)$ such that
\begin{equation}\label{eq:M_b_to_0}
\|\varphi_j\|_2 = 1 \quad \textrm{and}\quad \|(M_b - \Delta)\varphi_j\|_2 \to 0 \quad \textrm{ as } j\to \infty. 
\end{equation}
Since 
\[
\ip{(M_b - \Delta)\varphi_j}{\varphi_j} = \|\nabla \varphi_j\|_2^2 + \|\sqrt b \varphi_j\|^2_2 \geq \|\sqrt b \varphi_j\|^2_2 \geq \|b\|_\infty^{-1} \|b \varphi_j\|^2_2,
\]
by \eqref{eq:M_b_to_0} we have $\|M_b \varphi_j\|_2\to 0$ as $j\to \infty$.
Thus, for any $s\neq 0$,
\[
\|(is M_b - \Delta)\varphi_j\|_2
\leq
\|(M_b - \Delta)\varphi_j\|_2
+
|is - 1|\,\|M_b \varphi_j\|_2
\to 0
\quad \text{ as } j\to \infty,
\]
which contradicts the invertibility of $is M_b - \Delta$ that holds by  \eqref{hyp:bound_2}.

Thus, in particular, there exists $c>0$ such that for every $\varphi\in \Dom(\Delta)$ we have
\begin{equation}\label{eq:closed_range}
\ip{(M_b - \Delta)\varphi}{\varphi} \geq  c \|\varphi\|^2_2.
\end{equation}
Now let $\la \in \C$ with $\Re \la >0$. Since $M_b$ and $-\Delta$ are non-negative,  for all $\varphi\in \Dom(\la M_b - \Delta)$,
\begin{align*}
\|(\la M_b-\Delta)\varphi\|_2\|\varphi\|_2 & \geq | \ip{(\la M_b-\Delta)\varphi}{\varphi} | \geq \min(1, \Re \la) \ip{(M_b -\Delta)\varphi}{\varphi},
\end{align*}
Thus, by \eqref{eq:closed_range}, 
\begin{equation}\label{eq:from_below}
\|(\la M_b - \Delta) \varphi\|_2 \geq  c \min(1, \Re \la) \|\varphi\|_2
\end{equation}
for all $\varphi\in \Dom(\la M_b - \Delta)$. Since $(\la M_b - \Delta)^* = \bar{\la} M_b - \Delta$, the same estimate applies to the adjoint of $\la M_b - \Delta$, showing that $\la M_b - \Delta$ is invertible. Moreover, from \eqref{eq:from_below} we infer 
\[
\|R_D(\la) M_b\| \leq \|b\|_\infty \|R_D(\la)\| \leq c^{-1} \|b\|_\infty \max(1, 1/\Re \la),
\]
and \eqref{cond:Re_est_PD} follows.

Now, to establish \eqref{cond:almost_bound_2}, put $\varphi_\theta(t)=\min\{t,t^\theta\}$ and
\[
 F(\lambda):=\lambda(1+\lambda)^{\theta-1} R_D(\lambda)B,
 \qquad \operatorname{Re}\lambda>0,
\]
where the power is defined by the holomorphic branch on the right
half-plane. By \eqref{hyp:bound_2},
\[
\sup_{s\ne0}\|F(is)\| \le [R_D B]_\theta.
\]
Moreover, \eqref{cond:Re_est_PD} gives, with $r=\operatorname{Re}\lambda>0$,
\[
\|F(\lambda)\| \lesssim |\lambda|\max\{1,r^{-1}\} \lesssim (1+|\lambda|)^2r^{-1}.
\]
Since $F$ is holomorphic on the right half-plane and continuous at
its non-zero imaginary boundary points, Proposition~\ref{thm:MP-fudge-JEMS-correct} applies
(with $N=2$ and $K=1$ in \eqref{eq:MP-growth-JEMS-correct}), yielding
\[
\|F(\lambda)\|\le [R_D B]_\theta,
\qquad \operatorname{Re}\lambda>0.
\]
For $\operatorname{Re}\lambda>0$ we have
\[
\varphi_\theta(|\lambda|)
\le 2^{1-\theta}|\lambda|\,|1+\lambda|^{\theta-1},
\]
so $\varphi_\theta(|\lambda|)\|R_D(\lambda)B\|$ is bounded
throughout the right half-plane. Finally,
$H_\theta\cap\{\operatorname{Re}\lambda\le0\}\subset G_\theta$,
and Lemma~\ref{lem on pencil est}(ii) supplies the same bound there. This proves \eqref{cond:almost_bound_2}.
\end{proof}

\smallskip

We proceed with some auxiliary results needed for the proof of
Theorem~\ref{thm: pencil est}. First, in an abstract setting, we examine how the conditions $[R_D B]_\theta < \infty$ and $[R] < \infty$ in \eqref{pencils est} and \eqref{hyp:bound_2} are related. The following result shows that, once an a priori estimate for either $R_D(\lambda)$ or $R(\lambda)$ for $\la$ near $0$ is available, the two pencils exhibit exactly the same asymptotic behaviour as $\lambda\to0$; see \eqref{equ:equiv_on_U}.

\begin{lemma}\label{lem near 0} Let $A$ be a linear operator on a Banach space $X$ and $B\in \cL(X)$. Suppose that for some $\Lambda \subset \C$ we have $\Lambda \subset \rho(P_D)$ (resp. $\Lambda \subset \rho(P)$) and 
\[
\sup_{\la\in \Lambda} \| \la^2 R_D(\la)\| <1 \quad \textrm{(resp.}\quad \sup_{\la\in \Lambda} \| \la^2 R(\la)\| <1). 
\]
Then, $\Lambda \subset \rho(P)$ (resp. $\Lambda \subset \rho(P_D)$) and there exists a constant $c>0$ such that
\begin{equation}\label{equ:equiv_on_U}
{c^{-1}}\|R_D(\la) \|\leq \|R(\la)\| \leq c \|R_D (\la)\| \qquad \la \in \Lambda.
\end{equation}
\end{lemma}

\begin{proof}
For $\la \in \Lambda$, write $P(\la) = P_D(\la) + \la^2I$ and use $\|\la^2 R_D(\la)\| <1$ to infer that $\Lambda \subset \rho(P)$, $
R(\la) = R_D(\la) \sum_{n=0}^\infty (-1)^n \la^{2n} R_D(\la)^{n}$, 
and 
\[
\|R(\la)\|\leq  {\|R_D(\la)\|}/{(1-\|\la^2R_D(\la)\|)},
\]
which yields the second inequality in \eqref{equ:equiv_on_U}. The proof of the first inequality in \eqref{equ:equiv_on_U} is similar. 
\end{proof}



The inequality underlying the following lemma also arises in the form-theoretic treatment of unbounded damping and potentials in \cite[Lemma~2.3]{FrSiTr18}.

\begin{lemma}\label{lem: est via form}
Let $\Omega$ be an open subset of $\R^n$. Let $B = M_b$ be the multiplication operator by a nonzero function $b\in L_+^\infty(\Omega)$. Let $A = -\Delta$ be the negative Dirichlet Laplacian on $L^2(\Omega)$.  

Then for any $\omega \in (0, \pi)$ there exists a constant $c_\omega>0$ such that for any $\la\in \C\setminus \{0\}$ with $|\arg \la|\leq \omega$ and for every $\varphi\in \Dom(A)$
\begin{equation}\label{est: crucial}
|\la| \|\sqrt{b}\varphi\|_2^2 + \|\nabla \varphi \|_2^2  \leq c_\omega^2  \|\la B \varphi + A \varphi\|_2   \|\varphi\|_2.
\end{equation}
Furthermore, if $\la \in \rho(P_D)\setminus(-\infty,0]$ then
\begin{equation}\label{RdB est}
|\la| \|R_D(\la) B\|^2\leq \left(\cos(\arg \la/2)\right)^{-1} \|b\|_\infty \|R_D ( \bar{\la})\|.\\[1ex]
\end{equation}
\end{lemma}

\begin{proof} We begin with the following auxiliary estimate: for any $s, r\geq 0$ and $0\neq \la\in \C$ 
\begin{equation}\label{eq: fundamental est}
\cos^2(\arg \la/2) (|\la|r + s)^2 \leq |\la|^2 r^2 + s^2 + 2rs \Re \la.
\end{equation}

This trivially holds if $rs = 0$. Otherwise, setting $t = s/(|\la| r)$, we observe 
\[
\frac{\left|  \la r + s\right|^2 }{\left( |\la| r + s \right)^2} = \frac{1+t^2 + 2t \cos (\arg \la) }{(1+t)^2}.
\] 
The right-hand side is minimized at $t=1$, giving $\cos^2(\arg \la/ 2)$, thereby proving~\eqref{eq: fundamental est}.
We can now proceed as follows. Fix $\la \in \C\setminus (-\infty, 0]$  and   $\varphi\in \Dom(A)\subset W^{1,2}_0(\Omega)$. Set $f: = \la B \varphi + A \varphi = P_D(\la)\varphi$. 
We observe that   
\begin{equation}\label{coer 1}
|\langle f, \varphi \rangle|^2 =  |\la|^2 \|\sqrt{b}\, \varphi\|_2^4 + \| \nabla \varphi\|_2^4 + 2 \Re\la \|\sqrt{b} \varphi\|_2^2 \|\nabla \varphi\|_2^2.
\end{equation}
Taking $r: = \|\sqrt{b} \varphi\|_2^2$ and $s: = \|\nabla \varphi\|_2^2$ in \eqref{eq: fundamental est}, we obtain 
\begin{equation}\label{coer_2}
 |\la| \|\sqrt{b} \varphi\|_2^2 +  \|\nabla \varphi\|_2^2 \leq \cos(\arg \la/2)^{-1} |\langle f, \varphi \rangle| \leq  \cos(\arg \la/2)^{-1} \|f\|_2 \|\varphi\|_2.
\end{equation}
If, in addition, $\la \in \rho(P_D)$, then combining this with the estimate
\[
\|b\|_\infty^{-1} \|b \varphi \|_2^2 \leq  \|\sqrt{b} \varphi\|_2^2,
\]
we conclude that 
\[
\|b\|_\infty^{-1} |\la| \|B R_D(\la) f\|_2^2  + \|\nabla R_D(\la) f \|_2^2 \leq   \cos(\arg \la /2)^{-1} \| R_D(\la) f\|_2 \|f\|_2.
\]
Since $\|R_D(\la)B\| = \|(BR_D(\overline{\la}))^*\| = \|BR_D(\overline{\la})\|$, applying the preceding estimate to $\overline{\la}$ yields the desired
estimates \eqref{est: crucial} and \eqref{RdB est}.
\end{proof}

\begin{proof}[Proof of Theorem \ref{thm: pencil est}]
Fix $\omega \in (0,\pi)$ and $0 \neq \lambda \in \mathbb{C}$ with 
$|\arg \lambda| \le \omega$. 
Let $\varphi \in \Dom(A)$ and set 
$f := \lambda B \varphi + A \varphi = P_D(\lambda) \varphi$. 
By \eqref{est: crucial}, there exists a constant $c_\omega$ such that 
for any measurable set $D \subset \Omega$ satisfying 
$b(\xi) \ge \eta$ for a.e.\ $\xi \in D$, we have
\[
\|\varphi\|_{L^2(D)} \leq \frac{c_\omega}{\sqrt{\eta |\la|}} \|f\|_2^{1/2} \|\varphi\|_2^{1/2}.
\]

Below, we apply the estimate to the set $D$ given by
\[
D_\Omega = \Omega\setminus \{\xi: \chi(\xi) = 1\} \qquad \textrm{ and } \qquad D_U = U\setminus \{\xi: \chi(\xi) = 1\}.
\]

Also by \eqref{est: crucial} we infer that
\begin{equation}\label{est: imaginary part}
\|\nabla \varphi\|_2  \leq c_\omega \|f\|_2^{1/2} \|\varphi\|_2^{1/2}.
\end{equation}
Since $\mathrm{supp} (\chi) \subset U$ and $\chi(\xi) \equiv 1$ on $U'$,
\begin{align*}
\|\varphi\|_2 & \leq \|\chi \cdot \varphi\|_2 + \|(1-\chi)\cdot \varphi\|_2\\
& \leq \|\chi \cdot \varphi\|_{L^2(U)} + \| \varphi \|_{L^2(D_\Omega)}.
\end{align*} 
Further, since $\chi \cdot \varphi \in W^{1}_0 (U)$, applying the Poincar\'e property of $U$, we obtain 
\begin{align*}
\|\chi \varphi\|_{L^2(U)} & \leq c_P \|\nabla (\chi \varphi )\|_{L^2(U)}\\
& \leq c_P \|\varphi \nabla \chi  \|_{L^2(U)} + c_P \|\chi \nabla \varphi \|_{L^2(U)}\\
&\leq c_P \| \nabla \chi\|_{L^\infty(\Omega)} \|\varphi \|_{L^2(D_U)} + c_P \|\nabla \varphi \|_{L^2(U)}.
\end{align*}

Combining the above estimates yields 
\begin{equation}\label{est: main}
\|\varphi\|_2 \leq \frac{c_\omega}{\sqrt{\eta |\la|}} (c_P \| \nabla \chi\|_{L^\infty(\Omega)} + 1) \|f\|_2^{1/2} \|\varphi\|_2^{1/2} + c_P c_\omega \|f\|_2^{1/2} \|\varphi\|_2^{1/2}.
\end{equation}
We first consider the case $0 < |\lambda| < 1$, when 
\eqref{est: main} yields 
\begin{equation}\label{est:main_2}
\|\varphi\|_2 \leq \wt c_\omega |\la|^{-1} \|P_D(\la) \varphi \|_2 
\end{equation}
with $ \wt c_{\omega} =  \left( c_\omega \eta^{-1/2}(c_P \|\nabla \chi\|_{L^\infty(\Omega)} + 1) + c_P   c_\omega \right)^2 $. 
 Since $\varphi \in \Dom(A) = \Dom(P_D(\lambda))$ and $\lambda$ with 
$|\arg \lambda| \le \omega$ was arbitrary, we conclude that 
$P_D(\lambda)$ and $P_D(\overline{\lambda})$ are injective. 
Moreover, since 
$P_D(\lambda) = P_D(\overline{\lambda})^*$, 
the range $\Ran(P_D(\lambda))$ is dense in 
$L^2(\Omega)$. From \eqref{est:main_2}, it follows that $\lambda \in \rho(P_D)$ and that the first bound in \eqref{est: 1/t at 0} for $R_D(\lambda)=P_D(\lambda)^{-1}$ holds.
Since $\omega \in (0,\pi)$ was chosen arbitrarily, we conclude that
$
\{ \lambda \in \mathbb{C} : |\lambda| < 1 \}\setminus (-\infty, 0] 
\subset \rho(P_D).
$
Combining the first bound in \eqref{est: 1/t at 0} with Lemma \ref{lem near 0}  for  
\[
\Lambda : = \{ \la \in \C\setminus \{0\} : |\la|\leq \min(1/2, 1/(2\wt c_\omega) ) \textrm{ and } |\arg \la| \leq \omega\},
\]
we obtain the second bound in \eqref{est: 1/t at 0}.  

We begin the proof of \eqref{est: 1/sqrt at infty} by using \eqref{est: main} with $|\lambda|\geq 1$ to obtain
\[
\|\varphi\|_2 \leq \widetilde{c}_\omega \|f\|_2,
\]
where $\widetilde{c}_\omega$ is the constant defined above. The same argument as the one just applied for $0<|\la|<1$ yields  
$\{\la \in \C : |\arg \la|\leq \omega,\,  |\la|\geq 1 \} \subset \rho(P_D)$ and the first bound in \eqref{est: 1/sqrt at infty}. 
The second bound in \eqref{est: 1/sqrt at infty} follows now from the first bound in \eqref{est: 1/sqrt at infty} and \eqref{RdB est} in Lemma \ref{lem: est via form}. 
\end{proof}

\begin{proof}[{Proof of Proposition \ref{cond: 1/t at infty}}] To shorten the notation, let $\eta: = 2-1/\theta$. 
Fix $\omega\in (0, \pi)$. It is readily seen that for any $s,r \geq 0$  and $0\neq \la \in\C$ with $|\arg \la|\leq \omega$ we have that
\begin{equation}\label{est: main3}
|\la| s \leq \mu_\omega |\la s + r|, 
\end{equation}
where $\mu_\omega:=1$ for $\omega\in(0,\pi/2]$ and
$\mu_\omega:=1/\sin(\omega)$ for $\omega\in(\pi/2,\pi)$. 
Moreover, for any $\varphi \in \Dom(A)$ we have 
\begin{align*}
 \left|\la \|b \varphi \|_2^2 + \|\sqrt{b} \nabla \varphi\|_2^2\right|  & =  |\langle\la b \varphi -\Delta \varphi),  b \varphi\rangle  + \langle\Delta \varphi , b \varphi\rangle + \langle b \nabla \varphi, \nabla \varphi\rangle|\\
& = |\langle\la b \varphi -\Delta \varphi),  b \varphi\rangle  - \langle \nabla \varphi , b \nabla \varphi\rangle - \langle \nabla \varphi , \varphi \nabla b\rangle + \langle b \nabla \varphi, \nabla \varphi\rangle | \\
&\leq |\langle\la b\varphi -\Delta \varphi),  b \varphi\rangle| +  | \langle \nabla \varphi , \varphi \nabla b\rangle |\\
& \leq \|(\la B +A) \varphi\|_2 \|b\varphi\|_2 + \|\nabla \varphi\|_2 \|\varphi \nabla b\|_2
\end{align*}
and, using H\"older's inequality, 
\[
\|b \varphi \|_2  \leq \|b\|_\infty^{1-\eta} \|\varphi\|_2^{1-\eta} \|b \varphi\|_2^\eta.
\]
Assumption \eqref{cond: b'<b} yields
\[
\|\varphi \nabla b\|_2  \leq c \| \varphi b^\eta\|_2  \leq c  \|\varphi\|_2^{1-\eta} \|b \varphi\|_2^\eta.
\]
By \eqref{asum:RD_bound}, $\{\la \neq 0: |\arg \la | \leq \omega\} \subset \rho(P_D)$. 
Let $0\neq \la \in \C$ with $|\arg \la| \leq \omega$. Taking $\varphi\in \Dom(A)$ and $f := (\la B + A) \varphi$, i.e.,  $\varphi = R_D(\la) f$,  and applying \eqref{est: crucial}, we infer that
\[
\|\nabla \varphi\|_2 \leq c_\omega \|f\|_2^{1/2} \|\varphi\|_2^{1/2} \leq c_\omega \|R_D(\la)\|^{1/2} \|f\|_2.
\]
Consequently, combining \eqref{est: main3} with $s = \|b\varphi\|^2$ and $r = \|\sqrt{b}\nabla \varphi\|_2^2$, together with the estimates above, we proceed as follows:
\begin{align*}
|\lambda|\|b \varphi\|_2^2
&\leq \mu_\omega
\left(\|b\|_{\infty}^{1-\eta} + c c_\omega \|R_D(\lambda)\|^{1/2} \right)
\|f\| \, \|\varphi\|_2^{1-\eta} \, \|b\varphi\|_2^\eta \\
&\leq \mu_\omega
\left(\|b\|_{\infty}^{1-\eta} + c c_\omega \|R_D(\lambda)\|^{1/2} \right)
\|R_D(\lambda)\|^{1-\eta} \|f\|_2^{2-\eta} \|b \varphi\|_2^\eta,
\end{align*}
which yields
\begin{equation}\label{est:R_DB}
\|B R_D(\lambda)\| = \|R_D(\bar{\lambda}) B\|
\leq \nu_{\lambda}\, |\lambda|^{\frac{1}{\eta-2}},
\end{equation}
where
\[
 \nu_\lambda:=\left[
 \mu_\omega\bigl(\|b\|_\infty^{1-\eta}
    +cc_\omega\|R_D(\lambda)\|^{1/2}\bigr)
 \|R_D(\lambda)\|^{1-\eta}
 \right]^{1/(2-\eta)}.
\]
By \eqref{asum:RD_bound}, $\nu_\lambda$ is uniformly bounded whenever
$|\arg\lambda|\leq\omega$ and $|\lambda|\geq1$. Since $1/(2-\eta)=\theta$, we obtain the desired estimate \eqref{est: RD theta}.

By the first bound in \eqref{est: 1/sqrt at infty}, the last statement follows directly from the first part.
\end{proof}

\subsection{Proof of Propositions \ref{ex: of b}, \ref{prop:a_la_ui},
 \ref{prop:interm_freq}, and \ref{cor:fixed-scale}}\label{proof of Prop 'of b'}
A characterization of open and, necessarily,  quasi-cylindrical sets of $\R^n$ for which the Poincaré property holds can be found in \cite{EdEv18}.  
For our purposes, we rely only on the following geometric condition on the unbounded set $U$ which guarantees that $U$ has the Poincar\'e property.

\begin{lemma}[{\cite[Theorem 6.33]{AdFo03}}]\label{thm: P ineq} Let $U\subset \R^n$ be an unbounded open set. Suppose there exist constants $\kappa, N, h$ such that  $0 < \kappa $, $0 \leq N < \infty$, $0 < h < \infty$, and for every cube $Q \subset \mathbb{R}^n$ having the edge length $h$ and nonempty intersection with  $\{x \in U : \|x\| \ge N \}$ we have
\[
{\mu_{n-1}(Q,U)} \ge \kappa {h^{n-1}} ,
\]
where $\mu_{n-1}(Q,U)$ is the maximum of the $(n-1)$-measure of $\Pi(Q \setminus U)$  taken over all projections $\Pi$ onto $(n-1)$-dimensional faces of $Q$. Then $U$ has the Poincar\'e property. 
\end{lemma}

\begin{proof}[{Proof of Proposition \ref{ex: of b}}]
 (i) For any $p\in (r\ZZ)^n$ we write $Q(p) : = p + [0,r]^n$. By the regularity of Lebesgue's measure, for each lattice point $p \in (r\ZZ)^n$ such that $Q(p) \cap \Omega \cap \{\|\xi\|\geq r\ell\} \neq \emptyset$  
there exists an open set $U(p) \subset \Omega \cap \bigcup_{e\in \{-r,0,r\}^n} Q(p+e)$ such that 
\[
\left\{ \xi \in \Omega \cap Q(p):  b(\xi) < \varepsilon \right\} \subset U(p)
\]
and 
\begin{equation}\label{est:U(p)} 
\left|U(p)\setminus  \left\{\xi\in \Omega\cap Q(p):  b(\xi)< \varepsilon \right\}\right|\leq \varepsilon r^n /3^n.
\end{equation}
For every remaining lattice point $p\in(r\mathbb Z)^n$,
we put $U(p)=\varnothing$. Moreover, let $U_0 : = \{\|\xi\| < r\ell \} \cap \Omega$.
 
We claim that the set  $U: = U_0 \cup \bigcup_p U(p)$ satisfies the conditions stated in (i). 
 Indeed, let $\xi\in\Omega$ satisfy $b(\xi)<\varepsilon$.
If $|\xi|<r\ell$, then $\xi\in U_0$.
Otherwise, choose $p\in(r\mathbb Z)^n$ such that $\xi\in Q(p)$.
The condition defining \(U(p)\) is then satisfied, so \(\xi\in U(p)\subset U\).
If $U$ is bounded, it has the Poincar\'e property. We may therefore
assume that $U$ is unbounded.
Clearly, it suffices to show that the assumptions of Lemma \ref{thm: P ineq} holds for $U$ with $\kappa : = 4^{-n}\min(4,\varepsilon)$, $N :=r\ell+4r\sqrt n$, and $h := 4r$. 

Let $Q$ be any cube with edges parallel to the coordinate axes and length equal to $4r$ and non-empty intersection with $U\cap \{\|\xi\|\geq N\}$. Clearly, $Q\cap U_0 = \emptyset$.
(Indeed, a cube of edge length $4r$ has diameter $4r\sqrt n$.
If $Q$ meets $\{\xi\in U:\|\xi\|\ge N\}$,
then it cannot meet $U_0$, since otherwise two points of $Q$
would be more than $4r\sqrt n$ apart.)
Moreover, there exists a lattice point $p \in (r\ZZ)^n$ such that $U(p)\cap Q \neq \emptyset$. 
Therefore, either
\begin{itemize}
\item [(1)] there exists a vector $e\in \{-2r,-r,0,r, 2r\}^n$ such that $Q(p+e) \subset Q$ and $ Q(p+e)\cap U = \emptyset$, or
\item [(2)] for each $e\in \{-2r,-r,0,r, 2r\}^n$, if  $Q(p+e) \subset Q$ then  $Q(p+e) \cap U \neq~\emptyset$.
\end{itemize}  
In the first case, since $Q\setminus U \supset Q(p+e)$, we have
\[
\frac{\mu_{n-1}(Q, U)}{(4r)^{n-1}} \geq  \frac{r^{n-1}}{(4r)^{n-1}} = \frac{1}{4^{n-1}}.
\]
In the second case,  clearly
\begin{equation}\label{est:Q/U}
\frac{\mu_{n-1}(Q, U)}{(4r)^{n-1}} \geq  \frac{\mu_{n-1}(Q(p+e), U)}{(4r)^{n-1}}.
\end{equation}
Since $Q(p+e)\subset Q$ intersects $U$ and $Q\cap U_0=\varnothing$,
it also intersects $\Omega\cap\{\|\xi\|\ge r\ell\}$.
Hence the assumption of the proposition applies to $Q(p+e)$.

We first estimate $|Q(p+e)\setminus U|$. 
Since 
\[
Q(p+e) \cap U = \bigcup_{\wt e  \in \{-r,0,r\}^n} Q(p+ e)\cap U(p+ e+ \wt e).
\] 
and, using \eqref{est:U(p)} in the case when $U(p+e + \wt e)\neq \emptyset$, for any $\wt e \in  \{-r,0,r\}^n\setminus \{0\}$
\[
|U(p+e + \wt e) \cap Q(p+e)| \leq |U(p+e + \wt e) \setminus Q(p+e + \wt e)| \leq \varepsilon r^n /3^n,
\]
\begin{align*}
|U(p+e) \cap Q(p+e)| & \leq |U(p+e)\setminus (Q(p+e)\cap\{\xi \in \Omega: b(\xi)<\varepsilon \}| + | \{\xi \in \Omega: b(\xi)<\varepsilon \}|\\
& \leq \varepsilon r^n/3^n + (1 - 2\varepsilon)r^n ,
\end{align*}
we obtain 
\[
|Q(p+e)\cap U| \leq (3^n-1)\varepsilon r^n/3^n + \varepsilon r^n/3^n + (1 - 2\varepsilon)r^n = (1 - \varepsilon)r^n.
\]
Consequently, 
\[
|Q(p+e) \setminus U| \geq \varepsilon r^n.
\] 

Further, denoting $D: = Q(p+e) \setminus U$, we show that there exists an $s$-section of $D$, $\{\xi'\in \R^{n-1}: (\xi',s)\in D\}$ for some $s\in [p_n+e_n, p_n+e_n+r]$,  with the measure greater or equal to $\varepsilon r^{n-1}/4$. Then, $\mu_{n-1}(D) \geq |\{\xi'\in \R^{n-1}: (\xi',s)\in D\}|  \geq \varepsilon r^{n-1}/4$, and by \eqref{est:Q/U} we get the desired estimate.
By Fubini's theorem 
\[
\varepsilon r^n \leq  |D| = \int_{p_n+e_n}^{p_n+e_n+ r} \left( \int \mathbf{1}_D(\xi',\xi_n)\, \ud \xi' \right) \ud \xi_n.
\]
Therefore, there exists a set $D_n\subset [p_n+e_n, p_n+e_n+r]$ of positive measure such that for every $s \in D_n$ 
\[
\int \mathbf{1}_D (\xi', s) \ud \xi' = |\{ \xi' \in \R^{n-1}: (\xi',s) \in D\}|  \geq \varepsilon r^{n-1}/4.
\]
Thus, the proof of (i) is complete.

For (ii), let $U$ be as in (i).  Let $W: = \{\xi \in \Omega : b(\xi) < \varepsilon\}\subset U$
and $U': = \{ \xi\in \Omega: b(\xi)<\varepsilon/2\}$. 
If $U'=\emptyset$, then \(b\ge\varepsilon/2\) on \(\Omega\), and the conclusion follows directly from Theorem \eqref{thm: pencil est}. Hence we may assume that 
$U'\neq\emptyset$.
The uniform continuity of $b$ yields
\begin{equation}\label{def: 3rho}
  \rho: = \mathrm{dist}(\Omega\setminus U, \overline{U'}) \geq  \mathrm{dist}(\Omega \setminus W, \overline{U'}): =  \inf_{\xi\in \Omega \setminus  W, \zeta\in \overline{U'}}\|\xi-\zeta\| >0.
\end{equation}

 Choose $\eta\in C^\infty(\R)$ with $0\le\eta\le 1$, $\eta\equiv 0$ on
$(-\infty,1/2]$ and $\eta\equiv 1$ on $[1,\infty),$ and set
\[
\chi(\xi):=\eta\!\left(\rho^{-1}\mathrm{dist}(\xi, \Omega\setminus U)\right).
\]
Clearly, $\chi(\xi) = 1$ for any $\xi$ with $\dist(\xi, \Omega\setminus U)> \rho$ and $\chi(\xi) = 0$  for any $\xi$ with $\dist(\xi, \Omega\setminus U)\leq \rho/2$. 
Moreover, since $\mathrm{dist}(\cdot, \Omega\setminus U)$ is a $1$-Lipschitz function, Rademacher's theorem implies that it is differentiable a.e., with $
|\nabla \mathrm{dist}(\cdot, \Omega\setminus U)|\le 1$ a.e. Hence, by the chain rule, $\|\nabla\chi(\xi)\| \leq {\|\eta'\|_\infty}/\rho$ for a.e. $\xi\in\R^n$. Thus, $\chi$ has the regularity and support properties required in Theorem~\ref{thm: pencil est}, and the conclusions of that theorem apply.

Finally, (iii) follows from part (ii) and Proposition \ref{cond: 1/t at infty}. \end{proof}

\smallskip

\begin{proof}[Proof of Proposition \ref{prop:a_la_ui}]
(i) Let $s\neq 0$ and $\varphi\in \Dom(\Delta)$. We set $f: = (-s^2 + is b - \Delta)\varphi$ .  By~\eqref{ineq:b_*}, we have
\begin{align*}
|s| \|\varphi\|_2
&\leq c_s \bigl(\|f\|_2+|s|\|(b_*-b)\varphi\|_2\bigr),
\end{align*}
and therefore
\begin{align}\label{ineq:epsilon}
\nonumber \|\varphi\|_2^2
\nonumber &\leq 2c_s^2\bigl(s^{-2}\|f\|_2^2+\|(b_*-b)\varphi\|_2^2\bigr)\\
\nonumber &\leq 2c_s^2\bigl(s^{-2}\|f\|_2^2+(1+c_*)^2\|b\varphi\|_2^2\bigr)\\
&\leq 2c_s^2\bigl(s^{-2} \|f\|_2^2 + (1+c_*)^2\|b\|_\infty \|\sqrt{b}\varphi\|_2^2\bigr),
\end{align}
where the second inequality follows from the fact that
$0\leq b_*(\xi)\leq c_*b(\xi)$ for a.e.~$\xi\in \Omega$. 
Since $\Im \ip{f}{\varphi} = s \|\sqrt{b} \varphi\|_2^2$, we have 
\[
|s|\|\sqrt{b}\,\varphi\|_2^2
\leq \|f\|_2\|\varphi\|_2.
\]
Since for every $\varepsilon>0$ we have 
\[
|s|^{-1} \|f\|_2 \|\varphi\|_2 \leq (\varepsilon^{-2}/2)  |s|^{-2}\|f\|_2^2+ (\varepsilon^2/2) \, \|\varphi\|_2^2,
\]
choosing $\varepsilon$ such that $ 2 \varepsilon^2  c_s^2 (1+c_*)^2\|b\|_\infty = 1$, we can continue \eqref{ineq:epsilon} to obtain
\[
\|\varphi\|^2/2 \leq 2 c_s^2 \, \left( 1 + c_s^2(1+c_*)^4 \|b\|^2_\infty \right) \,  |s|^{-2} \|f\|_2^2,  
\] 
therefore showing \eqref{ineq:b}.

For the last statement in the proposition, assume that $b_*\in\cB(\Omega)$. Then \eqref{ineq:b_*} holds for every $s\neq 0$ with
\[
c_s \leq
\sup_{s\neq0}
|s| \bigl\|(-s^2+isM_{b_*}-\Delta)^{-1}\bigr\|<\infty.
\]
Consequently, \eqref{ineq:b} also holds for every
$s\neq 0$ with $\widetilde c_s$
uniformly bounded in $s\neq 0$. Since $P(is)^*=P(-is)$ for all $s\neq 0$, we conclude that
$i\mathbb{R}\setminus\{0\}\subset\rho(P)$, and \eqref{ineq:b} yields
$[R]<\infty$, i.e., $b\in \cB(\Omega)$.

\smallskip

(ii) Let $G:=\{\xi\in\Omega:b(\xi)>0\}$. Choose a sequence
$(\psi_j) \subset C_c^\infty(G)$ such that $0\le\psi_j\le1$ and 
$G=\bigcup_{j\ge1}\{\xi: \psi_j(\xi)>0\}$. Let $c_j: = \min \{b(\xi):  \xi \in \mathrm{supp}(\psi_j) \}$, $j\in \N$. Clearly, $c_j>0$ for all $j\in \N$.
Let $(p_\ell)_{\ell\ge1}$ be a sequence of seminorms defining the usual
Fr\'echet topology of $C^\infty(\Omega)$. Choose $\varepsilon_j>0$ so that
\[
\varepsilon_j\le 2^{-j}c_j,
\qquad \varepsilon_jp_\ell(\psi_j)\le2^{-j},
\quad 1\le \ell\le j.
\]
Then the series $a :=\sum_{j=1}^\infty\varepsilon_j\psi_j$
converges in $C^\infty(\Omega)$.

If $\psi_j(\xi)>0$, then $\xi\in\supp\psi_j$, and hence
$c_j\le b(\xi)$. Therefore,
\[
\begin{aligned}
a(\xi) &= \sum_{j=1}^\infty\varepsilon_j\psi_j(\xi)
 \le \sum_{\psi_j(\xi)>0}2^{-j}c_j\psi_j(\xi) \\
&\le b(\xi)\sum_{j=1}^\infty2^{-j} \le b(\xi).
\end{aligned}
\]
Moreover, $a(\xi)>0$ whenever $b(\xi)>0$, while all $\psi_j$ vanish
outside $G$. Thus $0\le a\le b$ and $\{\xi: a(\xi) > 0\}=\{\xi:b(\xi) >0\}$.

Choose radii $r_0<r_1<r_*$ and a radial cut-off
$\chi_\infty\in C^\infty(\mathbb R^n)$ with $0\le\chi_\infty\le1$,
$\chi_\infty \equiv 0$ on $\{|\xi|\le r_1\}$, and
$\chi_\infty \equiv 1$ on $\{|\xi|\ge r_*\}$.
In particular,
$\supp\chi_\infty\subset\{|\xi|>r_0\}$, where
$b\ge b_0$.
We define $b_*:=(1-\chi_\infty)a + \chi_\infty$. 
Then $b_*\in C^\infty(\Omega)\cap L_+^\infty(\Omega)$ and
$ b_*\equiv 1$ on $\{|\xi|\ge r_*\}$. 

Using $a\le b$ together with
$\chi_\infty\le b_0^{-1}\chi_\infty b$, we obtain
\[
\begin{aligned}
b_*
&\le (1-\chi_\infty)b+b_0^{-1}\chi_\infty b \\
&= \bigl(1-\chi_\infty+b_0^{-1}\chi_\infty\bigr)b \\
&\le \max\{1,b_0^{-1}\}\,b.
\end{aligned}
\]
Hence one may take $c_* =\max\{1,b_0^{-1}\}$.

Finally, if $b(\xi)=0$, then $a(\xi)=0$, while
$\chi_\infty(\xi)=0$ because $\supp\chi_\infty$ is contained in the
region where $b\ge b_0>0$. Thus $b_*(\xi)=0$.
Conversely, if $b(\xi)>0$, then $a(\xi)>0$. If
$\chi_\infty(\xi)<1$, then
$b_*(\xi)\ge(1-\chi_\infty(\xi))a(\xi)>0$, whereas if
$\chi_\infty(\xi)=1$, then $b_*(\xi)=1$. Consequently,
$\{\xi \in \Omega: b_*(\xi)>0\}=\{\xi \in \Omega: b(\xi)>0\}$.

Since the exterior geometric control condition depends only on the
positivity set, it is satisfied by $b_*$ if and only if it is satisfied
by $b$. \end{proof}

\smallskip

We begin the proof of Proposition~\ref{prop:interm_freq} with an auxiliary lemma. For a subset $G\subset\Omega$, $G\neq\Omega$, we write
$\|\cdot\|_{2,G}$ for the norm in $L^2(G)$, while $\|\cdot\|_2$ denotes the norm in $L^2(\Omega)$. We recall that $\langle\cdot,\cdot\rangle$ stands for the inner product in $L^2(\Omega)$.

\begin{lemma}\label{lem:interm_freq} Let $\Omega$, $U$, $U'$ and $\chi$ satisfy the assumptions of Proposition \ref{prop:interm_freq}. 
Then the operator $K: = K_{U, \chi}: W^{1,2}_0(\Omega) \rightarrow W^{-1,2}(\Omega)$ given by 
\[
 (K \varphi) (\psi)  : = \int_U \chi \varphi\,\overline{\psi} \ud \xi, \qquad \varphi, \psi \in W^{1,2}_0(\Omega) 
\]
is compact.
\end{lemma}

\begin{proof}
Let $(\varphi_j)$ be a bounded sequence in $Z:=W_0^{1,2}(\Omega)$. Then $(\chi \varphi_j)$ is a bounded sequence in $W^{1,2}_0(U)$, and using compactness of the embedding of $W_0^{1,2}(U)$ into $L^2(U)$, passing to a subsequence, we may assume that $(\chi\varphi_j)$ is convergent in $L^2(U)$. 
Let $\psi\in Z$ with $\|\psi\|_Z\le 1$, where $\|\psi\|_Z^2:=\|\nabla \psi\|_{2}^2+\|\psi\|_{2}^2$. 
Then
\[
|(K(\varphi_j - \varphi_\ell)) ( \psi) |
\le \|\chi \varphi_j - \chi \varphi_\ell\|_{2,U}\,\|\psi\|_{2,U}
\le \|\chi \varphi_j - \chi \varphi_\ell\|_{2,U}.
\]
Taking the supremum over $\|\psi\|_Z\le 1$ gives
\[
\|K(\varphi_j - \varphi_\ell)\|_{Z'}\le \|\chi\varphi_j - \chi\varphi_\ell\|_{2,U}.
\]
Therefore $(K \varphi_j)$ is Cauchy in $Z'$, so convergent. Consequently, $K$ is compact.
\end{proof}

\begin{proof}[Proof of Proposition \ref{prop:interm_freq}]
We set $Z:=W_0^{1,2}(\Omega)$. 
For $ s \neq 0$ we define the bounded sesquilinear form on $Z$ by 
\[
\mathfrak a_{is} (\varphi,\psi):=\langle \nabla \varphi , \nabla  \psi \rangle  + is \langle  b \varphi, \psi \rangle - s^2\langle \varphi,\psi \rangle,
\qquad \varphi,\psi\in Z.
\]
and its \emph{modification} 
\[
\widetilde {\mathfrak{a}}_{is}(\varphi,\psi):=\mathfrak a_{is} (\varphi,\psi) + (1 + s^2) (K\varphi)(\psi), \qquad \varphi,\psi\in Z,
\] where $K = K_{U,\chi}$ is defined as in Lemma \ref{lem:interm_freq}.
We claim that, for every $s\neq0$, $\widetilde{\mathfrak a}_{is}$ is coercive up to a rotation. To prove the claim, let $\varphi\in Z$. Since
$b(\xi)\ge\varepsilon$ for a.e. $\xi\in\Omega\setminus U'$, where
$\varepsilon>0$ is a constant, we have 
\[
\Im \widetilde {\mathfrak{a}}_{is}(\varphi,\varphi) = s \int_\Omega b |\varphi|^2\,d\xi
\quad\Rightarrow\quad
\sgn(s)\,\Im \widetilde {\mathfrak{a}}_{is}(\varphi,\varphi) \ge \varepsilon |s|\,\|\varphi\|_{2,\Omega\setminus U'}^2.
\]
Moreover, since $\chi \equiv 1$ on $U'$,
\[
\Re \widetilde {\mathfrak{a}}_{is}(\varphi,\varphi)
=\|\nabla \varphi\|_{2}^2
- s^2 \|\varphi\|_{2, \Omega\setminus U'}^2 + \|\varphi\|_{2,U'}^2 + (1+s^2) \|\sqrt{\chi} \varphi\|^2_{2, U\setminus U'}.
\]
Pick $\eta = \eta(\varepsilon, s) \in(0,\pi/2)$ so that $\varepsilon \sin\eta = 2|s|\cos \eta$.
Then $\varepsilon |s| \sin\eta - s^2 \cos\eta  =  \frac12 \varepsilon |s| \sin\eta > 0$, and for $s>0$, since $\cos (\eta) (1+s^2)\|\sqrt{\chi} \varphi\|^2_{2, U\setminus U'}\geq 0$,
\begin{align*}
\Re\!\bigl(e^{-i\eta}\widetilde{\mathfrak{a}}_{is}(\varphi,\varphi)\bigr)
&= \cos\eta\,\Re\widetilde {\mathfrak{a}}_{is}(\varphi,\varphi) + \sin\eta\,\,\Im\widetilde{\mathfrak{a}}_{is}(\varphi,\varphi)\\
&\ge \cos\eta\|\nabla \varphi\|_{2}^2 + \cos\eta\|\varphi\|_{2, U'}^2
+(\varepsilon |s| \sin\eta - |s|^2 \cos\eta )\|\varphi \|_{2,\Omega\setminus U'}^2\\
&\ge c(\varepsilon, s)  \|\varphi\|_{Z}^2,
\end{align*}
where $ c(s,\varepsilon): = \min(\cos \eta,  \frac12 \varepsilon |s| \sin\eta) >0$.
Thus the form $e^{-i\eta}\widetilde {\mathfrak{a}}_{is}$ is coercive on $Z$.
Similarly, if $s<0$ then the form $e^{i\eta}\widetilde {\mathfrak{a}}_{is}$ is coercive on $Z$, as claimed.

We define the (bounded) operators $\mathcal P(is), \widetilde{\mathcal P}(is):Z\to Z'$ by
\[
\langle \mathcal P(is) \varphi, \psi\rangle_{Z',Z} := {\mathfrak{a}}_{is}(\varphi,\psi),\qquad
\langle \widetilde{\mathcal P}(is) \varphi, \psi\rangle_{Z',Z} := \widetilde {\mathfrak{a}}_{is}(\varphi,\psi),
\qquad \varphi,\psi\in Z.
\]
Applying the Lax--Milgram theorem, we conclude that, for every $s\neq 0$, the operator $e^{-i\sgn(s)\eta}\widetilde{\mathcal P}(is)$ is an isomorphism from $Z$ onto $Z'$; consequently, so is $\widetilde{\mathcal P}(is)$.

Since
\[
\widetilde{\mathcal P}(is) = \mathcal P(is) + (1+ s^2)K,
\] 
thus, by Lemma~\ref{lem:interm_freq}, the operator $\mathcal P(is)$ is a compact perturbation of the isomorphism $\widetilde{\mathcal P}(is)$. 
In particular, $\mathcal P(is) : Z\to Z'$ is Fredholm of index $0$. 
Therefore, in order to show that $\cP(is)$ is an isomorphism it suffices to prove that it has a trivial kernel. Let $\varphi\in Z$ satisfy $\mathcal P(is) \varphi = 0$, i.e. ${\mathfrak{a}}_{is}(\varphi,\psi)=0$ for all $\psi\in Z$.
Taking $\psi=\varphi$ gives
\[
\|\nabla \varphi\|_{2}^2-s^2\|\varphi\|_{2}^2 + is\int_\Omega b|\varphi|^2\,d\xi=0.
\]
Taking imaginary parts yields $\int_\Omega b|\varphi|^2\ud \xi=0$. Since $b(\xi)\ge \varepsilon >0$ a.e. on $\Omega\setminus U$, then $\varphi=0$ a.e.\ on $\Omega\setminus U$. 
Moreover, from ${\mathfrak{a}}_{is}(\varphi,\psi)=0$ and $ b(\xi) \varphi(\xi) =0$ for a.e. $\xi\in\Omega$, we obtain 
\[
\langle \nabla \varphi, \nabla \psi \rangle =   s^2\langle \varphi,\psi\rangle, \qquad \psi\in Z,
\]
i.e., $\varphi$ is a weak solution of $(\Delta+s^2)\varphi = 0$ in $\Omega$. 
By elliptic regularity, $\varphi$ is real analytic in the connected domain $\Omega$. The zero set of any nontrivial real analytic function has Lebesgue measure zero. Since $\varphi$ vanishes on the set $\Omega\setminus U$, which has positive measure, it follows that $\varphi\equiv0$ in $\Omega$.  Consequently,
$\Ker\mathcal{P}(is)=\{0\}$, and thus $\cP(is)$ is an isomorphism.

Let $f\in L^2(\Omega)\subset Z'$. Since $\mathcal P(is)$ is an isomorphism $Z\to Z'$, there exists a unique $\varphi\in Z$ such that
\[
\mathfrak{a}_{is}(\varphi,\psi)=\langle f,\psi\rangle,\qquad \psi\in Z.
\]
Equivalently,
\[
\langle \nabla \varphi , \nabla \psi \rangle = \langle f + s^2 \varphi - i s b \varphi,\ \psi \rangle,\qquad \psi\in Z.
\]
Because $b\in L^\infty(\Omega)$ and $\varphi \in L^2(\Omega)$, the right-hand side is represented by the function
$g:=f + s^2 \varphi -i s b \varphi \in L^2(\Omega)$. By the definition of the Dirichlet Laplacian, this implies $\varphi \in \Dom(\Delta)$ and $-\Delta \varphi=g$. Hence, $(-s^2 + is M_b - \Delta) \varphi = f$, i.e., $P(is) \varphi =f$. Finally, since $\Dom(\Delta) \subset Z$, so $P(is): \Dom(\Delta)\to L^2(\Omega)$ is bijective. Since $P(is)$ is closed, its inverse is bounded by the closed graph theorem. 
\end{proof}

\begin{proof}[Proof of Proposition \ref{cor:fixed-scale}] 
We set
\[
\nu :=\sup_{\lambda\in\Sigma_{\omega+\pi/2}}
\min(|\lambda|,|\lambda|^\theta) \|R_D(\lambda)B\|<\infty.
\]
Assume (i), and set
\[
c_0:=\sup_{0<t\leq1}t^{1-\theta}\|T(t)\|,
\qquad
c_y:=\sup_{t\geq1}t^\eta\|T(t)y\|.
\]
Theorem \ref{thm:mild_DPE}(i) gives $c_0<\infty$,
while assumption (i) ensures $c_y<\infty$.
Using \eqref{eq:laplace_int}, for $0<\mu\leq1$ we obtain
\[
\begin{aligned}
\|R_D(\mu)By\|
&\leq c_0\|y\|\int_0^1 e^{-\mu t}t^{\theta-1}\,\ud t 
+c_y\int_1^\infty e^{-\mu t}t^{-\eta}\,\ud t\\
&\leq \left({c_0}\|y\|/{\theta} 
+ c_y \Gamma(1-\eta)\right)\mu^{\eta-1}.
\end{aligned}
\]
For $\lambda\in\Sigma_{\omega+\pi/2}$ with
$0<|\lambda|\leq1$, set $\mu:=|\lambda|$.
By the pseudo-resolvent identity (5.3) and
the definition of $\nu$,
\[
\begin{aligned}
\|R_D(\lambda)By\| 
&\leq \bigl(1+|\mu-\lambda|\|R_D(\lambda)B\|\bigr)
\|R_D(\mu)By\|\\
&\leq (1+2\nu) \left({c_0}\|y\|/{\theta} 
+ c_y \Gamma(1-\eta)\right) |\lambda|^{\eta-1}.
\end{aligned}
\]
This proves (ii).

Conversely, assume (ii) for some $\omega'\in(0,\omega)$,
and put
\[
\wt c_y:=\sup_{\substack{
\lambda\in\Sigma_{\omega'+\pi/2}\\
0<|\lambda|\leq1}}
|\lambda|^{1-\eta}\|R_D(\lambda)By\|<\infty.
\]
By \eqref{formula:semig}, for every $\rho>0$,
\[
T(t)y=\frac{1}{2\pi i}
\int_{\Gamma(\rho,\omega')}
e^{\lambda t}R_D(\lambda)By\,d\lambda.
\]
The contribution from the circular arc is
$O(\rho^\eta)$ as $\rho\downarrow0$.
Passing to the limit and using (ii) and
\eqref{DE-parab_est}, we obtain
\[
\begin{aligned}
\|T(t)y\|
&\leq \frac{1}{\pi}\left(
\wt c_y \int_0^1 e^{-tr\sin\omega'}r^{\eta-1}\,\ud r
+\nu \|y\|\int_1^\infty
e^{-tr\sin\omega'}r^{-\theta}\,\ud r
\right)\\
&\leq c t^{-\eta},\qquad t\geq1,
\end{aligned}
\]
where $c>0$ is independent of $t$. Hence (i) follows, completing the proof.
\end{proof}

\subsection{Proof of Theorem \ref{thm:diffusion_for_Delta}}\label{subs:diffusion_for_Delta} 
To verify that all assumptions of Theorems~\ref{thm:mild_DPE}--\ref{thm: diffusion} are satisfied, we first note that, in both cases (i) and (ii), combining Theorem~\ref{thm: DWE well-posed} with Proposition~\ref{prop:linear_growth}, we conclude that the corresponding (DWE) is well posed and that the $C_0$-semigroup $\mathcal{U}$ associated with the first-order problem exhibits at most linear growth as $t\to\infty$.  
Moreover, to verify the condition $B(\Dom(A))\subset \Dom\bigl((A+\omega)^{1/2}\bigr)$ in Theorem~\ref{thm: diffusion} (with $\omega = 0$, since $-\Delta$ is sectorial), in each case, $\Omega = \R^n$ or $\Omega = \R^n\setminus \overline{\cO}$,
since $b\in W^{1,\infty}(\Omega)$ and $\Dom(\Delta)\subset W^{1,2}_0(\Omega)$, we have 
\[
M_b(\Dom(\Delta)) \subset b\cdot\bigl(W^{1,2}_0(\Omega)\bigr)
\subset W^{1,2}_0(\Omega) = \Dom((-\Delta)^{1/2}).
\]

In the case of (i), since $b\in BUC(\R^n)$, Corollary~\ref{cor:agcc_implies_thm21} implies that $b\in \mathcal{B}(\R^n)\cap \mathcal{B}_D^{1/2}(\R^n)$,
and hence Hypothesis~\ref{hyp:main_2} with $\theta =\frac12$ holds. In particular, by Theorem~\ref{thm:mild_DPE}(ii), we have
$v_0\in \Delta^{-1}(\Ran(M_b))\subset Y_D^{1/2}=X_D^1$. Therefore,
$v_0$ satisfies the assumption on the initial datum corresponding to the problem (PE) in Theorem~\ref{thm: diffusion intro}(i).
 The remaining condition, namely $\Ran(M_b)\setminus\Ran(\Delta)\neq\emptyset$, required in Theorem~\ref{thm: diffusion}, is ensured by Corollary~\ref{cor:agcc_implies_r} (or, in the case $n\leq 4$, by Proposition~\ref{prop:RninR}(ii)).

For part (ii), Lemma~\ref{lem:B_infty_class}(ii)
gives the high-frequency hyperbolic estimate,
while Corollary~\ref{cor:examples_Omega_b}(ii)
supplies the low-frequency hyperbolic and
parabolic resolvent estimates.
Proposition~\ref{prop:correction_b}
(or Proposition~\ref{prop:ncc_implies_R})
gives the range condition. Together with
the well-posedness and compatibility
properties verified above, these imply
the assertion.

\providecommand{\bysame}{\leavevmode\hbox to3em{\hrulefill}\thinspace}

\bibliographystyle{abbrv}
\bibliographystyle{amsplain}

\end{document}